\documentclass[12pt,reqno]{amsart}
\usepackage{amsfonts}
\usepackage{amsmath}
\usepackage{amsthm}
\usepackage{tikz-cd}
\usepackage{mathtools}
\usepackage{amssymb}
\usepackage{mathrsfs}
\usepackage{dsfont}
\usepackage[shortlabels]{enumitem}
\usepackage[pagebackref=true]{hyperref}
\usepackage[dvipsnames]{xcolor}
\hypersetup{
    colorlinks=true,
    linkcolor=Blue,
    filecolor=red,      
    urlcolor=Plum,
    citecolor=teal,
}
\usepackage[alphabetic]{amsrefs}

\usepackage{setspace}
\newcounter{mainstepcounter}[subsubsection] 
\newcounter{altstepcounter}[subsection] 
\newcounter{theoremalph}[section] 

\theoremstyle{definition}
\newtheorem{Theorem}[subsubsection]{Theorem}
\newtheorem{Lemma}[subsubsection]{Lemma}
\newtheorem{Corollary}[subsubsection]{Corollary}
\newtheorem{Proposition}[subsubsection]{Proposition}
\newtheorem{Definition}[subsubsection]{Definition}

\newtheorem{Example}[subsubsection]{Example}
\newtheorem{Notation}[subsubsection]{Notation}
\newtheorem{Convention}[subsubsection]{Convention}
\newtheorem{Remark}[subsubsection]{Remark}
\newtheorem{Discussion}[subsubsection]{Discussion}
\newtheorem{Construction}[subsubsection]{Construction}
\newtheorem{Situation}[subsubsection]{Situation}

\newtheorem{step}[mainstepcounter]{Step}
\newtheorem{Step}[altstepcounter]{Step} 

\newtheorem{claim}{Claim}[subsubsection]

\newtheorem{IntroTheorem}{Theorem}[theoremalph]

\newtheorem{IntroCorollary}[IntroTheorem]{Corollary}

\makeatletter
\def\subsubsection{\@startsection{subsubsection}{3}%
  \z@{.5\linespacing\@plus.7\linespacing}{-.5em}%
  {\normalfont\bfseries}}
\makeatother

\newcounter{dummy}
\makeatletter
\newcommand\myitem[1][]{\item[#1]\refstepcounter{dummy}\def\@currentlabel{#1}}
\makeatother

\makeatletter
\DeclareRobustCommand{\cev}[1]{%
  \mathpalette\do@cev{#1}%
}
\newcommand{\do@cev}[2]{%
  \fix@cev{#1}{+}%
  \reflectbox{$\m@th#1\vec{\reflectbox{$\fix@cev{#1}{-}\m@th#1#2\fix@cev{#1}{+}$}}$}%
  \fix@cev{#1}{-}%
}
\newcommand{\fix@cev}[2]{%
  \ifx#1\displaystyle
    \mkern#23mu
  \else
    \ifx#1\textstyle
      \mkern#23mu
    \else
      \ifx#1\scriptstyle
        \mkern#22mu
      \else
        \mkern#22mu
      \fi
    \fi
  \fi
}
\makeatother

\makeatletter
\def\l@subsection{\@tocline{2}{0pt}{2.5pc}{5pc}{}}
\makeatother

\DeclareMathAlphabet\EuScript{U}{eus}{m}{n}
\SetMathAlphabet\EuScript{bold}{U}{eus}{b}{n}

\numberwithin{equation}{subsubsection}

\newcommand{\bC}{\mathbb{C}}

\newcommand{\bE}{\mathbb{E}}

\newcommand{\bN}{\mathbb{N}}

\newcommand{\bP}{\mathbb{P}}
\newcommand{\bQ}{\mathbb{Q}}
\newcommand{\bR}{\mathbb{R}}

\newcommand{\bV}{\mathbb{V}}

\newcommand{\bZ}{\mathbb{Z}}

\newcommand{\clA}{\mathcal{A}}

\newcommand{\clE}{\mathcal{E}}

\newcommand{\clG}{\mathcal{G}}

\newcommand{\clK}{\mathcal{K}}
\newcommand{\clL}{\mathcal{L}}

\newcommand{\clN}{\mathcal{N}}
\newcommand{\clO}{\mathcal{O}}

\newcommand{\clR}{\mathcal{R}}

\newcommand{\fj}{\mathfrak{j}}

\newcommand{\fm}{\mathfrak{m}}

\newcommand{\fp}{\mathfrak{p}}

\newcommand{\fs}{\mathfrak{s}}

\newcommand{\fC}{\mathfrak{C}}
\newcommand{\fD}{\mathfrak{D}}
\newcommand{\fE}{\mathfrak{E}}

\newcommand{\fL}{\mathfrak{L}}
\newcommand{\fM}{\mathfrak{M}}

\newcommand{\fP}{\mathfrak{P}}

\newcommand{\fS}{\mathfrak{S}}

\newcommand{\fX}{\mathfrak{X}}
\newcommand{\fY}{\mathfrak{Y}}

\newcommand{\sC}{\mathscr{C}}

\newcommand{\sF}{\mathscr{F}}

\newcommand{\esC}{\EuScript{C}}

\newcommand{\esO}{\EuScript{O}}

\newcommand{\esS}{\EuScript{S}}

\newcommand{\blambda}{\boldsymbol{\lambda}}
\newcommand{\bmu}{\boldsymbol{\mu}}
\newcommand{\bs}{\boldsymbol{s}}

\newcommand{\ul}{\underline}
\newcommand{\ol}{\overline}

\newcommand{\rrarrows}{\;\substack{\longrightarrow\\[-0.4ex]\longrightarrow}\;}
\newcommand{\taut}{\mathds{1}}
\DeclareMathOperator\Spec{Spec}

\DeclareMathOperator\relSpec{\bf Spec}

\DeclareMathOperator\coker{coker}
\DeclareMathOperator\img{im}
\DeclareMathOperator\Hom{Hom}

\DeclareMathOperator\Sym{Sym}

\DeclareMathOperator\Ext{Ext}
\DeclareMathOperator\Trop{Trop}

\DeclareMathOperator{\PreLog}{{\bf PreLog}}
\DeclareMathOperator{\Log}{{\bf Log}}
\DeclareMathOperator{\Sch}{{\bf Sch}}
\DeclareMathOperator{\LogSch}{{\bf LogSch}}
\DeclareMathOperator{\LogPt}{{\bf LogPt}}
\DeclareMathOperator{\SlbPt}{{\bf SlbPt}}
\DeclareMathOperator{\Zar}{Zar}
\DeclareMathOperator{\Mon}{{\bf Mon}}
\DeclareMathOperator{\Ab}{{\bf Ab}}

\DeclareMathOperator\eval{eval}

\newcommand{\plC}{\scalebox{0.8}[1.3]{$\sqsubset$}}

\begin{document}
    
    \title[What is a stable log map?]{What is a stable log map?}
    
    \author{Mohammad Farajzadeh-Tehrani}
    \address{Department of Mathematics, University of Iowa}
    \email{mohammad-tehrani@uiowa.edu}
    
    \author{Mohan Swaminathan}
    \address{School of Mathematics, Tata Institute of Fundamental Research, Mumbai}
    \email{mohans@math.tifr.res.in}
    
    \date{\today}

    \begin{abstract}
        Let $X$ be a smooth projective variety over $\bC$ with a simple normal crossings divisor $D\subset X$. We compare the notions of stable log maps to $(X,D)$ in algebraic geometry and symplectic topology. In particular, we prove an equivalence between fine (basic) algebraic log maps and symplectic log maps, and we define the symplectic analogue of fine saturated algebraic log maps by refining the notion of log Gromov convergence.
    \end{abstract}

    \maketitle

    \setcounter{tocdepth}{1}
    \tableofcontents
    \section{Introduction}\label{sec:intro}

Let $X$ be a smooth projective variety $X$ over $\bC$ equipped with a simple normal crossings divisor $D\subset X$ and consider the moduli space of maps from smooth curves to $X$ which have prescribed contact orders with the irreducible components of $D$. A natural compactification of this moduli space is provided by \emph{stable log maps} to $(X, D)$, following either the work of Abramovich--Chen--Gross--Siebert in the algebraic category (\cites{GS-log, Chen-logDF1, AC-logDF2}) or the work of the first author in the symplectic category (\cite{Tehrani-log}); see also the work of Parker in the category of exploded manifolds (\cite{Parker-exploded-vs-log}). When $D$ is a smooth divisor, a different compactification using \emph{stable relative maps} to $(X,D)$ was discovered earlier (\cites{Jun-Li-1, Jun-Li-2} in the algebraic category and \cites{Li-Ruan-relGW, Ionel-Parker-relGW} in the symplectic category).

The aim of this article is to clarify the relation between the algebraic and symplectic notions of stable log maps. The algebraic approach is elegant and natural, though its geometric content can feel somewhat hidden behind the abstract machinery of log structures. The symplectic perspective is more explicit and, importantly, extends naturally to the almost Kähler setting. In particular, it makes transparent the additional data required to lift an ordinary stable map to a stable log map. However, this perspective focuses primarily on individual stable log maps and leaves open the question of how to organize these objects into families. Consequently, the topology (and any further structure) on the moduli space must be constructed separately.

\subsection{Main results} 

We need to briefly recall the algebraic and symplectic approaches to defining stable log maps before stating our main results.

\subsubsection{Algebraic log maps}

The algebraic approach to stable log maps uses the language of logarithmic algebraic geometry. For the purpose of the introduction, let us provide an informal idea of log structures on schemes based on the two desiderata below; for more details, including the precise definition of log structures, see Subsection~\ref{subsec:log-schemes}.
\begin{itemize}
    \item If $(Y,E)$ is a pair consisting of a variety $Y$ and a closed codimension $1$ subset $E\subset Y$ having `mild singularities', then it determines a log structure $M_{(Y,E)}$ on $Y$.

    \item Given a log structure $M_Y$ on $Y$ and a morphism $f:Z\to Y$, there is a well-defined pullback log structure $f^*M_Y$ on $Z$. Any log structure on a scheme $Z$ is required to locally be the pullback of $M_{(Y,E)}$ for some pair $(Y,E)$ with `mild singularities'.
\end{itemize}

If we take `mild singularities' to mean that $(Y,E)$ is \emph{toroidal} (i.e., locally modeled on a toric variety and its boundary divisor), then we obtain the class of \emph{fine} log schemes. If we take `mild singularities' to mean that $(Y,E)$ is toroidal and $Y$ is normal, then we obtain the smaller subclass of \emph{fine saturated} (\emph{fs}) log schemes. In particular, endowing the variety $X$ with the log structure $M_{(X,D)}$ yields an fs log scheme, which we denote by $X^\dagger$.

An \emph{algebraic log map} to the log scheme $X^\dagger$ is a morphism of log schemes whose target is $X^\dagger$ and whose source is a nodal curve with marked points $(C,{\bf x})$ endowed with a log structure $M_C$ which is \emph{log smooth} over a log structure $M_S$ on the one-point scheme $S = \Spec\bC$; for more details, see Definition~\ref{def:alg-log-map}. Given an algebraic log map as above, the log structures $M_C$ and $M_S$ are not determined by the underlying morphism of schemes from $C$ to $X$. However, once we fix certain discrete invariants, there is a `universal' choice of $M_C$ and $M_S$. The meaning of `universal' is sensitive to which category of log schemes we work with: fine or fs. We refer to the `universal' choices in these two settings as \emph{fine-basic} and \emph{fs-basic} algebraic log maps respectively; for more details, see Definitions~\ref{def-part:fine-basic} and \ref{def-part:fs-basic}.

\subsubsection{Symplectic log maps} 

Enumerate the irreducible components of $D$ as $D_1,\ldots,D_r$ and write $[r] = \{1,\ldots,r\}$. A \emph{symplectic log map} to the pair $(X,D)$ consists of
\begin{itemize}
    \item a holomorphic map $f:C\to X$ from a nodal curve with marked points $(C,{\bf x})$, and
    
    \item a collection of integers $\bmu = \{\mu_{i,j},\mu_{\vec{e},j}\}$, where we have one integer for each pair consisting of a special point of the normalization $\tilde C$ (i.e., a marked point $x_i$ or an inverse image $y_{\vec{e}}$ of a node $y_e$) and an irreducible component $D_j$ of $D$,
\end{itemize}
subject to the following conditions.
\begin{itemize}    
    \item For each node $y_e$, denote the two special points lying over it by $y_{\vec{e}}$ and $y_{\cev{e}}$. The vectors ${\bmu}_{\vec{e}} = (\mu_{\vec{e},1},\ldots,\mu_{\vec{e},r})\in\bZ^r$ and ${\bmu}_{\cev{e}} = (\mu_{\cev{e},1},\ldots,\mu_{\cev{e},r})\in\bZ^{r}$ are negatives of one another.
    
    \item For each irreducible component $C_v$ of $C$, write $I_v$ for the set of $j$ with $f(C_v)\subset D_j$. Whenever $j\not\in I_v$, the map $f|_{\tilde C_v}$ meets $D_j$ only at the special points, with contact orders given by $\mu_{i,j}$ and $\mu_{\vec{e},j}$. Whenever $j\in I_v$, there exists a nonzero meromorphic section $\zeta_{v,j}$ of the line bundle $(f|_{\tilde C_v})^*\clN_{D_j/X}$ such that the zeros and poles of $\zeta_{v,j}$ occur only at the special points, and their orders are given by $\mu_{i,j}$ and $\mu_{\vec{e},j}$.

    \item There exist choices of vectors $\bs_v\in\bR_{>0}^{I_v}\times\{0\}^{[r]\setminus I_v}\subset\bR^r$ for each irreducible component $C_v$ of $C$ such that the vectors $\bs_{v'} - \bs_v\in\bR^r$ and $\bmu_{\vec{e}}\in\bR^r$ are positive scalar multiples of one another whenever we have $y_{\vec{e}}\in\tilde C_v$ and $y_{\cev{e}}\in\tilde C_{v'}$.

    \item There exist choices of meromorphic sections $\zeta_{v,j}$ as above together with choices of local holomorphic coordinates centred at each special point $y_{\vec{e}}$ such that the lowest order terms of $\zeta_{v,j}$ at $y_{\vec{e}}$ and $\zeta_{v',j}$ at $y_{\cev{e}}$ match whenever we have $y_{\vec{e}}\in\tilde C_v$ and $y_{\cev{e}}\in\tilde C_{v'}$.
\end{itemize}

In other words, a symplectic log map is a nodal curve mapping holomorphically to $X$ along with a prescription of contact orders at each special point of the curve with each irreducible component of $D$, subject to some combinatorial and analytic conditions. For more details and motivation, see Definition~\ref{def:symp-log-map}, Remark~\ref{rem:symp-log-map} and Discussion~\ref{disc:log-map-motivation}.

\subsubsection{Comparison} As before, write $X^\dagger$ for the log scheme obtained by endowing the smooth variety $X$ with the log structure $M_{(X,D)}$ determined by the simple normal crossings divisor $D\subset X$. We are now in a position to state our main results.

\begin{IntroTheorem}[{Theorem~\ref{thm:fine-log-map-comparison} \& Corollary~\ref{cor:fine-log-map-comparison}}]\label{intro-thm:fine-log-map-comparison}
    Any fine algebraic log map to $X^\dagger$ has a well-defined underlying symplectic log map to $(X,D)$. Conversely, any symplectic log map to $(X,D)$ lifts uniquely (up to unique isomorphism) to a fine-basic algebraic log map to $X^\dagger$.
\end{IntroTheorem}

The papers \cites{GS-log, Chen-logDF1, AC-logDF2} work exclusively with fs log schemes. Therefore, they do not consider the case of fine log schemes separately and, instead of our term \emph{fs-basic}, they use the term \emph{basic} (in the case of \cite{GS-log}) or \emph{minimal} (in the case of \cites{Chen-logDF1, AC-logDF2}). 

To find a symplectic analogue of fs algebraic log maps, we define a new notion of \emph{saturated symplectic log map} and a corresponding refinement of the notion of \emph{log Gromov convergence} from \cite{Tehrani-log}; see Definitions~\ref{def-part:sat-symp-log-map} and \ref{def:sat-log-gromov-conv}.

\begin{IntroTheorem}[{Theorem~\ref{thm:fs-log-map-comparison} \& Corollary~\ref{cor:fs-log-map-comparison}}]\label{intro-thm:fs-log-map-comparison}
    Any fs algebraic log map to $X^\dagger$ has a well-defined underlying saturated symplectic log map to $(X,D)$. Conversely, any saturated symplectic log map to $(X,D)$ lifts uniquely (up to unique isomorphism) to a fs-basic algebraic log map to $X^\dagger$, i.e., a \emph{basic log map} in the language of \cite{GS-log}.
\end{IntroTheorem}

A feature of our definition of saturated symplectic log maps is that any symplectic log map lifts to a saturated symplectic log map and, while this lift may not be unique, the number of distinct lifts is finite and is given by an explicit formula; see Lemma~\ref{lem:number-of-saturations}. Together with Theorem~\ref{intro-thm:fs-log-map-comparison}, this leads to the following result, which verifies \cite[Conjecture~1.1]{Tehrani-log}.

\begin{IntroCorollary}[{Corollary~\ref{cor:number-of-saturations} \& Remark~\ref{rem:number-of-saturations}}]\label{intro-cor:number-of-saturations}
    There exists a well-defined and surjective forgetful map from the points of the moduli space of fs-basic algebraic stable log maps to $X^\dagger$ constructed in \cites{GS-log, Chen-logDF1, AC-logDF2} to the points of the moduli space of symplectic stable log maps to $(X,D)$ constructed in \cite{Tehrani-log}. Moreover, for any symplectic stable log map, its inverse image under this forgetful map is a non-empty finite set whose cardinality is given by the formula in \cite[Equation~(5.15)]{Tehrani-log}.
\end{IntroCorollary}

We may summarize our main results by the following diagram, where the dashed vertical arrow on the left is the unique map which makes the diagram commute.

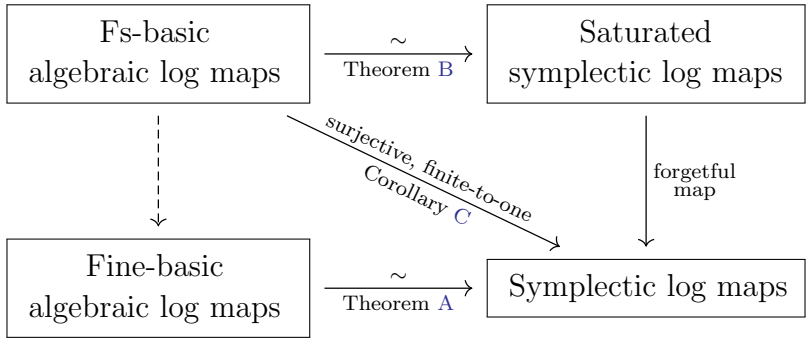
\begin{equation*}
\begin{tikzcd}
    \fbox{\begin{minipage}{3.8cm}{$\begin{array}{c} \text{Fs-basic} \\ \text{algebraic log maps} \end{array}$}\end{minipage}} \arrow[dd,dashed] \arrow[ddrr] \arrow[ddrr,phantom,shift left=1ex,"\mathbin{\rotatebox[origin=c]{334.5}{\tiny{surjective, finite-to-one}}}"] \arrow[ddrr,phantom,shift right=1.2ex,"\mathbin{\rotatebox[origin=c]{334.5}{\tiny{Corollary~\ref{intro-cor:number-of-saturations}}}}"] \arrow[rr,"\sim","\text{Theorem~\ref{intro-thm:fs-log-map-comparison}}"'] & & \fbox{\begin{minipage}{4cm}{$\begin{array}{c} \text{Saturated} \\ \text{symplectic log maps} \end{array}$}\end{minipage}} \arrow[dd,"\substack{\text{forgetful} \\ \text{map}}"] \\ \\
    \fbox{\begin{minipage}{3.8cm}{$\begin{array}{c} \text{Fine-basic} \\ \text{algebraic log maps} \end{array}$}\end{minipage}} \arrow[rr,"\sim","\text{Theorem~\ref{intro-thm:fine-log-map-comparison}}"'] & & \fbox{\begin{minipage}{4cm}{$\begin{array}{c} \text{Symplectic log maps} \end{array}$}\end{minipage}}
\end{tikzcd}
\end{equation*}

Let us elaborate a bit more on the relationship between fine-basic and fs-basic algebraic log maps. For any symplectic log map, consider the category of its lifts to a fine algebraic log map. This category is non-empty and, moreover, it has a final object (called a fine-basic algebraic log map). The category of lifts which are fs algebraic log maps is also non-empty but may no longer have a final object. Instead, it decomposes into finitely many connected components (i.e., pairwise disjoint full subcategories with no morphisms between them), each of which has a final object (called an fs-basic algebraic log map). On the level of moduli spaces, one may think of saturation as being similar to normalization: locally, a collection of irreducible components meeting at a point get separated out into a collection of normal connected components. Accordingly, the idea behind defining a saturated symplectic log map is to endow an ordinary symplectic log map with an additional geometric decoration (called a \emph{saturation datum}) which picks out one of these connected components; see Example~\ref{exa:simplest-sat-example} for a concrete instance of this when $D\subset X$ is a smooth divisor.

\subsubsection{Ingredients of the proof} 

Most of the work in this article goes towards establishing Theorem~\ref{intro-thm:fine-log-map-comparison}. Once this is done, Theorem~\ref{intro-thm:fs-log-map-comparison} and Corollary~\ref{intro-cor:number-of-saturations} are proved by largely formal arguments. The proof of Theorem~\ref{intro-thm:fine-log-map-comparison} has the following three main ingredients.
\begin{itemize}
    \item {\bf Systems of line bundles:} A log structure on a scheme naturally gives rise to a \emph{system of line bundles} on it; see Definitions~\ref{def:sys-line-bundle} and \ref{def:log-lb-vcd}. For a large class of log schemes (called \emph{Deligne--Faltings log schemes}), which includes the target $X^\dagger$ of our algebraic log maps, this system of line bundles completely encodes the log structure; see Definition~\ref{def:df-log} and Remark~\ref{rem:log-str-ghost-slb}. In particular, for any log scheme $S^\dagger$, a log morphism $f^\dagger:S^\dagger\to X^\dagger$ can be completely described by its underlying morphism $f:S\to X$ and a finite collection of line bundle isomorphisms on $S$ (one isomorphism corresponding to each irreducible component of $D$); see Proposition~\ref{prop:log-mor-to-snc-pair}.

    \item {\bf Classification of log curves over a log point:} The domain of an algebraic log map is a nodal curve with marked points, which is equipped with a log structure $M_C$ that is log smooth over some log structure $M_S$ on the point $S = \Spec\bC$. In Theorem~\ref{thm:log-curve-over-log-pt}, we classify such objects and use the classification to give a concrete description of the system of line bundles associated to the log structure $M_C$. To the best of the authors' knowledge, this result is not contained in the existing literature; see Remark~\ref{rem:log-curves-whats-new} for a more detailed discussion. Together with the previous point, this leads to an explicit description of the data needed to lift an ordinary stable map to an algebraic stable log map (over some given log structure on the point); see Proposition~\ref{prop:alg-log-map}.

    \item {\bf Log structures on a point via generators and relations:} We define the notion of a \emph{presentation for a system of line bundles} and identify a simple necessary and sufficient condition (called \emph{consistency}) for such a presentation to be \emph{realized} by an actual system of line bundles; see Definition~\ref{def:slb-presentation} and Proposition~\ref{prop:slb-presentation}. Crucially, the conditions in the definition of a symplectic log map together imply the consistency of a certain presentation; see Lemma~\ref{lem:basic-fine-slb-consistent-symplectic-log-map}. Using this observation, we obtain a surjective forgetful map from algebraic log maps to symplectic log maps, with a unique fine-basic element in each fibre; see the proof of Theorem~\ref{thm:fine-log-map-comparison} for more details.
\end{itemize}

\subsection{Future directions}

In future work, we plan to use the techniques developed in the present article to give a `symplectic' description of stable log maps in the central fibre of a \emph{log smooth degeneration} as in \cite{ACGS-decomp}; the case of semistable degenerations is worked out in \cite{Tehrani-semistable}. In a different direction, we also plan to construct \emph{global Kuranishi charts} (building on \cite{AMS21,HS24,AMS24}) for moduli spaces of (saturated) symplectic log maps and use them to define log Gromov--Witten theory in the symplectic setting.

\subsection{Structure of the article}

In Section~\ref{sec:log-prelims}, we give a detailed review of the concepts from logarithmic geometry which are needed in later sections; readers having some familiarity with log geometry can skip this section and refer back to it as necessary. In Section~\ref{sec:slb-prelims}, we introduce and discuss the notion of a \emph{system of line bundles}. In Section~\ref{sec:log-via-slb}, we show how to reformulate log geometry using the language of systems of line bundles introduced in the previous section. In Section~\ref{sec:log-curves}, we study families of log curves and classify log curves over a log point with a given underlying prestable curve. Finally, in Section~\ref{sec:log-maps}, we study log maps and prove Theorems~\ref{intro-thm:fine-log-map-comparison} and \ref{intro-thm:fs-log-map-comparison} as well as Corollary~\ref{intro-cor:number-of-saturations}. 


\subsection{Acknowledgments}

We are grateful to Ben Church, Mark Gross, Eleny Ionel and Sushmita Venugopalan for useful discussions and correspondences. We thank John Pardon, Dhruv Ranganathan and Bernd Siebert for comments on an earlier draft. The first author was supported by NSF-RTG grant DMS-2038103. The second author was supported by a Murty Science Fellowship during the final stages of writing.

\subsection{Conventions}\label{subsec:conventions}

Unless explicitly stated otherwise, we work in the category $\Sch$ of schemes which are locally of finite type over $\bC$ and the category $\LogSch$ of fine log schemes whose underlying schemes lie in $\Sch$. We use the word \emph{point} to mean $\bC$-valued point. When we use $\dim$, $\Hom$ and $\otimes$ for vector spaces, these are understood to be over $\bC$. When we use $\Hom$ and $\otimes$ for quasicoherent sheaves on a scheme $S$, these are understood to be over $\clO_S$. We mostly use $\Gamma(-,\sF)$ to denote sections of a sheaf $\sF$ on $S$; when $\sF$ is a quasicoherent sheaf, we sometimes use the notation $H^0(-,\sF)$ instead.

Throughout, we abbreviate the word \emph{logarithmic} to \emph{log}. We use the word \emph{morphism} for morphisms of ordinary schemes (or stacks) and the phrase \emph{log morphism} for morphisms of log schemes (or log stacks). We denote log schemes (resp. their underlying schemes) by notations such as $S^\dagger$ (resp. $S$), while we denote log morphisms (resp. their underlying morphisms) by notations such as $f^\dagger$ (resp. $f$). For a morphism $f:T\to S$ of schemes, we use the same notation $f^\#$ to denote both the natural maps $f^{-1}\clO_S\to\clO_T$ and $\clO_S\to f_*\clO_T$. The word \emph{variety} means irreducible, reduced, separated scheme of finite type over $\bC$. We use the abbreviation \emph{nc} (resp. \emph{snc}) to mean normal crossings (resp. simple normal crossings).

For $S\in\Sch$, we write $\Sch{\!/S}$ for the category of \emph{$S$-schemes}, i.e., schemes in $\Sch$ equipped with a morphism to $S$. Morphisms in $\Sch{\!/S}$ are morphisms of schemes which also respect the projection to $S$ and are called \emph{$S$-morphisms}. For a quasicoherent $\clO_S$-algebra $\clR$, the associated $S$-scheme is denoted by $\relSpec_S\clR$. For a quasicoherent sheaf $\sF$ on $S$ and a closed subscheme $T\subset S$, we use $\sF|_T$ to denote the quasicoherent sheaf on $T$ defined by $\sF\otimes_{\clO_S}\clO_T$.

Schemes are endowed with their \emph{Zariski topology}, but we will frequently need to work in the (finer) \emph{\'etale topology} since there is no Zariski-local inverse function theorem (e.g., finite unbranched covers typically do not have Zariski-local sections). The differential-geometrically minded reader may safely replace schemes (resp. the \'etale topology) throughout by complex analytic spaces (resp. the ``usual", i.e., complex analytic topology).

An \emph{\'etale neighborhood} in a scheme $S$ is an \'etale morphism $U\to S$. An \emph{\'etale neighborhood of a point} $s\in S$ is an \'etale neighborhood $U\to S$ together with a point $u\in U$ mapping to $s$; we typically don't explicitly mention $u$ if this omission is unlikely to cause confusion. An \emph{\'etale covering} of $S$ is a collection $\{U_i\to S\}_i$ of \'etale neighborhoods whose images together cover the whole of $S$. Descent theory allows us to regard quasicoherent $\clO_S$-modules as sheaves on $S$ with respect to the \'etale topology.
    \section{Log geometry background}\label{sec:log-prelims}

We review the relevant background on logarithmic structures which will be used throughout the paper. Our purpose in doing so is to make the article more accessible to readers who may not be familiar with logarithmic geometry. After recalling the notions of \emph{monoids}, \emph{log schemes} and \emph{log stacks}, we discuss two important classes of examples: \emph{Artin cones} and log structures associated to normal crossings divisors.

\begin{Notation}
    Throughout this section, we use $S$, $S'$, $T$ etc to denote objects of the category $\Sch$ of schemes which are locally of finite type over $\bC$.
\end{Notation}

\subsection{Monoids and toric geometry}

The main reference for the topics recalled in this subsection is \cite[Chapter~I]{ogus-log-book}.

Recall that a \emph{monoid} is a set $Q$ with a commutative associative binary operation $+$ which has an identity element $0$. Maps of monoids are maps of sets which respect $+$ and $0$. By formally adjoining inverses to a monoid $Q$, we get an abelian group $Q^\text{gp}$ called its \emph{associated abelian group}, with a natural map $Q\to Q^\text{gp}$.

\begin{Definition}\label{def:monoid-props}
    Let $Q$ be a monoid. We say that $Q$ is
    \begin{enumerate}[label = (\arabic*), itemsep = 0.3ex]
        \item \emph{integral} if the map $Q\to Q^\text{gp}$ is injective;
        
        \item \emph{fine} if $Q$ is finitely generated and integral;
        
        \item \emph{saturated} if it is integral and for any $m\in\bZ_{\ge 1}$ and $q\in Q^\text{gp}$ such that $mq\in Q$, we also have $q\in Q$;
        
        \item \emph{fs} if $Q$ is fine and saturated;
        
        \item \emph{sharp} if $Q$ has no invertible elements other than $0$;
        
        \item \emph{toric} if $Q$ is fs and $Q^\text{gp}$ is torsion-free (and thus, free of finite rank).
    \end{enumerate}
    Let $Q^\times\subset Q$ be the \emph{group of units}, i.e., elements which are invertible in $Q$. The \emph{sharpening} of $Q$ is defined to be the quotient $\ol Q := Q/Q^\times$ of the monoid $Q$ by the additive action of the group $Q^\times$. The monoid $\ol Q$ is sharp by construction.
\end{Definition}

\begin{Remark}\label{rem:sharp-fs-implies-no-torsion}
    If $Q$ is a sharp fs monoid, then it is toric. Indeed, if $q\in Q^\text{gp}$ is an element of finite order, then we must have $q\in Q$ since the monoid $Q$ is saturated. It follows that $q = 0$ since the monoid $Q$ is sharp and the element $q$ is invertible in $Q$.
\end{Remark}

\begin{Notation}
    Let $\Mon$ be the category of monoids. We use $\Ab$, $\Mon^\text{int}$, $\Mon^\text{f}$, $\Mon^\text{sat}$ and $\Mon^\text{fs}$ to denote the full subcategories of $\Mon$ consisting of abelian groups, integral monoids, fine monoids, saturated monoids and fs monoids respectively. These fit into the following diagram.
    \begin{equation*}
    \begin{tikzcd}
        & \Mon^\text{fs} \arrow[r,hookrightarrow] \arrow[d,hookrightarrow] & \Mon^\text{f} \arrow[d,hookrightarrow] & \\
        \Ab \arrow[r,hookrightarrow] & \Mon^\text{sat} \arrow[r,hookrightarrow] & \Mon^\text{int} \arrow[r,hookrightarrow] & \Mon
    \end{tikzcd}
    \end{equation*}
\end{Notation}

\begin{Lemma}\label{lem:mon-adj} 
    We have the following adjunctions.
    \begin{enumerate}[ref = \theLemma(\arabic*), label = (\arabic*), itemsep = 0.3ex]
        \item The functor $Q\mapsto Q^\text{gp}$ is left adjoint to $\Ab\hookrightarrow\Mon$.
        
        \item The functor $\Mon^\text{int}\hookrightarrow\Mon$ has a left adjoint, denoted $Q\mapsto Q^\text{int}$.
        
        \item The functor $\Mon^\text{sat}\hookrightarrow\Mon^\text{int}$ has a left adjoint, denoted $Q\mapsto Q^\text{sat}$.
        
        \item\label{lem-part:saturation-of-fine-is-fine} The restriction of $Q\mapsto Q^\text{sat}$ is left adjoint to $\Mon^\text{fs}\hookrightarrow\Mon^\text{f}$.
    \end{enumerate}
\end{Lemma}
\begin{proof}
    (1): This is immediate from the definition.

    (2): For $Q\in\Mon$, take $Q^\text{int}\in\Mon^\text{int}$ to be the image of $Q\to Q^\text{gp}$.

    (3):  For $Q\in\Mon^\text{int}$, take $Q^\text{sat}=\{q\in Q^\text{gp}\;|\;\exists\,m\in\bZ_{\ge 1}\text{ with }mq\in Q\}\in\Mon^\text{sat}$.
    
    (4): This follows from part~(3) and \cite[I.2.2.5]{ogus-log-book}.
\end{proof}

\begin{Definition}\label{def:integralization-saturation}
    Using the notation of Lemma~\ref{lem:mon-adj}, given $Q\in\Mon$ (resp. $\Mon^\text{int}$), we call $Q^\text{int}\in\Mon^\text{int}$ (resp. $Q^\text{sat}\in\Mon^\text{sat}$) the \emph{integralization} (resp. \emph{saturation}) of $Q$.
\end{Definition}

\begin{Lemma}\label{lem:monoid-colimit}
    Arbitrary colimits exist in the categories $\Ab$, $\Mon$, $\Mon^\text{int}$ and $\Mon^\text{sat}$. The functors $Q\mapsto Q^\text{gp}$, $Q\mapsto Q^\text{int}$ and $Q\mapsto Q^\text{sat}$ commute with colimits.
\end{Lemma}
\begin{proof}
    The second assertion follows from the general categorical fact that left adjoints preserve colimits. The existence of arbitary colimits is standard for $\Ab$. For $\Mon$, it follows from the discussion in \cite[Chapter~I, pages~2--5]{ogus-log-book}. By the adjunctions in Lemma~\ref{lem:mon-adj}, a colimit in $\Mon^\text{int}$ (resp. $\Mon^\text{sat}$) is computed by computing the same colimit first in $\Mon$ and then applying the functor $Q\mapsto Q^\text{int}$ (resp. $Q\mapsto (Q^\text{int})^\text{sat}$) to the result.
\end{proof}

\begin{Discussion}\label{disc:monoid-coeq}
    For us, the most important case of Lemma~\ref{lem:monoid-colimit} is the \emph{coequalizer} of two maps $f,g:P\rrarrows Q$ in $\Mon$, $\Mon^\text{int}$ or $\Mon^\text{sat}$. We describe each case below.
    \begin{enumerate}[label = (\arabic*), itemsep = 0.3ex]
        \item Let $\equiv$ be the equivalence relation on $Q$ generated by declaring $q + f(p)\equiv q + g(p)$ for all $p\in P$, $q\in Q$. The coequalizer of $f,g$ in $\Mon$ is the set of equivalence classes in $Q$ under $\equiv$, with the monoid operation descended from $Q$.
        
        \item The coequalizer of $f,g$ in $\Mon^\text{int}$ is the image of the natural map from $Q$ to the cokernel of the abelian group homomorphism $f^\text{gp}-g^\text{gp}:P^\text{gp}\to Q^\text{gp}$.
        
        \item The coequalizer of $f,g$ in $\Mon^\text{sat}$ is the saturation of the coequalizer of $f,g$ in $\Mon^\text{int}$.
    \end{enumerate}
\end{Discussion}

\begin{Definition}
    Let $P$, $Q$ and $Q'$ be monoids.
    \begin{enumerate}[label = (\arabic*), itemsep = 0.3ex]
        \item Given maps $f:P\to Q$ and $f':P\to Q'$ in $\Mon$, their \emph{pushout} $Q\oplus_P Q'$ is defined to be the coequalizer of the maps $(f,0),(0,f'):P\rrarrows Q\oplus Q'$ in $\Mon$. The pushout is sometimes also called the \emph{amalgamated sum}.

        \item When $P\subset Q$ is an inclusion of monoids, the \emph{quotient} $Q/P$ is defined to be the pushout $Q\oplus_P0$, where the map $P\to Q$ is given by inclusion.
    \end{enumerate}
\end{Definition}

\begin{Lemma}
\label{lem:integral-pushout}
    If $f:P\to Q$ and $f':P\to Q'$ are maps in $\Mon^\text{int}$ and one of $P$, $Q$, $Q'$ is an abelian group, then $Q\oplus_P Q'\in\Mon^\text{int}$. The same is true with $\Mon^\text{int}$ replaced by $\Mon^\text{f}$.
\end{Lemma}
\begin{proof}
    For $\Mon^\text{int}$, this follows from \cite[I.1.3.4]{ogus-log-book}. For $\Mon^\text{f}$, this follows by additionally observing that the natural map $Q\oplus Q'\to Q\oplus_P Q'$ is surjective. 
\end{proof}

\begin{Definition}\label{def:monoid-gen-rel}
    Let $Q$ be a fine monoid.
    \begin{enumerate}[label = (\arabic*), ref = \theDefinition(\arabic*), itemsep = 0.3ex]
        \item Let $I$, $J$ be finite sets and let $\{{\bf e}_i\}$ and $\{{\bf e}'_j\}$ denote the standard bases of $\bN^I$ and $\bN^J$ respectively. Given a pair of monoid maps $a,b:\bN^J\rrarrows\bN^I$, we define
        \begin{align}\label{eqn:generators-relations-int}
            \left\langle {\bf e}_i\,\middle|\,\sum a_{ij}{\bf e}_i = \sum b_{ij}{\bf e}_i\right\rangle^\text{int}
        \end{align}
        to be their coequalizer in $\Mon^\text{int}$, where $(a_{ij})$, $(b_{ij})$ are the matrices of $a$, $b$ respectively with respect to the standard bases. We refer to \eqref{eqn:generators-relations-int} as the \emph{fine monoid generated by the symbols ${\bf e}_i$ modulo the relations $\sum a_{ij}{\bf e}_i = \sum b_{ij}{\bf e}_i$.}

        \item\label{def-part:integral-presentation} An \emph{integral presentation} of $Q$ consists of finite sets $I$, $J$ and a diagram
        \begin{equation}\label{eqn:integral-presentation}
        \begin{tikzcd}
            \bN^J \arrow[r,shift left=.5ex,"a"] \arrow[r,shift right=.5ex,swap,"b"] & \bN^I\arrow [r,"\pi"] & Q
        \end{tikzcd}
        \end{equation}
        in $\Mon^\text{int}$ which expresses $Q$ as the coequalizer of the maps $a,b:\bN^J\rrarrows\bN^I$ in $\Mon^\text{int}$, i.e., the map $\pi$ descends to an isomorphism from the monoid \eqref{eqn:generators-relations-int} to $Q$.
    \end{enumerate}
\end{Definition}

\begin{Definition}
    For $Q\in\Mon^\text{f}$, we define its \emph{dual} to be the monoid
    \begin{align*}
        Q^\vee := \Hom_{\Mon}(Q,\bN).
    \end{align*}
\end{Definition}

\begin{Lemma}\label{lem:dual-monoid}
    Let $Q$ be a sharp fine monoid. 
    \begin{enumerate}[ref = \theLemma(\arabic*), label = (\arabic*), itemsep = 0.3ex]
        \item\label{lem-part:sharp-fine-interior-hom} There exists $\beta\in Q^\vee$ such that $\beta^{-1}(0) = \{0\}\subset Q$.
        
        \item The monoid $Q^\vee$ is sharp and fs. Thus, by Remark~\ref{rem:sharp-fs-implies-no-torsion}, $Q^\vee$ is toric.
        
        \item\label{lem-part:saturation-units-torsion} The group of units in $Q^\text{sat}$ is given by $Q^\text{gp}_\text{tor}\subset Q^\text{gp}$.
        
        \item\label{lem-part:fine-sharp-to-fs-sharp} The natural map $Q\to Q^\vee\phantom{}^\vee$ induces an isomorphism $\ol{Q^\text{sat}}\xrightarrow{\sim}Q^\vee\phantom{}^\vee$, where $\ol{Q^\text{sat}}$ is the sharpening of the monoid $Q^\text{sat}$.
    \end{enumerate}    
\end{Lemma}
\begin{proof}
    Follows from I.2.2.1, I.2.2.3(1), I.1.3.5(5) and I.2.2.3(3) respectively in \cite{ogus-log-book}.   
\end{proof}

We conclude this subsection with a brief discussion of the relation between fine monoids and toric geometry.

\begin{Definition}\label{def:scheme-from-monoid}
    Given a fine monoid $Q$, we obtain the algebra $\bC[Q^\text{gp}]$ and its sub-algebra $\bC[Q]$. Using this, define the affine schemes
    \begin{align*}
        V_Q := \Spec\bC[Q]\quad\text{and}\quad T_Q:=\Spec\bC[Q^\text{gp}].    
    \end{align*}
    We denote the monomial in $\bC[Q^\text{gp}]$ corresponding to $q\in Q^\text{gp}$ by $\chi^q$. 
    Let $Q^\text{gp}_\text{tor}\subset Q^\text{gp}$ denote the torsion subgroup. Using this, define the affine schemes
    \begin{align*}
        T^0_Q := \Spec\bC[Q^\text{gp}/Q^\text{gp}_\text{tor}]\quad\text{and}\quad\pi_0(T_Q) := \Spec\bC[Q^\text{gp}_\text{tor}].
    \end{align*}
\end{Definition}

\begin{Lemma}\label{lem:scheme-from-monoid}
    For a fine monoid $Q$, we have the following.
    \begin{enumerate}[ref = \theLemma(\arabic*), label = (\arabic*), itemsep = 0.3ex]
        \item The schemes $V_Q$, $T_Q$, $T^0_Q$ and $\pi_0(T_Q)$ are of finite type over $\bC$ and the natural morphism $T_Q\to V_Q$ is the inclusion of a dense Zariski open subset.
        
        \item The algebra maps given on monomials by $\chi^q\mapsto\chi^q\otimes\chi^q$ induce the structure of an algebraic group on $T_Q$, $T^0_Q$ and $\pi_0(T_Q)$ as well as an action of $T_Q$ on $V_Q$ extending the obvious left multiplication action of $T_Q$ on itself.
        
        \item\label{lem-part:alg-group-from-monoid} The subgroup $T^0_Q\subset T_Q$ is the identity component and is an algebraic torus, while the quotient $T_Q\twoheadrightarrow\pi_0(T_Q)$ is the group of connected components of $T_Q$.
        
        \item The scheme $V_Q$ is pure-dimensional, with dimension given by the rank of $Q^\text{gp}$, and has exactly $|Q^\text{gp}_\text{tor}|$ irreducible components.
        
        \item The morphism $V_{Q^\text{sat}}\to V_Q$ induced by $Q\subset Q^\text{sat}$ is the normalization of $V_Q$.
        
        \item Suppose $Q$ is also sharp. Then the action of $T_Q$ on $V_Q$ has a fixed point $0\in V_Q$, defined by the maximal ideal of $\bC[Q]$ generated by the monomials $\chi^q$ for all $0\ne q\in Q$. Moreover, all irreducible components of $V_Q$ pass through $0$.
    \end{enumerate}
\end{Lemma}
\begin{proof}
    (1): Since $Q$ is finitely generated as a monoid, it follows that $\bC[Q]$, $\bC[Q^\text{gp}]$, $\bC[Q^\text{gp}_\text{tor}]$ and $\bC[Q^\text{gp}/Q^\text{gp}_\text{tor}]$ are finitely generated as algebras over $\bC$. Specifically, if $q_1,\ldots,q_m\in Q$ are generators of $Q$, then $\bC[Q]$ is generated as a $\bC$-algebra by $\chi^{q_1},\ldots,\chi^{q_m}$ and $\bC[Q^\text{gp}]$ is obtained from $\bC[Q]$ by inverting the elements $\chi^{q_1},\ldots,\chi^{q_m}$. This shows that $T_Q\subset V_Q$ is Zariski open. Since $Q$ is an integral monoid, the natural map $\bC[Q]\to\bC[Q^\text{gp}]$ is injective and we deduce that the inclusion $T_Q\hookrightarrow V_Q$ is dominant, i.e., $T_Q$ is dense in $V_Q$.
    
    (2): The given algebra map induces morphisms 
    \begin{align*}
        T_Q\times T_Q\to T_Q,\quad T^0_Q\times T^0_Q\to T^0_Q,\quad \pi_0(T_Q)\times\pi_0(T_Q)\to\pi_0(T_Q),\quad\text{and}\quad V_Q\times T_Q\to V_Q.    
    \end{align*}
    It is straightforward to verify that the first three are group structures, the last one is a group action and that the inclusion $T_Q\hookrightarrow V_Q$ is a $T_Q$-equivariant morphism.
       
    (3): The result is immediate once we use the structure theorem for finitely generated abelian groups to obtain a direct sum decomposition 
    \begin{align}\label{eqn:str-thm-ab-gp-iso}
        Q^\text{gp}\xrightarrow{\sim} \textstyle\bigoplus_i\bZ\oplus\bigoplus_j(\bZ/d_j\bZ).    
    \end{align}
        
    (4): Since $T_Q\subset V_Q$ is dense, the connected components of the algebraic group $T_Q$ are in bijection with the irreducible components of $V_Q$. Using \eqref{eqn:str-thm-ab-gp-iso}, we now see that $V_Q$ has exactly $|Q^\text{gp}_\text{tor}|$ irreducible components, each of which has dimension equal to the rank of $Q^\text{gp}$.
        
    (5): This is a special case of \cite[I.3.4.1(2)]{ogus-log-book}.
        
    (6): Since $Q$ is sharp, we have a well-defined $\bC$-algebra map $\bC[Q]\to\bC$ given by $\chi^q\mapsto 0$ for all $0\ne q\in Q$. The kernel of this map is the maximal ideal defining $0\in V_Q$. It is clear from the definition that $0$ is fixed by the action of $T_Q$ on $V_Q$. Since $T_Q$ acts transitively on the irreducible components of $V_Q$, it follows that all of them contain $0$.
\end{proof}

\begin{Example}
    Let $a_1,a_2\ge 1$ be integers and let $d = \gcd(a_1,a_2)$. Define $Q\in\Mon^\text{f}$ by $Q = \langle {\bf e}_1,{\bf e}_2\,|\,a_1{\bf e}_1 = a_2{\bf e}_2\rangle^\text{int}$. The presentation of $Q$ yields an isomorphism
    \begin{align*}
        V_Q \xrightarrow{\sim}\Spec\left(\frac{\bC[z_1,z_2]}{(z_1^{a_1}-z_2^{a_2})}\right).
    \end{align*}
    The factorization $z_1^{a_1}-z_2^{a_2} = \prod_{\zeta}(z_1^{a_1/d}-\zeta z_2^{a_2/d})$, with $\zeta\in\bC$ ranging over the $d^\text{th}$ roots of unity, shows that $V_Q$ has exactly $d$ irreducible components, all meeting at $z_1 = z_2 = 0$. 
\end{Example}

\begin{Remark}
    When $Q$ is a toric monoid, the scheme $V_Q$ is a normal affine toric variety with dense algebraic torus $T_Q$.
\end{Remark}

\subsection{Log schemes}\label{subsec:log-schemes}

The main references for the topics recalled in this subsection are the papers \cite{KKato-log}, \cite{olsson-stack-of-log-str}, \cite{log-geometry-and-moduli} and the book \cite{ogus-log-book}.

\begin{Notation}\label{nota:monoid-and-scheme}
    Throughout this subsection, we use $Q$, $Q'$ etc to denote fine monoids. We denote the constant sheaf of monoids on a scheme $S$ with stalk $Q$ by the notation $\ul Q_S$ or simply $\ul Q$ when $S$ is clear from the context. 
\end{Notation}

\begin{Convention}
    All sheaves of monoids will be considered with respect to the \'etale topology; see \cite[\href{https://stacks.math.columbia.edu/tag/03N4}{Tag~03N4}]{stacks-project}. The structure sheaf $\clO_S$ of a scheme $S$ is regarded as a sheaf of monoids under multiplication; similarly for the subsheaf $\clO_S^\times\subset\clO_S$ of multiplicative units. Descent theory shows that $\clO_S$ and $\clO_S^\times$ are sheaves with respect to the \'etale topology.
\end{Convention}

\begin{Remark}\label{rem:pushout-monoid-sheaves}
    To compute quotients and pushouts in the category of sheaves of monoids, we first compute these at the level of presheaves of monoids (i.e., over each \'etale neighborhood, we use Lemma~\ref{lem:monoid-colimit} for $\Mon$) and then sheafify the resulting presheaf of monoids. 
\end{Remark}

\begin{Definition}\label{def:log-str-basics}
    Let $S$ be a scheme.
    \begin{enumerate}[label = (\arabic*), ref = \theDefinition(\arabic*), itemsep = 0.3ex]
        \item A \emph{pre-log structure} on $S$ is a sheaf of monoids $M$ with respect to the \'etale topology on $S$ and a map $\alpha:M\to\clO_S$ of sheaves of monoids. 
        
        A \emph{map of pre-log structures} between $\alpha:M\to\clO_S$ and $\alpha':M'\to\clO_S$ is a map $\varphi:M\to M'$ of sheaves of monoids such that the following diagram commutes.
        \begin{equation*}
        \begin{tikzcd}
            M \arrow[rr,"\varphi"] \arrow[dr,swap,"\alpha"] && M' \arrow[dl,"\alpha'"] \\
            & \clO_S &
        \end{tikzcd}
        \end{equation*}
        
        The \emph{category of pre-log structures} on $S$, denoted $\PreLog(S)$, has objects given by pre-log structures on $S$ and morphisms given by maps of pre-log structures. 
        
        \item\label{def-part:ghost-sheaf} The \emph{ghost sheaf} (or \emph{characteristic sheaf}) of a pre-log structure $\alpha:M\to\clO_S$ is the quotient $\ol M = M/(\alpha^{-1}(\clO_S^\times))$ in the category of sheaves of monoids on $S$. For each point $s\in S$, the stalk $\ol M_s$ of the ghost sheaf is a sharp monoid.

        A map $\varphi:M\to M'$ in $\PreLog(S)$ induces a map $\ol\varphi:\ol M\to\ol M\phantom{}'$ on ghost sheaves. For each point $s\in S$, the map $\ol\varphi_s:\ol M_s\to\ol M\phantom{}'_s$ on stalks satisfies $(\ol\varphi_s)^{-1}(0) = \{0\}$.

        \item A pre-log structure $\alpha:M\to\clO_S$ is called a \emph{log structure} if the induced map $\alpha^{-1}(\clO_S^\times)\to\clO_S^\times$ is an isomorphism. For a log structure $\alpha:M\to\clO_S$, we always regard $\clO_S^\times$ as a subsheaf of $M$ using $\alpha^{-1}$. 
        
        A \emph{map of log structures} is a map of the corresponding pre-log structures. The \emph{category of log structures} on $S$, denoted $\Log(S)$, has objects given by log structures on $S$ and morphisms given by maps of log structures.

        \item A \emph{log scheme} $S^\dagger$ consists of a scheme $S$ with a log structure $\alpha_S:M_S\to\clO_S$. 
        
        We usually omit $\alpha_S$ from the notation and denote a log scheme by $S^\dagger = (S,M_S)$ and its ghost sheaf by $\ol M_S$; we drop the subscripts when $S$ is clear from the context. We refer to $\alpha_S$ as the \emph{structure map} and $S$ as the \emph{underlying scheme} of $S^\dagger$.
    \end{enumerate} 
\end{Definition}

\begin{Convention}\label{conv:log-str-mult-ghost-add}
    We always write the operation on a pre-log structure $M$ (resp. its ghost sheaf $\ol M$) multiplicatively (resp. additively).
\end{Convention}

\begin{Definition}\label{def:assoc-log-str}
    Given a pre-log structure $\alpha:M\to\clO_S$ on the scheme $S$, the \emph{associated log structure} $\alpha^{\log}:M^{\log}\to\clO_S$ is defined via the following commutative diagram, where the square is a pushout diagram in the category of sheaves of monoids on $S$.
    \begin{equation}\label{eqn:assoc-log-str}
    \begin{tikzcd}
        \alpha^{-1}(\clO_S^\times) \arrow[d] \arrow[r,hookrightarrow] \arrow[dr,phantom,"\ulcorner",pos=0.99] & M \arrow[ddr,"\alpha",bend left] \arrow[d] & \\
         \clO_S^\times \arrow[drr,hookrightarrow,bend right] \arrow[r] & M^{\log} \arrow[dr,"\alpha^{\log}"] & \\
        & & \clO_S
    \end{tikzcd}
    \end{equation}
    By construction, the induced map $\ol M\to\ol M\phantom{}^\text{log}$ on ghost sheaves is an isomorphism.
    
    For any pre-log structure $M$ and log structure $M'$, composition with the natural map $M\to M^{\log}$ induces an isomorphism
    \begin{align}\label{eqn:prelog-log-adjunction}
        \Hom_{\Log(S)}(M^{\log},M')\xrightarrow{\sim}\Hom_{\PreLog(S)}(M,M').
    \end{align}
    In other words, $M\mapsto M^{\log}$ is the left adjoint to the inclusion $\Log(S)\subset\PreLog(S)$.
\end{Definition}

\begin{Lemma}\label{lem:pre-log-diff-unit}
    Let $\alpha:M\to\clO_S$ be a pre-log structure on the scheme $S$ and let $u:M\to\clO_S^\times$ be a map of sheaves of monoids; this yields a new pre-log structure $\alpha u:M\to\clO_S$. 
    \begin{enumerate}[label = (\arabic*), ref = \theLemma(\arabic*), itemsep = 0.3ex]
        \item\label{lem-part:same-ghost} The ghost sheaves of $\alpha$ and $\alpha u$ are both the same quotient of $M$.
        
        \item\label{lem-part:unit-induced-map} Let $\alpha^{\log}:M^{\log}_\alpha\to\clO_S$ and $(\alpha u)^{\log}:M^{\log}_{\alpha u}\to\clO_S$ be the log structures associated to $\alpha$ and $\alpha u$ respectively. Then the automorphism of $M\oplus\clO_S^\times$ given by the formula $(m,v)\mapsto (m,\textstyle\frac1{u(m)}v)$ descends to an isomorphism of log structures 
        \begin{align}\label{eqn:pre-log-diff-unit}
            \Psi_u: M_\alpha^{\log}\xrightarrow{\sim} M_{\alpha u}^{\log}
        \end{align}
        inducing the identity map on ghost sheaves.

        \item Given a log structure $\alpha':M'\to\clO_S$ and a map of pre-log structures $\varphi:M\to M'$, let $\varphi^{\log}$ and $(\varphi u)^{\log}$ be the maps corresponding to the maps $\varphi$ and $\varphi u$ respectively under the adjunction \eqref{eqn:prelog-log-adjunction}. Then the following diagram commutes.
        \begin{equation}\label{eqn:pre-log-diff-unit-map-compatible}
        \begin{tikzcd}
            M_\alpha^{\log} \arrow[dr,swap,"\varphi^{\log}"] \arrow[rr,"\Psi_u","\sim"'] && M_{\alpha u}^{\log} \arrow[dl,"(\varphi u)^{\log}"] \\
            & M' &
        \end{tikzcd}
        \end{equation}

        \item\label{lem-part:pre-log-diff-unit-cocyle} Given another map $u':M\to\clO_S^\times$ of sheaves of monoids, taking $M' = M_{\alpha uu'}^{\log}$ in \eqref{eqn:pre-log-diff-unit-map-compatible} yields the identity $\Psi_{u'}\circ\Psi_u = \Psi_{uu'}$ for maps $M_\alpha^{\log}\to M_{\alpha u}^{\log}\to M_{\alpha uu'}^{\log}$.
    \end{enumerate}
\end{Lemma}
\begin{proof}
    Since we have $(\alpha u)^{-1}(\clO_S^\times) = \alpha^{-1}(\clO_S^\times)\subset M$, the ghost sheaves of $\alpha$ and $u\alpha$ are both the same quotient of $M$. The rest is immediate from definitions.
\end{proof}

\begin{Definition}
    Let $\alpha:M\to\clO_S$ be a pre-log structure on the scheme $S$ and let $f:T\to S$ be a morphism of schemes.
    \begin{enumerate}[label = (\arabic*), itemsep = 0.3ex]
        \item The \emph{pullback pre-log structure} is the composition
        \begin{align}\label{eqn:pullback-pre-log}
            f^{-1}M\xrightarrow{f^{-1}\alpha} f^{-1}\clO_S\xrightarrow{f^\#}\clO_T,    
        \end{align}
        where $f^\#$ is the map on structure sheaves induced by $f$.

        \item When $\alpha:M\to\clO_S$ is a log structure, then the \emph{pullback log structure} is the log structure associated to the pre-log structure \eqref{eqn:pullback-pre-log} and is denoted by
        \begin{align*}
            f^*\alpha:f^*M\to\clO_T.    
        \end{align*}
        The ghost sheaf $\ol{f^*M}$ of the pullback log structure is the same as the pullback sheaf $f^{-1}\ol M$, which we will instead denote as $f^*\ol M$.
    \end{enumerate}
\end{Definition}

\begin{Definition}
    Let $S^\dagger = (S,M_S)$ and $T^\dagger = (T,M_T)$ be log schemes. 
    \begin{enumerate}[label = (\arabic*), itemsep = 0.3ex]
        \item A \emph{log morphism} $f^\dagger = (f,f^\flat):T^\dagger\to S^\dagger$ consists of a morphism of schemes $f:T\to S$ and a map of log structures $f^\flat:f^*M_S\to M_T$. Alternatively, a log morphism is a morphism $f:T\to S$ together with a commutative diagram
        \begin{equation}\label{eqn:log-mor-diagram}
        \begin{tikzcd}
            f^{-1}M_S \arrow[d,swap,"f^{-1}\alpha_S"] \arrow[r,"f^\flat_\text{pre}"] & M_T \arrow[d,swap,"\alpha_T"] \\
        f^{-1}\clO_S \arrow[r,"f^\#"] & \clO_T
        \end{tikzcd}
        \end{equation}
        where $f^\flat_\text{pre}$ is a map of sheaves of monoids and $f^\#$ is the map on structure sheaves induced by $f$. We call $f$ the \emph{underlying morphism} of $f^\dagger$. We also say that $f^\dagger$ is a \emph{lift} of $f$ to a log morphism.

        \item A log morphism $f^\dagger:T^\dagger\to S^\dagger$ is said to be \emph{strict} if the map $f^\flat:f^*M_S\to M_T$ of log structures is an isomorphism.
    \end{enumerate}
\end{Definition}

\begin{Example}\label{exa:div-log-str}
    Let $X$ be a pure-dimensional reduced scheme and let $D\subset X$ be a Zariski closed subset of pure codimension $1$. The \emph{divisorial log structure} $M_{(X,D)}\subset\clO_X$ is the subsheaf of monoids defined as follows: its sections over an \'etale neighborhood $U$ are the regular functions on $U$ which are nowhere vanishing on the inverse image of $X\setminus D$ in $U$.

    Let $Y$ be another pure-dimensional reduced scheme with a Zariski closed subset $E\subset Y$ of pure codimension $1$. Using the divisorial log structures associated to $D\subset X$ and $E\subset Y$, we get the log schemes $X^\dagger = (X,M_{(X,D)})$ and $Y^\dagger = (Y,M_{(Y,E)})$ respectively. Using \eqref{eqn:log-mor-diagram}, a morphism $f:Y\to X$ can be lifted to a log morphism $f^\dagger:Y^\dagger\to X^\dagger$ if and only if $f(Y\setminus E)\subset X\setminus D$. Moreover, in this situation, $f^\dagger$ is uniquely determined by $f$.
\end{Example}

Other than Example~\ref{exa:div-log-str}, the following construction provides the main source of log structures and log morphisms that arise in practice.

\begin{Construction}\label{cons:log-spec}
    Let $(Q,R)$ be a pair consisting of a 
    fine monoid $Q$ and a $\bC$-algebra $R$ of finite type. Let $Q\to R$ be a monoid map, where we take the multiplicative monoid structure on $R$. Define the log scheme
    \begin{align*}
        \Spec(Q\to R)    
    \end{align*}
    to be the affine scheme $\Spec R$, equipped with the log structure associated to the pre-log structure $\ul Q\to\clO_{\Spec R}$ induced by the map $Q\to R$.
    
    Given another such pair $(Q',R')$ and a map $Q'\to R'$, consider a monoid map $Q\to Q'$ and a $\bC$-algebra map $R\to R'$ such that the diagram
    \begin{equation*}
    \begin{tikzcd}
        Q \arrow[d] \arrow[r] & R \arrow[d] \\
        Q' \arrow[r] & R'
    \end{tikzcd}
    \end{equation*}
    commutes. Then we get an induced log morphism
    \begin{align*}
        \Spec(Q'\to R')\to\Spec(Q\to R).
    \end{align*}
\end{Construction}

\begin{Definition}\label{def:std-log-str-toric}
    For a fine monoid $Q$, the \emph{standard log structure} $M_{V_Q}$ on the scheme $V_Q = \Spec\bC[Q]$ from Definition~\ref{def:scheme-from-monoid} is defined as follows. We let
    \begin{align*}
        V_Q^\dagger = (V_Q,M_{V_Q}) := \Spec(Q\hookrightarrow\bC[Q]).
    \end{align*}    
    By construction, we have a natural map 
    \begin{align}\label{eqn:nat-map-Q-to-M-toric}
        Q\to\Gamma(V_Q,M_{V_Q}).
    \end{align}
\end{Definition}

\begin{Remark}
    When the monoid $Q$ is toric, it can be shown that $M_{V_Q}$ agrees with the divisorial log structure (Example~\ref{exa:div-log-str}) on the affine toric variety $V_Q$ determined by the toric boundary divisor $V_Q\setminus T_Q$; see \cite[11.6]{KKato-toric-singularities} or \cite[III.1.9.5]{ogus-log-book}.
\end{Remark}

\begin{Proposition}
\label{prop:log-map-to-toric}
    Fix a fine monoid $Q$. Then there is a natural bijection between the following two types of objects for any log scheme $S^\dagger = (S,M_S)$.
    \begin{enumerate}[itemsep = 0.3ex]
        \myitem[(TV$\phantom{}_1$)]\label{log-mor-to-toric} Log morphisms $S^\dagger\to V_Q^\dagger$.
        
        \myitem[(TV$\phantom{}_2$)]\label{monoid-map-to-log-str} Monoid maps $Q\to\Gamma(S,M_S)$. 
    \end{enumerate}
\end{Proposition}
\begin{proof}
    We define the maps \ref{log-mor-to-toric}~$\to$~\ref{monoid-map-to-log-str}~$\to$~\ref{log-mor-to-toric}, but omit the verification that their composites \ref{log-mor-to-toric}~$\to$~\ref{log-mor-to-toric} and \ref{monoid-map-to-log-str}~$\to$~\ref{monoid-map-to-log-str} are identity maps.

    \ref{log-mor-to-toric}~$\to$~\ref{monoid-map-to-log-str}: Given $f^\dagger=(f,f^\flat):S^\dagger\to V_Q^\dagger$, associate to it the composition 
    \begin{align*}
        Q\xrightarrow{\eqref{eqn:nat-map-Q-to-M-toric}}\Gamma(V_Q,M_{V_Q})\xrightarrow{f^*}\Gamma(S,f^*M_{V_Q})\xrightarrow{\Gamma(f^\flat)}\Gamma(S,M_S).    
    \end{align*}

    \ref{monoid-map-to-log-str}~$\to$~\ref{log-mor-to-toric}: Given $\varphi:Q\to\Gamma(S,M_S)$, let $\ul Q\to M_S$ be the induced map of sheaves of monoids. Then the composition $\ul Q\to M_S\to\clO_S$ defines an algebra map $\bC[Q]\to\Gamma(S,\clO_S)$, i.e., a morphism $f:S\to V_Q$. Moreover, $f^*M_{V_Q}$ is the log structure associated to the pre-log structure $\ul Q\to\clO_S$. Thus, $\ul Q\to M_S$ induces a map $f^\flat:f^*M_{V_Q}\to M_S$ of log structures via the adjunction \eqref{eqn:prelog-log-adjunction} and this lifts $f$ to a log morphism.
\end{proof}

\begin{Definition}\label{def:chart-fine-fs}
    Let $S^\dagger = (S,M_S)$ be a log scheme.
    \begin{enumerate}[label = (\arabic*), ref = \theDefinition(\arabic*), itemsep = 0.3ex]
        \item A monoid map $Q\to\Gamma(S,M_S)$, where $Q$ is a fine monoid, 
        is said to be a \emph{log chart} if the associated log morphism $S^\dagger\to V_Q^\dagger$ is strict.
        
        \item\label{def-part:chart-fine} The log scheme $S^\dagger$ is said to be \emph{fine} if every point of $S$ has an \'etale neighborhood over which the log structure admits a log chart.
        
        \item The log scheme $S^\dagger$ is said to be \emph{fine and saturated} (\emph{fs}) if it is fine and the monoids appearing in its \'etale local log charts can all be taken to be fs.
    \end{enumerate}
\end{Definition}

\begin{Remark}
    The log schemes $\Spec(Q\to R)$ from Construction~\ref{cons:log-spec} for a fine (resp. fs) monoid $Q$ are the \'etale local models of fine (resp. fs) log schemes.
\end{Remark}

\begin{Lemma}\label{lem:fine-log-str-integral}
    If $M\to\clO_S$ is a fine log structure, then $M\to M^\text{gp}$ and $\ol M\to \ol M\phantom{}^\text{gp}$ are injective while $M\to\ol M$ and $M^\text{gp}\to\ol M\phantom{}^\text{gp}$ are $\clO_S^\times$-torsors. In particular, the following is a pullback diagram in the category of sheaves of monoids on $S$.
    \begin{equation}\label{eqn:ghost-groupification-cartesian}
    \begin{tikzcd}
        M \arrow[r] \arrow[d] \arrow[dr,phantom,"\lrcorner",pos=0.001] & M^\text{gp} \arrow[d] \\
        \ol M \arrow[r] & \ol M\phantom{}^\text{gp}
    \end{tikzcd}
    \end{equation}
    Moreover, if $M'\to\clO_S$ is another fine log structure and $\varphi: M\to M'$ is a map of log structures on $S$, then the following are pullback diagrams in the category of sheaves of monoids on $S$.
    \begin{equation}\label{eqn:map-of-log-str-cartesian}
    \begin{tikzcd}
        M \arrow[r,"\varphi"] \arrow[d] \arrow[dr,phantom,"\lrcorner",pos=0.001] & M' \arrow[d] & & M^\text{gp} \arrow[r,"\varphi^\text{gp}"] \arrow[d] \arrow[dr,phantom,"\lrcorner",pos=0.001] & M\phantom{}'^{\text{gp}} \arrow[d] \\
        \ol M \arrow[r,"\ol\varphi"] & \ol M\phantom{}' & & \ol M\phantom{}^{\text{gp}} \arrow[r,"\ol\varphi^\text{gp}"] & \ol M\phantom{}'^{\text{gp}}
    \end{tikzcd}
    \end{equation}
\end{Lemma}
\begin{proof}
    Since the assertions are \'etale local on $S$, we may assume that $M$ is the log structure associated to a pre-log structure given by a constant sheaf $\ul Q$ whose stalk $Q$ is a fine monoid. The injectivity assertion now follows from Lemma~\ref{lem:integral-pushout} and the pushout diagram \eqref{eqn:assoc-log-str}. Since we have $\clO_S^\times\subset M\subset M^\text{gp}$, it follows that $\clO_S^\times$ acts freely on $M$ and $M^\text{gp}$. Thus, $M$ and $M^\text{gp}$ are $\clO_S^\times$-torsors over $\ol M$ and $\ol M\phantom{}^\text{gp}$ respectively. Finally, \eqref{eqn:ghost-groupification-cartesian} and \eqref{eqn:map-of-log-str-cartesian} are pullback diagrams since the vertical arrows are all $\clO_S^\times$-torsors.
\end{proof}

\begin{Lemma}[{\cite[2.3]{olsson-stack-of-log-str}}]\label{lem:local-lift-from-ghost}
    Let $M\to\clO_S$ be a fine log structure and let $Q$ be a fine monoid. Then any monoid map $Q\to\Gamma(S,\ol M)$ lifts to a monoid map $Q\to\Gamma(S,M)$ in some \'etale neighborhood of any given point of $S$.
\end{Lemma}
\begin{proof}
    (Our argument is essentially the same as the one in \cite{olsson-stack-of-log-str} and we include it for completeness.) From \eqref{eqn:ghost-groupification-cartesian}, we have $M = M^\text{gp}\times_{\ol M\phantom{}^\text{gp}}\ol M$. Thus, it suffices to prove the result after replacing $Q$, $M$ and $\ol M$ by their associated (sheaves of) abelian groups. Decompose the group $Q^\text{gp}$ as in \eqref{eqn:str-thm-ab-gp-iso}. Let $s\in S$ be a given point.
        
    For each free $\bZ$ summand of $Q^\text{gp}$ in \eqref{eqn:str-thm-ab-gp-iso}, lift the image of $1\in\bZ$ in $\ol M\phantom{}^\text{gp}$ to a section of $M^\text{gp}$ in an \'etale neighborhood of $s$.
    For each finite cyclic $\bZ/d\bZ$ summand of $Q^\text{gp}$ in \eqref{eqn:str-thm-ab-gp-iso}, lift the image of $1\in\bZ/d\bZ$ in $\ol M\phantom{}^\text{gp}$ to a section $\sigma$ of $M^\text{gp}$ in an \'etale neighborhood of $s$. Thus, we have $\sigma^d = u$ for a section $u$ of $\clO_S^\times$; recall that we regard $\clO_S^\times$ as a subsheaf of $M$. Passing to a further \'etale neighborhood of $s$ to get a section $v$ of $\clO_S^\times$ with $v^d = u$, replace the lift $\sigma$ by $v^{-1}\sigma$ to ensure that we now have $\sigma^d = 1$.
    Combining the lifts for each summand of $Q^\text{gp}$, we obtain the desired lift $Q^\text{gp}\to\Gamma(S,M^\text{gp})$ in an \'etale neighborhood of $s$.
\end{proof}

\begin{Notation}\label{nota:logsch-variants}
    Recall from Subsection~\ref{subsec:conventions} that $\Sch$ denotes the category of schemes which are locally of finite type over $\bC$, while $\Sch{\!/S}$ denotes the category of schemes in $\Sch$ lying over some fixed base scheme $S\in\Sch$.
    
    Let $\LogSch$ (resp. $\LogSch^\text{fs}$) denote the category of fine (resp. fs) log schemes, whose underlying schemes lie in $\Sch$. Let $\Log$ (resp. $\Log^\text{fs}$) denote the subcategory of $\LogSch$ (resp. $\LogSch^\text{fs}$) consisting of strict log morphisms. These fit into the following diagram.
    \begin{equation*}
    \begin{tikzcd}
        \Log^\text{fs} \arrow[d,hookrightarrow] \arrow[r,hookrightarrow] & \LogSch^\text{fs} \arrow[d,hookrightarrow] & \\
        \Log \arrow[r,hookrightarrow] & \LogSch \arrow[r] & \Sch
    \end{tikzcd}
    \end{equation*}
\end{Notation}

\begin{Convention}\label{conv:sch-logsch}
    For the rest of this article, unless otherwise specified, all schemes will be from $\Sch$ and all log schemes will be from $\LogSch$. For $S\in\Sch$, we will henceforth use $\Log(S)$ to denote the category of fine log structures on $S$.
\end{Convention}

\subsection{Descent for log structures}

\begin{Definition}\label{def:descent-datum-log-scheme}
    Let $p:U\to S$ be a surjective smooth morphism in $\Sch$. 
    
    For $i\in\{1,2\}$ and $jk\in\{12,13,23\}$, we have the natural projections
    \begin{align*}
        U\times_S U\times_S U \xrightarrow{p_{jk}} U\times_S U \xrightarrow{p_i} U.
    \end{align*}
    \begin{enumerate}[label = (\arabic*), itemsep = 0.3ex]
        \item A \emph{descent datum} $(M_U,\Psi)$ for a log structure with respect to $p:U\to S$ consists of 
        \begin{enumerate}[label = (\roman*), itemsep = 0.15ex]
            \item a log structure $M_U$ on $U$, and
        
            \item an isomorphism of log structures $\Psi:p_1^*M_U\xrightarrow{\sim} p_2^*M_U$ on $U\times_S U$
        \end{enumerate}
        which satisfy the \emph{cocyle condition} on $U\times_S U\times_S U$, i.e., we have
        \begin{align}\label{eqn:descent-cocycle-condition}
            p_{13}^*\Psi = p_{23}^*\Psi\circ p_{12}^*\Psi.    
        \end{align}

        \item For descent data $(M_U,\Psi)$ and $(M_U',\Psi')$ as above, a \emph{map of descent data} between these is a map $\varphi_U:M_U\to M_U'$ of log structures such that the following diagram commutes.
        \begin{equation}\label{eqn:descent-log-str-map}
        \begin{tikzcd}
            p_1^*M_U \arrow[d,"p_1^*\varphi_U"] \arrow[r,"\Psi","\sim"'] & p_2^*M_U \arrow[d,"p_2^*\varphi_U"] \\
            p_1^*M_U' \arrow[r,"\Psi'","\sim"'] & p_2^*M_U'
        \end{tikzcd}
        \end{equation}

        \item Descent data for log structures with respect to $p:U\to S$ and maps of descent data form a category, which we denote by $\Log(p:U\to S)$.
    \end{enumerate}
\end{Definition}

\begin{Proposition}[{\cite[Appendix~A]{olsson-stack-of-log-str}}]\label{prop:log-str-descent-scheme}
    In the setting of Definition~\ref{def:descent-datum-log-scheme}, the functor 
    \begin{align}\label{eqn:descent-equivalence-scheme}
        \Log(S) &\to \Log(p:U\to S) \\
        M_S &\mapsto (p^*M_S,\text{id})
    \end{align}
    is an equivalence of categories. Its inverse is given by mapping $(\alpha_U:M_U\to\clO_U,\Psi)$ to the log structure on $S$ given by the row wise equalizer of the diagram
    \begin{equation*}
    \begin{tikzcd}
        p_*M_U \arrow[rrr,shift left=.5ex,"\hat p_*\Psi\,\circ\,p_*\varepsilon_1"] \arrow[rrr,shift right=.5ex,swap,"p_*\varepsilon_2"] \arrow[d,"p_*\alpha_U"] &&& \hat p_*p_2^*M_U \arrow[d,"\hat p_*p_2^*\alpha_U"] \\
        p_*\clO_U \arrow[rrr,shift left=.5ex,"p_*p_1^\#"] \arrow[rrr,shift right=.5ex,swap,"p_*p_2^\#"] &&& \hat p_*\clO_{U\times_S U}    
    \end{tikzcd}
    \end{equation*}
    where $\hat p:U\times_S U\to S$ is the natural projection given by $\hat p = p\circ p_1 = p\circ p_2$ and, for $i=1,2$, the maps $\varepsilon_i:M_U\to {p_i}_*p_i^*M_U$ are induced by the pushforward-pullback adjunction.
\end{Proposition}
\begin{proof}
    It follows from \cite[A.5]{olsson-stack-of-log-str} that \eqref{eqn:descent-equivalence-scheme} is an equivalence. The explicit formula for the inverse of \eqref{eqn:descent-equivalence-scheme} is precisely the content of \cite[A.6]{olsson-stack-of-log-str}.
\end{proof}

\begin{Remark}
    In fact, \cite[Appendix~A]{olsson-stack-of-log-str} proves Proposition~\ref{prop:log-str-descent-scheme} in the much more general setting where $p:U\to S$ is an \emph{fppf covering}, i.e., a morphism which is faithfully flat and locally of finite presentation.
\end{Remark}

\subsection{Log stacks}

Recall that a \emph{stack} $\fX$ is a \emph{category fibred in groupoids} over $\Sch$, where we denote the fibration as $\fp:\fX\to\Sch$, 
satisfying the $2$-categorical analogue of the sheaf property for the fppf topology; see \cite[\href{https://stacks.math.columbia.edu/tag/0304}{Tag~0304}]{stacks-project}.

Given $S\in\Sch$, we regard it as a stack by identifying it with $\Sch{\!/S}\to\Sch$. Using the equivalence provided by the 2-Yoneda lemma \cite[\href{https://stacks.math.columbia.edu/tag/04SS}{Tag~04SS}]{stacks-project}, we always identify the groupoid of morphisms $S\to\fX$ with the groupoid $\fX(S)$ given by the fibre of $\fp$ over $S$.

The stack $\fX$ is said to be \emph{algebraic} if the diagonal morphism $\fX\to\fX\times_{\Sch}\fX$ is \emph{representable by algebraic spaces} and $\fX$ has a \emph{smooth  atlas}, i.e., a surjective smooth morphism $U\to\fX$ from a scheme $U\in\Sch$; see \cite[\href{https://stacks.math.columbia.edu/tag/026O}{Tag~026O}]{stacks-project}.

\begin{Definition}[{\cite[Definition~3.1]{FKato-logcurve}, \cite[Section~1]{gillam-minimal}}]\label{def:log-str-on-stack}
Let $\fp:\fX\to\Sch$ be a stack.
    \begin{enumerate}[label = (\arabic*), ref = \theDefinition(\arabic*), itemsep = 0.3ex]
        \item\label{def-part:log-str-on-stack} A \emph{fine log structure} on $\fX$ is a functor $M_\fX:\fX\to\Log$ lifting $\fp$. 
        
        \begin{enumerate}[label = (\roman*), itemsep = 0.15ex]
            \item For $\xi\in\fX$ lying over $S\in\Sch$, we denote the induced log structure on $S$ by $\xi^*M_\fX$.
            
            \item For a morphism $\psi:\eta\to\xi$ in $\fX$ lying over $f:T\to S$ in $\Sch$, we denote the induced isomorphism of log structures on $T$ by
            \begin{align*}
                M_\fX(\psi):f^*(\xi^*M_\fX)\xrightarrow{\sim}\eta^*M_\fX.
            \end{align*}
        \end{enumerate}

        Following Convention~\ref{conv:sch-logsch}, we will omit the adjective `fine' and simply call $M_\fX$ a \emph{log structure} on $\fX$. We call $M_\fX$ an \emph{fs log structure} if its image is contained in $\Log^\text{fs}$.
        
        \item Let $M_\fX$, $M'_\fX$ be log structures on $\fX$. A \emph{map of log structures} $\varphi:M_\fX\to M'_\fX$ is a natural transformation, where $M_\fX$, $M_\fX'$ are regarded as functors $\fX\to\LogSch$ (rather than $\fX\to\Log$), such that $\varphi$ lifts the identity natural transformation of $\fp$.
        
        \begin{enumerate}[label = (\roman*), itemsep = 0.15ex]
            \item For $\xi\in\fX$ lying over $S\in\Sch$, we denote the induced map of log structures on $S$ by $\xi^*\varphi:\xi^*M_\fX\to \xi^*M'_\fX$.
            
            \item Naturality of $\varphi$ states that, if $\psi:\eta\to\xi$ in $\fX$ lies over $f:T\to S$ in $\Sch$, then the following diagram commutes.
            \begin{equation*}
            \begin{tikzcd}
                f^*(\xi^*M_\fX) \arrow[rr,"M_\fX(\psi)"] \arrow[d,swap,"f^*(\xi^*\varphi)"] & & \eta^*M_\fX \arrow[d,swap,"\eta^*\varphi"] \\
                f^*(\xi^*M_\fX') \arrow[rr,"M_\fX'(\psi)"] & & \eta^*M_\fX'
            \end{tikzcd}
            \end{equation*}
        \end{enumerate}
        
        \item Log structures on $\fX$ and the maps between them form a category, denoted by $\Log(\fX)$. For $\xi\in\fX$ lying over $S\in\Sch$, we have the functor $\xi^*:\Log(\fX)\to\Log(S)$.
        
        \item A \emph{log stack} is a pair $\fX^\dagger = (\fX,M_\fX)$ consisting of a stack $\fX$ and a log structure $M_\fX$ on it. 
        We call $\fX^\dagger$ an \emph{fs log stack} if $M_\fX$ is an fs log structure.
        
        \item\label{def-part:log-mor-to-stack} Let $\fX^\dagger = (\fX,M_\fX)$ be a log stack.
        \begin{enumerate}[label = (\roman*), itemsep = 0.15ex]
            \item For a log scheme $S^\dagger = (S,M_S)\in\LogSch$, a \emph{log morphism} $f^\dagger = (f,f^\flat):S^\dagger\to\fX^\dagger$ is defined to be a morphism of stacks $f:S\to\fX$ together with a map of log structures $f^\flat:f^*M_\fX\to M_S$ on $S$.
        
            \item More generally, given another log stack $\fY^\dagger = (\fY,M_\fY)$, we define a log morphism $f^\dagger = (f,f^\flat):\fY^\dagger\to\fX^\dagger$ to be a morphism of stacks $f:\fY\to\fX$ together with a map $f^\flat:f^*M_\fX\to M_\fY$ in $\Log(\fY)$, where $f^*M_\fX:= M_\fX\circ f:\fY\to\Log$.
        \end{enumerate}
        
        In either case, the log morphism $f^\dagger$ is called \emph{strict} if $f^\flat$ is an isomorphism.

        \item For a log stack $\fX^\dagger = (\fX,M_\fX)$, we define the monoid
        \begin{align*}
            \Gamma(\fX,\ol M_\fX) = \displaystyle\lim_{\xi:S\to\fX}\Gamma(S,\xi^*\ol M_\fX),
        \end{align*}
        as a limit over the category $\fX$. By definition, for every $\xi:S\to\fX$ with $S\in\Sch$, we have a natural projection $\xi^*:\Gamma(\fX,\ol M_\fX)\to\Gamma(S,\xi^*\ol M_\fX)$. 
        
        Explicitly, specifying an element $m\in\Gamma(\fX,\ol M_\fX)$ is the same as specifying an element $\xi^*m\in\Gamma(S,f^*\ol M_\fX)$ for each $\xi:S\to\fX$ with $S\in\Sch$ such that, for any $\psi:\eta\to\xi$ in $\fX$ lying over $f:T\to S$ in $\Sch$, the map induced by $M_\fX(\psi)$ must send $f^*(\xi^*m)$ to $\eta^*m$.
    \end{enumerate}
\end{Definition}

We have the following analogues of Definition~\ref{def:descent-datum-log-scheme} and Proposition~\ref{prop:log-str-descent-scheme} in the setting of algebraic stacks.

\begin{Definition}\label{def:descent-datum-log-stack}
    Let $\fX$ be an algebraic stack. Let $\pi:U\to\fX$ with $U\in\Sch$ be a morphism which is representable by schemes and is a smooth atlas. 
    
    For $i\in\{1,2\}$ and $jk\in\{12,13,23\}$, we have the natural projections
    \begin{align*}
        U\times_\fX U\times_\fX U \xrightarrow{p_{jk}} U\times_\fX U \xrightarrow{p_i} U.
    \end{align*}
    The $2$-fibre product $U\times_\fX U$ comes with a canonical isomorphism $\psi:\pi\circ p_1\xrightarrow{\sim} \pi\circ p_2$ in $\fX$ lying over the identity morphism of $U\times_\fX U$, where we view $\pi\circ p_i$ as objects of $\fX$.
    \begin{enumerate}[label = (\arabic*), itemsep = 0.3ex]
        \item A \emph{descent datum} $(M_U,\Psi)$ for a log structure with respect to $\pi:U\to\fX$ consists of 
        \begin{enumerate}[label = (\roman*), itemsep = 0.15ex]
            \item a log structure $M_U$ on the scheme $U$, and
        
            \item an isomorphism of log structures $\Psi:p_1^*M_U\xrightarrow{\sim} p_2^*M_U$ on the scheme $U\times_\fX U$
        \end{enumerate}
        which satisfy the \emph{cocyle condition} \eqref{eqn:descent-cocycle-condition} on the scheme $U\times_\fX U\times_\fX U$.

        \item For descent data $(M_U,\Psi)$ and $(M_U',\Psi')$ as above, a \emph{map of descent data} between these is a map $\varphi_U:M_U\to M_U'$ of log structures such that \eqref{eqn:descent-log-str-map} commutes.

        \item Descent data for log structures with respect to $\pi:U\to \fX$ and maps of descent data form a category, which we denote by $\Log(\pi:U\to \fX)$.

        \item Given any log structure $M_\fX$ on $\fX$, we get a descent datum $(\pi^*M_\fX,M_\fX(\psi))$ with respect to $\pi:U\to\fX$, where $M_\fX(\psi):p_1^*(\pi^*M_\fX)\xrightarrow{\sim}p_2^*(\pi^*M_\fX)$ is the isomorphism of log structures induced by the isomorphism $\psi$ in $\fX$.
    \end{enumerate}
\end{Definition}

\begin{Proposition}
\label{prop:log-str-descent-stack}
    In the setting of Definition~\ref{def:descent-datum-log-stack}, the functor
    \begin{align}\label{eqn:descent-equivalence-stack}
        \Log(\fX)&\to\Log(\pi:U\to\fX)\\
        M_\fX&\mapsto(\pi^*M_\fX,M_\fX(\psi))
    \end{align}
    is an equivalence of categories.
    
    Moreover, if $M_\fX$ belongs to $\Log(\fX)$ and $\varphi_U:(\pi^*M_\fX,M_\fX(\psi))\xrightarrow{\sim}(M_U,\Psi)$ is an isomorphism in $\Log(\pi:U\to\fX)$, then the following is an equalizer diagram.
    \begin{equation}\label{eqn:global-sections-equalizer}
    \begin{tikzcd}
        \Gamma(\fX,\ol M_\fX) \arrow[r,"\ol\varphi_U\,\circ\,\pi^*"] & \Gamma(U,\ol M_U) \arrow[rr, shift left=.5ex,"\ol\Psi\,\circ\, p_1^*"] \arrow[rr, swap, shift right=.5ex,"p_2^*"] && \Gamma(U\times_\fX U,p_2^*\ol M_U).
    \end{tikzcd}
    \end{equation}
\end{Proposition}
\begin{proof}
    This is a formal consequence of Proposition~\ref{prop:log-str-descent-scheme}.
\end{proof}

We often use Proposition~\ref{prop:log-str-descent-stack} to construct log structures on algebraic stacks by descent from a suitable smooth atlas; for instance, see Definitions~\ref{def:std-log-str-artin-cone} and \ref{def:div-log-str-smooth-nc-stack} below.

\subsection{Log structures on Artin cones} For this subsection, we fix a fine monoid $Q$. We then have the associated scheme $V_Q = \Spec\bC[Q]$ and the algebraic group $T_Q = \Spec\bC[Q^\text{gp}]$ acting on $V_Q$; see Definition~\ref{def:scheme-from-monoid} and Lemma~\ref{lem:scheme-from-monoid}.

\begin{Definition}\label{def:artin-cone}
    The \emph{Artin cone} associated to $Q$ is the quotient stack
    \begin{align*}
        \clA_Q:=[V_Q/T_Q].
    \end{align*}
    
    Explicitly, for $S\in\Sch$, objects of $\clA_Q(S)$ are pairs $(P,f)$, where the scheme $P$ is an \'etale locally trivial principal $T_Q$-bundle over $S$ and $f:P\to V_Q$ is a $T_Q$-equivariant morphism. An arrow $(P',f')\to (P,f)$ in $\clA_Q$ lying over $S'\to S$ in $\Sch$ is a pullback diagram of principal $T_Q$-bundles which respects the morphisms to $V_Q$.
\end{Definition}

\begin{Remark}
    To show that $\clA_Q$ is indeed a stack with respect to the fppf topology, one uses the fact that the group $T_Q$ is smooth and affine; see \cite[4.5.6 and 4.5.7]{olsson-stacks-book}.
\end{Remark}

\begin{Remark}\label{rem:patching-map-to-quotient}
    To define a morphism $S\to\clA_Q$ for $S\in\Sch$, it suffices to specify an \'etale covering $\{S_i\to S\}_i$ together with morphisms $\{f_i:S_i\to V_Q\}_i$ and $\{g_{ij}:S_i\times_S S_j\to T_Q\}_{i,j}$ such that we have the compatibility $f_i\cdot g_{ij} = f_j$ over each double overlap $S_i\times_S S_j$ and the cocyle condition $g_{ij}g_{jk} = g_{ik}$ over each triple overlap $S_i\times_S S_j\times_S S_k$.
\end{Remark}

\begin{Definition}\label{def:quotient-artin-cone}
    The \emph{quotient morphism} $\pi_Q:V_Q\to\clA_Q$ is defined by the pair consisting of the trivial principal $T_Q$-bundle $V_Q\times T_Q\to V_Q$ and the action morphism $V_Q\times T_Q\to V_Q$.
\end{Definition}

The next statement records the standard fact that $\clA_Q$ is an algebraic stack, with a smooth atlas given by the quotient morphism $\pi_Q:V_Q\to\clA_Q$.

\begin{Proposition}
\label{prop:std-atlas-artin-cone}
    The stack $\clA_Q$ is algebraic. Given $S\in\Sch$ and $(P,f)\in\clA_Q(S)$, let $\pi:P\to S$ denote the bundle projection morphism. Define the isomorphism
    \begin{align*}
        \theta:f^*(V_Q\times T_Q) = P\times T_Q \xrightarrow{\sim} P\times_S P = \pi^*P
    \end{align*}
    of principal $T_Q$-bundles on $P$ by $(x,g)\mapsto(x,x\cdot g)$. 
    
    Then the following is a $2$-fibre product diagram
    \begin{equation}\label{eqn:pullback-from-artin-cone}
    \begin{tikzcd}
        P \arrow[r,"f"] \arrow[d,swap,"\pi"] \arrow[dr,phantom,"\lrcorner",pos=0.001] & V_Q \arrow[dl,phantom,pos=0.45, "\substack{\mathbin{\rotatebox[origin=c]{45}{$\Leftarrow$}} \\ \theta}"] \arrow[d,"\pi_Q"] \\
        S \arrow[r,swap,"{(P,f)}"] & \clA_Q
    \end{tikzcd}
    \end{equation}
    and induces an isomorphism $P\xrightarrow{\sim} S\times_{\clA_Q} V_Q$. In particular, the morphism $\pi_Q:V_Q\to\clA_Q$ is representable by schemes and is a smooth atlas for $\clA_Q$.
\end{Proposition}
\begin{proof}
    Algebraicity of the stack $\clA_Q$ follows from \cite[8.1.12]{olsson-stacks-book}. The remaining assertions are immediate from definitions.
\end{proof}

\begin{Corollary}\label{cor:std-atlas-artin-cone}
    In Proposition~\ref{prop:std-atlas-artin-cone}, take the scheme $S$ to be $V_Q$ and take $(P,f)\in\clA_Q(S)$ to be the pair corresponding to $\pi_Q:V_Q\to\clA_Q$. Then the isomorphism 
    \begin{align*}
        V_Q\times T_Q\xrightarrow{\sim} V_Q\times_{\clA_Q}V_Q    
    \end{align*}
    induced by \eqref{eqn:pullback-from-artin-cone} identifies the projections $p_1,p_2:V_Q\times_{\clA_Q}V_Q\to V_Q$ with the projection $(x,g)\mapsto x$ and the action morphism $(x,g)\mapsto x\cdot g$ respectively.
\end{Corollary}

\begin{Construction}
\label{cons:descent-datum-artin-cone}
    We define a descent datum for a log structure on the Artin cone $\clA_Q$, with respect to the quotient morphism $\pi_Q:V_Q\to \clA_Q$ (Definition~\ref{def:quotient-artin-cone}). By Proposition~\ref{prop:std-atlas-artin-cone}, the morphism $\pi_Q$ is representable by schemes and is a smooth atlas. 
    
    By Corollary~\ref{cor:std-atlas-artin-cone}, we have an identification $V_Q\times_{\clA_Q}V_Q \simeq V_Q\times T_Q$ under which the projections $p_1,p_2:V_Q\times_{\clA_Q} V_Q\to V_Q$ are given by the algebra maps 
    \begin{align*}
        p_1^\#:\chi^q\mapsto \chi^q\otimes 1\quad\text{and}\quad p_2^\#:\chi^q\mapsto\chi^q\otimes\chi^q.    
    \end{align*}
    This induces an identification $V_Q\times_{\clA_Q} V_Q\times_{\clA_Q} V_Q \simeq V_Q\times T_Q\times T_Q$ under which the projections $p_{12},p_{13},p_{23}:V_Q\times_{\clA_Q} V_Q\times_{\clA_Q} V_Q\to V_Q\times_{\clA_Q} V_Q$ are given by the algebra maps
    \begin{align*}
        p_{12}^\#:\chi^q\otimes\chi^{q'}&\mapsto \chi^q\otimes\chi^{q'}\otimes 1,\\ p_{13}^\#:\chi^q\otimes\chi^{q'}&\mapsto \chi^q\otimes\chi^{q'}\otimes\chi^{q'},\quad\text{and}\\
        p_{23}^\#:\chi^q\otimes\chi^{q'}&\mapsto \chi^q\otimes\chi^q\otimes\chi^{q'}.
    \end{align*}
    
    On $V_Q$, take the standard log structure $M_{V_Q}$ (Definition~\ref{def:std-log-str-toric}). The pullbacks $p_1^*M_{V_Q}$ and $p_2^*M_{V_Q}$ are seen to be the log structures on $V_Q\times T_Q$ associated to the pre-log structures $\alpha_1,\alpha_2:\ul Q\to\clO_{V_Q\times T_Q}$ respectively, which are defined by 
    \begin{align*}
        \alpha_1: q\mapsto \chi^q\otimes 1\quad\text{and}\quad \alpha_2: q\mapsto \chi^q\otimes\chi^q.    
    \end{align*}
    Define the map $u:\ul Q\to\clO_{V_Q\times T_Q}^\times$ of sheaves of monoids by
    \begin{align*}
        u:q\mapsto 1\otimes\chi^q.    
    \end{align*}
    Since $u$ satisfies $\alpha_2 = u\alpha_1$, it induces isomorphism of log structures
    \begin{align*}
        \Psi_Q:p_1^*M_{V_Q}\xrightarrow{\sim}p_2^*M_{V_Q}
    \end{align*} 
    via \eqref{eqn:pre-log-diff-unit}. For all $q\in Q$, we have the following identity in $\clO^\times_{V_Q\times T_Q\times T_Q}$.
    \begin{align*}
        (p_{13}^\#\circ u)(q) &= 1\otimes\chi^q\otimes\chi^q \\ &= (1\otimes 1\otimes\chi^q)\cdot(1\otimes\chi^q\otimes 1) = (p_{23}^\#\circ u)(q)\cdot(p_{12}^\#\circ u)(q).
    \end{align*}
    From this identity and Lemma~\ref{lem-part:pre-log-diff-unit-cocyle}, it follows that the pair $(M_{V_Q},\Psi_Q)$ satisfies the cocycle condition and thus, it defines a descent datum with respect to $\pi_Q:V_Q\to\clA_Q$.
\end{Construction}

\begin{Definition}\label{def:std-log-str-artin-cone}
    The \emph{standard log structure} $M_{\clA_Q}$ on the Artin cone $\clA_Q$ is defined to be the one which corresponds under \eqref{eqn:descent-equivalence-stack} to the descent datum in Construction~\ref{cons:descent-datum-artin-cone}. 
    
    Let $\clA_Q^\dagger = (\clA_Q,M_{\clA_Q})$ denote the resulting log stack. Via the natural identification $\pi_Q^*M_{\clA_Q} = M_{V_Q}$, the morphism $\pi_Q:V_Q\to\clA_Q$ lifts to a strict log morphism $\pi_Q^\dagger:V_Q^\dagger\to\clA_Q^\dagger$.
\end{Definition}

\begin{Lemma}\label{lem:artin-cone-global-sections-ghost-sheaf}
    The map $Q\to\Gamma(V_Q,\ol M_{V_Q})$ induced by \eqref{eqn:nat-map-Q-to-M-toric} factors uniquely through the map $\Gamma(\clA_Q,\ol M_{\clA_Q})\to\Gamma(V_Q,\ol M_{V_Q})$ induced by $\pi_Q^*M_{\clA_Q} = M_{V_Q}$, yielding a natural map
    \begin{align}\label{eqn:monoid-to-global-sections-ghost-artin-cone}
        Q\to\Gamma(\clA_Q,\ol M_{\clA_Q}).
    \end{align}
\end{Lemma}
\begin{proof}
    Using the diagram \eqref{eqn:global-sections-equalizer} for the descent datum $(M_{V_Q},\Psi_Q)$ from Construction~\ref{cons:descent-datum-artin-cone}, it suffices to show the equality of the two compositions in the following diagram.
    \begin{equation}\label{eqn:factor-through-equalizer}
    \begin{tikzcd}
        Q \arrow[r] & \Gamma(V_Q,\ol M_{V_Q}) \arrow[rr,shift left=.5ex,"\ol\Psi_Q\,\circ\, p_1^*"] \arrow[rr,shift right=.5ex,swap,"p_2^*"] & & \Gamma(V_Q\times T_Q,p_2^*\ol M_{V_Q})
    \end{tikzcd}
    \end{equation}
    
    Recall from Construction~\ref{cons:descent-datum-artin-cone} that $p_1^*M_{V_Q}$ (resp. $p_2^*M_{V_Q}$) is the log structure associated to the pre-log structure $\ul Q\to\clO_{V_Q\times T_Q}$ given by $q\mapsto\chi^q\otimes 1$ (resp. $q\mapsto \chi^q\otimes\chi^q$), while $\Psi_Q$ is defined via \eqref{eqn:pre-log-diff-unit} and the map $\ul Q\to\clO_{V_Q\times T_Q}^\times$ given by $q\mapsto 1\otimes\chi^q$. By Lemma~\ref{lem-part:same-ghost}, $p_1^*\ol M_{V_Q}$ and $p_2^*\ol M_{V_Q}$ are the \emph{same} quotient of $\ul Q$ and, by Lemma~\ref{lem-part:unit-induced-map}, $\ol\Psi_Q$ is the identity map on ghost sheaves. It follows that the two compositions in \eqref{eqn:factor-through-equalizer} are equal.
\end{proof}

\begin{Proposition}[{\cite[5.17]{olsson-stack-of-log-str}}]\label{prop:log-maps-to-artin-cone}
    Fix a fine monoid $Q$. Consider the following two types of objects for a log scheme $S^\dagger = (S,M_S)$. 
    \begin{enumerate}[itemsep = 0.3ex]
        \myitem[(AC$\phantom{}_1$)]\label{log-map-to-A_Q} Log morphism $S^\dagger\to\clA_Q^\dagger$.
        
        \myitem[(AC$\phantom{}_2$)]\label{Q-to-ghost-sheaf} Monoid map $Q\to\Gamma(S,\ol M_S)$.
    \end{enumerate}
    As the log scheme $S^\dagger$ varies, each of \ref{log-map-to-A_Q} and \ref{Q-to-ghost-sheaf} determines a category fibred in groupoids over $\LogSch$ using the obvious notion of pullbacks along log morphisms. 
    
    These two categories fibred in groupoids over $\LogSch$ are equivalent. 
    
    Moreover, under this equivalence, a log morphism $S^\dagger\to\clA_Q^\dagger$ is strict if and only if the corresponding monoid map $Q\to\Gamma(S,\ol M_S)$ lifts to a log chart in some \'etale neighborhood of any given point of $S$.
\end{Proposition}

\begin{proof}
    We define functors \ref{log-map-to-A_Q}~$\to$~\ref{Q-to-ghost-sheaf}~$\to$~\ref{log-map-to-A_Q}, but omit the verification that their composites \ref{log-map-to-A_Q}~$\to$~\ref{log-map-to-A_Q} and \ref{Q-to-ghost-sheaf}~$\to$~\ref{Q-to-ghost-sheaf} are naturally isomorphic to identity functors.

    \ref{log-map-to-A_Q}~$\to$~\ref{Q-to-ghost-sheaf}: Given $h^\dagger = (h,h^\flat):S^\dagger\to\clA_Q^\dagger$, associate to it the composition 
    \begin{equation*}
    \begin{tikzcd}
        Q \arrow[r,"\eqref{eqn:monoid-to-global-sections-ghost-artin-cone}"] & \Gamma(\clA_Q,\ol M_{\clA_Q}) \arrow[r,"h^*"] & \Gamma(S,h^*\ol M_{\clA_Q}) \arrow[r,"\Gamma(\ol h\phantom{}^\flat)"] & \Gamma(S,\ol M_S).
    \end{tikzcd}
    \end{equation*}
    
    \ref{Q-to-ghost-sheaf}~$\to$~\ref{log-map-to-A_Q}:
    Given $\beta:Q\to\Gamma(S,\ol M_S)$, associate a log morphism $h^\dagger:S^\dagger\to\clA_Q^\dagger$ to it as follows. For an informal explanation of the construction below, see Remark~\ref{rem:informal-log-mor-to-artin-cone}.

    We will work in the category $\Sch{\!/S}$ consisting of $S$-schemes and $S$-morphisms. For any $W\in\Sch{\!/S}$, let $M_W$ denote the pullback of $M_S$. In what follows, we refer to fppf sheaves simply as sheaves. The presheaf of sets $F$ on $\Sch{\!/S}$ defined by
    \begin{align*}
        F(W) := \left\{
        \begin{matrix}
            \tilde\beta_W:Q\to\Gamma(W,M_W)
        \end{matrix}
        \;\middle|\;
        \begin{matrix}
            \tilde\beta_W \text{ is a monoid map lifting the} \\
            \text{pullback of } \beta \text{ along } W\to S
        \end{matrix}
        \right\}
    \end{align*}
    is a sheaf by \cite[A.1]{olsson-stack-of-log-str}.
    From Lemmas~\ref{lem:fine-log-str-integral} and \ref{lem:local-lift-from-ghost}, we deduce that the sheaf $F$ is a torsor over the sheaf of abelian groups $G$ on $\Sch{\!/S}$ defined by 
    \begin{align*}
        G(W) := \Hom_{\Mon}(Q,\Gamma(W,\clO_W^\times)),    
    \end{align*}
    where the action $\mu:F\times G\to F$ is given by multiplication. The sheaf $G$ also acts by multiplication on the sheaf of sets $H$ on $\Sch{\!/S}$ defined by 
    \begin{align*}
        H(W) := \Hom_{\Mon}(Q,\Gamma(W,\clO_W))    
    \end{align*}
    and the map $F\to H$ induced by $M\to\clO_S$ is $G$-equivariant.
    
    The sheaf $G$ is represented by the group scheme $S\times T_Q = \relSpec_S\clO_S[Q^\text{gp}]$ over $S$, with the universal family $\eta_G\in G(S\times T_Q)$ given by $\eta_G(q) = 1\otimes\chi^q$. Similarly, $H$ is represented by the $S$-scheme $S\times V_Q = \relSpec_S\clO_S[Q]$, with the universal family $\eta_H\in H(S\times V_Q)$ given by $\eta_H(q) = 1\otimes\chi^q$. Since the morphism $S\times T_Q\to S$ is affine, \cite[4.5.6]{olsson-stacks-book} shows that the $G$-torsor $F$ is represented by a principal $T_Q$-bundle $\pi:P\to S$; denote the corresponding universal family in $F(P)$ by $\tilde\beta_P:Q\to\Gamma(P,M_P)$.
    
    Let $\pi_1,\pi_2:P\times T_Q\to P$ be the projection $(x,g)\mapsto x$ and the action morphism $(x,g)\mapsto x\cdot g$ respectively. Since the action $\mu:F\times G\to F$ corresponds to the action $\pi_2:P\times T_Q \to P$, the following diagram of sheaves on $\Sch{\!/S}$ commutes.
    \begin{equation}\label{eqn:key-comm-diag-for-descent}
    \begin{tikzcd}
        P\times_S(S\times T_Q) \arrow[d,"\tilde\beta_P\times\eta_G"] \arrow[r,"\pi_2"] & P \arrow[d,"\tilde\beta_P"] \\
        F \times G \arrow[r,"\mu"] & F
    \end{tikzcd}
    \end{equation}
    Here, we implicitly use the Yoneda embedding to identify objects of $\Sch{\!/S}$ as sheaves on $\Sch{\!/S}$. Concretely, \eqref{eqn:key-comm-diag-for-descent} translates to the following identity in $F(P\times T_Q)$.
    \begin{align}\label{eqn:key-identity-for-descent}
        \pi_1^*\tilde\beta_P \cdot (\pi\times\text{id})^*\eta_G = \pi_2^*\tilde\beta_P.
    \end{align}
    By Proposition~\ref{prop:log-map-to-toric}, $\tilde\beta_P$ yields a log morphism $f^\dagger = (f,f^\flat):P^\dagger\to V_Q^\dagger$. The underlying morphism $f:P\to V_Q$ is induced by $F\to H$ and is therefore $T_Q$-equivariant. Thus, the pair $(P,f)$ defines a morphism $h:S\to\clA_Q$. The log structure $h^*M_{\clA_Q}$ on $S$ is defined by the descent datum $f^*(M_{V_Q},\Psi_Q)\in\Log(\pi:P\to S)$. To complete the proof, we need the following claim.
    
    \begin{claim}\label{claim:artin-cone}
        The map $f^\flat:f^*M_{V_Q}\to M_P$ is in fact a map of descent data 
        \begin{align}\label{eqn:log-mor-to-artin-cone-descent-data-map}
            f^*(M_{V_Q},\Psi_Q)\to(M_P,\text{id}).    
        \end{align}
    \end{claim}
    \begin{proof}[Proof of claim]
        We need to show that the following diagram of log structures on $P\times T_Q$ commutes, where we use the identifications $P\times T_Q\simeq P\times_S P$ and $V_Q\times T_Q\simeq V_Q\times_{\clA_Q} V_Q$ from Proposition~\ref{prop:std-atlas-artin-cone} and  Corollary~\ref{cor:std-atlas-artin-cone} respectively.
        \begin{equation}\label{eqn:log-mor-to-artin-cone-descent}
        \begin{tikzcd}
            (f\circ\pi_1)^*M_{V_Q} \arrow[rr,"(f\times\text{id})^*\Psi_Q","\sim"'] \arrow[dr,swap,"\pi_1^*f^\flat"] & & (f\circ\pi_2)^*M_{V_Q} \arrow[dl,"\pi_2^*f^\flat"] \\
            & M_{P\times T_Q} &
        \end{tikzcd}
        \end{equation}
        Recall that $\Psi_Q$ is induced via \eqref{eqn:pre-log-diff-unit} by the map $u:\ul Q\to\clO^\times_{V_Q\times T_Q}$ given by $q\mapsto 1\otimes\chi^q$. Therefore, $(f\times\text{id})^*\Psi_Q$ is induced by $(f\times\text{id})^\#\circ u  = (\pi\times\text{id})^*\eta_G$. Similarly, since $f^\flat$ is induced by $\tilde\beta_P$, it follows that $\pi_i^*f^\flat$ is induced by $\pi_i^*\tilde\beta_P$ for $i=1,2$. Now, the identity \eqref{eqn:key-identity-for-descent} and the diagram \eqref{eqn:pre-log-diff-unit-map-compatible} together imply that \eqref{eqn:log-mor-to-artin-cone-descent} commutes.
    \end{proof}
    
    By Claim~\ref{claim:artin-cone}, $f^\flat$ descends to a map of log structures $h^\flat:h^*M_{\clA_Q}\to M_S$, which completes the definition of $h^\dagger = (h,h^\flat):S^\dagger\to\clA_Q^\dagger$. 
    
    Finally, note that $h^\flat$ is an isomorphism if and only if the same is true of $f^\flat$, i.e., if and only if $\tilde\beta_P:Q\to\Gamma(P,M_P)$ is a log chart. In other words, $h^\dagger$ is strict if and only if the universal lift $\tilde\beta_P$ of $\beta$ is a log chart, i.e., if and only if $\beta$ lifts to a log chart in some \'etale neighborhood of any given point of $S$.
\end{proof}

\begin{Remark}\label{rem:informal-log-mor-to-artin-cone}
    The informal idea behind the functor \ref{Q-to-ghost-sheaf} $\to$ \ref{log-map-to-A_Q} is the following. By Lemma~\ref{lem:local-lift-from-ghost}, we may find an \'etale covering $\{S_i\to S\}_i$ and lifts $\tilde\beta_i:Q\to\Gamma(S_i,M_{S_i})$ of the maps $\beta|_{S_i}$. Each $\tilde\beta_i$ induces a log morphism $f_i^\dagger: S_i^\dagger\to V_Q^\dagger$ by Proposition~\ref{prop:log-map-to-toric}. 
    
    By Lemma~\ref{lem:fine-log-str-integral}, $\tilde\beta_i$ and $\tilde\beta_j$ differ over the double overlaps $S_{ij}:=S_i\times_S S_j$ by a monoid map $\gamma_{ij}:Q\to\Gamma(S_{ij},\clO_{S_{ij}}^\times)$, i.e., a morphism $g_{ij}:S_{ij}\to T_Q$ with $\tilde\beta_i\gamma_{ij} = \tilde\beta_j$ and $f_i\cdot g_{ij} = f_j$.
    These satisfy the cocycle condition $g_{ij}g_{jk} = g_{ik}$ on the triple overlaps $S_{ijk}:=S_i\times_S S_j\times_S S_k$. As in Remark~\ref{rem:patching-map-to-quotient}, the data $\{f_i^\dagger,g_{ij}\}$ now patch together to define a log morphism $S^\dagger\to\clA_Q^\dagger$.
\end{Remark}

\begin{Remark}
    For fixed $Q$ and $S^\dagger$, the equivalence \ref{log-map-to-A_Q} $\simeq$ \ref{Q-to-ghost-sheaf} shows that the groupoid of log morphisms $S^\dagger\to\clA_Q^\dagger$ is equivalent to a set, i.e., every object has a trivial automorphism group. This should be contrasted with the groupoid of morphisms $S\to\clA_Q$, which is never equivalent to a set for non-trivial $Q$ and non-empty $S$.
\end{Remark}

\subsection{Log structures associated to nc divisors}\phantom{}

\begin{Definition}\label{def:nc-divisor-scheme}
    Let $X$ be a smooth scheme and $D\subset X$ be a reduced closed subscheme.
    \begin{enumerate}[label = (\arabic*), ref = \theDefinition(\arabic*), itemsep = 0.3ex]
        \item\label{def-part:snc} We call $D$ a \emph{simple normal crossings} (\emph{snc}) \emph{divisor} in $X$ if any point $x\in X$ has a Zariski open neighborhood $U$ and a smooth morphism $(z_1,\ldots,z_r):U\to\bC^r$ for some integer $r\ge 0$ such that $U\times_XD\subset U$ is defined by the equation $\prod z_i = 0$. 
        
        \item\label{def-part:nc} We call $D$ a \emph{normal crossings} (\emph{nc}) \emph{divisor} in $X$ if any point $x\in X$ has an \'etale neighborhood $U$ such that $U\times_X D\subset U$ is an snc divisor.
    \end{enumerate}
\end{Definition}

\begin{Lemma}\label{lem:pullback-of-nc-div-log}
    Let $X$ be a smooth scheme and let $D\subset X$ be an nc divisor. Let $M_{(X,D)}$ denote the associated divisorial log structure (Example~\ref{exa:div-log-str}). 
    \begin{enumerate}[label = (\arabic*), ref = \theLemma(\arabic*), itemsep = 0.3ex]
        \item\label{lem-part:nc-div-log-chart} Let $U\to X$ be an \'etale neighborhood and let $(z_1,\ldots,z_r):U\to\bC^r$ be a smooth morphism such that $U\times_XD\subset U$ is defined by $\prod z_i = 0$. Then the map
        \begin{align}\label{eqn:nc-div-log-chart}
            \bN^r\to\Gamma(U,M_{(X,D)}), \quad (a_1,\ldots,a_r)\mapsto \prod z_i^{a_i}.
        \end{align}
        is a log chart for $M_{(X,D)}$ over $U$. In particular, $M_{(X,D)}$ is a fine log structure.

        \item\label{lem-part:pullback-of-nc-div-log} If $\pi:X'\to X$ is a smooth morphism of schemes, then $D':=X'\times_XD\subset X'$ is an nc divisor. Moreover, the unique map of log structures $\pi^\flat:\pi^*M_{(X,D)}\to M_{(X',D')}$ provided by Example~\ref{exa:div-log-str} is an isomorphism.
    \end{enumerate}
\end{Lemma}
\begin{proof}
    (1): This is immediate from definitions.

    (2): The assertions may be checked \'etale locally on $X'$. Thus, we may assume that we have a smooth morphism $(z_1,\ldots,z_r):X\to\bC^r$ such that $D\subset X$ is defined by $\prod z_i = 0$. Since the composition $X'\to X\to\bC^r$ is smooth, the result follows from part (1).
\end{proof}

\begin{Definition}\label{def:branch}
    Let $D\subset X$ be an nc divisor in a smooth scheme and let $x\in D$.
    \begin{enumerate}[label = (\arabic*), ref = \theDefinition(\arabic*), itemsep = 0.3ex]
        \item Consider triples $(U\to X,u,\Delta)$, where $U\to X$ is an \'etale neighborhood such that $U\times_XD\subset U$ is an snc divisor, $u\in U$ is a point lying over $x$ and $\Delta$ is an irreducible (smooth) component of $U\times_XD$ passing through $u$.
        
        Consider two such triples $(U\to X,u,\Delta)$ and $(U'\to X,u',\Delta')$ to be equivalent if the pullbacks of the divisors $\Delta\subset U$ and $\Delta'\subset U'$ to the fibre product $U\times_X U'$ are equal in a Zariski open neighborhood of the point $(u,u')$.

        \item\label{def-part:branch} A \emph{branch} of $D$ at $x$ is defined to be an equivalence class of triples $(U\to X,u,\Delta)$ as above. The branch of $D$ at $x$ associated to a triple $(U\to X,u,\Delta)$ will usually be denoted simply by $\Delta$ when there is no risk of confusion.
            
        \item Each branch $\Delta$ of $D$ at $x$ has a well-defined \emph{tangent space} $T_{\Delta,x}\subset T_{X,x}$, which is a linear subspace described as follows: if $(U\to X,u,\Delta)$ is a triple representing the branch, then $T_{\Delta,x}$ is the image of $T_{\Delta,u}\hookrightarrow T_{U,u}\xrightarrow{\sim}T_{X,x}$. 
        
        \item\label{def-part:normal-space} The \emph{normal space} of a branch $\Delta$ of $D$ at $x$ is the quotient $\clN_{\Delta/X}|_x:=T_{X,x}/T_{\Delta,x}$.
    \end{enumerate}
\end{Definition}

\begin{Lemma}\label{lem:nc-ghost-stalk}
    Let $X$ be a smooth scheme and let $D\subset X$ be an nc divisor, with the associated divisorial log structure denoted by $M_{(X,D)}$. Let $x\in D$ be a point; we also view it as a morphism $x:\Spec\bC\to X$ which factors through $D$.
    \begin{enumerate}[label = (\arabic*), ref = \theDefinition(\arabic*), itemsep = 0.3ex]
        \item\label{lem-part:ghost-sheaf-stalk-generators} For each branch $\Delta$ of $x$, we have a well-defined element $\rho_\Delta\in x^*\ol M_{(X,D)}$ characterized by the following property: if the triple $(U\to X,u,\Delta)$ represents $\Delta$ and $z\in\clO_{U,u}$ is a local equation for $\Delta\subset U$ at $u$, then $\rho_\Delta$ is the image of $z$ in $x^*\ol M_{(X,D)}$.

        \item\label{lem-part:ghost-sheaf-stalk-free} The monoid $x^*\ol M_{(X,D)}$ is freely generated by the finite collection of elements $\{\rho_\Delta\}$, where $\Delta$ ranges over the branches of $D$ at $x$.

        \item\label{lem-part:branches-pullback} Let $X'\to X$ be a smooth morphism of schemes and define $D':=X'\times_XD$. Let $x'\in X'$ be a point lying over $x$. For each triple $(U\to X,u,\Delta)$ representing a branch of $\Delta$ at $x$, define $U':=X'\times_XU$ and let $\Delta'$ to be the unique irreducible component of the smooth divisor $X'\times_X\Delta\subset U'$ passing through the point $u':=(x',u)$. The map        
        \begin{align*}
            (U\to X,u,\Delta)\mapsto(U'\to X',u',\Delta')    
        \end{align*}
        induces a bijection from the branches of $D$ at $x$ to the branches of $D'$ at $x'$.
    \end{enumerate}
\end{Lemma}
\begin{proof}
    (1): The fact that $\rho_\Delta$ is independent of the choices of the triple $(U\to X,u,\Delta)$ and the local equation $z$ is immediate from definitions. 
    
    (2): After replacing $X$ by an \'etale neighborhood of $x$, the divisor $D$ is defined by $\prod z_i = 0$ for some smooth morphism $(z_1,\ldots,z_r):X\to\bC^r$ mapping $x$ to $0\in\bC^r$. The chart \eqref{eqn:nc-div-log-chart} now shows that the images of $z_1,\ldots,z_r$ freely generate the monoid $x^*\ol M_{(X,D)}$.

    (3): This is immediate from part (2) and Lemma~\ref{lem-part:pullback-of-nc-div-log}. 
\end{proof}

We now generalize the discussion of nc divisors and their associated divisorial log structures from smooth schemes to smooth algebraic stacks.

\begin{Definition}\label{def:div-log-str-smooth-nc-stack}
    Let $\fX$ be a smooth algebraic stack. 
    
    Consider an nc divisor $\fD\subset\fX$, i.e., a reduced closed substack such that, for any smooth morphism $\xi:S\to\fX$ with $S\in\Sch$, the pullback $\fD_\xi:=S\times_\fX\fD\subset S$ is an nc divisor. Let $M_{(S,\fD_\xi)}$ denote the associated divisorial log structure on $S$ (Example~\ref{exa:div-log-str}).

    A \emph{divisorial log structure} associated to $\fD\subset\fX$ is a log structure $M_{(\fX,\fD)}:\fX\to\Log$ such that, for any smooth morphism $\xi:S\to\fX$ with $S\in\Sch$, the structure map $\xi^*M_{(\fX,\fD)}\to\clO_S$ identifies $\xi^*M_{(\fX,\fD)}$ with the subsheaf $M_{(S,\fD_\xi)}\subset\clO_S$.
\end{Definition}

\begin{Proposition}\label{prop:div-log-str-smooth-nc-stack}
    In the setting of Definition~\ref{def:div-log-str-smooth-nc-stack}, a divisorial log structure associated to $\fD\subset\fX$ exists and is unique up to unique isomorphism.
\end{Proposition}
\begin{proof}
    The proposition follows from Lemma~\ref{lem-part:pullback-of-nc-div-log} when $\fX$ is a scheme. The strategy is to `bootstrap' from this to the case where $\fX$ is an algebraic space and then to the general case. In fact, we will only give the details of how to go from algebraic spaces to algebraic stacks; going from schemes to algebraic spaces is very similar and slightly simpler. 
    
    Accordingly, assume that the proposition is true for the case of algebraic spaces. Let $\fD\subset\fX$ be as in Definition~\ref{def:div-log-str-smooth-nc-stack}. We now construct a divisorial log structure $M_{(\fX,\fD)}$.
    
    \begin{step}
        Consider an object $\xi\in\fX$ lying over $S\in\Sch$. We must define the log structure $\xi^*M_{(\fX,\fD)}$ on $S$. For this, we need the following claim.
        \begin{claim}\label{claim:div-log-str-smooth-nc-stack}
            There exist $U,S'\in\Sch$ and a 2-commutative diagram
            \begin{equation}\label{eqn:div-log-str-smooth-nc-stack-1}
            \begin{tikzcd}
                S' \arrow[r] \arrow[d] & U \arrow[d] \arrow[dl,phantom,pos=0.35, "\substack{\mathbin{\rotatebox[origin=c]{45}{$\Leftarrow$}}}"] \\
                S \arrow[r,"\xi"] & \fX
            \end{tikzcd}
            \end{equation}
            such that $U\to\fX$ is smooth and $S'\to S$ is surjective and smooth.
        \end{claim}
        \begin{proof}[Proof of claim]
            Choose $U\to\fX$ to be a smooth atlas for $\fX$. Then choose $S'\to U\times_\fX S$ to be a smooth atlas for the algebraic space $U\times_\fX S$.
        \end{proof}

        Choose a diagram \eqref{eqn:div-log-str-smooth-nc-stack-1} as in Claim~\ref{claim:div-log-str-smooth-nc-stack}. Since the algebraic spaces $U$, $U':=U\times_\fX U$ and $U'':=U\times_\fX U\times_\fX U$ are all smooth over $\fX$, pulling back $\fD\subset\fX$ gives rise to nc divisors $\fD_U\subset U$, $\fD_{U'}\subset U'$ and $\fD_{U''}\subset U''$ respectively. The proposition is true for algebraic spaces and so, we get associated divisorial log structures $M_{(U,\fD_U)}$, $M_{(U',\fD_{U'})}$ and $M_{(U'',\fD_{U''})}$. Write 
        \begin{align*}
            p_1,p_2:U'\to U,\quad p_{12},p_{13},p_{23}:U''\to U'\quad\text{and}\quad \hat p_1,\hat p_2,\hat p_3:U''\to U
        \end{align*}
        for the natural projections. As $p_i$ is smooth, $p_i^*M_{(U,\fD_U)}$ is a divisorial log structure associated to $\fD_{U'}\subset U'$; similar considerations apply to $p_{jk}$ and $\hat p_l$. The proposition is true for algebraic spaces and so, we have \emph{unique} isomorphisms of log structures $p_i^\flat$, $p_{jk}^\flat$ and $\hat p_l^\flat$ which lift $p_i$, $p_{jk}$ and $\hat p_l$ to strict log morphisms. Uniqueness of these isomorphisms implies that
        \begin{align*}
            \hat p_1^\flat &= p_{12}^\flat\circ(p_{12}^*p_1^\flat) = p_{13}^\flat\circ(p_{13}^*p_1^\flat),\\
            \hat p_2^\flat &= p_{12}^\flat\circ(p_{12}^*p_2^\flat) = p_{23}^\flat\circ(p_{23}^*p_1^\flat),\quad\text{and}\\
            \hat p_3^\flat &= p_{13}^\flat\circ(p_{13}^*p_2^\flat) = p_{23}^\flat\circ(p_{23}^*p_2^\flat).
        \end{align*}
        Define the isomorphism of log structures $\Psi := (p_2^\flat)^{-1}\circ p_1^\flat: p_1^*M_{(U,\fD_U)}\xrightarrow{\sim} p_2^*M_{(U,\fD_U)}$. The above identities imply that the cocycle condition $p_{13}^*\Psi = p_{23}^*\Psi\circ p_{12}^*\Psi$ holds.
        
        Pulling back $(M_{(U,\fD_U)},\Psi)$ along $S'\to U$ and $S'\times_S S'\to U\times_\fX U$ yields a descent datum belonging to $\Log(S'\to S)$. Define $\xi^*M_{(\fX,\fD)}$ by applying Proposition~\ref{prop:log-str-descent-scheme} to this descent datum. One verifies that a different choice of diagram \eqref{eqn:div-log-str-smooth-nc-stack-1} yields a canonically isomorphic log structure. Thus, $\xi^*M_{(\fX,\fD)}$ is well-defined up to unique isomorphism.        
        
        In particular, when $\xi:S\to\fX$ is smooth, we may choose \eqref{eqn:div-log-str-smooth-nc-stack-1} to be the diagram in which $S'\to S$ and $S'\to U$ are both the identity morphism $\text{id}_S$. It follows that the structure map $\xi^*M_{(\fX,\fD)}\to\clO_S$ identifies $\xi^*M_{(\fX,\fD)}$ with the subsheaf $M_{(S,\fD_\xi)}\subset\clO_S$.
    \end{step}

    \begin{step}
        Consider a morphism $\psi:\eta\to\xi$ in $\fX$ lying over $f:T\to S$ in $\Sch$. We must define the isomorphism of log structures $M_{(\fX,\fD)}(\psi):f^*(\xi^*M_{(\fX,\fD)})\xrightarrow{\sim}\eta^*M_{(\fX,\fD)} $.
        
        To begin, choose a $2$-commutative diagram \eqref{eqn:div-log-str-smooth-nc-stack-1} as in Claim~\ref{claim:div-log-str-smooth-nc-stack} and complete it to the following diagram with $T':=S'\times_S T\in\Sch$.
        \begin{equation}\label{eqn:div-log-str-smooth-nc-stack-2}
        \begin{tikzcd}
            T'\arrow [r] \arrow[dr,phantom,"\lrcorner",pos=0.001] \arrow[d] & S' \arrow[r] \arrow[d] & U \arrow[d] \arrow[dl,phantom,pos=0.35, "\substack{\mathbin{\rotatebox[origin=c]{45}{$\Leftarrow$}}}"] \\
            T \arrow[r,"f"] & S \arrow[r,"\xi"] & \fX
        \end{tikzcd}
        \end{equation}
        Use the right square (resp. the composition of the two squares) in \eqref{eqn:div-log-str-smooth-nc-stack-2} to define the log structure $\xi^*M_{(\fX,\fD)}$ (resp. $\eta^*M_{(\fX,\fD)}$). The descent datum in $\Log(T'\to T)$ corresponding to $\eta^*M_{(\fX,\fD)}$ is the pullback of the descent datum in $\Log(S'\to S)$ corresponding to $\xi^*M_{(\fX,\fD)}$. This defines $M_{(\fX,\fD)}(\psi)$ via Proposition~\ref{prop:log-str-descent-scheme}. One verifies that this isomorphism is well-defined and functorial.
    \end{step}
    This proves the existence of a divisorial log structure. The proof that it is unique up to unique isomorphism is very similar and will be omitted.
\end{proof}

\begin{Remark}
    In view of Proposition~\ref{prop:div-log-str-smooth-nc-stack}, we may unambiguously refer to \emph{the} divisorial log structure $M_{(\fX,\fD)}$ associated to an nc divisor $\fD$ in a smooth algebraic stack $\fX$.
\end{Remark}

\begin{Example}\label{exa:div-log-str-smooth-nc-stack-mor}
    Let $\fX$, $\fY$ be smooth algebraic stacks and let $\fD\subset\fX$, $\fE\subset\fY$ be nc divisors. Consider a morphism $\varphi:\fY\to\fX$ representable by algebraic spaces such that its restriction to the open substack $\fY\setminus\fE\subset\fY$ factors through the open substack $\fX\setminus\fD\subset\fX$.
    
    Define the log stacks $\fX^\dagger = (\fX,M_{(\fX,\fD)})$ and $\fY^\dagger = (\fY,M_{(\fY,\fE)})$. Using the case of schemes discussed in Example~\ref{exa:div-log-str} and arguments very similar to Proposition~\ref{prop:div-log-str-smooth-nc-stack}, one finds that $\varphi$ lifts uniquely to a log morphism $\varphi^\dagger:\fY^\dagger\to\fX^\dagger$.
\end{Example}
    \section{Systems of line bundles}\label{sec:slb-prelims}

We discuss the notion of a \emph{system of line bundles} on a scheme, which will be used in Section~\ref{sec:log-via-slb} to give concrete descriptions of log structures and log morphisms. Very closely related notions are discussed in \cite[Section~3]{borne-vistoli-DF} and \cite[Section~2.4]{RSPW}. We also discuss \emph{presentations} for systems of line bundles (Subsection~\ref{subsec:slb-presentation}); this is a key ingredient in the comparison of algebraic and symplectic log maps carried out in Section~\ref{sec:log-maps}.

\begin{Notation}
    Throughout this section, we use $S$, $S'$, $T$ etc to denote schemes belonging to $\Sch$ and we use $P,P',Q,Q'$ etc to denote integral monoids; see Definition~\ref{def:monoid-props}.
\end{Notation}

\subsection{Preliminaries}

\begin{Discussion}\label{disc:principal-vs-line-bundle}
    For a scheme $S$, recall the natural one-to-one correspondence between 
    \begin{enumerate}[itemsep = 0.3ex]
        \item \emph{principal $\bC^\times$-bundles} on $S$, i.e., $\clO_S^\times$-torsors, and
        
        \item \emph{line bundles} on $S$, i.e., locally free $\clO_S$-modules of rank $1$.
    \end{enumerate}

    Given a rank $1$ locally free $\clO_S$-module $L$, the associated $\clO_S^\times$-torsor $L^\times$ is obtained as the sheaf of nowhere vanishing sections of $L$. Conversely, given an $\clO_S^\times$-torsor $L^\times$, the associated rank $1$ locally free $\clO_S$-module $L$ is obtained as 
    \begin{align}\label{eqn:line-bundle-from-torsor}
        L = \clO_S\oplus_{\clO_S^\times}L^\times,    
    \end{align}
    i.e., the quotient of $\clO_S\times L^\times$ under the equivalence relation generated by $(fg,\sigma)\sim(f,g\sigma)$ for local sections $f$, $g$ and $\sigma$ of $\clO_S$, $\clO_S^\times$ and $L^\times$ respectively. 
\end{Discussion}

\begin{Discussion}\label{disc:sheaf-cartier-div}
    Let $\clK_S$ be the \emph{sheaf of rational functions} on $S$, as in \cite{kleiman-rational-functions} or \cite[page~140]{Har77}, and let $\clK_S^\times\subset\clK_S$ be its subsheaf of multiplicative units. 
    The sheaf of \emph{Cartier divisors} on $S$ is the quotient $\text{CaDiv}_S = \clK_S^\times/\clO_S^\times$ and we write the operation on it additively. A Cartier divisor is \emph{effective} if its local lifts to $\clK_S^\times$ lie within the subsheaf $\clO_S\subset\clK_S$.
    
    Any Cartier divisor $D$ on $S$ yields a rank $1$ locally free $\clO_S$-submodule $\clO_S(D)\subset\clK_S$. Its associated $\clO_S^\times$-torsor, denoted by $\clO_S(D)^\times\subset\clK_S^\times$, can be described as the inverse image of $-D\in \Gamma(S,\text{CaDiv}_S)$ under the quotient map $\clK_S^\times\to\clK_S^\times/\clO_S^\times$.
    If $D'$ is another Cartier divisor such that the difference $D-D'$ is effective, then we have an inclusion of $\clO_S$-submodules $\clO_S(D')\subset\clO_S(D)\subset\clK_S$. In particular, when $D$ itself is effective, taking $D' = 0$ yields a canonical global section $\taut_D\in H^0(S,\clO_S(D))$ which is obtained as the image of the constant function $1\in H^0(S,\clO_S)$ under the inclusion $\clO_S\subset\clO_S(D)$. More generally, $\taut_D$ is defined over the dense open subset $U\subset S$ where $D|_U$ is effective and may therefore be regarded as a \emph{rational section} of $\clO_S(D)$ even when $D$ is not effective.

    For any Cartier divisors $D$ and $D'$ on $S$, multiplication in $\clK_S$ induces an isomorphism $\clO_S(D)\otimes\clO_S(D')\xrightarrow{\sim}\clO_S(D+D')$ of rank $1$ locally free $\clO_S$-modules. Moreover, under this isomorphism, the rational section $\taut_{D}\otimes\taut_{D'}$ maps to the rational section $\taut_{D+D'}$.
\end{Discussion}

\begin{Definition}\label{def:vcd}
    A \emph{virtual Cartier divisor} on $S$ is defined to be a pair $(L,\sigma)$ consisting of a line bundle $L$ on $S$ together with a global section $\sigma\in H^0(S,L)$; in particular, the section $\sigma$ is allowed to be identically zero. An \emph{isomorphism of virtual Cartier divisors} $(L,\sigma)\xrightarrow{\sim}(L',\sigma')$ on $S$ is an isomorphism of line bundles $L\xrightarrow{\sim} L'$ which maps $\sigma$ to $\sigma'$.

    Given a virtual Cartier divisor $(L,\sigma)$ on $S$, its \emph{pullback} along a morphism $f:S'\to S$ of schemes is defined to be the virtual Cartier divisor $f^*(L,\sigma) := (f^*L,f^*\sigma)$ on $S'$.
\end{Definition}

\begin{Example}
\label{exa:eff-div-to-vdiv}
    If $D$ is an effective Cartier divisor on $S$, then $\clO_S(-D)\subset\clO_S$ is a sheaf of locally principal ideals, and we use this to conflate $D$ with the corresponding codimension $1$ closed subscheme of $S$. Dualizing the inclusion $\clO_S(-D)\subset\clO_S$ yields the virtual Cartier divisor $(\clO_S(D),\taut_D)$ on $S$, which was described in Discussion~\ref{disc:sheaf-cartier-div}.  
    
    Conversely, if $(L,\sigma)$ is a virtual Cartier divisor on $S$ such that $\sigma$ is a non-zerodivisor in any local trivialization of $L$, then the scheme-theoretic zero locus $D\subset S$ of $\sigma$ is an effective Cartier divisor whose ideal sheaf is the image of the map $L^\vee\to\clO_S$ induced by $\sigma$. This provides an isomorphism of virtual Cartier divisors $(L,\sigma)\xrightarrow{\sim}(\clO_S(D),\taut_D)$.
\end{Example}

\subsection{Key definitions and properties}

\begin{Definition}\label{def:sys-line-bundle}
    A \emph{system of line bundles} (\emph{slb}) on a scheme $S$ indexed by an integral monoid $Q$ is the data $\clL = \{L_q,\Phi_{q_1,q_2},\sigma_q\}$ consisting of the following.
    \begin{enumerate}[label = (\arabic*), ref = \theDefinition(\arabic*), itemsep = 0.3ex]
        \item For each $q\in Q$, a line bundle $L_q$ on $S$. 
        
        \item\label{def-part:slb-Phi-conditions} For all $q_1,q_2\in Q$, a line bundle isomorphism $\Phi_{q_1,q_2}:L_{q_1}\otimes L_{q_2}\xrightarrow{\sim} L_{q_1 + q_2}$ satisfying the following properties.
        \begin{enumerate}[label = (\roman*), itemsep = 0.15ex]
            \item We have $\Phi_{q_1,q_2} = \Phi_{q_2,q_1}$ under the natural identification $L_{q_1}\otimes L_{q_2} = L_{q_2}\otimes L_{q_1}$.
            
            \item For all $q_1,q_2,q_3\in Q$, the following diagram commutes.
            \begin{equation*}
            \begin{tikzcd}
                L_{q_1}\otimes L_{q_2}\otimes L_{q_3} \arrow[d,"\normalfont\text{id}_{L_{q_1}}\otimes\Phi_{q_2,q_3}","\mathbin{\rotatebox[origin=c]{90}{$\sim$}}"'] \arrow[rrr,"\Phi_{q_1,q_2}\otimes\normalfont\text{id}_{L_{q_3}}","\sim"'] & & & L_{q_1+q_2}\otimes L_{q_3} \arrow[d,"\Phi_{q_1+q_2,q_3}","\mathbin{\rotatebox[origin=c]{90}{$\sim$}}"'] \\
                L_{q_1}\otimes L_{q_2 + q_3} \arrow[rrr,"\Phi_{q_1,q_2+q_3}","\sim"'] & & & L_{q_1 + q_2 + q_3}
            \end{tikzcd}   
            \end{equation*}
        \end{enumerate}
        In particular, we have a trivialization $\Phi_0:L_0\xrightarrow{\sim}\clO_S$ defined by $\Phi_0\otimes\text{id}_{L_0} = \Phi_{0,0}$.
        
        \item For each $q\in Q$, a section $\sigma_q\in H^0(S,L_q)$ satisfying the following properties.
        \begin{enumerate}[label = (\roman*), itemsep = 0.15ex]
            \item The trivialization $\Phi_0$ maps $\sigma_0$ to the element $1\in H^0(S,\clO_S)$.
            
            \item For all $q_1,q_2\in Q$, the isomorphism $\Phi_{q_1,q_2}$ maps $\sigma_{q_1}\otimes\sigma_{q_2}$ to $\sigma_{q_1+q_2}$.
        \end{enumerate}
    \end{enumerate}
\end{Definition}

\begin{Discussion}\label{disc:extending-slb-to-gp}
    In the setting of Definition~\ref{def:sys-line-bundle}, we may naturally extend the definitions of $L_q$ and $\Phi_{q_1,q_2}$ (but not $\sigma_q$) to $q,q_1,q_2\in Q^\text{gp}$ as described below, so that the properties in Definition~\ref{def-part:slb-Phi-conditions} continue to hold.

    \begin{enumerate}[label = (\arabic*), itemsep = 0.3ex] 
        \item Given $q\in Q^\text{gp}$, let $I_q$ be the set of pairs $(a,b)\in Q\times Q$ with $q=a-b$. Define a sheaf $\fL_q$ on $S$ as follows: a local section $t\in\Gamma(U,\fL_q)$ is a tuple of the form
        \begin{align*}
            t = (t_{a,b}\in\Gamma(U,L_a\otimes L_b^\vee))_{(a,b)\,\in\, I_q}    
        \end{align*}
        such that for all $(a,b),(c,d)\in I_q$, the isomorphism $L_a\otimes L_b^\vee\xrightarrow{\sim} L_c\otimes L_d^\vee$ induced by $\Phi_{b,c}^{-1}\otimes\Phi_{a,d}$ sends $t_{a,b}$ to $t_{c,d}$. It follows from Definition~\ref{def-part:slb-Phi-conditions} that the projection $t\mapsto t_{a,b}$ defines an isomorphism $\Pi_{a,b}: \fL_q\xrightarrow{\sim}L_a\otimes L_b^\vee$ for each $(a,b)\in I_q$. 
        
        In particular, $\fL_q$ is a locally free $\clO_S$-module of rank $1$, i.e., a line bundle. Moreover, for each $q\in Q$, the maps $\Pi_{q,0}$ and $\Phi_0$ induce a canonical isomorphism $\fL_q\xrightarrow{\sim}L_q$. In view of this, we will henceforth write $L_q$ instead of $\fL_q$ for all $q\in Q^\text{gp}$.
        
        \item Given $q_1,q_2\in Q^\text{gp}$, define $\Phi_{q_1,q_2}:L_{q_1}\otimes L_{q_2}\xrightarrow{\sim} L_{q_1+q_2}$ via the following commutative diagram after choosing $(a_i,b_i)\in I_{q_i}$ for $i=1,2$. It follows from Definition~\ref{def-part:slb-Phi-conditions} that this is well-defined independent of the choices of $(a_i,b_i)$. As before, one checks that the new definition of $\Phi_{q_1,q_2}$ agrees with the old one if $q_1,q_2\in Q$.
        \begin{equation*}
        \begin{tikzcd}
            L_{q_1}\otimes L_{q_2} \arrow[rrr,"\Phi_{q_1,q_2}","\sim"'] \arrow[d,"\Pi_{a_1,b_1}\otimes\Pi_{a_2,b_2}","\mathbin{\rotatebox[origin=c]{90}{$\sim$}}"'] &&& L_{q_1+q_2} \arrow[d,"\Pi_{a_1+a_2,b_1+b_2}","\mathbin{\rotatebox[origin=c]{90}{$\sim$}}"'] \\
            L_{a_1}\otimes L_{b_1}^\vee\otimes L_{a_2}\otimes L_{b_2}^\vee \arrow[rrr,"\Phi_{a_1,a_2}\otimes(\Phi_{b_1,b_2}^\vee)^{-1}","\sim"'] &&& L_{a_1+a_2}\otimes L_{b_1+b_2}^\vee
        \end{tikzcd}
        \end{equation*}
    \end{enumerate}
\end{Discussion}

\begin{Definition}
    Let $\clL = \{L_q,\Phi_{q_1,q_2},\sigma_q\}$ and $\hat\clL = \{\hat L_q,\hat\Phi_{q_1,q_2},\hat\sigma_q\}$ be two systems of line bundles on a scheme $S$, indexed by the integral monoid $Q$. 
    
    An \emph{slb isomorphism} $\Psi:\clL\xrightarrow{\sim}\hat\clL$ is the data $\Psi = \{\Psi_q\}$ of a line bundle isomorphism $\Psi_q:L_q\xrightarrow{\sim}\hat L_q$, for each $q\in Q$, satisfying the following properties.
    \begin{enumerate}[label = (\arabic*), itemsep = 0.3ex]
        \item For all $q_1,q_2\in Q$, the following diagram commutes.
        \begin{equation*}
        \begin{tikzcd}
            L_{q_1}\otimes L_{q_2} \arrow[rr,"\Psi_{q_1}\otimes\Psi_{q_2}","\sim"'] \arrow[d,"\Phi_{q_1,q_2}","\mathbin{\rotatebox[origin=c]{90}{$\sim$}}"'] & & \hat L_{q_1}\otimes\hat L_{q_2} \arrow[d,"\hat\Phi_{q_1,q_2}","\mathbin{\rotatebox[origin=c]{90}{$\sim$}}"']\\
            L_{q_1 + q_2} \arrow[rr,"\Psi_{q_1 + q_2}","\sim"'] & & \hat L_{q_1 + q_2}
        \end{tikzcd}
        \end{equation*}
        
        \item For all $q\in Q$, the isomorphism $\Psi_q$ maps $\sigma_q$ to $\hat\sigma_q$.
    \end{enumerate}
\end{Definition}

\begin{Example}
\label{exa:df-rank-r-sys-line-bundle}
    Consider $Q = \bN^r$. An slb on $S$ indexed by $\bN^r$ is specified uniquely up to unique isomorphism by a collection of virtual Cartier divisors $(L_i,\sigma_i)$ on $S$, for $1\le i\le r$. Explicitly, for $(k_1,\ldots,k_r)\in\bN^r$, we take
    \begin{align*}
        \begin{matrix}
            L_{(k_1,\ldots,k_r)} := L_1^{\otimes k_1}\otimes\cdots\otimes L_r^{\otimes k_r}
            \quad\text{and}\quad
            \sigma_{(k_1,\ldots,k_r)} := \sigma_1^{\otimes k_1}\otimes\cdots\otimes\sigma_r^{\otimes k_r},
        \end{matrix}
    \end{align*} 
    while the isomorphisms $\Phi_{q_1,q_2}$ are taken to be the natural identifications. We refer to this as the slb \emph{freely generated} by the virtual Cartier divisors $(L_i,\sigma_i)$ for $1\le i\le r$.
\end{Example}

\begin{Example}\label{exa:n-free-cart-div}
    Combining Examples~\ref{exa:eff-div-to-vdiv} and \ref{exa:df-rank-r-sys-line-bundle}, any collection of $r$ effective Cartier divisors $D_1,\ldots,D_r\subset S$ determines an slb on $S$ indexed by $Q = \bN^r$.
\end{Example}

\begin{Example}\label{exa:trivial-slb}
    Let $Q$ be an integral monoid and $S$ be a scheme, with $\clO_S$ viewed as a sheaf of monoids under multiplication. Consider a map of monoids 
    \begin{align*}
        \tau:Q\to \Gamma(S,\clO_S).    
    \end{align*}
    This data gives rise to an slb $\{L_q,\Phi_{q_1,q_2},\sigma_q\}$ on $S$ indexed by $Q$, called the \emph{trivial slb} associated to $\tau$, which is defined as follows. For all $q\in Q$, we take $L_q := \clO_S$ and $\sigma_q := \tau(q)$ and, for all $q_1,q_2\in Q$, we take the isomorphisms $\Phi_{q_1,q_2}$ to be the identity maps.
\end{Example}

\begin{Example}\label{exa:toric-slb}
    Consider a fine monoid $Q$. We then obtain an slb on $V_Q = \Spec\bC[Q]$ indexed by $Q$ by applying Example~\ref{exa:trivial-slb} to the natural inclusion $Q\hookrightarrow\bC[Q]=\Gamma(V_Q,\clO_{V_Q})$.
\end{Example}

\begin{Definition}
    Let $\clL=\{L_q,\Phi_{q_1,q_2},\sigma_q\}$ be an slb on $S$ indexed by $Q$.
    \begin{enumerate}[label = (\arabic*), itemsep = 0.3ex]
        \item If $\varphi:Q'\to Q$ is a homomorphism of integral monoids, then the \emph{pullback of $\clL$ along $\varphi$} is the slb on $S$ indexed by $Q'$ defined as
        \begin{align*}
            \varphi^*\clL := \{L_{\varphi(q')},\Phi_{\varphi(q'_1),\varphi(q'_2)},\sigma_{\varphi(q')}\}.
        \end{align*}
        
        \item If $T$ is another scheme and $f:T\to S$ is a morphism of schemes, then the \emph{pullback of $\clL$ along $f$} is the slb on $T$ indexed by $Q$ defined as
        \begin{align*}
            f\phantom{}^*\clL := \{f\phantom{}^*L_q,f\phantom{}^*\Phi_{q_1,q_2},f\phantom{}^*\sigma_q\}.
        \end{align*}
        We sometimes write $\clL|_{T}$ instead of $f^*\clL$ provided $f$ is clear from the context.
    \end{enumerate}
\end{Definition}

\begin{Lemma}
\label{lem:slb-local-triv}
    Let $\clL$ be an slb on $S$ indexed by a fine monoid $Q$.
    
    Then any point $s\in S$ has an \'etale neighborhood $U$ such that $\clL|_U$ is isomorphic to the trivial slb (Example~\ref{exa:trivial-slb}) associated to some monoid map $\tau:Q\to\Gamma(U,\clO_U)$.
\end{Lemma}
\begin{proof}
    Write $\clL = \{L_q,\Phi_{q_1,q_2},\sigma_q\}$ and recall from Discussion~\ref{disc:extending-slb-to-gp} that we can make sense of $L_q$ and $\Phi_{q_1,q_2}$ for $q,q_1,q_2\in Q^\text{gp}$ as well. Use the structure theorem for finitely generated abelian groups to choose a direct sum decomposition
    \begin{align}\label{eqn:str-thm-ab-gp-iso-2}
        Q^\text{gp}\xrightarrow{\sim} \textstyle\bigoplus_{i}\bZ\oplus\bigoplus_{j}(\bZ/d_j\bZ).
    \end{align}
    For each $i$, let $g_i\in Q^\text{gp}$ be a generator of the corresponding $\bZ$ summand and, for each $j$, let $h_j\in Q^\text{gp}$ be a generator of the corresponding $\bZ/d_j\bZ$ summand in \eqref{eqn:str-thm-ab-gp-iso-2}. Since we have the relation $d_j\cdot h_j = 0$ in $Q^\text{gp}$, iterating the isomorphisms $\Phi_{q_1,q_2}$ and using the trivialization $\Phi_0$ of $L_0$ from Definition~\ref{def:sys-line-bundle} yields a global trivialization of $L_{h_j}^{\otimes d_j}$ which we denote as
    \begin{align*}
        \Phi_{(d_j,h_j)}:L_{h_j}^{\otimes d_j}\xrightarrow{\sim}\clO_S
    \end{align*}
    for each $j$. Given $s\in S$, pick a Zariski open neighborhood $U$ of $s$ and trivializations 
    \begin{align}\label{eqn:slb-local-triv-generators}
        \Psi_{g_i}:L_{g_i}|_U\xrightarrow{\sim}\clO_U \quad\text{and}\quad \Psi_{h_j}:L_{h_j}|_U\xrightarrow{\sim}\clO_U
    \end{align}
    for each $i,j$. Define $u_j\in\Gamma(U,\clO_U^\times)$ for each $j$ by
    \begin{align*}
        \Phi_{(d_j,h_j)}|_U = u_j\cdot\Psi_{h_j}^{\otimes d_j}.
    \end{align*}
    Replacing $U$ by an \'etale neighborhood of $s$, find $v_j\in\Gamma(U,\clO_U^\times)$ such that $v_j^{d_j} = u_j$. Replace $\Psi_{h_j}$ by $v_j\cdot\Psi_{h_j}$ so that we now have $u_j = 1$. From \eqref{eqn:str-thm-ab-gp-iso-2}, it is immediate that the new collection \eqref{eqn:slb-local-triv-generators} extends uniquely to a collection of trivializations $\Psi_q:L_q|_U\xrightarrow{\sim}\clO_U$, for all $q\in Q^\text{gp}$, which are compatible with the isomorphisms $\Phi_{q_1,q_2}|_U$ for all $q_1,q_2\in Q^\text{gp}$.
    
    For $q\in Q$, let $\tau(q)\in\Gamma(U,\clO_U) = \Hom(\clO_U,\clO_U)$ be the image of the section $\sigma_q|_U$ under the trivialization $\Psi_q$. One checks that $\tau:Q\to\Gamma(U,\clO_U)$ is a monoid map and that the collection $\Psi = \{\Psi_q\}$ defines an isomorphism from $\clL|_U$ to the trivial slb associated to $\tau$. 
\end{proof}

\subsection{Slb presentations}\label{subsec:slb-presentation}

In this subsection, we discuss a generalization of Example~\ref{exa:df-rank-r-sys-line-bundle} which allows us to construct systems of line bundles by specifying finitely many virtual Cartier divisors (generators) and finitely many line bundle isomorphisms among them (relations). This will be used in Section~\ref{sec:log-maps} to compare algebraic and symplectic log maps.

\begin{Definition}\label{def:slb-presentation}
    Fix a scheme $S$ and a fine monoid $Q$. Let
    \begin{equation}\label{eqn:slb-monoid-presentation}
    \begin{tikzcd}
        \bN^J \arrow[r,shift left=.5ex,"a"] \arrow[r,shift right=.5ex,swap,"b"] & \bN^I \arrow [r,"\pi"] & Q
    \end{tikzcd}
    \end{equation}
    be an integral presentation of $Q$, as in Definition~\ref{def-part:integral-presentation}. 

    Let $(a_{ij})$, $(b_{ij})$ be the matrices of the maps $a,b:\bN^J\rrarrows\bN^I$ respectively with respect to the standard bases $\{{\bf e}_i\}$ of $\bN^I$ and $\{{\bf e}_j'\}$ of $\bN^J$. For each $i\in I$, define $\gamma_i:=\pi({\bf e}_i)\in Q$.
    \begin{enumerate}[ref = \theDefinition(\arabic*), label = (\arabic*), itemsep = 0.3ex]
        \item\label{def-part:slb-presentation-defined} An \emph{slb presentation} on $S$, with respect to the integral presentation \eqref{eqn:slb-monoid-presentation} of $Q$, is the data $\fP = \{\fL_i,\fs_i,\varphi_j\}$ consisting of the following.
        
        \begin{enumerate}[label = (\roman*), itemsep = 0.15ex]
            \item For each $i\in I$, a virtual Cartier divisor $(\fL_i,\fs_i)$ on $S$.
        
            \item For each $j\in J$, a line bundle isomorphism 
            \begin{align*}
                \varphi_j:\textstyle\bigotimes_i \fL_i^{\otimes a_{ij}}\xrightarrow{\sim}\bigotimes_i \fL_i^{\otimes b_{ij}}.
            \end{align*}
        \end{enumerate}

        \item\label{def-part:slb-presentation-realization-category} Given an slb presentation $\fP = \{\fL_i,\fs_i,\varphi_j\}$ on $S$, with respect to \eqref{eqn:slb-monoid-presentation}, we define its \emph{category of realizations}, denoted by $\|\fP\|$, as follows.

        \begin{enumerate}[label = (\roman*), itemsep = 0.15ex]
            \item Objects of $\|\fP\|$ are tuples $(P,\psi,\clL,\{\Psi_i\})$, where
            \begin{enumerate}[label = (\alph*), itemsep = 0.15ex]
                \item $P$ is a fine monoid,
                
                \item $\psi:Q\to P$ is a map of monoids,
                
                \item $\clL = \{L_p,\Phi_{p_1,p_2},\sigma_p\}$ is an slb on $S$ indexed by $P$, and

                \item for each $i\in I$,
                \begin{align*}
                    \Psi_i:(\fL_i,\fs_i)\xrightarrow{\sim}(L_{\psi(\gamma_i)},\sigma_{\psi(\gamma_i)})  
                \end{align*}
                is an isomorphism of virtual Cartier divisors on $S$
            \end{enumerate}
            such that, for each $j\in J$, the isomorphism $\bigotimes_i L_{\psi(\gamma_i)}^{\otimes a_{ij}}\xrightarrow{\sim}\bigotimes_i L_{\psi(\gamma_i)}^{\otimes b_{ij}}$ induced by the maps $\Phi_{p_1,p_2}$ and the identity $\sum_i a_{ij}\psi(\gamma_i) = \sum_i b_{ij}\psi(\gamma_i)\in P$ agrees with
            \begin{align*}
                \textstyle
                \left(\bigotimes_i\Psi_i^{\otimes b_{ij}}\right)\circ\varphi_j\circ\left(\bigotimes_i\Psi_i^{\otimes a_{ij}}\right)^{-1}:\bigotimes_i L_{\psi(\gamma_i)}^{\otimes a_{ij}}\xrightarrow{\sim}\bigotimes_i L_{\psi(\gamma_i)}^{\otimes b_{ij}}.
            \end{align*}

            \item Morphisms $(P,\psi,\clL,\{\Psi_i\})\to(P',\psi',\clL',\{\Psi'_i\})$ in $\|\fP\|$ are pairs $(\omega,\Omega)$, where
            \begin{enumerate}[label = (\alph*), itemsep = 0.15ex]
                \item $\omega:P\to P'$ is a monoid map, and
                
                \item $\Omega = \{\Omega_p\}:\clL\xrightarrow{\sim}\omega^*\clL'$ is an slb isomorphism
            \end{enumerate}
            such that $\psi' = \omega\circ\psi$ and, for each $i\in I$, we have $\Psi_i' = \Omega_{\psi(\gamma_i)}\circ\Psi_i$.
        \end{enumerate}

        \item\label{def-part:consistent-slb-presentation} An slb presentation $\fP = \{\fL_i,\fs_i,\varphi_j\}$ on $S$, with respect to \eqref{eqn:slb-monoid-presentation}, is said to be \emph{consistent} if, for each point $s\in S$, there exist an \'etale neighborhood $U$ of $s$ and trivializations $\chi_i:\fL_i|_U\xrightarrow{\sim}\clO_U$ for each $i\in I$ which satisfy the following conditions.
        \begin{enumerate}[ref = \theenumi(\roman*), label = (\roman*), itemsep = 0.3ex]
            \item\label{def-subpart:relation-compatible-trivialization} For each $j\in J$, we have $\textstyle\bigotimes_i\chi_i^{\otimes a_{ij}} = \bigotimes_i\chi_i^{\otimes b_{ij}}\circ\varphi_j|_U$.
        
            \item\label{def-subpart:section-compatible-trivialization} For all ${\bf x},{\bf y}\in\bN^I$ with $\pi({\bf x})=\pi({\bf y})\in Q$, we have $\textstyle\prod_i(\chi_i(\fs_i|_U))^{x_i} = \prod_i(\chi_i(\fs_i|_U))^{y_i}$ as elements of $\Gamma(U,\clO_U)$, where $x_i,y_i$ are the coordinates of ${\bf x},{\bf y}$ respectively.
    \end{enumerate} 
    \end{enumerate}
\end{Definition}

\begin{Notation}
    In the setting of Definition~\ref{def:slb-presentation}, we use the following abbreviations. For ${\bf x}\in\bZ^I$, we write $\fL^{\otimes\bf x} := \bigotimes_i\fL_i^{\otimes x_i}$, $L_\psi^{\otimes{\bf x}} := \bigotimes_i L_{\psi(\gamma_i)}^{\otimes x_i}$ and $\Psi^{\otimes{\bf x}} := \bigotimes_i\Psi_i^{\otimes x_i}$, where $x_i$ are the coordinates of ${\bf x}$. If ${\bf x}\in\bN^I$, we also write $\fs^{\otimes{\bf x}} := \bigotimes_i\fs_i^{\otimes x_i}$. Similarly, for ${\bf z}\in\bZ^J$, we write $\varphi^{\otimes\bf z} := \bigotimes_j\varphi_j^{\otimes z_j}$, where $z_j$ are the coordinates of ${\bf z}$.
\end{Notation}

\begin{Proposition}\label{prop:slb-presentation}
    Consider an slb presentation $\fP = \{\fL_i,\fs_i,\varphi_j\}$ on $S$, with respect to \eqref{eqn:slb-monoid-presentation}, as in Definition~\ref{def:slb-presentation}. Then the following statements are equivalent.
    \begin{enumerate}[label = (\arabic*), itemsep = 0.3ex]
        \item The category $\|\fP\|$ has an initial object.

        \item The category $\|\fP\|$ is non-empty.

        \item The slb presentation $\fP$ is consistent.
    \end{enumerate}
    
    Moreover, if any (and therefore each) of these three statements is true, then an object $(P,\psi,\clL,\{\Psi_i\})$ of $\|\fP\|$ is an initial object if and only if $\psi:Q\to P$ is an isomorphism.
\end{Proposition}
\begin{proof}
    (1) $\Rightarrow$ (2): This is clear.

    (2) $\Rightarrow$ (3): Choose an object $(P,\psi,\clL,\{\Psi_i\})$ of $\|\fP\|$. Given a point $s\in S$, Lemma~\ref{lem:slb-local-triv} yields an \'etale neighborhood $U$ of $s$, a monoid map $\tau:P\to\Gamma(U,\clO_U)$ and an slb isomorphism from $\clL|_U$ to the trivial slb associated to $\tau$. This slb isomorphism, after composition with $\Psi_i|_U$, induces trivializations $\chi_i:\fL_i|_U\to\clO_U$ for each $i\in I$. One checks that these trivializations satisfy the conditions in Definition~\ref{def-part:consistent-slb-presentation}.

    (3) $\Rightarrow$ (1): By assumption, $\fP$ is consistent. This has the following consequence.

    \begin{claim}\label{claim:consistent-slb-presentation}
        If ${\bf x},{\bf y}\in\bN^I$ and ${\bf z}\in\bZ^J$ satisfy ${\bf x} - {\bf y} = (a^\text{gp}-b^\text{gp})({\bf z})$, then
        \begin{align*}
            \varphi^{\otimes{\bf z}}\in \Hom(\fL^{\otimes\, a^\text{gp}({\bf z})},\fL^{\otimes\, b^\text{gp}({\bf z})})\xrightarrow{\sim}\Hom(\fL^{\otimes{\bf x}},\fL^{\otimes{\bf y}})
        \end{align*}
        maps $\fs^{\otimes{\bf x}}$ to $\fs^{\otimes{\bf y}}$. In particular, if $(a^\text{gp}-b^\text{gp})({\bf z}) = {\bf 0}$, then $\varphi^{\otimes{\bf z}}$ is the identity map.
    \end{claim}
    \begin{proof}[Proof of claim]
        Both assertions can be checked \'etale locally on $S$. Given a point $s\in S$, choose an \'etale neighborhood $U$ of $s$ and trivializations $\chi_i:\fL_i|_U\xrightarrow{\sim}\clO_U$ as in Definition~\ref{def-part:consistent-slb-presentation}. Under these trivializations, $\varphi^{\otimes{\bf z}}$ corresponds to the identity map on $\clO_U$ and the two sections $\fs^{\otimes{\bf x}}$ and $\fs^{\otimes{\bf y}}$ correspond to the same element of $\Gamma(U,\clO_U)$. This proves the first assertion. To deduce the second assertion from the first, take ${\bf x} = {\bf y} = {\bf 0}\in\bN^I$ and note that $\fs^{\otimes{\bf 0}} = 1$.
    \end{proof}

    The remainder of the argument is divided into four steps.

    \begin{step}\label{proof-step:slb-presentation-prelim}
        Given $q\in Q$ and ${\bf x},{\bf y}\in\pi^{-1}(q)\subset\bN^I$, define $\varphi_{{\bf x},{\bf y}}:\fL^{\otimes {\bf x}}\xrightarrow{\sim}\fL^{\otimes{\bf y}}$ by the formula $\varphi_{{\bf x},{\bf y}} := \varphi^{\otimes{\bf z}}$ for any ${\bf z}\in\bZ^J$ satisfying ${\bf x}-{\bf y}=(a^\text{gp}-b^\text{gp})({\bf z})$. This is well-defined since, by the second assertion in Claim~\ref{claim:consistent-slb-presentation}, the element $\varphi_{{\bf x},{\bf y}}\in\Hom(\fL^{\otimes{\bf x}},\fL^{\otimes{\bf y}})=\Gamma(X,\fL^{\otimes({\bf y}-{\bf x})})$ is independent of the choice of ${\bf z}$ and depends only on ${\bf x}-{\bf y}$. By the same reasoning, we have the compatibility relation $\varphi_{{\bf x},{\bf y}}\circ\varphi_{{\bf w},{\bf x}} = \varphi_{{\bf w},{\bf y}}$ whenever ${\bf w},{\bf x},{\bf y}\in\pi^{-1}(q)$. Moreover, the first assertion in Claim~\ref{claim:consistent-slb-presentation} shows that $\varphi_{{\bf x},{\bf y}}$ maps $\fs^{\otimes{\bf x}}$ to $\fs^{\otimes{\bf y}}$ whenever ${\bf x},{\bf y}\in\pi^{-1}(q)$. 
    \end{step}

    \begin{step}\label{proof-step:slb-presentation-init-obj-def}
        In this step, we construct an object of $\|\fP\|$ of the form $(Q,\text{id}_Q,\hat\clL,\{\hat\Psi_i\})$ by an argument which is very similar to Discussion~\ref{disc:extending-slb-to-gp}.
        \begin{enumerate}[label = (\alph*), itemsep = 0.3ex]
            \item For $q\in Q$, define the sheaf $\hat L_q$ on $S$ as follows: a local section $t\in\Gamma(U,\hat L_q)$ is a tuple of the form $t=(t_{\bf x}\in\Gamma(U,\fL^{\otimes{\bf x}}))_{{\bf x}\,\in\,\pi^{-1}(q)}$ with $\varphi_{{\bf x},{\bf y}}(t_{\bf x}) = t_{\bf y}$ for all ${\bf x},{\bf y}\in\pi^{-1}(q)$. 
        
            \item By Step~\ref{proof-step:slb-presentation-prelim}, the projection $t\mapsto t_{\bf x}$ defines an isomorphism $\Pi_{{\bf x},q}:\hat L_q\xrightarrow{\sim}\fL^{\otimes{\bf x}}$ for each ${\bf x}\in\pi^{-1}(q)$. Moreover, if ${\bf y}\in\pi^{-1}(q)$ is another element, then $\varphi_{{\bf x},{\bf y}} = \Pi_{{\bf y},q}\circ\Pi_{{\bf x},q}^{-1}$. In particular, $\hat L_q$ is a locally free $\clO_S$-module of rank $1$, i.e., a line bundle.
        
            \item For $q_1,q_2\in Q$, define the isomorphism $\hat\Phi_{q_1,q_2}:\hat L_{q_1}\otimes\hat L_{q_2}\xrightarrow{\sim}\hat L_{q_1+q_2}$ via the following commutative diagram, after choosing ${\bf x}_k\in\pi^{-1}(q_k)$ for $k=1,2$. By Step~\ref{proof-step:slb-presentation-prelim}, the isomorphism $\Phi_{q_1,q_2}$ is independent of the choices of ${\bf x}_1,{\bf x}_2$.
            \begin{equation*}
            \begin{tikzcd}
                \hat L_{q_1}\otimes\hat L_{q_2} \arrow[d,"\Pi_{{\bf x}_1,q_1}\otimes \Pi_{{\bf x}_2,q_2}","\mathbin{\rotatebox[origin=c]{90}{$\sim$}}"'] \arrow[rr,"\hat\Phi_{q_1,q_2}","\sim"'] & & \hat L_{q_1+q_2}\arrow[d,"\Pi_{{\bf x},q}","\mathbin{\rotatebox[origin=c]{90}{$\sim$}}"'] \\ 
                \fL^{\otimes{\bf x}_1}\otimes\fL^{\otimes{\bf x}_2} \arrow[rr,equal] & & \fL^{\otimes({\bf x}_1+{\bf x}_2)}
            \end{tikzcd}
            \end{equation*}
        
            \item For $q\in Q$, define the section $\hat\sigma_q$ of $\hat L_q$ to be the image of $\fs^{\otimes{\bf x}}$ under $\Pi_{{\bf x},q}^{-1}$, after choosing ${\bf x}\in\pi^{-1}(q)\subset\bN^I$. By Step~\ref{proof-step:slb-presentation-prelim}, the section $\hat\sigma_q$ is independent of the choice of ${\bf x}$.

            \item For $i\in I$, define the isomorphism $\hat\Psi_i:\fL_i\xrightarrow{\sim}\hat L_{\gamma_i}$ by $\hat\Psi_i := \Pi_{{\bf e}_i,\gamma_i}^{-1}$, where we recall that ${\bf e}_i\in\bN^I$ is a standard basis vector and $\gamma_i = \pi({\bf e}_i)$.
        \end{enumerate}
        
        Write $\hat\clL:=\{\hat L_q,\hat\Phi_{q_1,q_2},\hat\sigma_q\}$. One checks that $(Q,\text{id}_Q,\hat\clL,\{\hat\Psi_i\})$ is indeed an object of $\|\fP\|$.
    \end{step}

    \begin{step}\label{proof-step:slb-presentation-init-obj-check}
        In this step, we show that $(Q,\text{id}_Q,\hat\clL,\{\hat\Psi_i\})$ from Step~\ref{proof-step:slb-presentation-init-obj-def} is an initial object of $\|\fP\|$. 
        
        Let $(P,\psi,\clL,\{\Psi_i\})$ be a given object of $\|\fP\|$ and write $\clL = \{L_p,\Phi_{p_1,p_2},\sigma_p\}$. We will show that there is exactly one morphism $(\omega,\Omega):(Q,\text{id}_Q,\hat\clL,\{\hat\Psi_i\})\to (P,\psi,\clL,\{\Psi_i\})$
        in $\fP$.

        \noindent\emph{Uniqueness}. 
        
        \noindent Since $\psi = \omega\circ\text{id}_Q$, we must have $\omega = \psi$. For each $i\in I$, since $\Psi_i = \Omega_{\gamma_i}\circ\hat\Psi_i$, we must have $\Omega_{\gamma_i} = \Psi_i\circ\hat\Psi_i^{-1}$. Since the elements $\gamma_i$ generate $Q$, it follows that $\Omega_q$ is uniquely determined for all $q\in Q$. Explicitly, given $q\in Q$ and ${\bf x}\in\pi^{-1}(q)\subset\bN^I$, the isomorphism $\Omega_q:\hat L_q\xrightarrow{\sim}L_{\psi(q)}$ fits into the following commutative diagram, where the right vertical isomorphism is induced by the maps $\Phi_{p_1,p_2}$ and the identity $\psi(q) = \sum_i x_i\psi(\gamma_i)\in P$.
        
        \begin{equation}\label{eqn:slb-init-obj-map}
        \begin{tikzcd}
            \hat L_q \arrow[d,"\Pi_{{\bf x},q}","\mathbin{\rotatebox[origin=c]{90}{$\sim$}}"'] \arrow[rrr,"\Omega_q","\sim"'] & & & L_{\psi(q)} \arrow[d,"\mathbin{\rotatebox[origin=c]{90}{$\sim$}}"'] \\
            \fL^{\otimes {\bf x}} \arrow[rrr,"\Psi^{\otimes{\bf x}}","\sim"'] & & & L_\psi^{\otimes{\bf x}}
        \end{tikzcd}
        \end{equation}

        \noindent\emph{Existence}. 
        
        \noindent Define $\omega:=\psi$. For each $q\in Q$, define the isomorphism $\Omega_q:\hat L_q\xrightarrow{\sim}L_{\psi(q)}$ by the commutative diagram \eqref{eqn:slb-init-obj-map}. For this to be well-defined, we must show that, for any ${\bf x},{\bf y}\in\pi^{-1}(q)$, the following diagram commutes, where the right vertical isomorphism is induced by the maps $\Phi_{p_1,p_2}$ and the identity $\sum_i x_i\psi(\gamma_i) = \sum_i y_i\psi(\gamma_i)\in P$.

        \begin{equation}\label{eqn:slb-init-obj-well-defined}
        \begin{tikzcd}
            \fL^{\otimes{\bf x}} \arrow[d,"\Pi_{{\bf y},q}\,\circ\,\Pi_{{\bf x},q}^{-1}","\mathbin{\rotatebox[origin=c]{90}{$\sim$}}"'] \arrow[rrr,"\Psi^{\otimes {\bf x}}","\sim"'] & & &  L_\psi^{\otimes{\bf x}} \arrow[d,"\mathbin{\rotatebox[origin=c]{90}{$\sim$}}"'] \\
            \fL^{\otimes{\bf y}} \arrow[rrr,"\Psi^{\otimes {\bf y}}","\sim"'] & & & L_\psi^{\otimes{\bf y}}
        \end{tikzcd}
        \end{equation}
        To see this, we argue as follows. From Step~\ref{proof-step:slb-presentation-init-obj-def}, we have $\Pi_{{\bf y},q}\circ\Pi_{{\bf x},q}^{-1} = \varphi_{{\bf x},{\bf y}} = \varphi^{\otimes{\bf z}}$ for any choice of ${\bf z}\in\bZ^J$ such that ${\bf x} - {\bf y} = (a^\text{gp} - b^\text{gp})({\bf z})$. By writing ${\bf z}$ as a $\bZ$-linear combination of standard basis vectors ${\bf e}'_j$, it suffices to check the commutativity of \eqref{eqn:slb-init-obj-well-defined} in the case where ${\bf z} = {\bf e}'_j$ is a standard basis vector. In this special case, Definition~\ref{def-part:slb-presentation-realization-category} shows that \eqref{eqn:slb-init-obj-well-defined} commutes. One checks that $(\omega,\Omega)$ is indeed a morphism in $\|\fP\|$.  
    \end{step}

    \begin{step}
        In this final step, we characterize initial objects of $\|\fP\|$. 
        
        Note that an object $(P,\psi,\clL,\{\Psi_i\})$ is an initial object if and only if the unique morphism $(\omega,\Omega):(Q,\text{id}_Q,\hat\clL,\{\hat\Psi_i\})\to(P,\psi,\clL,\{\Psi_i\})$ is an isomorphism. It is clear that $(\omega,\Omega)$ is an isomorphism if and only if $\omega$ is an isomorphism. Moreover, from Step~\ref{proof-step:slb-presentation-init-obj-check}, have $\omega = \psi$.
    \end{step}
    Thus, $(P,\psi,\clL,\{\Psi_i\})$ is an initial object of $\|\fP\|$ if and only if $\psi$ is an isomorphism.
\end{proof}
    \section{Log geometry via systems of line bundles}\label{sec:log-via-slb}

We explain how to use systems of line bundles to describe \emph{Deligne--Faltings log structures} on a scheme as well as log morphisms whose target is such a log scheme; an essentially equivalent discussion appears in \cite[Section~3]{borne-vistoli-DF}. As an illustration, we describe two important classes of log schemes from this viewpoint: \emph{log points} and \emph{smooth snc pairs}.

\begin{Notation}
    Throughout this section, we use $S$, $S'$, $T$ etc to denote schemes belonging to $\Sch$ and we use $Q,Q'$ etc to denote sharp fine monoids; see Definition~\ref{def:monoid-props}. All log structures $M$ that we consider are assumed to be fine; see Definition~\ref{def:chart-fine-fs}. Recall from Convention~\ref{conv:log-str-mult-ghost-add} that we write the operation on $M$ (resp. $\ol M$) multiplicatively (resp. additively).
\end{Notation}

\subsection{Global system of line bundles on a log scheme}\label{subsec:ghost-to-slb}

\begin{Definition}\label{def:log-lb-vcd}
    Let $S^\dagger = (S,M)$ be a log scheme in $\LogSch$, i.e., $S$ is locally of finite type over $\bC$ and $M$ is a fine log structure. Denote the structure map by $\alpha:M\to\clO_S$.
    \begin{enumerate}[label = (\arabic*), ref = \theDefinition(\arabic*), itemsep = 0.3ex]
        \item Given a global section $\normalfont m\in\Gamma(S,\ol M\phantom{}^\text{gp})$, Lemma~\ref{lem:fine-log-str-integral} shows that the inverse image of $m$ under $\normalfont M^\text{gp}\to\ol M\phantom{}^\text{gp}$ is an $\clO_S^\times$-torsor, which we denote by $\clO_{S^\dagger}(-m)^\times\subset\normalfont M^\text{gp}$.
        
        The \emph{line bundle $\clO_{S^\dagger}(m)$ associated to} the global section $m\in\Gamma(S,\ol M\phantom{}^\text{gp})$ is defined to be the dual of the line bundle determined by $\clO_{S^\dagger}(-m)^\times$ via \eqref{eqn:line-bundle-from-torsor}.

        \item For any $m_1,m_2\in \Gamma(S,\ol M\phantom{}^\text{gp})$, the group operation on $M^\text{gp}$ yields a multiplication map 
        \begin{align*}
            \clO_{S^\dagger}(-m_1)^\times\times \clO_{S^\dagger}(-m_2)^\times\to \clO_{S^\dagger}(-(m_1+m_2))^\times   
        \end{align*}
        and this induces a line bundle isomorphism
        \begin{align*}
            \Phi_{S^\dagger,m_1,m_2}: \clO_{S^\dagger}(m_1)\otimes \clO_{S^\dagger}(m_2)\xrightarrow{\sim} \clO_{S^\dagger}(m_1+m_2).
        \end{align*}
        
        \item If $m$ is a section of the subsheaf $\ol M\subset\ol M\phantom{}^\text{gp}$, then we have $\clO_{S^\dagger}(-m)^\times\subset M$. The restriction of $\alpha:M\to\clO_S$ induces a global section $\sigma_{S^\dagger,m}\in H^0(S,\clO_{S^\dagger}(m))$.
        
        The \emph{virtual Cartier divisor associated to} the global section $m\in\Gamma(S,\ol M)$ is defined to be the pair $(\clO_{S^\dagger}(m),\sigma_{S^\dagger,m})$.

        \item\label{def-part:global-slb} The \emph{global system of line bundles} on the log scheme $S^\dagger$ is defined to be the slb on $S$ indexed by the sharp integral monoid $\Gamma(S,\normalfont\ol M)$, consisting of the line bundles $\clO_{S^\dagger}(m)$, the isomorphisms $\Phi_{S^\dagger,m_1,m_2}$ and the global sections $\sigma_{S^\dagger,m}$.
        One verifies that $\clO_{S^\dagger}(m)$, $\Phi_{S^\dagger,m_1,m_2}$ and $\sigma_{S^\dagger,m}$ satisfy the conditions of Definition~\ref{def:sys-line-bundle}.
    \end{enumerate}
\end{Definition}

The global slb on a log scheme is functorial in the following sense.

\begin{Lemma}\label{lem:map-of-log-str-slb}
    Let $f^\dagger = (f,f^\flat):S'^\dagger\to S^\dagger$ be a log morphism between two log schemes $S'^\dagger = (S',M')$ and $S^\dagger = (S,M)$ in $\LogSch$. 
    
    Then the map $f^\flat:f^*M\to M'$ induces an isomorphism between the following two systems of line bundles on $S'$ indexed by the partially ordered abelian group $\Gamma(S,\ol M)$.
    \begin{enumerate}[label = (\arabic*), itemsep = 0.3ex]
        \item The pullback of the global slb on the log scheme $S^\dagger$ along the morphism $f:S'\to S$.
        
        \item The pullback of the global slb on the log scheme $S'^\dagger$ along the map of sharp integral monoids $\Gamma(S,\ol M)\to\Gamma(S',\ol M\phantom{}')$ induced by $f^\flat$.
    \end{enumerate}
\end{Lemma}
\begin{proof}
    This is a formal consequence of the pullback diagrams \eqref{eqn:map-of-log-str-cartesian}.
\end{proof}

\begin{Discussion}\label{disc:log-lb-vcd-triv-dual}
    Consider a log scheme $S^\dagger = (S,M)$ in $\LogSch$. We make the following series of elementary observations about its global slb.
    \begin{enumerate}[label = (\arabic*), ref = \theDiscussion(\arabic*), itemsep = 0.3ex]
        \item\label{disc-part:log-lb-vcd-triv} The line bundle $\clO_{S^\dagger}(0)$ is trivialized by the section $\sigma_{S^\dagger,0}$. Equivalently, the $\clO_S^\times$-torsor $\clO_{S^\dagger}(0)^\times = \alpha^{-1}(\clO_S^\times)\subset M$ is trivialized by the structure map $\alpha$. This corresponds to the trivialization $\Phi_0$ appearing in the definition of an slb (Definition~\ref{def:sys-line-bundle}).
        
        \item\label{disc-part:log-lb-vcd-dual} For a section $m\in\Gamma(S,\ol M\phantom{}^\text{gp})$, the line bundles $\clO_{S^\dagger}(\pm m)$ are dual to each other via the isomorphism $\Phi_{S^\dagger,m,-m}$ and the trivialization of $\clO_{S^\dagger}(0)$ in part (1).

        \item\label{disc-part:log-line-bundle-trivialized} Consider a section $m\in\Gamma(S,\ol M\phantom{}^\text{gp})$ and a morphism of schemes $f:S'\to S$ such that the pullback $f^*m\in\Gamma(S',f^*\ol M\phantom{}^\text{gp})$ vanishes. From part (1) and Lemma~\ref{lem:map-of-log-str-slb} applied to $S'^\dagger = (S',f^*M)$ and $f^\dagger = (f,\text{id})$, we get a trivialization $\tau:\clO_{S'}\xrightarrow{\sim}f^*(\clO_{S^\dagger}(m))$. 
        
        \item\label{disc-part:log-line-bundle-trivialized-by-section} If we additionally have $m\in\Gamma(S,\ol M)$ in part (3), then Lemma~\ref{lem:map-of-log-str-slb} shows that the trivialization $\tau$ agrees with the map induced by the section $f^*\sigma_{S^\dagger,m}$.

        \item\label{disc-part:log-lb-vcd-lift-induced-trivialization} Let $m\in\Gamma(S,\ol M)$ be a section and suppose $\tilde m\in\Gamma(S,M)$ is a lift of $m$ under the map $M\to\ol M$. The section $\tilde m$, viewed as a nowhere vanishing section of $\clO_{S^\dagger}(m)^\vee$, induces isomorphism of virtual Cartier divisors $(\clO_{S^\dagger}(m),\sigma_{S^\dagger,m})\xrightarrow{\sim}(\clO_S,\alpha(\tilde m))$.

        \item\label{disc-part:log-lb-vcd-lift-induced-triv-product} If $m_1,m_2\in\Gamma(S,\ol M)$ are sections with respective lifts $\tilde m_1,\tilde m_2\in\Gamma(S,M)$, then the isomorphisms from part (5) for $m_1$, $m_2$ and $m_1 + m_2$ (induced by $\tilde m_1$, $\tilde m_2$ and $\tilde m_1\tilde m_2$ respectively) are compatible with the isomorphism $\Phi_{S^\dagger,m_1,m_2}$.
    \end{enumerate}
\end{Discussion}

\begin{Example}
\label{exa:global-slb-affine-toric}
    For a sharp fine monoid $Q$, consider the standard log structure $M_{V_Q}$ on $V_Q = \Spec\bC[Q]$, as in Definition~\ref{def:std-log-str-toric}. The natural map $Q\to\Gamma(V_Q,\ol M_{V_Q})$ induced by \eqref{eqn:nat-map-Q-to-M-toric} can be shown to be an isomorphism. Under this isomorphism, the slb described in Example~\ref{exa:toric-slb} is precisely the global slb on the log scheme $V_Q^\dagger = (V_Q,M_{V_Q})$.
\end{Example}

\begin{Remark}
\label{rem:log-str-ghost-slb}
    By the discussions in \cite[Section~3]{borne-vistoli-DF} and \cite[Section~2.4]{RSPW}, a fine log structure $M\to\clO_S$ is equivalent to the data of a sheaf of monoids $\ol M$ on $S$, which is fine in the sense of \cite[3.20]{borne-vistoli-DF}, and an slb on $S$ indexed by $\ol M\phantom{}$ which has trivial kernel in the sense of \cite[2.3]{borne-vistoli-DF}. For a detailed proof of this equivalence, see \cite[3.6]{borne-vistoli-DF}.
    
    In one direction, this equivalence is easy to describe. Explicitly, given $M\to\clO_S$, one repeats the construction in Definition~\ref{def-part:global-slb} for each \'etale neighborhood in $S$ and keeps track of all restriction maps to get the desired slb on $S$ indexed by $\ol M$. We will explain the inverse of this equivalence in the case of Deligne--Faltings log structures; see Definition~\ref{def:df-log}. Fine log structures are \'etale locally of Deligne--Faltings type.
    
    The description of a log structure $M$ in terms of an slb indexed by $\ol M$ has the virtue of clearly separating the information contained in the log structure into its discrete (i.e., combinatorial) and continuous (i.e., geometric) constituents. It allows us to informally regard $\ol M$ as the sheaf of `boundary Cartier divisors' and $M\to\ol M$ as the $\clO_S^\times$-torsor keeping track of the `local equations' of these boundary Cartier divisors.
\end{Remark}

\subsection{Deligne--Faltings log structures}\label{subsec:DF-log}
In this subsection, we describe the class of log structures which are completely specified by a single slb.

\begin{Proposition}[{\cite[5.14]{olsson-stack-of-log-str}}, {\cite[3.21]{borne-vistoli-DF}}, {\cite[Section~6]{log-geometry-and-moduli}}]\label{prop:equiv-DF-log}
    Let $Q$ be a sharp fine monoid. For $S\in\Sch$, consider the following three types of objects.
    \begin{enumerate}[itemsep = 0.3ex]
        \myitem[(DF$\phantom{}_1$)]\label{log-str-Q-chart} Log structure $\alpha: M\to\clO_S$ on $S$ along with monoid map $\beta: Q\to\Gamma(S,\ol M)$ which lifts to a log chart $Q\to\Gamma(S,M)$ in some \'etale neighborhood of any given point of $S$.
        
        \myitem[(DF$\phantom{}_2$)]\label{slb-Q} Slb on $S$ indexed by $Q$.
        
        \myitem[(DF$\phantom{}_3$)]\label{map-to-A_Q} Morphism of stacks from $S$ to the Artin cone $\clA_Q = [V_Q/T_Q]$.
    \end{enumerate}
    As the scheme $S$ varies, each of \ref{log-str-Q-chart}, \ref{slb-Q} and \ref{map-to-A_Q} determines a category fibred in groupoids over $\Sch$ using the obvious notion of pullbacks along morphisms in $\Sch$.
    
    These categories fibred in groupoids over $\Sch$ are mutually equivalent. Explicitly, there are functors \ref{log-str-Q-chart}~$\to$~\ref{slb-Q}, \ref{slb-Q}~$\to$~\ref{map-to-A_Q} and \ref{map-to-A_Q}~$\to$~\ref{log-str-Q-chart} such that the resulting composites \ref{log-str-Q-chart}~$\to$~\ref{log-str-Q-chart}, \ref{slb-Q}~$\to$~\ref{slb-Q} and \ref{map-to-A_Q}~$\to$~\ref{map-to-A_Q} are naturally isomorphic to identity functors.
\end{Proposition}

\begin{proof}
    We define functors \ref{log-str-Q-chart}~$\to$~\ref{slb-Q}~$\to$~\ref{map-to-A_Q}~$\to$~\ref{log-str-Q-chart}, but omit the verification that their composites \ref{log-str-Q-chart}~$\to$~\ref{log-str-Q-chart}, \ref{slb-Q}~$\to$~\ref{slb-Q} and \ref{map-to-A_Q}~$\to$~\ref{map-to-A_Q} are naturally isomorphic to identity functors.
    
    \ref{log-str-Q-chart}~$\to$~\ref{slb-Q}: Given $\alpha:M\to\clO_S$ and $\beta:Q\to\Gamma(S,\ol M)$, we map it to the pullback of the global slb on the log scheme $S^\dagger = (S,M)$ along the monoid map $\beta$.

    \ref{slb-Q}~$\to$~\ref{map-to-A_Q}: Given an slb $\{L_q,\Phi_{q_1,q_2},\sigma_q\}$ on $S$ indexed by $Q$, we associate to it a morphism $S\to\clA_Q$ as follows. Consider the $Q^\text{gp}$-graded algebra map
    \begin{align}\label{eqn:graded-C[Q]-algebra-map}
        \bC[Q]\to\clR:=\bigoplus_{q\in Q^\text{gp}} L_q,\qquad
        \chi^q\mapsto\sigma_q\quad\text{for }q\in Q.
    \end{align}
    Here, the multiplication on the quasicoherent $Q^\text{gp}$-graded $\clO_S$-algebra $\clR$ is induced by the isomorphisms $\Phi_{q_1,q_2}$; to make sense of $L_q$ and $\Phi_{q_1,q_2}$ for $q,q_1,q_2\in Q^\text{gp}$, see Discussion~\ref{disc:extending-slb-to-gp}. Note that $P = \relSpec_S\clR\to S$ is a principal $T_Q$-bundle (\'etale local triviality follows from Lemma~\ref{lem:slb-local-triv}), and that it comes with a $T_Q$-equivariant morphism $f:P\to V_Q$ induced by the graded algebra map \eqref{eqn:graded-C[Q]-algebra-map}. The pair $(P,f)$ defines a morphism $S\to \clA_Q$.

    \ref{map-to-A_Q}~$\to$~\ref{log-str-Q-chart}: Given $g:S\to\clA_Q$, lift it to a strict log morphism 
    \begin{align*}
        g^\dagger:S^\dagger = (S,M)\to\clA_Q^\dagger    
    \end{align*}
    by taking $M = g^*M_{\clA_Q}$ and $g^\flat = \text{id}$. Define the map $\beta:Q\to\Gamma(S,\ol M)$ by applying the functor \ref{log-map-to-A_Q} $\to$ \ref{Q-to-ghost-sheaf} from Proposition~\ref{prop:log-maps-to-artin-cone} to $g^\dagger$. Since the log morphism $g^\dagger$ is strict, Proposition~\ref{prop:log-maps-to-artin-cone} also implies that $\beta$ lifts to a log chart in some \'etale neighborhood of any given point of $S$.
\end{proof}

\begin{Definition}\label{def:df-log}
    A \emph{Deligne--Faltings} (\emph{DF}) \emph{log structure} on the scheme $S$, modeled on the sharp fine monoid $Q$, is an object of one of the equivalent groupoids \ref{log-str-Q-chart}, \ref{slb-Q} or \ref{map-to-A_Q}.
\end{Definition}

\begin{Remark}
    Viewpoint \ref{slb-Q} on DF log structures is the simplest to handle from the differential-geometric perspective, since it only involves (maps of) line bundles.
\end{Remark}

\begin{Discussion}\label{disc:equiv-DF-log}
    It is useful to have explicit descriptions of the functors \ref{slb-Q}~$\to$~\ref{log-str-Q-chart}, \ref{log-str-Q-chart}~$\to$~\ref{map-to-A_Q} and \ref{map-to-A_Q}~$\to$~\ref{slb-Q}, defined by composing the functors appearing in the proof of Proposition~\ref{prop:equiv-DF-log}.

    \ref{slb-Q}~$\to$~\ref{log-str-Q-chart}: Let $\{L_q,\Phi_{q_1,q_2},\sigma_q\}$ be a given slb. Define $\alpha: M\to\clO_S$ to be the log structure on $S$ which is associated to the pre-log structure
    \begin{align}\label{eqn:slb-to-prelog}
        \coprod_{q\in Q}\sigma_q:\coprod_{q\in Q} (L_q^\vee)^\times\to\clO_S.
    \end{align}
    The left side of \eqref{eqn:slb-to-prelog} is obtained by sheafifying $U\mapsto\coprod_{q\in Q}\Gamma(U,(L_q^\vee)^\times)$, and it carries the monoid structure defined by the one on $Q$ and the isomorphisms $\Phi_{q_1,q_2}$. Recall from Discussion~\ref{disc:principal-vs-line-bundle} that $(L_q^\vee)^\times$ denotes the sheaf of nowhere vanishing sections of $L_q^\vee$. 
    
    Given $s\in S$, use Lemma~\ref{lem:slb-local-triv} to choose an \'etale neighborhood $U$ of $s$ and trivializations of the line bundles $\{L_q\}_{q\in Q^\text{gp}}$ over $U$ which are compatible with the isomorphisms $\Phi_{q_1,q_2}$. (Here, we are using Discussion~\ref{disc:extending-slb-to-gp} to make sense of $L_q$ and $\Phi_{q_1,q_2}$ for $q,q_1,q_2\in Q^\text{gp}$.) 
    
    These trivializations yield a splitting of the natural projection $\coprod_{q\in Q}\Gamma(U,(L_q^\vee)^\times)\to Q$, i.e., a log chart $Q\to\Gamma(U,M_U)$, where $M_U$ denotes the restriction of $M$ to $U$. The induced maps $Q\to\Gamma(U,\ol M_U)$ patch together to define the global map $\beta:Q\to\Gamma(S,\ol M)$.
    
    \ref{log-str-Q-chart}~$\to$~\ref{map-to-A_Q}: Given $M$ and $\beta$, write $S^\dagger = (S,M)$. Apply the functor \ref{Q-to-ghost-sheaf} $\to$ \ref{log-map-to-A_Q} from Proposition~\ref{prop:log-maps-to-artin-cone} to get a strict log morphism $S^\dagger\to\clA_Q^\dagger$ and then pass to the underlying morphism of stacks to get $S\to\clA_Q$.

    \ref{map-to-A_Q}~$\to$~\ref{slb-Q}: Let a morphism $S\to\clA_Q$ be given, i.e., a pair $(P,f)$ consisting of a principal $T_Q$-bundle $\pi:P\to S$ and a $T_Q$-equivariant map $f:P\to V_Q$. In particular, $\pi$ is an affine morphism and we get a quasicoherent $Q^\text{gp}$-graded $\bC[Q]$-algebra $\clR = \pi_*\clO_P$ on $S$. 
    
    For $q\in Q^\text{gp}$, let $L_q$ be the $q$-graded piece of $\clR$; it is a locally free $\clO_S$-module of rank $1$. Multiplication in $\clR$ induces isomorphisms $\Phi_{q_1,q_2}:L_{q_1}\otimes L_{q_2}\xrightarrow{\sim} L_{q_1+q_2}$ for all $q_1,q_2\in Q^\text{gp}$. Given any $q\in Q$, define $\sigma_q\in H^0(S,L_q)$ to be the image of the monomial $\chi^q$ under the algebra map $\bC[Q]\to\clR$ corresponding to $f:P\to V_Q$. This defines the slb $\{L_q,\Phi_{q_1,q_2},\sigma_q\}$.
\end{Discussion}

\subsection{Log morphisms to Deligne--Faltings targets}\label{subsec:morphism-to-DF-log}

For this subsection, fix a scheme $S$ in $\Sch$ equipped with a DF log structure modeled on a sharp fine monoid $Q$. From this, we obtain the following equivalent structures on $S$ using viewpoints \ref{log-str-Q-chart}, \ref{slb-Q} and \ref{map-to-A_Q}.
\begin{enumerate}[itemsep = 0.3ex]
    \myitem[(DF$\phantom{}_1$)] A log scheme $S^\dagger = (S,M_S)$ and a monoid map $\beta:Q\to\Gamma(S,\ol M_S)$ which admits \'etale local lifts to a log chart $Q\to\Gamma(S,M_S)$.
    
    \myitem[(DF$\phantom{}_2$)] A system of line bundles $\clL_{S^\dagger,\beta}$ on $S$ indexed by $Q$. Letting $\clL_{S^\dagger}$ denote the global slb on the log scheme $S^\dagger$ in the sense of Definition~\ref{def-part:global-slb}, the slb $\clL_{S^\dagger,\beta}$ is seen to be the pullback of $\clL_{S^\dagger}$ along the monoid map $\beta$.
    
    \myitem[(DF$\phantom{}_3$)] A morphism of stacks $\hat\beta:S\to\clA_Q$.
\end{enumerate}

\begin{Definition}\label{def:tropicalization}
    Given a log scheme $T^\dagger = (T,M_T)\in\LogSch$ along with a log morphism $f^\dagger=(f,f^\flat):T^\dagger\to S^\dagger$, the composition
    \begin{align*}
        Q\xrightarrow{\beta}\Gamma(S,\ol M_S)\xrightarrow{f^*}\Gamma(T,f^*\ol M_S)\xrightarrow{\Gamma(\ol f\phantom{}^\flat)}\Gamma(T,\ol M_T)
    \end{align*}
    is called the \emph{tropicalization} of $f^\dagger$ and is denoted by $\Trop(f^\dagger)$.
\end{Definition}

\begin{Remark}\label{rem:tropicalization}
    In the setting of Definition~\ref{def:tropicalization}, if the monoids $Q$ and $\Gamma(T,\ol M_T)$ are fs, then there is no loss of information in replacing $\Trop(f^\dagger)$ by the dual map of \emph{cones} 
    \begin{align*}
        \Hom_{\Mon}(\Gamma(T,\ol M_T),\bR_{\ge 0})\to\Hom_{\Mon}(Q,\bR_{\ge 0}).
    \end{align*}
    
    When the log schemes $T^\dagger$ and $S^\dagger$ are both fs but $S^\dagger$ is not necessarily of Deligne--Faltings type, the tropicalization of a log morphism $f^\dagger:T^\dagger\to S^\dagger$ can be defined as a map of \emph{cone complexes}; see \cite[Appendix~B]{GS-log}.
\end{Remark}

The next two results give a classification of log morphisms $f^\dagger$ for which either the underlying morphism $f$ or tropicalization $\Trop(f^\dagger)$ is already prescribed.

\begin{Proposition}
\label{prop:log-mor-to-DF-target}
    Let $T^\dagger = (T,M_T)$ be an object in $\LogSch$ and let $f:T\to S$ be a morphism in $\Sch$. Denote the global slb on the log scheme $T^\dagger$ by $\clL_{T^\dagger}$. Then there is a natural bijection between the following two types of objects.
    \begin{enumerate}[itemsep = 0.3ex]
        \myitem[(LM$\phantom{}_1$)]\label{log-mor-lift} Log morphism $f^\dagger=(f,f^\flat):T^\dagger\to S^\dagger$ with underlying morphism $f:T\to S$.
        
        \myitem[(LM$\phantom{}_2$)]\label{ghost-plus-slb-iso} Monoid map $\psi:Q\to\Gamma(T,\ol M_T)$ and slb isomorphism $\Psi:f^*\clL_{S^\dagger,\beta}\xrightarrow{\sim}\psi^*\clL_{T^\dagger}$.
    \end{enumerate}
    
    Moreover, under this bijection, we have $\psi = \Trop(f^\dagger)$.
\end{Proposition}

\begin{proof}
    We show how to define maps \ref{log-mor-lift}~$\to$~\ref{ghost-plus-slb-iso}~$\to$~\ref{log-mor-lift} but omit the verification that the composites \ref{log-mor-lift}~$\to$~\ref{log-mor-lift} and \ref{ghost-plus-slb-iso}~$\to$~\ref{ghost-plus-slb-iso} are identity maps.

    \ref{log-mor-lift}~$\to$~\ref{ghost-plus-slb-iso}: 
    Given a log morphism $f = (f,f^\flat):T^\dagger\to S^\dagger$ with underlying morphism $f:T\to S$, define $\psi = \Trop(f^\dagger)$.

    By Lemma~\ref{lem:map-of-log-str-slb}, $f^\flat$ induces an slb isomorphism
    from the pullback of the global slb $\clL_{S^\dagger}$ along $f:T\to S$
    to the pullback of the global slb $\clL_{T^\dagger}$ along the map $\Gamma(S,\ol M\phantom{}_S)\to \Gamma(T,\ol M\phantom{}_T)$ induced by $f^\flat$.
    Define $\Psi$ to be the pullback of this slb isomorphism along $\beta$.

    \ref{ghost-plus-slb-iso}~$\to$~\ref{log-mor-lift}: Given a monoid map $\psi:Q\to\Gamma(T,\ol M_T)$ and an slb isomorphism $\Psi:f^*\clL_{S^\dagger,\beta}\xrightarrow{\sim}\psi^*\clL_{T^\dagger}$, we must define a map of log structures $f^\flat:f^*M_S\to M_T$. For simplicity, abbreviate the line bundles appearing in $\clL_{S^\dagger,\beta}$ as $L_q = \clO_{S^\dagger}(\beta(q))$ for $q\in Q$.
    
    Using the explicit formula \eqref{eqn:slb-to-prelog} for $f^*M_S$, defining $f^\flat$ amounts to defining a map
    \begin{align*}
        f^\flat_\text{pre}:\coprod_{q\in Q}
        (f^*L_q^\vee)^\times\to M_T
    \end{align*}
    of pre-log structures on $T$; recall the adjunction \eqref{eqn:prelog-log-adjunction}. Given a local section $\tau$ of $(f^*L_q^\vee)^\times$ for some $q\in Q$, the line bundle isomorphism $\Psi_q:f^*L_q\xrightarrow{\sim}\clO_{T^\dagger}(\psi(q))$ converts it into a local section $\tau'$ of $\clO_{T^\dagger}(-\psi(q))^\times\subset M_T$. In this situation, we define $f^\flat_\text{pre}(\tau) = \tau'$. Using the fact that $\Psi$ is an slb isomorphism, it is immediate that $f^\flat_\text{pre}$ is a map of pre-log structures on $T$.
\end{proof}

\begin{Proposition}
\label{prop:log-mor-to-DF-target-2}
    Let $T^\dagger = (T,M_T)$ be an object in $\LogSch$ and let $\psi:Q\to\Gamma(T,\ol M_T)$ be a monoid map. Let $\hat\psi:T\to\clA_Q$ be the result of applying the functor \ref{Q-to-ghost-sheaf} $\to$ \ref{log-map-to-A_Q} from Proposition~\ref{prop:log-maps-to-artin-cone} to $\psi$ and then passing to the underlying morphism. Define the algebraic space $S_\psi$ by the following $2$-fibre product diagram.
    \begin{equation*}
    \begin{tikzcd}
        S_\psi \arrow[r] \arrow[d] \arrow[dr,phantom,"\lrcorner",pos=0.001] & S \arrow[d,"\hat\beta"] \arrow[dl,phantom,pos=0.45, "\substack{\mathbin{\rotatebox[origin=c]{45}{$\Leftarrow$}}}"] \\
        T \arrow[r,swap,"\hat\psi"] & \clA_Q
    \end{tikzcd}    
    \end{equation*}
    
    Then there is a natural bijection between the following two types of objects.
    \begin{enumerate}[itemsep = 0.3ex]
        \myitem[(LM$\phantom{}'_1$)]\label{log-mor-with-given-trop} Log morphism $f^\dagger:T^\dagger\to S^\dagger$ with tropicalization $\psi:Q\to \Gamma(T,\ol M_T)$.
        
        \myitem[(LM$\phantom{}'_2$)]\label{target-fibration-sections} Section $\sigma:T\to S_\psi$ of the projection morphism $S_\psi\to T$.
    \end{enumerate}
\end{Proposition}
\begin{proof}
    We need the following claim.
    \begin{claim}\label{claim:target-fibration}
        For any $V$ in $\Sch$, morphisms $V\to S_\psi$ are naturally in bijection with the data of a pair of morphisms $V\to S$ and $V\to T$, together with an slb isomorphism between the pullbacks to $V$ along these morphisms of $\clL_{S^\dagger,\beta}$ and $\psi^*\clL_{T^\dagger}$.
    \end{claim}
    \begin{proof}[Proof of claim] 
        This is immediate from the universal property of 2-fibre products and the fact that the equivalence \ref{slb-Q}~$\simeq$~\ref{map-to-A_Q} from Proposition~\ref{prop:equiv-DF-log} identifies $\psi^*\clL_{T^\dagger}$ with $\hat\psi$.
    \end{proof}

    We now show how to define maps \ref{log-mor-with-given-trop}~$\to$~\ref{target-fibration-sections}~$\to$~\ref{log-mor-with-given-trop} but omit the verification that the composites \ref{log-mor-with-given-trop}~$\to$~\ref{log-mor-with-given-trop} and \ref{target-fibration-sections}~$\to$~\ref{target-fibration-sections} are identity maps.
    
    \ref{log-mor-with-given-trop}~$\to$~\ref{target-fibration-sections}: Given a log morphism $f^\dagger = (f,f^\flat):T^\dagger\to S^\dagger$ with $\Trop(f^\dagger) = \psi$, the map \ref{log-mor-lift} $\to$ \ref{ghost-plus-slb-iso} from Proposition~\ref{prop:log-mor-to-DF-target} yields an slb isomorphism $\Psi:f^*\clL_{S^\dagger,\beta}\xrightarrow{\sim}\psi^*\clL_{T^\dagger}$.By Claim~\ref{claim:target-fibration} applied to $V = T$, we find that the morphisms $f:T\to S$, $\text{id}_T:T\to T$ and the slb isomorphism $\Psi$ together define a section $\sigma$ of the projection $S_\psi\to T$.

    \ref{target-fibration-sections}~$\to$~\ref{log-mor-with-given-trop}: Given a section $\sigma$ of $S_\psi\to T$, Claim~\ref{claim:target-fibration} determines a morphism $f:T\to S$ and an slb isomorphism $\Psi:f^*\clL_{S^\dagger,\beta}\xrightarrow{\sim}\psi^*\clL_{T^\dagger}$. The map \ref{ghost-plus-slb-iso} $\to$ \ref{log-mor-lift} from Proposition~\ref{prop:log-mor-to-DF-target} now yields a log morphism $f^\dagger:T^\dagger\to S^\dagger$ with $\Trop(f^\dagger) = \psi$.
\end{proof}

\subsection{Log points}\label{subsec:log-points}

\begin{Definition}\label{def:log-point}
    A \emph{log point} is a log scheme $S^\dagger = (S,M_S)$ in $\LogSch$ with $S = \Spec\bC$. 
    
    The full subcategory of $\LogSch$ consisting of log points is denoted by $\LogPt$.
\end{Definition}

\begin{Proposition}
\label{prop:cat-of-log-pts}
    Consider the category $\SlbPt$ in which
    \begin{enumerate}[label = (\arabic*), itemsep = 0.3ex]
        \item objects are pairs $(Q,\clL)$, where $Q$ is a sharp fine monoid and $\clL=\{L_q,\Phi_{q_1,q_2},\sigma_q\}$ is an slb on $S$ indexed by $Q$ such that $\sigma_{q} = 0$ for all $0\ne q\in Q$, and
        
        \item morphisms $(Q,\clL)\to(Q',\clL')$ are pairs $(\psi,\Psi)$, where $\psi:Q'\to Q$ is a monoid map such that $\psi^{-1}(0) = \{0\}$ and $\Psi:\clL'\xrightarrow{\sim}\psi^*\clL$ is an slb isomorphism.
    \end{enumerate}
    
    Then the functor $\LogPt\to\SlbPt$ which sends a log point $S^\dagger = (S,M_S)$ to the sharp fine monoid $\ol M_S$ and the global slb on $S^\dagger$ is an equivalence.
\end{Proposition}
\begin{proof}
    We only explain why the functor $\LogPt\to\SlbPt$ in the statement is well-defined. The fact that it is an equivalence is an easy consequence of Propositions~\ref{prop:equiv-DF-log} and \ref{prop:log-mor-to-DF-target}.
    
    Given a log structure $\alpha:M\to\clO_S = \bC$, write $S^\dagger = (S,M)$. The monoid $\ol M$ is sharp by Definition~\ref{def-part:ghost-sheaf}. For any $0\ne m\in\ol M$ and any lift $\tilde m\in M$ of $m$, we have $\alpha(\tilde m)\not\in\clO_S^\times = \bC^\times$, i.e., we have $\alpha(\tilde m) = 0$. This shows that $\sigma_{S^\dagger,m} = 0$ for all $0\ne m\in\ol M$. 
    
    Given another log structure $\alpha':M'\to\clO_S$, write $S'^\dagger = (S,M')$. Then a log morphism $f^\dagger:S^\dagger\to S'^\dagger$ in $\LogPt$ is the same as a map of log structures $f^\flat:M'\to M$, which explains the reversal of arrows in the definition of $\SlbPt$. The induced map $\ol f\phantom{}^\flat:\ol M\phantom{}'\to\ol M$ on ghost sheaves satisfies $(\ol f\phantom{}^\flat)^{-1}(0) = \{0\}$ by Definition~\ref{def-part:ghost-sheaf}.
\end{proof}

\begin{Definition}
    Let $Q$ be a sharp fine monoid. We then have a natural map $Q\to\bC$, which is defined by $0\mapsto 1$ and $q\mapsto 0$ for all $0\ne q\in Q$.
    
    The \emph{split log point with monoid $Q$} is the log scheme $\Spec(Q\to\bC) = (\Spec\bC,\bC^\times\oplus Q)$, as in Construction~\ref{cons:log-spec}. The split log point with monoid $\bN$ is often also referred to as the \emph{standard log point}.
\end{Definition}

\begin{Remark}
\label{rem:splitting-log-pt}
    For a log point $S^\dagger = (S,M_S)$, an isomorphism $S^\dagger\simeq \Spec(\ol M_S\to\bC)$ of log schemes which induces the identity on $\ol M_S$ is equivalent to a section of $M_S\to\ol M_S$. The set of such sections is non-empty by Lemma~\ref{lem:local-lift-from-ghost} and it forms a $\Hom_{\Ab}(\ol M\phantom{}_S^\text{gp},\bC^\times)$-torsor.
\end{Remark}

\subsection{Smooth snc pairs}\label{subsec:smooth-snc}

\begin{Definition}\label{def:smooth-snc-pair}
    A \emph{smooth snc pair} $(X,D)$ consists of a smooth variety $X$ equipped with an snc divisor $D\subset X$; see Definition~\ref{def:nc-divisor-scheme}.
\end{Definition}

For the rest of this subsection, fix a smooth snc pair $(X,D)$ and an ordering $D_1,\ldots,D_r$ of the (smooth) irreducible components of $D$. Using this, we decompose $D$ as 
\begin{align}\label{eqn:snc-div-decomp}
    D = D_1\cup\cdots\cup D_r.    
\end{align}
For $1\le j\le r$, the inclusion $\clO_X\subset \clO_X(D_j)$ induces a global section $\taut_{D_j}\in H^0(X,\clO_X(D_j))$; see Discussion~\ref{disc:sheaf-cartier-div}. The linearization of the section $\taut_{D_j}$ along its vanishing locus $D_j\subset X$ yields a canonical isomorphism
\begin{align}\label{eqn:normal-bundle}
    \clN_{D_j/X}:= T_X|_{D_j}/T_{D_j}\xrightarrow{\sim}\clO_X(D_j)|_{D_j}
\end{align}
between the normal bundle of $D_j\subset X$ and the restriction of $\clO_X(D_j)$ to $D_j$.

\begin{Discussion}\label{disc:snc-pair-log-str}
    The decomposition \eqref{eqn:snc-div-decomp} of the snc divisor $D$ endows $X$ with a DF log structure modeled on $\bN^r$, which we may describe as follows using the three equivalent viewpoints \ref{log-str-Q-chart}, \ref{slb-Q} and \ref{map-to-A_Q}.

    \begin{enumerate}[itemsep = 0.3ex]
        \myitem[(DF$\phantom{}_1$)] We have the divisorial log structure $M_{(X,D)}\subset\clO_X$, as in Example~\ref{exa:div-log-str}, whose local sections are the regular functions on $X$ which are nowhere vanishing on $X\setminus D$. We write the resulting log scheme as $X^\dagger = (X,M_{(X,D)})$.
    
        Moreover, for $1\le j\le r$, we have a global section $s_j\in\Gamma(X,\ol M_{(X,D)})$ characterized by the property that each of its local lifts to $M_{(X,D)}$ gives a local equation for $D_j\subset X$. The sections $s_1,\ldots,s_r$ define a monoid map $\beta:\bN^r\to\Gamma(X,\ol M_{(X,D)})$, which lifts to a log chart near any point since $D\subset X$ is an snc divisor.

        \myitem[(DF$\phantom{}_2$)] The divisors $D_1,\ldots,D_r\subset X$ determine an slb indexed by $\bN^r$ on $X$ (Example~\ref{exa:n-free-cart-div}). 
    
        More explicitly, for ${\bf a} = (a_1,\ldots,a_r)\in\bN^r$, we take the corresponding line bundle to be $L_{\bf a} = \bigotimes_j\clO_X(D_j)^{\otimes a_j}$ and the global section of $L_{\bf a}$ to be $\sigma_{\bf a} = \bigotimes_j\taut_{D_j}^{\otimes a_j}$. The isomorphisms $\Phi_{{\bf a}_1,{\bf a}_2}$ are all taken to be natural identifications.
        
        For all ${\bf a}\in\bN^r$, multiplication of rational functions on $X$ induces an identification of $L_{\bf a}$ with $\clO_X(\sum a_jD_j)$; see Discussion~\ref{disc:sheaf-cartier-div}. This identification maps the section $\sigma_{\bf a}$ to the section $\taut_{\sum a_jD_j}$ induced by the inclusion $\clO_X\subset\clO_X(\sum a_jD_j)$.

        \myitem[(DF$\phantom{}_3$)] For $1\le j\le r$, interpreting the inclusion $\clO_X(-D_j)^\times\subset\clO_X$ as a principal $\bC^\times$-bundle over $X$ mapping $\bC^\times$-equivariantly to $\bC$, we get a morphism $\hat\beta_j:X\to[\bC/\bC^\times] = \clA_\bN$. Putting these together for $1\le j\le r$, we get a morphism
        \begin{align}\label{eqn:mor-snc-pair-to-artin-cone}
            \hat\beta = (\hat\beta_1,\ldots,\hat\beta_r):X\to[\bC/\bC^\times]^r = \clA_{\bN^r}.
        \end{align}
    \end{enumerate}
\end{Discussion}

It is easy to explicitly describe log morphisms to $X^\dagger = (X,M_{(X,D)})$ using the discussion from the preceding subsections.

\begin{Proposition}
\label{prop:log-mor-to-snc-pair}
    Let $S^\dagger = (S,M_S)$ be a log scheme. Then there is a natural bijection between the following two types of objects.
    \begin{enumerate}[itemsep = 0.3ex]
        \myitem[(LM$\phantom{}_1$)]\label{abstract-log-mor-to-snc-pair} Log morphism $f^\dagger:S^\dagger\to X^\dagger$.
        
        \myitem[(LM$\phantom{}_2$)]\label{concrete-log-mor-to-snc-pair} Morphism $f:S\to X$ and, for each $1\le j\le r$, a section $m_j\in\Gamma(S,\ol M_S)$ together with an isomorphism of virtual Cartier divisors
        \begin{align*}
            \Psi_j:f^*(\clO_X(D_j),\taut_{D_j})\xrightarrow{\sim}(\clO_{S^\dagger}(m_j),\sigma_{S^\dagger,m_j}).
        \end{align*}
    \end{enumerate}

    Under this bijection, the underlying morphism of $f^\dagger:S^\dagger\to X^\dagger$ is given by $f:S\to X$ and its tropicalization $\Trop(f^\dagger):\bN^r\to\Gamma(S,\ol M_S)$ is given by ${\bf e}_j\mapsto m_j$ for $1\le j\le r$.
\end{Proposition}
\begin{proof}
    This follows from the special case $Q = \bN^r$ of Proposition~\ref{prop:log-mor-to-DF-target}.
\end{proof}

\begin{Definition}
    Given a subset $I\subset\{1,\ldots,r\}$, the corresponding \emph{closed stratum} $D_I\subset X$ and the \emph{locally closed stratum} $D_I^\circ\subset X$ of the smooth snc pair $(X,D)$ are defined to be
    \begin{align*}
        D_I := \bigcap_{j\in I}D_j\quad\text{and}\quad D_I^\circ := \bigcap_{j\in I}D_j\cap\bigcap_{j\not\in I}(X\setminus D_j).
    \end{align*}
\end{Definition}

\begin{Example}
\label{exa:snc-log-str-on-stratum}
    Let $f:S\to X$ be a morphism which factors through a locally closed stratum $D_I^\circ$. For instance, this happens if $S = \Spec\bC$ and $f$ is the inclusion of a point. We may view $f$ as a morphism to the smaller smooth snc pair $(X,D)\setminus\bigcup_{j\not\in I}D_j$. 
    
    It follows that $f^*M_{(X,D)}$ comes from a DF log structure on $S$ modeled on $\bN^I$, which is described in viewpoint \ref{slb-Q} by the slb freely generated by the collection $\{(f^*\clN_{D_j/X},0)\}_{j\in I}$ of virtual Cartier divisors on $S$; see Example~\ref{exa:df-rank-r-sys-line-bundle}. Here, we have used \eqref{eqn:normal-bundle} to identify the restriction $\clO_X(D_j)|_{D_j}$ with the normal bundle $\clN_{D_j/X}$.

    By the explicit formula \eqref{eqn:slb-to-prelog} for converting an slb to a log structure, we find that the ghost sheaf $f^*\ol M_{(X,D)}$ is canonically isomorphic to the constant sheaf on $S$ with stalk $\bN^I$.
\end{Example}

We conclude with a useful consequence of Example~\ref{exa:snc-log-str-on-stratum} for the case of nc divisors.

\begin{Proposition}\label{prop:nc-log-str-at-point}
    Let $E\subset Y$ be an nc divisor in a smooth scheme $Y$ and let $M_{(Y,E)}$ be the associated divisorial log structure. Let $y\in Y$ be a point and let $\{\Delta_j\}_{j\in I}$ be the collection of distinct branches of $E$ at $y$; see Definition~\ref{def-part:branch}. We take $I = \emptyset$ if $y\not\in E$.

    Let $\{\rho_{\Delta_j}\}_{j\in I}$ be the basis of the free monoid $y^*\ol M_{(Y,E)}$ provided by Lemma~\ref{lem-part:ghost-sheaf-stalk-free}. Then, for each $j\in I$, the virtual Cartier divisor associated to $\rho_{\Delta_j}$ is given by $(\clN_{\Delta_j/Y}|_y,0)$.
\end{Proposition}
\begin{proof}
    We may replace $Y$ by an \'etale neighborhood of $y$ to reduce to the case $(Y,E)$ is a smooth snc pair. Now, the result is immediate from Example~\ref{exa:snc-log-str-on-stratum}.
\end{proof}
    \section{Log curves}\label{sec:log-curves}

We discuss the notion of a \emph{family of log curves} over a log scheme. We then give a complete classification and explicit description of log curve families over a log point (Theorem~\ref{thm:log-curve-over-log-pt}). Our treatment of this topic is slightly non-standard but nevertheless equivalent to the original treatment in \cite{FKato-logcurve}; this is explained in Remarks~\ref{rem:log-curve-defs} and \ref{rem:log-curves-whats-new} below.

\begin{Notation}
    For this section, we fix integers $g,n\ge 0$.
\end{Notation}

\subsection{Prestable curves}

Recall that a scheme $C\in\Sch$ is called a \emph{prestable curve} if it is projective, connected, reduced, $1$-dimensional and has at worst nodal singularities. 

An \emph{$n$-pointed prestable genus $g$ curve} $(C,{\bf x})$ consists of a prestable curve $C$ satisfying $g = \dim H^1(C,\clO_C)$, together with an $n$-tuple ${\bf x}=(x_1,\ldots,x_n)$ of distinct smooth points of $C$. We refer to the points $x_i$ for $1\le i\le n$ as \emph{marked points}.

\begin{Notation}\label{nota:prestable-curve}
    Given a prestable curve $(C,{\bf x})$, we have the following associated objects.
    \begin{enumerate}[label = (\arabic*), itemsep = 0.3ex]
        \item The normalization of $C$ is denoted by $\nu:\tilde C\to C$. Given $p\in C$ which is not a node, we denote its unique inverse image in $\tilde C$ also by $p$ by a minor abuse of notation.

        \item The complement of the marked points and nodes is denoted by $C_\text{gen}\subset C$.

        \item The combinatorial structure of $C$ is encoded in its \emph{dual graph} $\plC$.
        \begin{enumerate}[label = (\roman*), itemsep = 0.15ex]
            \item The vertices of $\plC$, denoted by $v$, index the connected components $\tilde C_v$ of $\tilde C$ or equivalently, the irreducible components $C_v$ of $C$. 
        
            \item The edges of $\plC$, denoted by $e$, index the nodes $y_e\in C$. The oriented edges of $\plC$, denoted by $\vec{e}$, index the points $y_{\vec{e}}\in\tilde C$ which map to nodes under $\nu$. For an oriented edge $\vec{e}$, the same edge $e$ with the opposite orientation is denoted by $\cev{e}$. 
        
            \item We write $\vec{e}: v\to v'$ to indicate that an oriented edge $\vec{e}$ goes from the vertex $v$ to the vertex $v'$, i.e., $y_{\vec{e}}\in\tilde C_v$ and $y_{\cev{e}}\in\tilde C_{v'}$.
        \end{enumerate}

        \item We have the following two $\bC$-vector spaces associated to each node $y_e\in C$.
        \begin{enumerate}[label = (\roman*), itemsep = 0.15ex]
            \item The \emph{space of smoothing parameters} at $y_e$ is defined as $\Lambda_{C,y_e}:=T_{\tilde C,y_{\vec{e}}}\otimes T_{\tilde C,y_{\cev{e}}}$.
            
            \item The \emph{Zariski tangent space} of $C$ at $y_e$ is defined as $T^{\Zar}_{C,y_e}:=(\fm_{C,y_e}/\fm_{C,y_e}^2)^\vee$, where $\fm_{C,y_e}$ is the maximal ideal of $\clO_{C,y_e}$. The normalization morphism $\nu$ induces a canonical direct sum decomposition $T^{\Zar}_{C,y_e} = T_{\tilde C,y_{\vec{e}}}\oplus T_{\tilde C,y_{\cev{e}}}$.
            
            This decomposition induces a canonical projection map $\Sym^2T^{\Zar}_{C,y_e}\to \Lambda_{C,y_e}$.
        \end{enumerate}
    \end{enumerate}
\end{Notation}

\subsection{Stacks of prestable curves} 

\begin{Definition}
\label{def:prestable-curves-stack}
    Let $\fM_{g,n}$ be the \emph{stack of $n$-pointed prestable genus $g$ curves} over $\Sch$. Explicitly, for $S\in\Sch$, an object of $\fM_{g,n}(S)$ is a tuple $(\pi:C\to S,{\bf x})$, where
    \begin{enumerate}[label = (\arabic*), itemsep = 0.3ex]
        \item $C$ is an \emph{algebraic space} over $\bC$ as defined in \cite[\href{https://stacks.math.columbia.edu/tag/025Y}{Tag~025Y}]{stacks-project},
        
        \item $\pi:C\to S$ is a proper, flat morphism of finite presentation, and
        
        \item ${\bf x}=(x_1,\ldots,x_n)$ is an $n$-tuple consisting of sections $x_i:S\to C$ of $\pi$ for $1\le i\le n$,
    \end{enumerate}
    such that, for each point $s\in S$, the fibre $C_s := \pi^{-1}(s)$ equipped with ${\bf x}(s) = (x_1(s),\ldots,x_n(s))$ is an $n$-pointed prestable genus $g$ curve. A morphism in $\fM_{g,n}$ is a pullback diagram respecting the $n$ marked point sections. Let $\fM_{g,n}^\times\subset\fM_{g,n}$ denote the open substack consisting of those families for which $\pi:C\to S$ is a smooth morphism.

    Let $\fC_{g,n}$ be the \emph{universal $n$-pointed prestable genus $g$ curve}, which is a stack over $\Sch$. Explicitly, for $S\in\Sch$, an object of $\fC_{g,n}(S)$ consists of an object $(\pi:C\to S,{\bf x})$ of $\fM_{g,n}(S)$ together with an additional section $x_0:S\to C$, which is not required to be disjoint from the relative singular locus of $\pi$ or the images of any of $x_1,\ldots,x_n$. A morphism in $\fC_{g,n}$ is a pullback diagram respecting the $n+1$ sections. Let $\fC_{g,n}^\times\subset\fC_{g,n}$ denote the open substack consisting of those families for which $\pi:C\to S$ is a smooth morphism and, for each point $s\in S$ and $1\le i\le n$, we have $x_0(s) \ne x_i(s)$.

    We define the \emph{universal family of $n$-pointed prestable genus $g$ curves} to be the natural morphism $\Pi_{g,n}:\fC_{g,n}\to\fM_{g,n}$ given by forgetting the section $x_0$. For each $1\le i\le n$, we have a section $x_{g,n;i}$ of $\Pi_{g,n}$ given by setting $x_0 = x_i$. We write ${\bf x}_{g,n} = (x_{g,n;1},\ldots,x_{g,n;n})$.
\end{Definition}

\begin{Notation}
    For $S\in\Sch$, we usually denote a family in $\fM_{g,n}(S)$ by using notation such as $(\pi:C\to S,{\bf x})$. When $\pi$ is clear from the context, we shorten this to $(C/S,{\bf x})$.
\end{Notation}

\begin{Remark}
    \label{rem:why-alg-space}
    In Definition~\ref{def:prestable-curves-stack}, if we considered only families $(C/S,{\bf x})$ over $S\in\Sch$ for which $C$ is a scheme, then $\fM_{g,n}$ would fail to be a stack; see \cite[\href{https://stacks.math.columbia.edu/tag/0DMJ}{Tag~0DMJ}]{stacks-project}. Allowing $C$ to be an algebraic space solves this problem; see \cite[\href{https://stacks.math.columbia.edu/tag/04U0}{Tag~04U0}]{stacks-project}.
    
    However, it is true that there always exists an \'etale covering $\{S_i\to S\}_i$ such that each $C\times_S S_i\to S_i$ is a projective morphism of schemes; see \cite[\href{https://stacks.math.columbia.edu/tag/0E6F}{Tag~0E6F}]{stacks-project}. 
    Thus, for \'etale local considerations on $S$, we may assume that $C$ is a scheme.
\end{Remark}

\begin{Remark}
\label{rem:prestable-curves-algebraic-stack}
    By \cite[5.1]{olsson-log-twisted-curves} or \cite[A.1]{GS-log}, the stack $\fM_{g,n}$ is algebraic. 
    
    Given $S\in\Sch$ and a family $\xi = (\pi:C\to S,{\bf x})\in\fM_{g,n}(S)$,
    let $\hat\xi\in\fC_{g,n}(C)$ be the family defined by the pullback $\pi^*\xi\in\fM_{g,n}(C)$, together with the section $x_0:C\to C\times_SC$ given by the diagonal morphism. Then we have the following 2-fibre product diagram.
    \begin{equation}\label{eqn:classifying-map-curve-family}
    \begin{tikzcd}
        C \arrow[d,swap,"\pi"] \arrow[r,"\hat\xi"] \arrow[dr,phantom,"\lrcorner",pos=0.001] & \fC_{g,n} \arrow[dl,phantom,pos=0.45, "\substack{\mathbin{\rotatebox[origin=c]{45}{$\Leftarrow$}} \\ \text{id}}"] \arrow[d,"\Pi_{g,n}"] \\
        S \arrow[r,swap,"\xi"] & \fM_{g,n}
    \end{tikzcd}
    \end{equation}
    Thus, the morphism $\Pi_{g,n}:\fC_{g,n}\to\fM_{g,n}$ is representable by algebraic spaces and it follows from \cite[\href{https://stacks.math.columbia.edu/tag/05UM}{Tag~05UM}]{stacks-project} that the stack $\fC_{g,n}$ is also algebraic.
\end{Remark}

\begin{Definition}
    Let $\Delta_{\fM_{g,n}}\subset\fM_{g,n}$ be the reduced closed substack $\fM_{g,n}\setminus\fM_{g,n}^\times$. Similarly, let $\Delta_{\fC_{g,n}}\subset\fC_{g,n}$ be the reduced closed substack $\fC_{g,n}\setminus\fC_{g,n}^\times$.
\end{Definition}

The next result describes the local structure of a versal family of prestable curves. To understand the statement and its proof, recall Notation~\ref{nota:prestable-curve} and the notion of \emph{branches} of an nc divisor from Definition~\ref{def-part:branch}.

\begin{Proposition}\label{prop:versal-curve}
    Consider $S\in\Sch$ and let $\xi = (\pi:C\to S,{\bf x})\in\fM_{g,n}(S)$ be a versal family, i.e., $\xi:S\to\fM_{g,n}$ is smooth. Using $\xi$, define the family $\hat\xi\in\fC_{g,n}(C)$ as in Remark~\ref{rem:prestable-curves-algebraic-stack}; by \eqref{eqn:classifying-map-curve-family}, $\hat\xi:C\to\fC_{g,n}$ is also smooth. We then have the following properties.
    \begin{enumerate}[label = (\arabic*), ref = \theProposition(\arabic*), itemsep = 0.3ex]
        \item\label{prop-part:versal-curve-smooth-nc} Both $S$ and $C$ are smooth. 
        The pullbacks $\Delta_\xi := S\times_{\fM_{g,n}}\Delta_{\fM_{g,n}}$ and $\Delta_{\hat\xi} := C\times_{\fC_{g,n}}\Delta_{\fC_{g,n}}$ are nc divisors in $S$ and $C$ respectively. Moreover, we have an equality of divisors 
        \begin{align}\label{eqn:versal-curve-boundary-divisor}
            \Delta_{\hat\xi} = \pi^{-1}(\Delta_\xi)+\sum_ix_i(S).
        \end{align}
        
        \item\label{prop-part:versal-curve-basic-log} Given a point $s\in S$, define $C_s:=\pi^{-1}(s)$. Then there is a unique bijection 
        \begin{align}\label{eqn:versal-edge-branch-bijection}
            e\mapsto\Delta_{\xi,e}    
        \end{align}
        from the edges of the dual graph $\plC_s$ of the curve $C_s$ to the branches of the nc divisor $\Delta_\xi\subset S$ at $s$, characterized by the following properties for each edge $e$.
        \begin{enumerate}[label = (\roman*), ref = \theenumi(\roman*), itemsep = 0.15ex]            
            \item The first derivative $(D\pi)_{y_e}: T_{C,y_e}\to T_{S,s}$ of $\pi$ at $y_e$
            satisfies
            \begin{align*}
                T^{\Zar}_{C_s,y_e} = \ker(D\pi)_{y_e}\quad\text{and}\quad
                \clN_{\Delta_{\xi,e}/S}|_s = \coker(D\pi)_{y_e},
            \end{align*}
            where $\clN_{\Delta_{\xi,e}/S}|_s$ is the normal space as in Definition~\ref{def-part:normal-space}.
            
            \item The second derivative of $\pi$ at $y_e$ is well-defined as a linear map 
            \begin{align*}
                (D^2\pi)_{y_e}: \Sym^2\ker(D\pi)_{y_e} \to \coker(D\pi)_{y_e}
            \end{align*}
            and it descends, via the projection $\Sym^2T^{\Zar}_{C_s,y_e} \to \Lambda_{C_s,y_e}$, to an isomorphism
            \begin{align}\label{eqn:versal-smoothing-parameter-iso}
                \Theta_{\xi,e}:\Lambda_{C_s,y_e}\xrightarrow{\sim} \clN_{\Delta_{\xi,e}/S}|_s.
            \end{align}

            \item\label{prop-subpart:versal-curve-1-goes-to-(1,1)} 
            Let $p\in C_s$ be a point and let $(U\to X,u,\Delta_{\xi,e})$ be a triple representing the branch $\Delta_{\xi,e}$ as in Definition~\ref{def-part:branch}. Then $C\times_S\Delta_{\xi,e}\subset C\times_SU$ is an nc divisor which is smooth at $(p,u)$ if $p\ne y_e$ and has exactly two branches at $(p,u)$ if $p = y_e$.
        \end{enumerate}
        
        \item\label{prop-part:versal-curve-naturality} The bijection \eqref{eqn:versal-edge-branch-bijection} and the isomorphisms \eqref{eqn:versal-smoothing-parameter-iso} are compatible with pullbacks along morphisms $\xi'\to\xi$ in $\fM_{g,n}$ which lie over smooth morphisms $S'\to S$ in $\Sch$.
    \end{enumerate}
\end{Proposition}
\begin{proof}
    (We only sketch the proof since this is a well-known result.)    
    Let $s\in S$ be a given point. Enumerate the edges of $\plC_s$ as $e_1,\ldots,e_r$. Using \cite[\href{https://stacks.math.columbia.edu/tag/0CBY}{Tag~0CBY}]{stacks-project}, replace $S$ by an \'etale neighborhood of $s$ and, for each $1\le i\le r$, choose
    \begin{enumerate}[label = (\alph*), itemsep = 0.3ex]
        \item an \'etale neighborhood $U_i\to C$ together with a point $\hat y_i\in U_i$ lying over $y_{e_i}$,
        
        \item a regular function $t_i$ on $S$ vanishing at $s$, and
        
        \item an \'etale $S$-morphism
        \begin{align}\label{eqn:etale-local-model-curve-family}
            U_i\to\relSpec_S\left(\frac{\clO_S[z,w]}{(zw-t_i)}\right)
        \end{align}
        mapping $\hat y_i$ to the point lying over $s$ given by $z = w = 0$.
    \end{enumerate}

    \begin{claim}\label{claim:versal-smoothing-parameters-coords}
        The vector-valued function $(t_1,\ldots,t_r):S\to\bC^r$ is a smooth morphism at $s$.    
    \end{claim}
    \begin{proof}[Proof of claim]
        (The argument relies on some facts from the deformation theory of local complete intersections; see \cite[4.4]{vistoli-lci} for the relevant background.) Since $(C/S,{\bf x})$ is a versal family, so is $C/S$. Combining the versality of $C/S$ with the fact that the obstructions to deforming $C_s$ lie in $\Ext^2(\Omega^1_{C_s},\clO_{C_s}) = 0$, we conclude that $S$ is smooth at $s$. 
        
        Using the \'etale morphisms $U_i\to C$ and \eqref{eqn:etale-local-model-curve-family}, we may identify the differential of the vector-valued function $(t_1,\ldots,t_r)$ at the point $s$ with the composition
        \begin{align*}
            T_{S,s}\to\Ext^1(\Omega^1_{C_s},\clO_{C_s})\to\bigoplus_{1\le i\le r}\left(\clE xt^1(\Omega^1_{C_s},\clO_{C_s})\right)_{y_i},
        \end{align*}
        where the first arrow is the Kodaira--Spencer map of $C/S$ at $s$, which is surjective by versality, while the second arrow is the natural restriction map, which is surjective by the local-to-global Ext spectral sequence. It follows that $dt_1,\ldots,dt_r$ are linearly independent at $s$.
    \end{proof}

    With Claim~\ref{claim:versal-smoothing-parameters-coords} in hand, we complete the proof as follows. For each $1\le i\le r$, define $\Delta_{\xi,e_i}$ by $t_i = 0$; this defines the bijection \eqref{eqn:versal-edge-branch-bijection}. Using the \'etale morphism \eqref{eqn:etale-local-model-curve-family}, we see that the map \eqref{eqn:versal-smoothing-parameter-iso} induced by the second derivative of $\pi$ at $y_{e_i}$ is given by
    \begin{align*}
        \textstyle\left.\frac{\partial}{\partial z}\right|_{y_{e_i}}\otimes\left.\frac{\partial}{\partial w}\right|_{y_{e_i}}\in\Lambda_{C_s,y_{e_i}}\mapsto \left.\frac{\partial}{\partial t_i}\right|_s\in\clN_{\Delta_{\xi,e_i}/S}|_{s}
    \end{align*}
    and thus, it is a linear isomorphism.
\end{proof}

\begin{Corollary}\label{cor:nodal-curves-nc}
    The stacks $\fM_{g,n}$ and $\fC_{g,n}$ are smooth over $\bC$. Moreover, the substacks $\Delta_{\fM_{g,n}}$ and $\Delta_{\fC_{g,n}}$ are nc divisors in $\fM_{g,n}$ and $\fC_{g,n}$ respectively.
\end{Corollary}

\subsection{Families of log curves}\label{subsec:log-curve-families}

\begin{Definition}\label{def:log-stack-of-curves}
    The log stacks 
    \begin{align*}
        \fM_{g,n}^\dagger = (\fM_{g,n},M_{\fM_{g,n}})
        \quad\text{and}\quad
        \fC_{g,n}^\dagger = (\fC_{g,n},M_{\fC_{g,n}})
    \end{align*}
    are defined by endowing the smooth stacks $\fM_{g,n}$ and $\fC_{g,n}$ with the divisorial log structure Definition~\ref{def:div-log-str-smooth-nc-stack} associated to their respective nc divisors $\Delta_{\fM_{g,n}}$ and $\Delta_{\fC_{g,n}}$.
    
    Since the restriction of $\Pi_{g,n}:\fC_{g,n}\to\fM_{g,n}$ to the open substack $\fC_{g,n}^\times\subset\fC_{g,n}$ factors through the open substack $\fM_{g,n}^\times\subset\fM_{g,n}$, Example~\ref{exa:div-log-str-smooth-nc-stack-mor} shows that $\Pi_{g,n}$ lifts uniquely to a log morphism
    \begin{align*}
        \Pi_{g,n}^\dagger = (\Pi_{g,n},\Pi_{g,n}^\flat):\fC_{g,n}^\dagger\to\fM_{g,n}^\dagger.
    \end{align*}
\end{Definition}

\begin{Definition}\label{def:family-of-log-curves}
    A \emph{family $\xi^\dagger = (\xi,\xi^\flat)$ of $n$-pointed log curves of genus $g$} over a log scheme $S^\dagger = (S,M_S)$ in $\LogSch$ is defined to be a log morphism $\xi^\dagger=(\xi,\xi^\flat):S^\dagger\to\fM_{g,n}^\dagger$. 
    
    Explicitly, using Definition~\ref{def-part:log-mor-to-stack}, this means $\xi^\dagger$ consists of a family of prestable curves $\xi = (C/S,{\bf x})\in\fM_{g,n}(S)$, called the \emph{underlying family} of $\xi^\dagger$, together with a map of log structures $\xi^\flat:\xi^*M_{\fM_{g,n}}\to M_S$.
    
    We call $\xi^\dagger$ a \emph{basic} family of log curves if the map $\xi^\flat$ is an isomorphism of log structures, i.e., if $\xi^\dagger:S^\dagger\to\fM_{g,n}^\dagger$ is a strict log morphism.
\end{Definition}

\begin{Remark}\label{rem:log-curve-defs}
    In the setting of fs log schemes, it is more standard to define a family of log curves as a log smooth, integral morphism $\pi^\dagger:C^\dagger\to S^\dagger$ of fs log schemes such that the fibres of the underlying morphism $\pi:C\to S$ are reduced connected $1$-dimensional schemes; see \cite[Definition~1.1]{FKato-logcurve} or \cite[4.5]{log-geometry-and-moduli}. 
    
    The standard definition of log curve families is equivalent to Definition~\ref{def:family-of-log-curves} in the setting of fs log schemes; see \cite[Appendix~A]{GS-log} or \cite[Appendix~B]{Chen-logDF1}. The key point is that both these approaches to families of log curves have a notion of \emph{basic family} (sometimes the term \emph{minimal family} is used instead) which satisfies the following two properties.
    \begin{enumerate}[label = (\arabic*), itemsep = 0.3ex]
        \item Any family $\xi$ of prestable curves uniquely lifts to a basic family of log curves.
        \item Any family $\xi^\dagger = (\xi,\xi^\flat)$ of log curves over a log scheme $S^\dagger = (S,M_S)$ is canonically the pullback of the basic family corresponding to $\xi$, along a log morphism whose underlying morphism is given by $\text{id}_S:S\to S$.
    \end{enumerate}
\end{Remark}

\begin{Definition}\label{def:family-of-log-curves-total-space}
    Let $\xi^\dagger = (\xi,\xi^\flat)$ be a family of $n$-pointed log curves of genus $g$ over a log scheme $S^\dagger = (S,M_S)$ in $\LogSch$, with $\xi = (\pi:C\to S,{\bf x})$ being its underlying family. Using \eqref{eqn:classifying-map-curve-family}, we define the log structure $M_C$ and the maps of log structures
    \begin{equation*}
    \begin{tikzcd}
        \pi^*M_S \arrow[r,"\pi^\flat"] & M_C & \arrow[l,swap,"\hat\xi^\flat"] \hat\xi^*M_{\fC_{g,n}}
    \end{tikzcd}
    \end{equation*}
    by the following pushout diagram in the category of log structures on $C$. Explicitly, we first compute the pushout in the category of sheaves of monoids (Remark~\ref{rem:pushout-monoid-sheaves}) to get a pre-log structure and then pass to the associated log structure (Definition~\ref{def:assoc-log-str}).
    \begin{equation}\label{eqn:family-of-log-curves-total-space}
    \begin{tikzcd}
        \pi^*(\xi^*M_{\fM_{g,n}}) \arrow[r,equals,"\eqref{eqn:log-curve-implicit-identification}"] \arrow[d,"\pi^*\xi^\flat"] & \hat\xi^*(\Pi_{g,n}^*M_{\fM_{g,n}}) \arrow[r,"\hat\xi^*\Pi_{g,n}^\flat"] \arrow[dr,phantom,"\ulcorner",pos=0.99] & \hat\xi^*M_{\fC_{g,n}} \arrow[d,"\hat\xi^\flat"] \\
        \pi^*M_S \arrow[rr,"\pi^\flat"] & & M_C
    \end{tikzcd}
    \end{equation}
    The identification in the top row of \eqref{eqn:family-of-log-curves-total-space} uses the isomorphism 
    \begin{align}\label{eqn:log-curve-implicit-identification}
        \pi^*(\xi^*M_{\fM_{g,n}})\xrightarrow{\sim}\hat\xi^*(\Pi_{g,n}^*M_{\fM_{g,n}})    
    \end{align}
    of log structures induced by the functor $M_{\fM_{g,n}}:\fM_{g,n}\to\Log$ and 
    the natural projection morphism $\Pi_{g,n}(\hat\xi)=\pi^*\xi\to\xi$ in $\fM_{g,n}$ lying over $\pi:C\to S$ in $\Sch$; see Definition~\ref{def-part:log-str-on-stack}. 
    
    We show in Proposition~\ref{prop:log-curve-family-total-space-fine} below that $M_C$ is a fine log structure.
    
    We call $M_C$ the \emph{log structure on the total space} of the family $\xi^\dagger$ and denote the resulting log scheme by $C^\dagger = (C,M_C)$. Additionally, the maps $\pi^\flat$ and $\hat\xi^\flat$ define log morphisms 
    \begin{align*}
        \pi^\dagger = (\pi,\pi^\flat):C^\dagger\to S^\dagger\quad\text{and}\quad\hat\xi^\dagger=(\hat\xi,\hat\xi^\flat):C^\dagger\to\fC_{g,n}^\dagger.
    \end{align*}
\end{Definition}

\begin{Proposition}\label{prop:log-curve-family-total-space-fine}
    In the setting of Definition~\ref{def:family-of-log-curves-total-space}, the log structure $M_C$ is fine.
\end{Proposition}
\begin{proof}
    Let $s\in S$ be a given point. Write $C_s:=\pi^{-1}(s)$ for the fibre over $s$. We will construct log charts for $M_C$ in an \'etale neighborhood of any given point in $C_s$, which will show that $M_C$ is fine. The argument is divided into three steps.
    \begin{step}\label{proof-step:emb-into-versal}
        After replacing $S$ by an \'etale neighborhood of $s$, choose a morphism $(\iota,\hat\iota):\xi\to\xi'$ in the category $\fM_{g,n}$ such that $\xi' = (\pi':C'\to S',{\bf x}')$ is a versal family lying over some $S'\in\Sch$. Write $s':=\iota(s)\in S'$ and let $C'_{s'}$ be the fibre over $s'$. We may assume that $C'$ is a scheme after replacing $S'$ by an \'etale neighborhood of $s'$; see Remark~\ref{rem:why-alg-space}. Define the families $\hat\xi\in\fC_{g,n}(C)$ and $\hat\xi'\in\fC_{g,n}(C')$ as in Remark~\ref{rem:prestable-curves-algebraic-stack} using $\xi$ and $\xi'$ respectively. 
        
        The situation is summarized by the following $2$-commutative diagram, in which each square is a $2$-fibre product diagram. The $2$-isomorphisms on the faces of this diagram are the evident ones; they have been left implicit for the sake of readability.

        \begin{equation}\label{eqn:versal-family-emb}
        \begin{tikzcd}
            & C \arrow[dd,"\pi"] \arrow[rrr,"\hat\xi"] \arrow[dr,"\hat\iota"] & & & \fC_{g,n} \arrow[dd,"\Pi_{g,n}"] \\
            & & C' \arrow[urr,"\hat\xi'"] & & \\
            \phantom{} \arrow[r,phantom,"s\,\in"] & S \arrow[rrr,"\xi" near end] \arrow[dr,"\iota"] \arrow[uu,"{\bf x}",bend left] & & & \fM_{g,n} \arrow[uu,"{\bf x}_{g,n}",bend left]\\
            & \phantom{} \arrow[r,phantom,"s'\,\in"] & S' \arrow[urr,"\xi'"] \arrow[from = uu,crossing over,"\pi'" near start] \arrow[uu,crossing over,"{\bf x}'",bend left,near end] & &
        \end{tikzcd}
        \end{equation}
        
        By Proposition~\ref{prop-part:versal-curve-smooth-nc}, the schemes $S'$ and $C'$ are smooth and the pullbacks 
        \begin{align*}
            \Delta_{\xi'} := S'\times_{\fM_{g,n}}\Delta_{\fM_{g,n}}\subset S' \quad\text{and} \quad\Delta_{\hat\xi'} := C'\times_{\fC_{g,n}}\Delta_{\fC_{g,n}}\subset C'    
        \end{align*}
        are nc divisors satisfying $C'\times_{S'}\Delta_{\xi'}\subset\Delta_{\hat\xi'}$. From Example~\ref{exa:div-log-str}, we get divisorial log structures $M_{(S',\Delta_{\xi'})}\subset\clO_{S'}$, $M_{(C',\Delta_{\hat\xi'})}\subset\clO_{C'}$ and a unique map of log structures
        \begin{align}\label{eqn:versal-nc-divisorial}
            \pi'\hspace{0.05em}\phantom{}^\flat: \pi'^*M_{(S,\Delta_{\xi'})}\to M_{(C',\Delta_{\hat\xi'})}.
        \end{align} 
        \begin{claim}\label{claim:basic-log-str-via-versal}
            The diagram \eqref{eqn:versal-family-emb} induces isomorphisms of log structures
            \begin{align}\label{eqn:basic-log-str-via-versal-1}
                \Psi_{\xi,\iota}:\xi^*M_{\fM_{g,n}}\xrightarrow{\sim}\iota^*M_{(S',\Delta_{\xi'})} \quad\text{and}\quad \Psi_{\hat\xi,\hat\iota}:\hat\xi^*M_{\fC_{g,n}}\xrightarrow{\sim}\hat\iota^*M_{(C',\Delta_{\hat\xi'})}
            \end{align}
            such that the following diagram commutes.
            \begin{equation}\label{eqn:basic-log-str-via-versal-2}
            \begin{tikzcd}
                \pi^*(\xi^*M_{\fM_{g,n}}) \arrow[r,equals,"\eqref{eqn:log-curve-implicit-identification}"] \arrow[d,"\pi^*\Psi_{\xi,\iota}","\mathbin{\rotatebox[origin=c]{90}{$\sim$}}"'] & \hat\xi^*(\Pi_{g,n}^*M_{\fM_{g,n}}) \arrow[rr,"\hat\xi^*\Pi_{g,n}^\flat"] & & \hat\xi^*M_{\fC_{g,n}} \arrow[d,"\Psi_{\hat\xi,\hat\iota}","\mathbin{\rotatebox[origin=c]{90}{$\sim$}}"'] \\
                \pi^*(\iota^*M_{(S',\Delta_{\xi'})}) \arrow[r,equals] & \hat\iota^*(\pi'^*M_{(S',\Delta_{\xi'})}) \arrow[rr,"\hat\iota^*\pi'\hspace{0.01em}\phantom{}^\flat"] & & \hat\iota^*M_{(C',\Delta_{\hat\xi'})}
            \end{tikzcd}
            \end{equation}
        \end{claim}
        \begin{proof}[Proof of claim]
            Since $\xi'$ is versal, both $\xi':S'\to\fM_{g,n}$ and $\hat\xi':C'\to\fC_{g,n}$ are smooth. Thus, the structure maps $\xi'^*M_{\fM_{g,n}}\to\clO_{S'}$ and $\hat\xi'^*M_{\fC_{g,n}}\to\clO_{C'}$ induce isomorphisms
            \begin{align}\label{eqn:basic-log-str-on-versal-1}
                \Psi_{\xi'}:\xi'^*M_{\fM_{g,n}}\xrightarrow{\sim} M_{(S',\Delta_{\xi'})} \quad\text{and}\quad
                \Psi_{\hat\xi'}:\hat\xi'^*M_{\fC_{g,n}}\xrightarrow{\sim} M_{(C',\Delta_{\hat\xi'})}
            \end{align}
            by definition; see Definitions~\ref{def:log-stack-of-curves} and \ref{def:div-log-str-smooth-nc-stack}. Since \eqref{eqn:versal-nc-divisorial} is the unique map of log structures from $\pi'^*M_{(S',\Delta_{\xi'})}$ to $M_{(C',\Delta_{\hat\xi'})}$, the following diagram must commute.
            \begin{equation}\label{eqn:basic-log-str-on-versal-2}
            \begin{tikzcd}
                \pi'^*(\xi'^*M_{\fM_{g,n}}) \arrow[r,equals,"\eqref{eqn:log-curve-implicit-identification}"] \arrow[d,"\pi'^*\Psi_{\xi'}","\mathbin{\rotatebox[origin=c]{90}{$\sim$}}"'] & \hat\xi'^*(\Pi_{g,n}^*M_{\fM_{g,n}}) \arrow[rr,"\hat\xi'^*\Pi_{g,n}^\flat"] && \hat\xi'^*M_{\fC_{g,n}} \arrow[d,"\Psi_{\hat\xi'}","\mathbin{\rotatebox[origin=c]{90}{$\sim$}}"'] \\
                \pi'^*M_{(S',\Delta_{\xi'})} \arrow[rrr,"\pi'\hspace{0.01em}\phantom{}^\flat"] &&& M_{(C',\Delta_{\hat\xi'})}
            \end{tikzcd}
            \end{equation}
            Pulling back \eqref{eqn:basic-log-str-on-versal-1} along $(\iota,\hat\iota):\xi\to\xi'$ defines the isomorphisms \eqref{eqn:basic-log-str-via-versal-1}. Finally, \eqref{eqn:basic-log-str-via-versal-2} commutes since it is the pullback of \eqref{eqn:basic-log-str-on-versal-2} along $(\iota,\hat\iota):\xi\to\xi'$.
        \end{proof}
    \end{step}

    \begin{step}\label{proof-step:new-pushout-expression} 
        Define the maps of log structures
        \begin{align*}
            \iota^\flat &:= \xi^\flat\circ\Psi_{\xi,\iota}^{-1}:\iota^*M_{(S',\Delta_{\xi'})}\to M_S,\quad\text{and}\\
            \hat\iota^\flat &:= \hat\xi^\flat\circ\Psi_{\hat\xi,\hat\iota}^{-1}:\hat\iota^*M_{(C',\Delta_{\hat\xi'})}\to M_C
        \end{align*}
        using the isomorphisms from \eqref{eqn:basic-log-str-via-versal-1}. Combining \eqref{eqn:family-of-log-curves-total-space} with \eqref{eqn:basic-log-str-via-versal-2}, we obtain the following pushout diagram in the category of log structures on $C$.
        \begin{equation}\label{eqn:general-log-str-via-versal}
        \begin{tikzcd}
            \pi^*(\iota^*M_{(S',\Delta_{\xi'})}) \arrow[r,equals] \arrow[d,"\pi^*\iota^\flat"] & \hat\iota^*(\pi'^*M_{(S',\Delta_{\xi'})}) \arrow[r,"\hat\iota^*\pi'\hspace{0.01em}\phantom{}^\flat"] \arrow[dr,phantom,"\ulcorner",pos=0.99] & \hat\iota^*M_{(C',\Delta_{\hat\xi'})} \arrow[d,"\hat\iota^\flat"] \\
            \pi^*M_S \arrow[rr,"\pi^\flat"] & & M_C
        \end{tikzcd}
        \end{equation}
    \end{step}
    
    \begin{step}\label{proof-step:log-curve-charts} 
        After replacing $S$ by an \'etale neighborhood of $s$, choose a fine monoid $Q$ and a log chart $\chi: Q\to\Gamma(S,M_S)$. This is possible since $S^\dagger$ belongs to $\LogSch$; see Notation~\ref{nota:logsch-variants} and Definition~\ref{def-part:chart-fine}. By \cite[III.1.2.7(2)]{ogus-log-book}, we may ensure that the log chart $\chi$ induces an isomorphism of monoids $\ol\chi_s:Q\xrightarrow{\sim}\ol M_{S,s}$. In particular, $Q$ is sharp.
        
        Let $p\in C_s$ be a given point and write $p':=\hat\iota(p)\in C'_{s'}$. We will construct a log chart for $M_C$ in an \'etale neighborhood of $p$. During the construction, we will replace $C$ (resp. $C'$) repeatedly by \'etale neighborhoods of $p$ (resp. $p'$) but, for convenience, we will not introduce any new notation for them. We consider three cases depending on whether $p$ is a general (i.e., smooth non-marked) point, marked point or node of $C_s$.
        
        \begin{enumerate}[label = (\alph*), ref = \thestep(\alph*), left = 0pt]
            \item\label{proof-step-part:general-point}
            {\bf General point.} After replacing $C'$ by a Zariski open neighborhood of $p'$, we may use Proposition~\ref{prop-part:versal-curve-smooth-nc} to ensure that $\pi'$ is smooth and $\Delta_{\hat\xi'} = C'\times_{S'}\Delta_{\xi'}$. 
            
            By Lemma~\ref{lem-part:pullback-of-nc-div-log}, the map $\pi'\hspace{0.05em}\phantom{}^\flat$ is an isomorphism and, by pulling back along $\hat\iota$, the top row of \eqref{eqn:general-log-str-via-versal} is also an isomorphism. It follows that $\pi^\flat:\pi^*M_S\to M_C$ is an isomorphism of log structures and we have a log chart given by
            \begin{align*}
                \hat\chi_\text{gen}:Q&\to\Gamma(C,M_C),\\
                q&\mapsto\pi^\flat(\chi(q)).
            \end{align*}
            
            \item\label{proof-step-part:marked-point} 
            {\bf Marked point.} Let $x:S\to C$ (resp. $x':S'\to C'$) be the marked point section containing $p$ (resp. $p'$). After replacing $C'$ by a Zariski open neighborhood of $p'$, we may use Proposition~\ref{prop-part:versal-curve-smooth-nc} to ensure that $\pi'$ is smooth and $\Delta_{\hat\xi'} = C'\times_{S'}\Delta_{\xi'} + x'(S')$. 
            Choose a local equation $z$ for the smooth divisor $x'(S')\subset C'$. Note that $z$ is a section of $M_{(C',\Delta_{\hat\xi'})}\subset\clO_{C'}$. We have a pre-log structure $\ul\bN\to\clO_{C'}$ given by $a\mapsto z^a$. 
            
            Since $\pi'$ is smooth and $\Delta_{\hat\xi'} = C'\times_{S'}\Delta_{\xi'} + x'(S')$ is an nc divisor, combining the two parts of Lemma~\ref{lem:pullback-of-nc-div-log} shows that the map of pre-log structures
            \begin{align*}
               \pi'^*M_{(S',\Delta_{\xi'})}\oplus\ul\bN\to M_{(C',\Delta_{\hat\xi'})},\quad (\sigma,a)\mapsto\pi'\hspace{0.05em}\phantom{}^\flat(\sigma)\cdot z^a 
            \end{align*}
            induces an isomorphism after passing to the associated log structures. Using \eqref{eqn:general-log-str-via-versal}, it now follows that we have a log chart given by
            \begin{align*}
                \hat\chi_x:Q\oplus\bN&\to\Gamma(C,M_C),\\
                (q,a)&\mapsto\pi^\flat(\chi(q))\cdot \hat\iota^\flat(z^a).
            \end{align*}
            The log charts $\hat\chi_x$ and $\hat\chi_\text{gen}$ are related in the following way. Define the regular function $Z:=\hat\iota^\#(z)$ on $C$. Note that $Z$ is a local equation for the Cartier divisor $x(S)\subset C$. On the complement of the closed subscheme $Z = 0$, we have
            \begin{align}\label{eqn:log-curve-chart-marked-generize}
                \hat\chi_x(q,a) = \hat\chi_\text{gen}(q)\cdot Z^a.
            \end{align}
            
            \item\label{proof-step-part:node} 
            {\bf Node.} Following Notation~\ref{nota:prestable-curve}, let us change the notation from $p$ (resp. $p'$) to $y$ (resp. $y'$). After replacing $C'$ by a Zariski open neighborhood of $y'$, we may use Proposition~\ref{prop-part:versal-curve-smooth-nc} to ensure that $\Delta_{\hat\xi'} = C'\times_{S'}\Delta_{\xi'}$. Now, using \cite[\href{https://stacks.math.columbia.edu/tag/0CBY}{Tag~0CBY}]{stacks-project}, replace $S'$ (resp. $C'$) by an \'etale neighborhood of $s'$ (resp. $y'$) and choose a regular function $t$ on $S'$ vanishing at $s'$ together with an \'etale $S'$-morphism
            \begin{align}\label{eqn:basic-log-curve-node-local-model}
                C' \to \relSpec_{S'}\left(\frac{\clO_{S'}[z,w]}{(zw-t)}\right)
            \end{align}
            mapping the node $y'\in C'_{s'}$ to the point lying over $s'$ given by $z = w = 0$. Denote the images of $z$, $w$ in $\clO_{C'}$ by the same letters. Note that $t$ is a section of $M_{(S',\Delta_{\xi'})}\subset\clO_{S'}$, and $z$, $w$ are sections of $M_{(C',\Delta_{\hat\xi'})}\subset\clO_{C'}$. We have pre-log structures
            \begin{align*}
                \ul\bN\to\clO_{S'},\quad a\mapsto t^a\quad\text{and}\quad
                \ul\bN^2\to\clO_{C'},\quad (a,b)\mapsto z^aw^b.
            \end{align*}
            
            Since $\Delta_{\hat\xi'} = C'\times_{S'}\Delta_{\xi'}$ is an nc divisor, combining Proposition~\ref{prop-subpart:versal-curve-1-goes-to-(1,1)} and Lemma~\ref{lem-part:nc-div-log-chart} shows that the map of pre-log structures
            \begin{align*}
                \pi'^*M_{(S',\Delta_{\xi'})}\oplus_{\ul\bN}\ul\bN^2\to M_{(C',\Delta_{\hat\xi'})},\quad [\sigma,(a,b)]\mapsto\pi'\hspace{0.05em}\phantom{}^\flat(\sigma)\cdot z^a\cdot w^b
            \end{align*}
            induces an isomorphism after passing to the associated log structures, where we define the amalgamated sum using the maps $a\mapsto\pi'\hspace{0.05em}\phantom{}^\flat(t^a)$ and $a\mapsto (a,a)$.

            After replacing $S$ by an \'etale neighborhood of $s$, write $\iota^\flat(t) = \chi(\rho)\cdot\gamma^{-1}$ for some $\rho\in Q\xrightarrow{\sim}\ol M_{S,s}$ and some section $\gamma$ of $\clO_S^\times$. Since $t$ vanishes at $s'$, we have $\rho\ne 0$.

            \begin{claim}\label{claim:node-monoid-fine}
                Consider the monoid maps 
                \begin{align*}
                    \bN\to Q,\quad 1\mapsto\rho\quad\text{and}\quad\bN\to\bN^2,\quad 1\mapsto(1,1).
                \end{align*}
                Then the resulting amalgamated sum $Q\oplus_\bN\bN^2$ is a fine monoid.
            \end{claim}
            \begin{proof}[Proof of claim]
                Since $Q$ is finitely generated, so is $Q\oplus_\bN\bN^2$. We are left to show that $Q\oplus_\bN\bN^2$ is integral. This follows from the equivalence of conditions (1) and (3) in \cite[I.4.6.2]{ogus-log-book} for the monoid map $\bN\to\bN^2$ given by $1\mapsto (1,1)$.
            \end{proof}
            Choose sections $\lambda$, $\mu$ of $\clO_C^\times$ such that $\lambda\mu = \pi^\#(\gamma)$, i.e., first choose $\lambda$ arbitrarily and then take $\mu:=\pi^\#(\gamma)\cdot \lambda^{-1}$. Using \eqref{eqn:general-log-str-via-versal} and the fine monoid $Q\oplus_\bN\bN^2$ from Claim~\ref{claim:node-monoid-fine}, it now follows that we have a well-defined log chart given by
            \begin{align*}
                \hat\chi_y:Q\oplus_\bN\bN^2 &\to\Gamma(C,M_C),\\
                [q,(a,b)] &\mapsto\pi^\flat(\chi(q))\cdot(\hat\iota^\flat(z)\lambda)^a\cdot(\hat\iota^\flat(w)\mu)^b.
            \end{align*}
            The log charts $\hat\chi_y$ and $\hat\chi_\text{gen}$ are related in the following way. Define the regular functions $Z:=\hat\iota^\#(z)\lambda$ and $W:=\hat\iota^\#(w)\mu$ on $C$. Note that $ZW$ equals the image of $\chi(\rho)$ in $\clO_S$. On the complement of the closed subscheme $Z = W = 0$, we have
            \begin{align}\label{eqn:log-curve-chart-node-generize}
                \hat\chi_y([q,(a,b)]) = 
                \begin{cases}
                    \hat\chi_\text{gen}(q+a\rho)\cdot W^{b-a} & \text{ if }W\ne 0, \\
                    \hat\chi_\text{gen}(q+b\rho)\cdot Z^{a-b} & \text{ if }Z\ne 0.
                \end{cases}
            \end{align}
            Although it is not necessary for the present proof, we record the following claim here because of its very close relation to the formulas appearing in \eqref{eqn:log-curve-chart-node-generize}.
            \begin{claim}\label{claim:log-curve-node-ghost-stalk}
                The monoid map 
                \begin{align*}
                    Q\oplus_\bN\bN^2\to Q\times Q,\quad [q,(a,b)]\mapsto(q+a\rho,q+b\rho)    
                \end{align*}
                is well-defined and injective. Its image consists of all pairs $(q,q')\in Q$ such that there exists a (necessarily unique) $c\in\bZ$ satisfying $q' = q + c\rho$.
            \end{claim}
            \begin{proof}[Proof of claim]
                This is exactly the content of \cite[Remark~1.2(2)]{GS-log}.
            \end{proof}
        \end{enumerate}
    \end{step}
    This completes the proof that $M_C$ is a fine log structure.
\end{proof}

\begin{Remark}
    The log charts for $M_C$ constructed in the proof of Proposition~\ref{prop:log-curve-family-total-space-fine} recover the \'etale local description of a family of log curves given in \cite[page~222]{FKato-logcurve}.
\end{Remark}

\begin{Remark}
    In the setting of Proposition~\ref{prop:log-curve-family-total-space-fine}, suppose that $S = \Spec\bC$ is a point and that the dual graph $\plC$ of $C$ has a loop, i.e., an edge both of whose endpoints are the same vertex. It can then be shown that it is impossible to completely specify the log structure $M_C$ by a single slb on $C$, unlike the DF log structures considered in Subsection~\ref{subsec:DF-log}. Instead, $M_C$ is specified by an slb on $C$ indexed by the sheaf $\ol M_C$; see Remark~\ref{rem:log-str-ghost-slb}.
\end{Remark}

\subsection{Log curves over log points}
For this subsection, we take $S = \Spec\bC$ to be a point and fix an $n$-pointed prestable genus $g$ curve, i.e., a family $\xi=(\pi:C\to S,{\bf x})\in\fM_{g,n}(S)$.

The following result gives a concrete description and classification of families $\xi^\dagger$ of log curves over log points $S^\dagger = (S,M_S)$ whose underlying family is $\xi$. We freely use Notation~\ref{nota:prestable-curve} in its statement and proof.

\begin{Theorem}\label{thm:log-curve-over-log-pt}
    Fix $S = \Spec\bC$ and $\xi = (\pi:C\to S,{\bf x})\in\fM_{g,n}(S)$ as above. For a fine log structure $M_S$ on $S$, write $S^\dagger = (S,M_S)$ and consider the following two types of objects.
    \begin{enumerate}[itemsep = 0.3ex]
        \myitem[(LC$\phantom{}_1$)]\label{abstract-log-curve} Family $\xi^\dagger$ of $n$-pointed log curves of genus $g$ over $S^\dagger$ with underlying curve $\xi$.
        
        \myitem[(LC$\phantom{}_2$)]\label{concrete-log-curve} Collection of elements $0\ne \rho_e\in \ol M_S$ and linear isomorphisms $\Theta_e:\Lambda_{C,y_e}\xrightarrow{\sim}\clO_{S^\dagger}(\rho_e)$, where $e$ ranges over the edges of the dual graph $\plC$ of $C$.
    \end{enumerate}
    As the log structure $M_S$ varies, each of \ref{abstract-log-curve} and \ref{concrete-log-curve} determines a category fibred in groupoids over $\LogPt$, using the obvious notion of pullbacks along log morphisms. 
    
    These two categories fibred in groupoids over $\LogPt$ are equivalent.

    Given a family $\xi^\dagger$ over $S^\dagger = (S,M_S)$ as in \ref{abstract-log-curve}, let $\{\rho_e,\Theta_e\}$ be the collection as in \ref{concrete-log-curve} corresponding to it under this equivalence. Using $\xi^\dagger$, define the log scheme $C^\dagger = (C,M_C)$ and the log morphism $\pi^\dagger=(\pi,\pi^\flat):C^\dagger\to S^\dagger$ as in Definition~\ref{def:family-of-log-curves-total-space}.
    Then we have the following.
    \begin{enumerate}[label = (\arabic*), ref = \theTheorem(\arabic*), itemsep = 0.3ex]
        \item\label{thm-part:basic-log-curve} The family $\xi^\dagger$ is basic if and only if the collection $\{\rho_e\}$ freely generates $\ol M_S$.
        
        \item\label{thm-part:log-curve-gen-point} The map $\ol\pi^\flat:\pi^*\ol M_S\to\ol M_C$ is an isomorphism over $C_\text{gen}\subset C$. In particular, we have a well-defined restriction map $\Gamma(C,\ol M_C)\to\Gamma(C_v\cap C_\text{gen},\ol M_C) = \ol M_S$ for each $v$.
        
        \item\label{thm-part:log-curve-marked-point} We have a canonical identification $\ol M_{C,x_i} = \ol M_S\oplus\bN$ for each $i$, under which the map $\ol\pi\phantom{}^\flat_{x_i}:\ol M_S\to \ol M_{C,x_i}$ is given by the inclusion $\ol M_S\subset\ol M_S\oplus\bN$. In particular, we have a well-defined restriction map $\Gamma(C,\ol M_C)\to\ol M_{C,x_i}/\ol M_S = \bN$ for each $i$.
        
        \item\label{thm-part:ghost-sections-log-curve} The restriction maps in (2) and (3) together define an isomorphism
        \begin{align}\label{eqn:ghost-sections-log-curve}
            \Gamma(C,\ol M_C) \xrightarrow{\sim} \left\{\begin{matrix}
                \{m_v\in\ol M_S\}_v \\
                \{m_i\in\bN\}_i
            \end{matrix}\,\middle|\,\begin{matrix}
                \forall\,\vec{e}:v\to v',\exists!\,m_{\vec{e}}\in\bZ\\
                \text{with }m_{v'} = m_v + m_{\vec{e}}\rho_e
            \end{matrix}\right\}.
        \end{align} 
        
        \item\label{thm-part:global-slb-log-curve} Given $m\in\Gamma(C,\ol M_C)$, apply \eqref{eqn:ghost-sections-log-curve} to $m$ to obtain $\{m_v\}_v$, $\{m_i\}_i$ and $\{m_{\vec{e}}\}_{\vec{e}}$. 
        \begin{enumerate}[label = (\roman*), ref = \theenumi(\roman*), itemsep = 0.15ex]
            \item\label{thm-subpart:log-curve-slb-on-irr-comps} 
            For each $v$, we have a Cartier divisor $D_{m,v}$ on $\tilde C_v$, a line bundle $L_{m,v}$ on $\tilde C_v$ and a well-defined global section $\sigma_{m,v}\in H^0(\tilde C_v,L_{m,v})$ given by the formulas
            \begin{align*}
                D_{m,v} &:= \textstyle\sum_{x_i\in \tilde C_v}m_i\cdot x_i+\sum_{y_{\vec{e}}\in\tilde C_v}m_{\vec{e}}\cdot y_{\vec{e}},\\
                L_{m,v} &:= \clO_{\tilde C_v}(D_{m,v})\otimes(\pi\circ\nu|_{\tilde C_v})^*\clO_{S^\dagger}(m_v),\quad\text{and}\\
                \sigma_{m,v} &:= \taut_{D_{m,v}}\otimes (\pi\circ\nu|_{\tilde C_v})^*\sigma_{S^\dagger,m_v}.
            \end{align*}
            In particular, $D_{m,v}\subset\tilde C_v$ is effective if $m_v = 0$ and we have $\sigma_{m,v} = 0$ if $m_v\ne 0$.
            
            \item\label{thm-subpart:log-curve-slb-at-nodes} 
            For each $\vec{e}:v\to v'$, define $\psi_{m,\vec{e}}:L_{m,v}|_{y_{\vec{e}}}\xrightarrow{\sim} L_{m,v'}|_{y_{\cev{e}}}$ to be the linear isomorphism corresponding to $\Theta_e^{\otimes m_{\vec{e}}}$ under the following chain of identifications                
            \begin{align*}
                \Hom\left(L_{m,v}|_{y_{\vec{e}}}, L_{m,v'}|_{y_{\cev{e}}}\right) &\xrightarrow{\sim} \Hom\left(T_{\tilde C,y_{\vec{e}}}^{\otimes m_{\vec{e}}}\otimes\clO_{S^\dagger}(m_v),T_{\tilde C,y_{\cev{e}}}^{\otimes m_{\cev{e}}}\otimes\clO_{S^\dagger}(m_{v'})\right) \\
                &\xrightarrow{\sim} \Hom\left(\Lambda_{C,y_e}^{\otimes m_{\vec{e}}},\clO_{S^\dagger}(m_{v'}-m_v)\right) \\
                &\xrightarrow{\sim} \Hom\left(\Lambda_{C,y_e},\clO_{S^\dagger}(\rho_e)\right)^{\otimes m_{\vec{e}}}.
            \end{align*}
            These identifications are induced by the relations $m_{\vec{e}}+m_{\cev{e}} = 0$, $m_{v'} = m_v + m_{\vec{e}}\rho_e$ and the isomorphisms appearing in the global slb on $S^\dagger$.
            
            \item Using the isomorphisms $\psi_{m,\vec{e}}$ to glue the line bundles $L_{m,v}$ on $\tilde C_v$ at the nodes $y_e$, we obtain a line bundle $L_m$ on $C$. Moreover, the sections $\sigma_{m,v}$ together yield a well-defined global section $\sigma_m\in H^0(C,L_m)$.

            \item Repeating this for any given pair $m_1,m_2\in\Gamma(C,\ol M_C)$, we get an isomorphism $\Phi_{m_1,m_2}:L_{m_1}\otimes L_{m_2}\xrightarrow{\sim}L_{m_1+m_2}$ induced by the isomorphisms appearing in the global slb on $S^\dagger$ and multiplication of rational functions on $C$.
        \end{enumerate}
        \vspace{1ex}
        The resulting slb $\{L_m,\Phi_{m_1,m_2},\sigma_m\}$ is naturally isomorphic to the global slb on $C^\dagger$.
    \end{enumerate}
\end{Theorem}

\begin{Remark}\label{rem:curve-divisor-fibre}
    For the first map in the sequence of identifications in Theorem~\ref{thm-subpart:log-curve-slb-at-nodes}, we use the canonical isomorphism $\clO_{\tilde C}(p)|_p\xrightarrow{\sim}T_{\tilde C,p}$ for any point $p\in\tilde C$; see \eqref{eqn:normal-bundle}.
\end{Remark}

\begin{Remark}\label{rem:log-curves-whats-new}
    We explain the relation of Theorem~\ref{thm:log-curve-over-log-pt} to the existing literature. For a family of log curves $\xi^\dagger$ over an fs log point $S^\dagger = (S,M_S)$, its \emph{tropicalization} $\Trop(\xi^\dagger)$ can be defined as a family of \emph{abstract tropical curves} (i.e., metric graphs) parametrized by the dual cone $\Hom_{\Mon}(\ol M_S,\bR_{\ge 0})$. The tropicalization encodes the same information as the collection $\{\rho_e\}$. Thus, the first four parts of Theorem~\ref{thm:log-curve-over-log-pt} may be rephrased in terms of $\Trop(\xi^\dagger)$. For details on tropicalizing a family of log curves, see \cite[Section~7.2]{moduli-trop-curves}.
    
    To the best of the authors' knowledge, the equivalence \ref{abstract-log-curve} $\simeq$ \ref{concrete-log-curve} has not appeared before in the literature, though it is possibly known to experts. The key new insight is that the collection of isomorphisms $\{\Theta_e\}$ is exactly the additional data needed to reconstruct the log curve $\xi^\dagger$ from its underlying curve $\xi$ and its tropicalization $\Trop(\xi^\dagger)$. A similar statement appears in \cite[Theorem~4.4]{FRTU-chip-firing}; however, it is incomplete (\cite{Ranganathan-personal-comm}) as it claims that $\xi^\dagger$ can be reconstructed from just $\xi$ and $\Trop(\xi^\dagger)$. The collection $\{\rho_e,\Theta_e\}$ allows us to give a completely explicit description of the global slb on the total space of a log curve $\xi^\dagger$; this is the content of Theorem~\ref{thm-part:global-slb-log-curve}.
\end{Remark}

\subsection{Proof of Theorem~\ref{thm:log-curve-over-log-pt}} 

We are given a prestable curve $\xi = (\pi:C\to S,{\bf x})\in\fM_{g,n}(S)$ over $S = \Spec\bC$. We will draw heavily on the proof of Proposition~\ref{prop:log-curve-family-total-space-fine} and so, the reader is advised to review that first. The argument is divided into six steps.

\begin{Step}\label{proof-step:emb-into-versal-point-case}
    \stepcounter{subsubsection}
    This is the same as Step~\ref{proof-step:emb-into-versal} in the proof of Proposition~\ref{prop:log-curve-family-total-space-fine} for the case $S = \Spec\bC$, but we provide a short summary in the next paragraph for the reader's convenience. 
    
    Choose a morphism $(\iota,\hat\iota):\xi\to\xi'$ in $\fM_{g,n}$ such that $\xi' = (\pi':C'\to S',{\bf x}')$ is a versal family with $S',C'\in\Sch$. Let $s$ denote the unique point of $S$, let $s'\in S'$ be the image of $s$ under $\iota$ and let $C'_{s'}$ be the fibre of $C'\to S'$ over $s'$. Define the families $\hat\xi\in\fC_{g,n}(C)$ and $\hat\xi'\in\fC_{g,n}(C')$ as in Remark~\ref{rem:prestable-curves-algebraic-stack} using $\xi$ and $\xi'$ respectively. 
    From this setup, we obtain the $2$-commutative diagram \eqref{eqn:versal-family-emb}, in which each square is a $2$-fibre product diagram. As before, we also obtain nc divisors $\Delta_{\xi'}\subset S'$, $\Delta_{\hat\xi'}\subset C'$ and the unique map $\pi'\hspace{0.05em}^\flat$ of log structures as in \eqref{eqn:versal-nc-divisorial}. Claim~\ref{claim:basic-log-str-via-versal} applies verbatim to our situation.
\end{Step}

\begin{Step}\label{proof-step:equivalence-log-curve}
    \stepcounter{subsubsection}
    In this step, we prove the equivalence \ref{abstract-log-curve}~$\simeq$~\ref{concrete-log-curve}. To produce the equivalence, we will go via an intermediate category fibred in groupoids over $\LogPt$, denoted by \ref{semi-concrete-log-curve}, whose objects over a log point $S^\dagger = (S,M_S)$ are of the following type.
    \begin{enumerate}
        \myitem[(LC$\phantom{}_3$)]\label{semi-concrete-log-curve} Map of log structures $\iota^\flat:\iota^*M_{(S',\Delta_{\xi'})}\to M_S$.
    \end{enumerate}
    
    For each edge $e$ of the dual graph $\plC$, denote the image of the node $y_e\in C$ under $\hat\iota:C\to C'$ by $y_e'\in C'$. Using $\hat\iota$, we identify $C$ with the fibre $C'_{s'}$. In particular, the graph $\plC$ (resp. the vector space $\Lambda_{C,y_e}$) is identified with the dual graph of $C'_{s'}$ (resp. the vector space $\Lambda_{C'_{s'},y'_e}$). We make the following preliminary remarks.
    \begin{enumerate}[label = R\arabic*., ref = R\arabic*, itemsep = 0.3ex]
        \item Proposition~\ref{prop-part:versal-curve-basic-log} provides a bijection 
        \begin{align}\label{eqn:edge-smoothing-parameter}
            e\mapsto\Delta_{\xi',e}    
        \end{align}
        from the edges of $\plC$ to the branches of $\Delta_{\xi'}$ at $s'$.
        
        \item For each edge $e$, the cokernel of the first derivative $(D\pi')_{y_e'}$ is the normal space $\clN_{\Delta_{\xi',e}/S'}|_{s'}$, and the second derivative $(D^2\pi')_{y_e'}$ induces an isomorphism
        \begin{align}\label{eqn:normal-direction-smoothing}
            \Theta_{\xi',e}:\Lambda_{C'_{s'},y_e'} \xrightarrow{\sim} \clN_{\Delta_{\xi',e}/S'}|_{s'}.
        \end{align}

        \item\label{versal-basic-log-basis} Lemma~\ref{lem-part:ghost-sheaf-stalk-free} applied to the bijection \eqref{eqn:edge-smoothing-parameter} provides a basis $\{\rho_{\Delta_{\xi',e}}\}_e$ for the free monoid $\iota^*\ol M_{(S',\Delta_{\xi'})}$, where $e$ ranges over the edges of $\plC$. 
        
        \item\label{versal-basic-log-vcd} By Proposition~\ref{prop:nc-log-str-at-point}, for each edge $e$, the virtual Cartier divisor associated to the basis vector $\rho_{\Delta_{\xi',e}}\in\iota^*\ol M_{(S',\Delta_{\xi'})}$ is given by $(\clN_{\Delta_{\xi',e}/S'}|_{s'},0)$.
    \end{enumerate}
    
    Using these remarks, we now define the desired equivalences as follows.

    \ref{abstract-log-curve}~$\to$~\ref{semi-concrete-log-curve}: Given a family $\xi^\dagger = (\xi,\xi^\flat)$ of log curves over $S^\dagger = (S,M_S)$ with underlying family $\xi$, we get a map of log structures $\xi^\flat:\xi^*M_{\fM_{g,n}}\to M_S$. Using the isomorphism $\Psi_{\xi,\iota}$ in \eqref{eqn:basic-log-str-via-versal-1}, we get a map $\iota^\flat:\iota^*M_{(S',\Delta_{\xi'})}\to M_S$ of log structures defined by $\iota^\flat:=\xi^\flat\circ\Psi_{\xi,\iota}^{-1}$. 

    The assignment $\xi^\dagger\mapsto\iota^\flat$ is clearly functorial and invertible.

    \ref{semi-concrete-log-curve}~$\to$~\ref{concrete-log-curve}: Given a log point $S^\dagger = (S,M_S)$ and a map $\iota^\flat:\iota^* M_{(S',\Delta_{\xi'})}\to M_S$ of log structures, define the element $\rho_e$ and the isomorphism $\Theta_e$ for each edge $e$ as follows.

    Define $\rho_e:=\ol\iota^\flat(\rho_{\Delta_{\xi',e}})\in\ol M_S$, where $\rho_{\Delta_{\xi',e}}$ is the basis vector from \ref{versal-basic-log-basis}. We have $\rho_e\ne 0$ since $(\ol\iota^\flat)^{-1}(0) = \{0\}$; see Definition~\ref{def-part:ghost-sheaf}. By Lemma~\ref{lem:map-of-log-str-slb} and \ref{versal-basic-log-vcd}, the map $\iota^\flat$ induces an isomorphism $\clN_{\Delta_{\xi',e}/S'}|_{s'} \xrightarrow{\sim} \clO_{S^\dagger}(\rho_e)$. Define $\Theta_e:\Lambda_{C,y_e}\xrightarrow{\sim}\clO_{S^\dagger}(\rho_e)$ to be the composition
    \begin{align*}
        \Lambda_{C,y_e} \xrightarrow{\sim}\Lambda_{C'_{s'},y'_e} \xrightarrow{\sim} \clN_{\Delta_{\xi',e}/S'}|_{s'} \xrightarrow{\sim} \clO_{S^\dagger}(\rho_e),
    \end{align*}
    where the first map is the identification induced by $\hat\iota:C\to C'$, the second map is the isomorphism \eqref{eqn:normal-direction-smoothing} and the third map is the isomorphism induced by $\iota^\flat$.

    The assignment $\iota^\flat\mapsto\{\rho_e,\Theta_e\}$ is clearly functorial. It is invertible by Proposition~\ref{prop:cat-of-log-pts}.
\end{Step}

\begin{Step}
    \stepcounter{subsubsection}
    In this step, we prove part (1). By definition, a family $\xi^\dagger = (\xi,\xi^\flat)$ of log curves is basic if and only if $\xi^\flat:\xi^*M_{\fM_{g,n}}\to M_S$ is an isomorphism, i.e., if and only if $\xi^\dagger$ is a final object of the category \ref{abstract-log-curve}. On the other hand, a collection $\{\rho_e,\Theta_e\}$ is a final object of the category \ref{concrete-log-curve} if and only if $\ol M_S$ is freely generated by the collection $\{\rho_e\}$. We are done since final objects correspond to final objects under any equivalence of categories.
\end{Step}

\begin{Step}\label{proof-step:log-curve-tropical-data}
    \stepcounter{subsubsection}
    In this step, we prove parts (2), (3) and (4). Fix an object $\xi^\dagger = (\xi,\xi^\flat)$ of \ref{abstract-log-curve} lying over a log point $S^\dagger = (S,M_S)$. Let $\{\rho_e,\Theta_e\}$ and $\iota^\flat$ denote the corresponding objects of \ref{concrete-log-curve} and \ref{semi-concrete-log-curve} respectively under the chain of equivalences \ref{abstract-log-curve}~$\simeq$~\ref{semi-concrete-log-curve}~$\simeq$~\ref{concrete-log-curve} defined in Step~\ref{proof-step:equivalence-log-curve}. Define the map of $\hat\iota^\flat:\hat\iota^*M_{(C',\Delta_{\hat\xi'})}\to M_C$ log structures by $\hat\iota^\flat := \hat\xi^\flat\circ\Psi_{\hat\xi,\hat\iota}^{-1}$. 
    
    From Step~\ref{proof-step:new-pushout-expression} in the proof of Proposition~\ref{prop:log-curve-family-total-space-fine}, we get the diagram \eqref{eqn:general-log-str-via-versal} expressing the log structure $M_C$ as a pushout. We now prove parts (2), (3) and (4) as follows.

    (2): By Step~\ref{proof-step-part:general-point} in the proof of Proposition~\ref{prop:log-curve-family-total-space-fine}, $\pi^\flat:\pi^*M_S\to M_C$ is an isomorphism over $C_\text{gen}\subset C$. Passing to ghost sheaves, we conclude that $\ol\pi^\flat$ is an isomorphism over $C_\text{gen}$.

    (3): This is immediate from Step~\ref{proof-step-part:marked-point} in the proof of Proposition~\ref{prop:log-curve-family-total-space-fine}.

    (4): This is a consequence of Step~\ref{proof-step-part:node} in the proof of Proposition~\ref{prop:log-curve-family-total-space-fine}. In more detail, using parts (2) and (3) together with \eqref{eqn:log-curve-chart-node-generize}, we see that the map \eqref{eqn:ghost-sections-log-curve} is well-defined. Using Claim~\ref{claim:log-curve-node-ghost-stalk}, we conclude that the map \eqref{eqn:ghost-sections-log-curve} is bijective.
\end{Step}

\begin{Step}\label{proof-step:log-curve-slb-description}
    \stepcounter{subsubsection}
    In this step, we prove part (5). We continue with exactly the same setup and notation as in Step~\ref{proof-step:log-curve-tropical-data}. Let $m\in\Gamma(C,\ol M_C)$ be a given global section and let $m_v\in\ol M_S$, $m_i\in\bN$ and $m_{\vec{e}}\in\bZ$ be the elements obtained by applying \eqref{eqn:ghost-sections-log-curve} to $m$. From $m$, we also get the associated virtual Cartier divisor $(\clO_{C^\dagger}(m),\sigma_{C^\dagger,m})$.
    
    For each vertex $v$ of $\plC$, consider the difference $\check m_v := m - \ol\pi^\flat(m_v)\in\Gamma(C,\ol M\phantom{}^\text{gp}_C)$ and the associated line bundle $\clO_{C^\dagger}(\check m_v)$. Using the isomorphisms appearing in the global slb on $C^\dagger$, we will identify $\clO_{C^\dagger}(\check m_v)$ with the line bundle $\clO_{C^\dagger}(m)\otimes\pi^*\clO_{S^\dagger}(m_v)^\vee$. By definition, $\check m_v$ vanishes on $C_v\cap C_\text{gen}$ and thus, from Discussion~\ref{disc-part:log-line-bundle-trivialized}, we get a trivialization
    \begin{align}\label{eqn:global-slb-log-curve-interior-triv}
        \clO_{C_v\cap C_\text{gen}}\xrightarrow{\sim}\clO_{C^\dagger}(\check m_v)|_{C_v\cap C_\text{gen}}.
    \end{align}
    View $C_v\cap C_\text{gen}$ as a Zariski open subset of $\tilde C_v$, so that \eqref{eqn:global-slb-log-curve-interior-triv} defines a rational section of the line bundle $(\nu|_{\tilde C_v})^*\clO_{C^\dagger}(\check m_v)$; define $E_{m,v}$ to be the divisor of zeros and poles of this rational section. Then $E_{m,v}$ is a (not necessarily effective) Cartier divisor on $\tilde C_v$ supported on the complement of $C_v\cap C_\text{gen}$. By definition, \eqref{eqn:global-slb-log-curve-interior-triv} extends to an isomorphism
    \begin{align}\label{eqn:global-slb-log-curve-twisted-triv}
        \clO_{\tilde C_v}(E_{m,v})\xrightarrow{\sim}(\nu|_{\tilde C_v})^*\clO_{C^\dagger}(\check m_v)
    \end{align}
    which maps $\taut_{E_{m,v}}$ to the rational section defined by \eqref{eqn:global-slb-log-curve-interior-triv}; see Discussion~\ref{disc:sheaf-cartier-div}. 
    
    \begin{claim}\label{claim:log-curve-slb-description}
        Recall the divisor $D_{m,v}$, the line bundle $L_{m,v}$, the section $\sigma_{m,v}$ and the isomorphism $\psi_{m,\vec{e}}$ from the statement of part (5). These have the following properties.
        \begin{enumerate}[label = P\arabic*., ref = P\arabic*, itemsep = 0.3ex]
            \item\label{log-curve-vertex-iso} For each vertex $v$ of $\plC$, we have the equality $E_{m,v} = D_{m,v}$ of Cartier divisors on $\tilde C_v$. In particular, \eqref{eqn:global-slb-log-curve-twisted-triv} induces an isomorphism of line bundles
            \begin{align*}
                \Psi_{m,v}: L_{m,v} \xrightarrow{\sim} (\nu|_{\tilde C_v})^*\clO_{C^\dagger}(m).
            \end{align*}
            
            \item\label{log-curve-edge-iso} For each oriented edge $\vec{e}:v\to v'$ of $\plC$, the following diagram commutes.
            \begin{equation*}
            \begin{tikzcd}
                L_{m,v}|_{q_{\vec{e}}} \arrow[rr,"\Psi_{m,v}|_{q_{\vec{e}}}","\sim"'] \arrow[dd,"\psi_{m,\vec{e}}","\mathbin{\rotatebox[origin=c]{90}{$\sim$}}"'] & & \clO_{C^\dagger}(m)|_{\nu(q_{\vec{e}})} \arrow[d,equals]\\
                & & \clO_{C^\dagger}(m)|_{q_e} \\
                L_{m,v'}|_{q_{\cev{e}}} \arrow[rr,"\Psi_{m,v}|_{q_{\cev{e}}}","\sim"'] & & \clO_{C^\dagger}(m)|_{\nu(q_{\cev{e}})} \arrow[u,equals]
            \end{tikzcd}
            \end{equation*}

            \item\label{log-curve-section-compatible} For each vertex $v$ of $\plC$, the isomorphism $\Psi_{m,v}$ maps $\sigma_{m,v}$ to $(\nu|_{\tilde C_v})^*\sigma_{C^\dagger,m}$.
        \end{enumerate}
    \end{claim}
    \begin{proof}[Proof of claim]
        We will check \ref{log-curve-vertex-iso} (resp. \ref{log-curve-vertex-iso} and \ref{log-curve-edge-iso}) over a Zariski open (resp. an \'etale) neighborhood of each marked point (resp. node) of $C$. Note that \ref{log-curve-section-compatible} follows from \ref{log-curve-vertex-iso} since $\Psi_{m,v}$ maps $\sigma_{m,v}$ to $(\nu|_{\tilde C_v})^*\sigma_{C^\dagger,m}$ over the dense open subset $C_v\cap C_\text{gen}$ of $\tilde C_v$.
        
        To check \ref{log-curve-vertex-iso} and \ref{log-curve-edge-iso}, we will use the log charts from Step~\ref{proof-step:log-curve-charts} in the proof of Proposition~\ref{prop:log-curve-family-total-space-fine}. Accordingly, fix a fine monoid $Q$ and a log chart $\chi:Q\to M_S$ which induces an isomorphism of monoids $\ol\chi:Q\xrightarrow{\sim}\ol M_S$. For the argument below, we will work over a Zariski open (resp. an \'etale) neighborhood of the marked point (resp. node) under consideration; for convenience, we will not introduce any new notation for this neighborhood and simply call it $C$.
        
        \begin{enumerate}[label = (\alph*), itemsep = 2ex, left = 0pt]
            \item {\bf Marked point.} Let $x\in C$ be a given marked point. Let $x':S'\to C'$ denote the marked point section containing $\hat\iota(x)\in C'_{s'}$. Choose a local equation $z$ for $x'(S')\subset C'$ and define the regular function $Z:=\hat\iota^\#(z)$ on $C$. Note that $Z$ is a local coordinate at $x\in C$. For the remainder, we will work in a Zariski open neighborhood of $x\in C$ over which the log chart $\hat\chi_x$ from Step~\ref{proof-step-part:marked-point} in the proof of Proposition~\ref{prop:log-curve-family-total-space-fine} is defined.
            
            Let $v$ be the vertex of $\plC$ for which $x\in C_v$ and let $1\le i\le n$ be such that $x = x_i$. Express $m\in\Gamma(C,\ol M_C)$ as the image of the section $\hat\chi_x(q,a)\in\Gamma(C,M_C)$ for some choice of $(q,a)\in Q\oplus\bN$. By \eqref{eqn:log-curve-chart-marked-generize} and the definition of $m_v\in\ol M_S$, we have $\ol\chi(q) = m_v$. Moreover, by the definition of $m_i\in\bN$, we have $a = m_i$. 
            
            Therefore, $\check m_v = m - \ol\pi^\flat(m_v)$ is the image of            
            \begin{align}\label{eqn:log-curve-marked-triv}
                \frac{\hat\chi_x(q,a)}{\hat\chi_x(q,0)} = \hat\chi_x(0,a) = (\hat\iota^\flat(z))^{m_i}
            \end{align}
            under $M\phantom{}^\text{gp}_C\to\ol M\phantom{}^\text{gp}_C$.
            It follows that $\check m_v$ is a section of the subsheaf $\ol M_C\subset\ol M\phantom{}^\text{gp}_C$. From this, we obtain a well-defined section $\sigma_{C^\dagger,\check m_v}$ of the line bundle $\clO_{C^\dagger}(\check m_v)$.
            \\[1ex]
            \emph{Checking \ref{log-curve-vertex-iso} at $x$.}
            \\[1ex]
            By Discussion~\ref{disc-part:log-line-bundle-trivialized-by-section}, the trivialization \eqref{eqn:global-slb-log-curve-interior-triv} is induced by the restriction of $\sigma_{C^\dagger,\check m_v}$ to the open subset $C_v\cap C_\text{gen}$. On the other hand, by Discussion~\ref{disc-part:log-lb-vcd-lift-induced-trivialization}, the section \eqref{eqn:log-curve-marked-triv} induces an isomorphism of virtual Cartier divisors
            \begin{align*}
                (\clO_{C^\dagger}(\check m_v),\sigma_{C^\dagger,\check m_v})\xrightarrow{\sim}(\clO_C,Z^{m_i}).    
            \end{align*} 
            From this, we conclude that the coefficient of the marked point $x = x_i$ in the divisor $E_{m,v}$ is $m_i$. This completes the proof that \ref{log-curve-vertex-iso} holds at $x$.
            
            \item {\bf Node.} Given a node $y\in C$, define $y':=\hat\iota(y)\in C'_{s'}$. After replacing $S'$ (resp. $C'$) by an \'etale neighborhood of $s'$ (resp. $y'$), choose an \'etale local model \eqref{eqn:basic-log-curve-node-local-model} for $\pi':C'\to S'$ given by regular functions $z$, $w$ (resp. $t$) on $C'$ (resp. $S'$) vanishing at $y'$ (resp. $s'$) and satisfying $zw = t$. Let $0\ne\rho\in Q$ and $\gamma\in\bC^\times$ be such that we have 
            \begin{align}\label{eqn:smoothing-parameter-decomp}
                \iota^\flat(t) = \chi(\rho)\cdot\gamma^{-1}.
            \end{align}
            Choose $\lambda,\mu\in\bC^\times$ such that $\lambda\mu = \gamma$ and define the regular functions $Z:=\hat\iota^\#(z)\lambda$ and $W:=\hat\iota^\#(w)\mu$ on $C$. For the remainder, we will work in an \'etale neighborhood of $y$ over which the log chart $\hat\chi_y$ from Step~\ref{proof-step-part:node} in the proof of Proposition~\ref{prop:log-curve-family-total-space-fine} is defined.

            Let $e$ be the edge of $\plC$ for which $y = y_e$. Lift $e$ to an oriented edge $\vec{e}:v\to v'$. Since $m_{\vec{e}} + m_{\cev{e}} = 0$, after reversing the orientation of $\vec{e}$ if necessary, we may assume that the integer $m_{\vec{e}}$ is nonnegative. After swapping the labels $z$, $w$ if necessary, we may also assume that $\tilde C_v$ (resp. $\tilde C_{v'}$) is given by $Z = 0$ (resp. $W = 0$). Note that $W$ (resp. $Z$) restricts to a local coordinate at $y_{\vec{e}}\in\tilde C_v$ (resp. $y_{\cev{e}}\in\tilde C_{v'}$). 
            
            Since the branch $\Delta_{\xi',e}$ of the nc divisor $\Delta_{\xi'}\subset S'$ is defined by the equation $t = 0$, it follows from \eqref{eqn:smoothing-parameter-decomp} that we have $\ol\chi(\rho) = \rho_e$. By Discussion~\ref{disc-part:log-lb-vcd-lift-induced-trivialization}, the element $\chi(\rho)\in M_S$ induces a linear isomorphism 
            \begin{align*}
                \clO_{S^\dagger}(\rho_e)\xrightarrow{\sim}\clO_S = \bC.    
            \end{align*}
            Define $\omega\in\clO_{S^\dagger}(\rho_e)$ to be the inverse image of $1\in\bC$ under this isomorphism. By \eqref{eqn:smoothing-parameter-decomp}, the isomorphism $\clN_{\Delta_{\xi',e}/S'}|_{s'}\xrightarrow{\sim}\clO_{S^\dagger}(\rho_e)$ induced by $\iota^\flat$ is given by $\left.\frac{\partial}{\partial t}\right|_{s'}\mapsto\gamma\omega$. Recall from the proof of Proposition~\ref{prop:versal-curve} that the isomorphism \eqref{eqn:normal-direction-smoothing} is given by $\left.\frac{\partial}{\partial z}\right|_{y'}\otimes\left.\frac{\partial}{\partial w}\right|_{y'}\mapsto\left.\frac{\partial}{\partial t}\right|_{s'}$. Since $\left.\frac{\partial}{\partial Z}\right|_y$ and $\left.\frac{\partial}{\partial W}\right|_y$ map to $\lambda^{-1}\left.\frac{\partial}{\partial z}\right|_{y'}$ and $\mu^{-1}\left.\frac{\partial}{\partial w}\right|_{y'}$ respectively under $\hat\iota:C\to C'$, it follows that  $\Theta_e:\Lambda_{C,y_e}\xrightarrow{\sim}\clO_{S^\dagger}(\rho_e)$ is given by
            \begin{align}\label{eqn:log-curve-theta-iso-coords}
                \Theta_e\left(\textstyle\left.\frac{\partial}{\partial Z}\right|_y\otimes\left.\frac{\partial}{\partial W}\right|_y\right) = \omega.
            \end{align}
            
            Express $m\in\Gamma(C,\ol M_C)$ as the image of the section $\hat\chi_y([q,(a,b)])\in\Gamma(C,M_C)$ for some choice of $[q,(a,b)]\in Q\oplus_\bN\bN^2$. By \eqref{eqn:log-curve-chart-node-generize} and the definitions of $m_v$, $m_{v'}\in\ol M_S$, we have $\ol\chi(q+a\rho) = m_v$ and $\ol\chi(q+b\rho) = m_{v'}$. Moreover, by the definition of $m_{\vec{e}}\in\bZ$ and our choice of orientation on the edge $e$, we have $b - a = m_{\vec{e}}\ge 0$. 
            
            Therefore, $\check m_v = m - \ol\pi^\flat(m_v)$ and $-\check m_{v'} = \ol\pi^\flat(m_{v'}) - m$ are the images of
            \begin{align}\label{eqn:log-curve-node-triv-1}
                \frac{\hat\chi_y([q,(a,b)])}{\hat\chi_y([q+a\rho,(0,0)])} = \hat\chi_y([0,(0,b-a)]) &= (\hat\iota^\flat(w)\mu)^{m_{\vec{e}}},\quad\text{and} \\[1.5ex]
                \label{eqn:log-curve-node-triv-2}
                \frac{\hat\chi_y([q+b\rho,(0,0)])}{\hat\chi_y([q,(a,b)])} = \hat\chi_y([0,(b-a,0)]) &= (\hat\iota^\flat(z)\lambda)^{m_{\vec{e}}}
            \end{align}
            respectively under $M\phantom{}_C^\text{gp}\to\ol M\phantom{}_C^\text{gp}$. It follows that $\check m_v$ and $-\check m_{v'}$ are sections of the subsheaf $\ol M_C\subset\ol M\phantom{}^\text{gp}_C$. From this, we obtain well-defined sections $\sigma_{C^\dagger,\check m_v}$ and $\sigma_{C^\dagger,-\check m_{v'}}$ of the line bundles $\clO_{C^\dagger}(\check m_v)$ and $\clO_{C^\dagger}(-\check m_{v'})$ respectively. 
            \\[1ex]
            \noindent\emph{Checking \ref{log-curve-vertex-iso} at $y$.} 
            \\[1ex]
            By Discussion~\ref{disc-part:log-line-bundle-trivialized-by-section}, the trivialization \eqref{eqn:global-slb-log-curve-interior-triv} is induced by the restriction of $\sigma_{C^\dagger,\check m_v}$ to the open subset $C_v\cap C_\text{gen}$. Similarly, we also have the corresponding statement for the trivialization of $\clO_{C^\dagger}(-\check m_{v'})$ over $C_{v'}\cap C_\text{gen}$ defined by dualizing the analogue of \eqref{eqn:global-slb-log-curve-interior-triv} for the vertex $v'$. On the other hand, by Discussion~\ref{disc-part:log-lb-vcd-lift-induced-trivialization}, the sections \eqref{eqn:log-curve-node-triv-1} and \eqref{eqn:log-curve-node-triv-2} induce isomorphisms of virtual Cartier divisors
            \begin{align}
                \label{eqn:vertex-v-triv}
                (\clO_{C^\dagger}(\check m_v),\sigma_{C^\dagger,\check m_v})&\xrightarrow{\sim}(\clO_C,W^{m_{\vec{e}}}),\quad\text{and} \\
                \label{eqn:vertex-v'-triv}
                (\clO_{C^\dagger}(-\check m_{v'}),\sigma_{C^\dagger,-\check m_{v'}})&\xrightarrow{\sim}(\clO_C,Z^{m_{\vec{e}}})
            \end{align}
            respectively. From this, we conclude that the coefficient of $y_{\vec{e}}$ (resp. $y_{\cev{e}}$) in the divisor $E_{m,v}$ (resp. $E_{m,v'}$) is $m_{\vec{e}}$ (resp. $m_{\cev{e}}$). This completes the proof that \ref{log-curve-vertex-iso} holds at $y$.
            \\[1ex]
            \emph{Checking \ref{log-curve-edge-iso} at $y$.}
            \\[1ex]
            Consider the isomorphism \eqref{eqn:global-slb-log-curve-twisted-triv} and the dual of its analogue for the vertex $v'$. By restricting these to $y_{\vec{e}}\in\tilde C_v$ and $y_{\cev{e}}\in\tilde C_{v'}$ respectively, we obtain isomorphisms
            \begin{align}\label{eqn:log-curve-section-leading-term}
                \varphi_{m,\vec{e}}:T_{\tilde C,y_{\vec{e}}}^{\otimes m_{\vec{e}}}\xrightarrow{\sim}\clO_{C^\dagger}(\check m_v)|_{y_e} \quad\text{and}\quad
                \varphi_{m,\cev{e}}:T_{\tilde C,y_{\cev{e}}}^{\otimes m_{\vec{e}}}\xrightarrow{\sim}\clO_{C^\dagger}(-\check m_{v'})|_{y_e}.
            \end{align}
            The tensor product of the maps in \eqref{eqn:log-curve-section-leading-term} yields an isomorphism
            \begin{align}\label{eqn:log-curve-section-leading-term-product}
                \varphi_{m,\vec{e}}\otimes\varphi_{m,\cev{e}}:\Lambda_{C,y_e}^{\otimes m_{\vec{e}}}\xrightarrow{\sim}(\clO_{C^\dagger}(\check m_v)\otimes \clO_{C^\dagger}(-\check m_{v'}))|_{y_e}.    
            \end{align}
            For \ref{log-curve-edge-iso} to hold at $y$, the maps \eqref{eqn:log-curve-section-leading-term-product} and $\Theta_e^{\otimes m_{\vec{e}}}$ must agree under the identification
            \begin{align}\label{eqn:log-curve-tensor-edge-iso}
                (\clO_{C^\dagger}(\check m_v)\otimes \clO_{C^\dagger}(-\check m_{v'}))|_{y_e}\xrightarrow{\sim}\clO_{S^\dagger}(\rho_e)^{\otimes m_{\vec{e}}}
            \end{align}
            induced by the relation $\check m_v - \check m_{v'} = m_{\vec{e}}\cdot\ol\pi^\flat(\rho_e)$ and the isomorphisms appearing in the global slb on $C^\dagger$. We check this by a direct computation in the next paragraph.
            
            Define the vectors $\Omega\in\clO_{C^\dagger}(\check m_v)|_{y_e}$ and $\Omega'\in\clO_{C^\dagger}(-\check m_{v'})|_{y_e}$ to be the inverse images of $1\in\clO_C|_{y_e} = \bC$ under the isomorphisms \eqref{eqn:vertex-v-triv} and \eqref{eqn:vertex-v'-triv} respectively. Then the two isomorphisms in \eqref{eqn:log-curve-section-leading-term} are given by 
            \begin{align*}
                \varphi_{m,\vec{e}}\left(\textstyle(\left.\frac{\partial}{\partial W}\right|_{y})^{\otimes m_{\vec{e}}}\right) = \Omega \quad\text{and}\quad \varphi_{m,\cev{e}}\left(\textstyle(\left.\frac{\partial}{\partial Z}\right|_{y})^{\otimes m_{\vec{e}}}\right) = \Omega'.    
            \end{align*}
            Looking at the expression \eqref{eqn:log-curve-theta-iso-coords} for $\Theta_e$ found earlier, we see that it suffices to show that the identification \eqref{eqn:log-curve-tensor-edge-iso} is given by $\Omega\otimes\Omega'\mapsto\omega^{\otimes m_{\vec{e}}}$. However, this is immediate from Discussion~\ref{disc-part:log-lb-vcd-lift-induced-triv-product}, the relation $(\hat\iota^\flat(z)\lambda)\cdot(\hat\iota^\flat(w)\mu) = \pi^\flat(\chi(\rho))$ in $\Gamma(C,M_C)$ and the definitions of $\Omega$, $\Omega'$ and $\omega$. This completes the proof that \ref{log-curve-edge-iso} holds at $y$.
        \end{enumerate}

        This completes the proof of Claim~\ref{claim:log-curve-slb-description}.
    \end{proof}

    For each section $m\in\Gamma(C,\ol M_C)$, the isomorphisms $\Psi_{m,v}$ provided by Claim~\ref{claim:log-curve-slb-description} glue together at the nodes to define a canonical isomorphism $\Psi_m:(L_m,\sigma_m)\xrightarrow{\sim}(\clO_{C^\dagger}(m),\sigma_{C^\dagger}(m))$ of virtual Cartier divisors on $C$. The proof that these are compatible with the isomorphisms $\Phi_{m_1,m_2}:L_{m_1}\otimes L_{m_2}\xrightarrow{\sim}L_{m_1 + m_2}$ is straightforward and will be omitted.
\end{Step}

\begin{Step}
    \stepcounter{subsubsection}
    In this final (optional) step, we show that the equivalence \ref{abstract-log-curve} $\simeq$ \ref{concrete-log-curve} constructed in Step~\ref{proof-step:equivalence-log-curve} is in fact independent of the morphism $\xi\to \xi'$ to a versal family chosen in Step~\ref{proof-step:emb-into-versal-point-case}.

    Consider another morphism $\xi\to\xi''$ in $\fM_{g,n}$ lying over $S\to S''$ in $\Sch$ such that $\xi''$ is a versal family. By considering a smooth atlas of the algebraic space $S'\times_{\fM_{g,n}} S''$, we may reduce to the case where $\xi\to\xi'$ admits a factorization $\xi\to\xi''\to\xi'$ in $\fM_{g,n}$ such that the corresponding morphism $S''\to S'$ in $\Sch$ is smooth. In this special case, Lemma~\ref{lem-part:pullback-of-nc-div-log} and Proposition~\ref{prop-part:versal-curve-naturality} together show that the equivalences \ref{abstract-log-curve} $\simeq$ \ref{concrete-log-curve} constructed using the two choices $\xi\to\xi''$ and $\xi\to\xi'$ agree.
\end{Step}

This concludes the proof of Theorem~\ref{thm:log-curve-over-log-pt}.\qed
    \section{Log maps}\label{sec:log-maps}

We discuss (\emph{pre})\emph{stable log maps} to a smooth snc pair from the viewpoints of algebraic geometry (\cite{GS-log}, \cite{Chen-logDF1}, \cite{AC-logDF2}) and symplectic topology (\cite{Tehrani-log}) and compare the two versions. To make the comparison, it turns out to be necessary to distinguish between the case of fine (Theorem~\ref{thm:fine-log-map-comparison}) and fs log schemes (Theorem~\ref{thm:fs-log-map-comparison}); see Remark~\ref{rem:why-fine-vs-fs-basic}. We also briefly recall the notion of \emph{log Gromov convergence} (\cite[Definition~3.7]{Tehrani-log}) and explain how to extend it to \emph{saturated symplectic stable log maps}; see Subsection~\ref{subsec:sat-log-gromov-conv}.

\begin{Notation}
    For this section, we fix a smooth snc pair $(X,D)$; see Definition~\ref{def:smooth-snc-pair}. We also fix integers $g,n\ge 0$, a homology class $\beta\in H_2(X,\bZ)$ and an ordering $D_1,\ldots,D_r$ of the irreducible components of $D$. Denote the divisorial log structure associated to $D\subset X$ (Example~\ref{exa:div-log-str}) by $M_{(X,D)}$ and denote the resulting log scheme by $X^\dagger = (X, M_{(X,D)})$.
\end{Notation}

\subsection{Prestable maps} Recall that a \emph{prestable map} to the scheme $X$ is a tuple $(C,{\bf x},f)$, where $(C,{\bf x})$ is a prestable curve and $f:C\to X$ is a morphism. We say that it is of \emph{class} $(g,n,\beta)$ if $(C,{\bf x})$ is an $n$-pointed prestable genus $g$ curve and we have $f_*[C] = \beta$. We say that it is a \emph{stable map} if the group of automorphisms of $C$ preserving ${\bf x}$ and $f$ is finite.

\begin{Convention}
    For the rest of this section, prestable maps $(C,{\bf x},f)$ are understood to be mapping to $X$ and to be of class $(g,n,\beta)$. Similarly, (families of) prestable (or log) curves are understood to be $n$-pointed and of genus $g$.
\end{Convention}

\begin{Notation}\label{nota:prestable-map}
    Given a prestable map $(C,{\bf x},f)$, we consistently use the following notations.
    \begin{enumerate}[ref = \theNotation(\arabic*), label = (\arabic*), itemsep = 0.3ex]
        \item We follow Notation~\ref{nota:prestable-curve} for the prestable curve $(C,{\bf x})$. In particular,
        \begin{enumerate}[label = (\roman*), itemsep = 0.15ex]
            \item $\plC$ denotes the dual graph of $C$,

            \item $C_\text{gen}\subset C$ denotes the complement of the nodes and marked points,
            
            \item $v$, $e$ and $\vec{e}$ denote vertices, edges and oriented edges of $\plC$ respectively,
            
            \item $C_v$, $\tilde C_v$, $y_e$ and $y_{\vec{e}}$ denote irreducible components of $C$, connected components of $\tilde C$, nodes of $C$ and inverse images of nodes under $\tilde C\to C$ respectively, and

            \item $\Lambda_{C,y_e} = T_{\tilde C,y_{\vec{e}}}\otimes T_{\tilde C,y_{\cev{e}}}$ denotes the space of smoothing parameters at a node $y_e\in C$.
        \end{enumerate}

        \item We use $i$ (resp. $j$) to index the marked points $x_i\in C$ (resp. the irreducible components $D_j$ of $D$). In particular, $i$ and $j$ range over $\{1,\ldots,n\}$ and $\{1,\ldots,r\}$ respectively.
        
        \item\label{nota-part:restriction-degeneracy-data} Associated to each vertex $v$ of $\plC$, we have the following objects. 
        \begin{enumerate}[label = (\roman*), itemsep = 0.15ex]
            \item The pullback of $f:C\to X$ to $C_v$ (resp. $\tilde C_v$) is denoted by $f_v$ (resp. $\tilde f_v$).
            
            \item We let $I_v\subset\{1,\ldots,r\}$ denote the set of $j$ such that $f(C_v)\subset D_j$.
        \end{enumerate}
    \end{enumerate}
\end{Notation}

In the following subsections, we study \emph{stable log maps} to the pair $(X,D)$. The discussion below, which is elementary and makes no reference to log geometry, provides the motivation behind the definition of these objects; see also \cite[Appendix~A]{Chen-Church-Zhao}.

\begin{Discussion}\label{disc:log-map-motivation}
    Consider a degenerating $1$-parameter family of smooth curves which map to $X$ and have prescribed contact orders with the irreducible components of $D$ at their marked points. More precisely, consider a setup consisting of the following objects.
    \begin{enumerate}[label = (\alph*), itemsep = 0.3ex]
        \item An irreducible smooth curve $\esS$ with a distinguished point $s\in\esS$.
        
        \item A family $(\pi:\esC\to \esS,{\bf x})$ of prestable curves such that its restriction to $\esS\setminus\{s\}$ has smooth fibres. Define $C:=\pi^{-1}(s)$ and the closed subset $\Delta:=C\cup\bigcup_i x_i(\esS)\subset\esC$.

        \item A morphism $F:\esC\to X$ which maps $\esC\setminus\Delta$ to $X\setminus D$. Define $f:C\to X$ to be the restriction of the morphism $F$ to the subset $C\subset\esC$.
    \end{enumerate}    

    We freely use Notation~\ref{nota:prestable-map} for the prestable map $(C,{\bf x}(s),f)$.
    
    Choose a local coordinate $t$ at $s\in\esS$. For each edge $e$ of $\plC$, there is a unique $\rho_e\in\bN\setminus\{0\}$ such that $\esC\to\esS$ is modeled on the $A_{\rho_e-1}$ singularity $zw = t^{\rho_e}$ in an \'etale neighborhood of the node $y_e\in C\subset\esC$. Moreover, the family $\esC\to\esS$ induces a canonical isomorphism
    \begin{align*}
        \Theta_e:\Lambda_{C,y_e}\xrightarrow{\sim}T_{\esS,s}^{\otimes\rho_e},    
    \end{align*}
    which is given in coordinates by $\frac{\partial}{\partial z}\otimes\frac{\partial}{\partial w}\mapsto(\frac{\partial}{\partial t})^{\otimes\rho_e}$. 
    
    Fix $1\le j\le r$ and consider the section $\sigma_j := F^*\taut_{D_j}$ of the line bundle $F^*(\clO_X(D_j))$ on $\esC$. By assumption, the section $\sigma_j$ has no zeros in $\esC\setminus\Delta$. Note that each irreducible component of $\Delta$ is either $x_i(\esS)$ for some $1\le i\le n$ or $C_v$ for some vertex $v$ of $\plC$.

    \begin{enumerate}[label = (\arabic*), itemsep = 0.3ex]
        \item For each $1\le i\le n$, let $\mu_{i,j}\in\bN$ be the generic order of vanishing of $\sigma_j$ along $x_i(\esS)\subset\esC$. Given any point $s'\in\esS\setminus\{s\}$, if we write $\esC_{s'}:=\pi^{-1}(s')$, then we see that $\mu_{i,j}$ is the contact order of $F|_{\esC_{s'}}$ with $D_j$ at the marked point $x_i(s')$.

        \item For each vertex $v$ of $\plC$, let $\delta_{v,j}\in\bN$ be the generic order of vanishing of $\sigma_j$ along $C_v\subset\esC$. We have $\delta_{v,j} = 0$ if and only if $f(C_v)\not\subset D_j$, i.e., if and only if $j\not\in I_v$. 
        
        If we expand $\sigma_j$ in the normal direction to $C_v\cap C_\text{gen}\subset\esC$ as a power series, then the lowest order term of this expansion yields a well-defined rational section $0\ne\sigma_{v,j}$ of the line bundle  $(T_{\esS,s}^{\otimes\delta_{v,j}})^\vee\otimes_\bC\tilde f_v^*(\clO_X(D_j))$ on $\tilde C_v$. In terms of the chosen local coordinate $t$ at $s\in\esS$, we have the explicit expression
        \begin{align*}
            \sigma_{v,j} = (dt)^{\otimes\delta_{v,j}}\otimes(t^{-\delta_{v,j}}\sigma_j)|_{\tilde C_v}.
        \end{align*}
        \begin{enumerate}[label = (\roman*), itemsep = 0.2ex]
            \item By construction, $\sigma_{v,j}$ has no zeros or poles in $C_v\cap C_\text{gen}\subset\tilde C_v$.
            \item For each $1\le i\le n$ with $x_i(s)\in\tilde C_v$, the order of vanishing of $\sigma_{v,j}$ at the point $x_i(s)$ is $\mu_{i,j}\in\bN$. Indeed, if $w = 0$ is a local equation for $x_i(\esS)\subset\esC$, then $(w,t)$ are local coordinates at $x_i(s)\in\esC$ and we have $\sigma_{v,j} = t^{\delta_{v,j}}w^{\mu_{i,j}}u$ for some local section $u$ of $F^*(\clO_X(D_j))$ which is regular and non-vanishing at $x_i(s)$.

            \item For each oriented edge $\vec{e}:v\to v'$ of $\plC$, let $\mu_{\vec{e},j}\in\bZ$ be the order of vanishing of $\sigma_{v,j}$ at the point $y_{\vec{e}}\in\tilde C_v$. If we expand $\sigma_{v,j}$ as a Laurent series centred at $y_{\vec{e}}$, then the lowest order term of this expansion yields a well-defined element 
            \begin{align*}
                0\ne \sigma_{\vec{e},j}\in(T_{\tilde C,y_{\vec{e}}}^{\otimes\mu_{\vec{e},j}})^\vee\otimes(T_{\esS,s}^{\otimes\delta_{v,j}})^\vee\otimes\clO_X(D_j)|_{f(y_e)}.   
            \end{align*}
            If $w$ is a local coordinate at $y_{\vec{e}}\in\tilde C_v$, then we have the explicit expression
            \begin{align*}
                \sigma_{\vec{e},j} = (dw)^{\otimes\mu_{\vec{e},j}}\otimes (w^{-\mu_{\vec{e},j}}\sigma_{v,j})|_{y_{\vec{e}}}.
            \end{align*}
        \end{enumerate}

        \item Fix an oriented edge $\vec{e}:v\to v'$ of $\plC$. Choose a local model $zw = t^{\rho_e}$ for an \'etale neighborhood of $y_e\in\esC$. After swapping $z$, $w$ if necessary, we may assume that $\tilde C_v$ (resp. $\tilde C_{v'}$) is given by $z = t = 0$ (resp. $w = t = 0$). In particular, $w$ (resp. $z$) restricts to a local coordinate at $y_{\vec{e}}\in\tilde C_v$ (resp. $y_{\cev{e}}\in\tilde C_{v'}$). 
        
        Since $\sigma_j$ has no zeros on $\esC\setminus\Delta$, we must have $\sigma_j = t^az^bw^ch$ for some $a,b,c\in\bN$ and some local section $h$ of $F^*(\clO_X(D_j))$ which is regular and non-vanishing at $y_e$. By substituting $z = t^{\rho_e}w^{-1}$, we obtain the alternative expression $\sigma_j = t^{a+b\rho_e}w^{c-b}h$, from which we deduce that
        \begin{align*}
            \begin{array}{ll}
                \delta_{v,j} = a + b\rho_e,\quad & \sigma_{v,j} = (dt)^{\otimes(a+b\rho_e)}\otimes (w^{c-b}h)|_{z\,=\,t\,=\,0}, \\
                \mu_{\vec{e},j} = c-b,\quad & \sigma_{\vec{e},j} = (dw)^{\otimes(c-b)}\otimes(dt)^{\otimes(a+b\rho_e)}\otimes h(y_e).
            \end{array}
        \end{align*}
        Similarly, the alternative expression $\sigma_j = t^{a+c\rho_e}z^{b-c}h$ shows that 
        \begin{align*}
            \begin{array}{ll}
                \delta_{v',j} = a + c\rho_e,\quad & \sigma_{v',j} = (dt)^{\otimes(a+c\rho_e)}\otimes (z^{b-c}h)|_{w\,=\,t\,=\,0}, \\
                \mu_{\cev{e},j} = b-c,\quad & \sigma_{\cev{e},j} = (dz)^{\otimes(b-c)}\otimes(dt)^{\otimes(a+c\rho_e)}\otimes h(y_e).
            \end{array}
        \end{align*}
        These computations lead to the following conclusions.
        \vspace{0.3ex}
        \begin{enumerate}[label = (\roman*), itemsep = 0.2ex]
            \item We have the identities $\delta_{v',j} = \delta_{v,j} + \mu_{\vec{e},j}\rho_e$ and $\mu_{\vec{e},j} + \mu_{\cev{e},j} = 0$.
            \item The lowest order terms $\sigma_{\vec{e},j}$ and $\sigma_{\cev{e},j}$ agree under the identification
            \begin{align*}
                T_{\tilde C,y_{\vec{e}}}^{\otimes\mu_{\vec{e},j}}\otimes T_{\esS,s}^{\otimes\delta_{v,j}} \xrightarrow{\sim} T_{\tilde C,y_{\cev{e}}}^{\otimes\mu_{\cev{e},j}}\otimes T_{\esS,s}^{\otimes\delta_{v',j}}
            \end{align*}
            induced by the isomorphism $\Theta_e^{\otimes\mu_{\vec{e},j}}$.
        \end{enumerate}
    \end{enumerate}

    Finally, we examine what happens under a base change. Let $\esS'$ be an irreducible smooth curve with a distinguished point $s'\in\esS'$ and let $\esS'\to\esS$ be a morphism under which the scheme-theoretic inverse image of $s$ is given by $k\cdot s'$ for some $k\in\bN\setminus\{0\}$. If we pull back the entire setup along $\esS'\to\esS$ and repeat the preceding discussion, then is not hard to see that $\delta_{v,j}$ and $\rho_e$ get multiplied by $k$, while $\mu_{i,j}$ and $\mu_{\vec{e},j}$ remain unchanged.
\end{Discussion}

The remainder of this subsection is motivated by Discussion~\ref{disc:log-map-motivation}.

\begin{Definition}\label{def:contact-data}
    For a prestable map $(C,{\bf x},f)$ to $X$, a \emph{contact datum} with respect to $D$ is a collection $\bmu = \{\mu_{i,j}\in\bN,\mu_{\vec{e},j}\in\bZ\}$ such that $\mu_{\vec{e},j}+\mu_{\cev{e},j} = 0$ for all $\vec{e}$ and $j$, where $i$, $j$ and $\vec{e}$ range over $\{1,\ldots,n\}$, $\{1,\ldots,r\}$ and the oriented edges of $\plC$ respectively.
\end{Definition}

\begin{Convention}
    We omit the qualifier `with respect to $D$' and simply say `contact datum' if there is no risk of confusion. Moreover, given a contact datum $\bmu$ for a prestable map $(C,{\bf x},f)$, its components will always be denoted by $\mu_{i,j}$ and $\mu_{\vec{e},j}$.
\end{Convention}

\begin{Construction}\label{cons:fine-basic-monoid}
    Consider a prestable map $(C,{\bf x},f)$ and a contact datum $\bmu$ for it. We construct the \emph{fine monoid associated to} $(C,{\bf x},f,\bmu)$, denoted by $\ol M(C,{\bf x},f,\bmu)$, as follows. 
    
    \begin{enumerate}[ref = \theConstruction(\arabic*), label = (\arabic*), itemsep = 0.3ex]
        \item Let $\clG$ be the free abelian group generated by the symbols
        \begin{align*}
            \begin{array}{rl}
                m_{v,j} & \text{for each vertex } v \text{ of }\plC \text{ and } 1\le j\le r,\text{ and}  \\
                m_e & \text{for each edge } e \text{ of }\plC 
            \end{array}
        \end{align*}
        and let $\clR$ be the subgroup of $\clG$ generated by the elements
        \begin{align*}
            \begin{array}{rl}
               m_{v,j}  & \text{for each } v \text{ and } j\not\in I_v, \text{ i.e., } f(C_v)\not\subset D_j,\text{ and}\\
               m_{v',j} - m_{v,j} - \mu_{\vec{e},j}m_e  & \text{for each } \vec{e}:v\to v' \text{ and }j.
            \end{array}
        \end{align*}

        \item Define $\ol M(C,{\bf x},f,\bmu)\subset\clG/\clR$ to be the submonoid generated by (the images of) the elements $m_{v,j}$ (for all $v$, $j$) and $m_e$ (for all $e$). It is a fine monoid by construction.
        
        \item\label{cons-part:reduced-gen-set} In fact, the monoid $\ol M(C,{\bf x},f,\bmu)$ is generated by the \emph{reduced generating set} obtained from the generating set in part (2) by omitting $m_{v,j}$ whenever $j\not\in I_v$.
    \end{enumerate}
\end{Construction}

\begin{Remark}\label{rem:fine-basic-monoid}
    In terms of Definition~\ref{def-part:integral-presentation}, we have an integral presentation
    \begin{align}\label{eqn:fine-basic-monoid}
        \ol M(C,{\bf x},f,\bmu) := \left\langle
        m_{v,j}, m_e\;
        \middle|\!
        \begin{array}{rcl}
            j\not\in I_v & \!\!\!\Rightarrow & \!\!\! m_{v,j} = 0 \\
            \vec{e}:v\to v' & \!\!\!\Rightarrow & \!\!\! m_{v',j} = m_{v,j} + \mu_{\vec{e},j}m_e
        \end{array}
        \!\!
        \right\rangle^\text{int}.
    \end{align}
    Moreover, using the notation from Construction~\ref{cons:fine-basic-monoid}, the inclusion $\ol M(C,{\bf x},f,\bmu)\hookrightarrow\clG/\clR$ induces an isomorphism $\ol M(C,{\bf x},f,\bmu)^\text{gp}\xrightarrow{\sim}\clG/\clR$ of abelian groups.
\end{Remark}

\begin{Lemma}\label{lem:fine-basic-monoid-sharp}
    In the setting of Construction~\ref{cons:fine-basic-monoid}, the following statements are equivalent.
    \begin{enumerate}
        \item The monoid $\ol M(C,{\bf x},f,\bmu)$ is sharp and, in addition, each element of the reduced generating set from Construction~\ref{cons-part:reduced-gen-set} has nonzero image in $\ol M(C,{\bf x},f,\bmu)$.

        \item There exists a monoid map $\ol M(C,{\bf x},f,\bmu)\to \bN$ under which each element of the reduced generating set maps to $\bN\setminus\{0\}$.

        \item There exists a sharp monoid $Q$ and a monoid map $\ol M(C,{\bf x},f,\bmu)\to Q$ under which each element of the reduced generating set maps to $Q\setminus\{0\}$.
    \end{enumerate}
\end{Lemma}
\begin{proof}
    (1) $\Rightarrow$ (2): This follows from Lemma~\ref{lem-part:sharp-fine-interior-hom}.

    The implications (2) $\Rightarrow$ (3) $\Rightarrow$ (1) are clear.
\end{proof}

\begin{Definition}
    Consider a prestable map $(C,{\bf x},f)$ and a contact datum $\bmu$ for it.
    \begin{enumerate}[ref = \theDefinition(\arabic*), label = (\arabic*), itemsep = 0.3ex]
        \item\label{def-part:pullback-twisted-by-contact-divisor} For each vertex $v$ of $\plC$ and $1\le j\le r$, define a Cartier divisor $D_{v,j}$ on $\tilde C_v$ and a line bundle $L_{v,j}$ on $\tilde C_v$ by the formulas
        \begin{align*}
            D_{v,j}&:=\textstyle\sum_{x_i\in\tilde C_v}\mu_{i,j}\cdot x_i + \sum_{y_{\vec{e}}\in\tilde C_v}\mu_{\vec{e},j}\cdot y_{\vec{e}}, \quad\text{and}\\
            L_{v,j}&:=\tilde f_v^*(\clO_X(D_j))\otimes\clO_{\tilde C_v}(-D_{v,j}).
        \end{align*}
        
        \item\label{def-part:loc-realizable} We say that $\bmu$ is \emph{locally realizable} if, for each vertex $v$ of $\plC$ and $1\le j\le r$,
        
        \begin{enumerate}[ref = \theenumi(\roman*), label = (\roman*), itemsep = 0.2ex]
            \item\label{def-subpart:loc-realizable-noncanonically-trivial} the line bundle $L_{v,j}$ is isomorphic to $\clO_{\tilde C_v}$, and
            
            \item\label{def-subpart:loc-realizable-canonically-trivial} if $j\not\in I_v$, then the Cartier divisor $D_{v,j}$ is effective and the scheme-theoretic zero locus of the section $\tilde f_v^*\taut_{D_j}$ of $\tilde f_v^*(\clO_X(D_j))$ is given by $D_{v,j}\subset\tilde C_v$.
        \end{enumerate}
    \end{enumerate}
\end{Definition}

\begin{Remark}\label{rem:non-degenerate-contact-order}
    Condition~\ref{def-subpart:loc-realizable-canonically-trivial} may be rephrased as follows: if we have $f(C_v)\not\subset D_j$ then, for each $x_i\in\tilde C_v$ (resp. $y_{\vec{e}}\in \tilde C_v$), the number $\mu_{i,j}\in\bN$ (resp. $\mu_{\vec{e},j}\in\bN$) is the contact order of $\tilde f_v:\tilde C_v\to X$ with the divisor $D_j$ at the point $x_i$ (resp. $y_{\vec{e}}$).
\end{Remark}

\begin{Lemma}\label{lem:loc-realizable}
    Consider a prestable map $(C,{\bf x},f)$ and a contact datum $\bmu$ for it. If $\bmu$ is locally realizable, then we have the following properties.
    \begin{enumerate}[ref = \theLemma(\arabic*), label = (\arabic*), itemsep = 0.3ex]
        \item For each vertex $v$ of $\plC$ and $1\le j\le r$, we have the following.
        \begin{enumerate}[ref = \theenumi(\roman*), label = (\roman*), itemsep = 0.15ex]
            \item The $\bC$-vector space $H^0(\tilde C_v,L_{v,j})$ is $1$-dimensional.

            \item\label{lem-part:twisted-pullback-eval-iso} If $\vec{e}$ is an oriented edge of $\plC$ such that $y_{\vec{e}}\in\tilde C_v$, then evaluating global sections of the line bundle $L_{v,j}$ at the point $y_{\vec{e}}$ induces an isomorphism
            \begin{align}\label{eqn:pullback-twisted-by-contact-divisor-eval}
                \eval_{\vec{e},j}:H^0(\tilde C_v,L_{v,j})\xrightarrow{\sim} L_{v,j}|_{y_{\vec{e}}}.
            \end{align}

            \item\label{lem-subpart:rational-section-is-section} The rational section $\tilde f_v^*\taut_{D_j}\otimes\taut_{-D_{v,j}}$ of $L_{v,j}$ is in fact a global section. If $j\not\in I_v$, then it is a nowhere vanishing section and there exists a unique isomorphism
            \begin{align}\label{eqn:fine-basic-slb-non-degeneracy}
                \varphi_{v,j}:H^0(\tilde C_v,L_{v,j})\xrightarrow{\sim}\bC
            \end{align}
            which maps $\tilde f_v^*\taut_{D_j}\otimes\taut_{-D_{v,j}}\in H^0(\tilde C_v,L_{v,j})$ to $1\in\bC$.
        \end{enumerate}

        \item\label{disc-part:fine-basic-slb-edge-iso} For each oriented edge $\vec{e}:v\to v'$ of $\plC$ and $1\le j\le r$, there exists a unique isomorphism 
        \begin{align}\label{eqn:fine-basic-slb-edge-iso}
            \varphi_{\vec{e},j}:H^0(\tilde C_{v'},L_{v',j})\xrightarrow{\sim} H^0(\tilde C_v,L_{v,j})\otimes\Lambda_{C,y_e}^{\otimes\mu_{\vec{e},j}}
        \end{align}
        which makes the following diagram commute.
        \begin{equation}\label{eqn:fine-basic-slb-edge-iso-diagram}
        \begin{tikzcd}
            H^0(\tilde C_{v'},L_{v',j}) \arrow[r,"\varphi_{\vec{e},j}","\sim"'] \arrow[d,"\eval_{\cev{e},j}","\mathbin{\rotatebox[origin=c]{90}{$\sim$}}"'] & H^0(\tilde C_v,L_{v,j})\otimes\Lambda_{C,y_e}^{\otimes\mu_{\vec{e},j}} \arrow[d,"\eval_{\vec{e},j}\,\otimes\,\text{id}","\mathbin{\rotatebox[origin=c]{90}{$\sim$}}"'] \\
            L_{v',j}|_{y_{\cev{e}}} \arrow[d, equal] & L_{v,j}|_{y_{\vec{e}}}\otimes \Lambda_{C,y_e}^{\otimes\mu_{\vec{e},j}} \arrow[d, equal] \\
            \clO_X(D_j)|_{f(y_e)}\otimes T_{\tilde C,y_{\cev{e}}}^{\otimes(-\mu_{\cev{e},j})} \arrow[r,equal] & \clO_X(D_j)|_{f(y_e)}\otimes T_{\tilde C,y_{\vec{e}}}^{\otimes(-\mu_{\vec{e},j})}\otimes\Lambda_{C,y_e}^{\otimes\mu_{\vec{e},j}}
        \end{tikzcd}
        \end{equation}
        Here, the vertical identifications are induced by Remark~\ref{rem:curve-divisor-fibre}, while the horizontal identification is induced by $\mu_{\vec{e},j} + \mu_{\cev{e},j} = 0$.
    \end{enumerate}
\end{Lemma}
\begin{proof}
    (1): Assertions~(i) and (ii) follow from Condition~\ref{def-subpart:loc-realizable-noncanonically-trivial}, i.e., the triviality of $L_{v,j}$. 
    Assertion~(iv) is clear for $j\in I_v$ and follows from Condition~\ref{def-subpart:loc-realizable-canonically-trivial} for $j\not\in I_v$.

    (2): This is immediate from the fact the maps \eqref{eqn:pullback-twisted-by-contact-divisor-eval} are isomorphisms.
\end{proof}

\begin{Construction}\label{cons:fine-basic-slb-presentation}
    Consider a prestable map $(C,{\bf x},f)$ and a locally realizable contact datum $\bmu$ for it. We construct an slb presentation $\fP(C,{\bf x},f,\bmu)$ on $S:=\Spec\bC$, with respect to the integral presentation \eqref{eqn:fine-basic-monoid} of $\ol M(C,{\bf x},f,\bmu)$, as follows; see Definition~\ref{def-part:slb-presentation-defined}. 
    \begin{enumerate}[label  = (\arabic*), itemsep = 0.3ex]
        \item We define the virtual Cartier divisors $(\fL_{v,j},\fs_{v,j})$ and $(\fL_e,\fs_e)$ on $S$ corresponding to the generators $m_{v,j}$ and $m_e$ appearing in \eqref{eqn:fine-basic-monoid} as follows.

        \begin{enumerate}[label = (\roman*), itemsep = 0.15ex]
            \item For each vertex $v$ of $\plC$ and $1\le j\le r$, we define
            \begin{align*}
                \fL_{v,j}:= H^0(\tilde C_v,L_{v,j}) \quad\text{and}\quad \fs_{v,j}:=\tilde f_v^*\taut_{D_j}\otimes\taut_{-D_{v,j}}.
            \end{align*}

            \item For each edge $e$ of $\plC$, we define
            \begin{align*}
                \fL_e:=\Lambda_{C,y_e} \quad\text{and}\quad \fs_e:=0.
            \end{align*}
        \end{enumerate}

        \item We define the line bundle isomorphisms $\varphi_{v,j}$ and $\varphi_{\vec{e},j}$ on $S$ corresponding to the relations $m_{v,j} = 0$ and $m_{v',j} = m_{v,j}+\mu_{\vec{e},j}\rho_e$ appearing in \eqref{eqn:fine-basic-monoid} as follows. 

        \begin{enumerate}[label = (\alph*), itemsep = 0.15ex]
            \item For each vertex $v$ of $\plC$ and $j\not\in I_v$, we use \eqref{eqn:fine-basic-slb-non-degeneracy} to define
            \begin{align*}
                \varphi_{v,j}:\fL_{v,j}\xrightarrow{\sim}\clO_S = \bC.
            \end{align*}

            \item For each oriented edge $\vec{e}:v\to v'$ and $1\le j\le r$, we use \eqref{eqn:fine-basic-slb-edge-iso} to define
            \begin{align*}
                \varphi_{\vec{e},j}:\fL_{v',j}\xrightarrow{\sim}\fL_{v,j}\otimes\fL_e^{\otimes\mu_{\vec{e},j}}.
            \end{align*}
        \end{enumerate}
    \end{enumerate}
    
    We refer to $\fP(C,{\bf x},f,\bmu)$ as the \emph{slb presentation associated to} $(C,{\bf x},f,\bmu)$.
\end{Construction}

\subsection{Symplectic log maps}

\begin{Definition}\label{def:symp-log-map}
    Let $(C,{\bf x},f)$ be a prestable map to $X$ and let $\bmu$ be a locally realizable contact datum for it with respect to $D$. Recall the monoid $\ol M(C,{\bf x},f,\bmu)$ and its reduced generating set from Construction~\ref{cons:fine-basic-monoid}, the line bundles $L_{v,j}$ on $\tilde C_v$ from Definition~\ref{def-part:pullback-twisted-by-contact-divisor}, and the isomorphisms $\varphi_{v,j}$ and $\varphi_{\vec{e},j}$ from Lemma~\ref{lem:loc-realizable}.
    
    In this setting, the tuple $(C,{\bf x},f,\bmu)$ is called a \emph{symplectic} (\emph{prestable}) \emph{log map} to the smooth snc pair $(X,D)$ if it satisfies the following conditions.
    \begin{enumerate}[ref = \theDefinition(\arabic*), label = (\arabic*), itemsep = 0.3ex]
        \item\label{def-part:tropical-map-existence-condition}
        There exists a monoid map $\ol M(C,{\bf x},f,\bmu)\to\bR_{\ge 0}$ under which each element of the reduced generating set maps to $\bR_{>0}$, where we view $\bR_{\ge 0}$ as an additive monoid.        
        
        \item\label{def-part:analytic-condition}
        There exist sections $0\ne\sigma_{v,j}\in H^0(\tilde C_v,L_{v,j})$, for each vertex $v$ of $\plC$ and $1\le j\le r$, and elements $0\ne\theta_e\in\Lambda_{C,y_e}$, for each edge $e$ of $\plC$, with the following properties.
        \begin{enumerate}[ref = \theenumi(\alph*), label = (\alph*), itemsep = 0.15ex]
            \item\label{def-subpart:analytic-non-degeneracy-condition} For each $v$ and $j\not\in I_v$, we have $\varphi_{v,j}(\sigma_{v,j}) = 1$.

            \item\label{def-subpart:analytic-matching-condition} For each $\vec{e}:v\to v'$ and $j$, we have $\varphi_{\vec{e},j}(\sigma_{v',j}) = \sigma_{v,j}\otimes\theta_e^{\otimes\mu_{\vec{e},j}}$.
        \end{enumerate}
    \end{enumerate}
    
    We call $(C,{\bf x},f,\bmu)$ a \emph{symplectic stable log map} if $(C,{\bf x},f)$ is a stable map.
\end{Definition}

\begin{Remark}\label{rem:symp-log-map}
    We explain why Definition~\ref{def:symp-log-map} is in fact equivalent to the original definition of \emph{log holomorphic maps} given in \cite[Definition~2.8]{Tehrani-log}. Though the definition in \cite{Tehrani-log} is more general, we only consider the case of nonnegative contact orders at marked points. For the discussion below, fix a prestable map $(C,{\bf x},f)$ to $X$.
    
    We first recall the notion of \emph{prelog holomorphic maps}, in the sense of \cite[Definition~2.4]{Tehrani-log}. Lifting $(C,{\bf x},f)$ to a prelog holomorphic map amounts to choosing a contact datum $\bmu$ for it with the following two properties. First, for each vertex $v$ of $\plC$ and $j\not\in I_v$, if we have $x_i\in\tilde C_v$ (resp. $y_{\vec{e}}\in\tilde C_v$) then the contact order of $\tilde f_v:\tilde C_v\to X$ with $D_j$ at this point is given by $\mu_{i,j}$ (resp. $\mu_{\vec{e},j})$. Second, for each vertex $v$ of $\plC$ and $j\in I_v$, there exists a rational (i.e., meromorphic) section $0\ne\zeta_{v,j}$ of $\tilde f_v^*\clN_{D_j/X}$ whose divisor of zeros and poles is given by $D_{v,j}$ from Definition~\ref{def-part:pullback-twisted-by-contact-divisor}; note that $\zeta_{v,j}$, if it exists, is unique up to $\bC^\times$-scaling.

    These two properties are readily seen to be equivalent to the local realizability of $\bmu$. For $j\not\in I_v$, this is the content of Remark~\ref{rem:non-degenerate-contact-order}. For $j\in I_v$, \eqref{eqn:normal-bundle} provides a canonical identification $\tilde f_v^*\clN_{D_j/X}\xrightarrow{\sim}\tilde f_v^*(\clO_X(D_j))$ and so, $\zeta_{v,j}$ gives rise to a nowhere vanishing global section $\sigma_{v,j}$ of the line bundle $L_{v,j}$ from Definition~\ref{def-part:pullback-twisted-by-contact-divisor}. Thus, lifting $(C,{\bf x},f)$ to a prelog holomorphic map is the same as choosing a locally realizable contact datum $\bmu$ for it. For the remainder, we will find it convenient to view $\zeta_{v,j}$ for $j\in I_v$ as rational sections of $\tilde f_v^*(\clO_X(D_j))$, rather than $\tilde f_v^*\clN_{D_j/X}$, with the divisor of zeros and poles being $D_{v,j}$. (This viewpoint is exclusive to the holomorphic setting and does not admit any obvious generalization to the pseudo-holomorphic setting, since it is hard to define and work with the sheaf $\clO_X(D_j)$ for non-integrable almost complex structures on $X$.) Using this new viewpoint, we extend the definition of $\zeta_{v,j}$ by setting $\zeta_{v,j}:=\tilde f_v^*\taut_{D_j}$ for $j\not\in I_v$; these are global sections of $\tilde f_v^*(\clO_X(D_j))$ whose divisor of zeros is given by $D_{v,j}$. As before, we get nowhere vanishing global sections
    \begin{align}\label{eqn:rescaled-limits-for-analytic-matching-condition}
        0\ne \sigma_{v,j}:=\zeta_{v,j}\otimes\taut_{-D_{v,j}}\in H^0(\tilde C_v,L_{v,j}).
    \end{align} 
    Observe that the section $\zeta_{v,j}$ is defined `on the nose' for $j\not\in I_v$, whereas it is determined only up to $\bC^\times$-scaling for $j\in I_v$.
    
    The prelog holomorphic map defined by the tuple $(C,{\bf x},f,\bmu)$ is a \emph{log holomorphic map} according to the original definition precisely when it satisfies Conditions (1) and (2) in \cite[Definition~2.8]{Tehrani-log}. The first of these conditions is identical to Condition~\ref{def-part:tropical-map-existence-condition}. By the paragraph preceding \cite[Definition~2.8]{Tehrani-log}, the second of these conditions is equivalent to the existence of rational sections $0\ne\zeta_{v,j}$ as above and local coordinates $z_{\vec{e}}$ centred at $y_{\vec{e}}\in\tilde C_v$ such that, for each oriented edge $\vec{e}:v\to v'$ of $\plC$ and $1\le j\le r$, the lowest order coefficients $\eta_{\vec{e},j}$ (resp. $\eta_{\cev{e},j}$) in the power series expansions of $\zeta_{v,j}$ (resp. $\zeta_{v',j}$) centred at $y_{\vec{e}}$ (resp. $y_{\cev{e}}$) with respect to the coordinates $z_{\vec{e}}$ (resp. $z_{\cev{e}}$) are equal as elements of $\clO_X(D_j)|_{f(y_e)}$.
    
    Condition~\ref{def-subpart:analytic-non-degeneracy-condition} simply states that, for $j\not\in I_v$, we have $\sigma_{v,j} = \tilde f_v^*\taut_{D_j}\otimes\taut_{-D_{v,j}}$, i.e., $\zeta_{v,j} = \tilde f_v^*\taut_{D_j}$. Using the notation of Lemma~\ref{lem-part:twisted-pullback-eval-iso}, we have the identities
    \begin{align*}
        \eval_{\vec{e},j}(\sigma_{v,j}) = (dz_{\vec{e}})^{\otimes\mu_{\vec{e},j}}\otimes\eta_{\vec{e},j} \quad\text{and}\quad \eval_{\cev{e},j}(\sigma_{v',j}) = (dz_{\cev{e}})^{\otimes\mu_{\cev{e},j}}\otimes\eta_{\cev{e},j}.
    \end{align*}
    Now, Condition~\ref{def-subpart:analytic-matching-condition} is seen to be equivalent to $\eta_{\vec{e},j} = \eta_{\cev{e},j}$ by taking
    \begin{align}\label{eqn:smoothing-parameter-for-analytic-matching-condition}
        \textstyle 0\ne\theta_e:=\left.\frac{\partial}{\partial z_{\vec{e}}}\right|_{y_{\vec{e}}}\otimes\left.\frac{\partial}{\partial z_{\cev{e}}}\right|_{y_{\cev{e}}}\in\Lambda_{C,y_e}.
    \end{align}
    Thus, Condition (2) in \cite[Definition~2.8]{Tehrani-log} is equivalent to Condition~\ref{def-part:analytic-condition}. We conclude that Definition~\ref{def:symp-log-map} is indeed equivalent to \cite[Definition~2.8]{Tehrani-log}.
\end{Remark}

\begin{Lemma}\label{lem:basic-fine-slb-consistent-symplectic-log-map}
    Consider a prestable map $(C,{\bf x},f)$ and a locally realizable contact datum $\bmu$ for it. Then the following statements are equivalent.
    \begin{enumerate}
        \item The tuple $(C,{\bf x},f,\bmu)$ is a symplectic log map.

        \item The slb presentation $\fP(C,{\bf x},f,\bmu)$ from Construction~\ref{cons:fine-basic-slb-presentation} is consistent in the sense of Definition~\ref{def-part:consistent-slb-presentation}.
    \end{enumerate}
\end{Lemma}
\begin{proof}
    Recall the monoid $\ol M(C,{\bf x},f,\bmu)$ and its reduced generating set from Construction~\ref{cons:fine-basic-monoid}.

    (1) $\Rightarrow$ (2): Choose elements $\sigma_{v,j}$ and $\theta_e$ as in Condition~\ref{def-part:analytic-condition}. Using these choices, define the trivializations $\chi_{v,j}:\fL_{v,j}\xrightarrow{\sim}\bC$ and $\chi_e:\fL_e\xrightarrow{\sim}\bC$ by $\sigma_{v,j}\mapsto 1$ and $\theta_e\mapsto 1$ respectively. We verify that these trivializations satisfy the two conditions in Definition~\ref{def-part:consistent-slb-presentation}. 
    \begin{enumerate}[label = (\roman*), itemsep = 0.3ex]
        \item We verify the compatibility of $\chi_{v,j}$ and $\chi_e$ with $\varphi_{v,j}$ and $\varphi_{\vec{e},j}$.
        \begin{enumerate}[label = (\alph*), itemsep = 0.15ex]
            \item For each $v$ and $j\not\in I_v$, we have $\chi_{v,j} = \varphi_{v,j}$ since both sides map $\sigma_{v,j}$ to $1$.

            \item For each $\vec{e}:v\to v'$ and $j$, we have $\chi_{v',j} = (\chi_{v,j}\otimes\chi_e^{\otimes\mu_{\vec{e},j}})\circ\varphi_{\vec{e},j}$ since both sides map $\sigma_{v',j}$ to $1$.
        \end{enumerate}

        \item Given any element $q\in\ol M(C,{\bf x},f,\bmu)$ and a choice of coefficients $x_{v,j},x_e\in\bN$ such that $\sum x_{v,j}m_{v,j} + \sum x_em_e$ maps to $q$, we verify that $\prod(\chi_{v,j}(\fs_{v,j}))^{x_{v,j}}\prod(\chi_e(\fs_e))^{x_e}$ depends only on $q$ and not on the choice of $x_{v,j},x_e$. It suffices to show that we have
        \begin{align}\label{eqn:sections-in-split-log-pt-hom}
            \prod(\chi_{v,j}(\fs_{v,j}))^{x_{v,j}}\prod(\chi_e(\fs_e))^{x_e} = \begin{cases}
                1 & \text{if }q = 0, \\
                0 & \text{if }q\ne 0.
            \end{cases}
        \end{align}
        The identity \eqref{eqn:sections-in-split-log-pt-hom} follows from the next two observations. First, we have \begin{align}\label{eqn:sections-in-split-log-pt-0-or-1}
            \chi_e(\fs_e) = 0\text{ for all }e, \quad\text{and}\quad \chi_{v,j}(\fs_{v,j}) = \begin{cases}
                1 & \text{if }j\not\in I_v, \\
                0 & \text{if }j\in I_v.
            \end{cases}
        \end{align}
        Second, by Condition~\ref{def-part:tropical-map-existence-condition} and Lemma~\ref{lem:fine-basic-monoid-sharp}, $\ol M(C,{\bf x},f,\bmu)$ is sharp and each element of its reduced generating set has nonzero image in $\ol M(C,{\bf x},f,\bmu)$.
    \end{enumerate}

    (2) $\Rightarrow$ (1): Choose trivializations $\chi_{v,j}:\fL_{v,j}\xrightarrow{\sim}\bC$ and $\chi_e:\fL_e\xrightarrow{\sim}\bC$ as in Definition~\ref{def-part:consistent-slb-presentation}. Using these choices, define $0\ne \sigma_{v,j}\in\fL_{v,j}$ (resp. $0\ne \theta_e\in\fL_e$) to be the inverse image of $1\in\bC$ under $\chi_{v,j}$ (resp. $\chi_e$). We verify that $(C,{\bf x},f,\bmu)$ is a symplectic log map. \begin{enumerate}[label = (\roman*), itemsep = 0.3ex]
        \item We verify Condition~\ref{def-part:tropical-map-existence-condition}. For $j\not\in I_v$, we have $\chi_{v,j} = \varphi_{v,j}$ by Condition~\ref{def-subpart:relation-compatible-trivialization} and thus, $\chi_{v,j}(\fs_{v,j}) = 1$. Moreover, $\fs_{v,j} = 0$ for $j\in I_v$ and $\fs_e = 0$ for all $e$. Thus, we again have \eqref{eqn:sections-in-split-log-pt-0-or-1}. Condition~\ref{def-subpart:section-compatible-trivialization} implies that the assignment 
        \begin{align*}
            m_{v,j} \mapsto \chi_{v,j}(\fs_{v,j}),\quad m_e\mapsto\chi_e(\fs_e)
        \end{align*}
        gives a well-defined monoid map $\ol M(C,{\bf x},f,\bmu)\to\{0,1\}$, where we endow $\{0,1\}$ with the multiplicative monoid structure. Note that the identity element of the monoid $\{0,1\}$ is $1$ rather than $0$. The implication (3) $\Rightarrow$ (2) in Lemma~\ref{lem:fine-basic-monoid-sharp}, with $(\{0,1\},1)$ in place of $(Q,0)$, shows the existence of a monoid map $\ol M(C,{\bf x},f,\bmu)\to\bN\subset\bR_{\ge 0}$ under which each element of the reduced generating set has a nonzero image.
       
        \item We verify Condition~\ref{def-part:analytic-condition}. From Condition~\ref{def-subpart:relation-compatible-trivialization}, we get the following.
        \begin{enumerate}[label = (\alph*), itemsep = 0.15ex]
            \item For each $v$ and $j\not\in I_v$, we have $\chi_{v,j} = \varphi_{v,j}$. Evaluating both sides at $\sigma_{v,j}$, we deduce that $\varphi_{v,j}(\sigma_{v,j}) = 1$.

            \item For each $\vec{e}:v\to v'$ and $j$, we have $\chi_{v',j} = (\chi_{v,j}\otimes\chi_e^{\otimes\mu_{\vec{e},j}}) \circ\varphi_{\vec{e},j}$. Evaluating both sides at $\sigma_{v',j}$, we deduce that $\varphi_{\vec{e},j}(\sigma_{v',j}) = \sigma_{v,j}\otimes\theta_e^{\otimes\mu_{\vec{e},j}}$.
        \end{enumerate}
    \end{enumerate}
    
    This completes the proof.
\end{proof}

\subsection{Algebraic log maps}

\begin{Definition}\label{def:alg-log-map}
    An \emph{algebraic} (\emph{prestable}) \emph{log map} to the log scheme $X^\dagger = (X,M_{(X,D)})$ is a tuple $(\xi^\dagger/S^\dagger,f^\dagger)$ consisting of the following data.
    \begin{enumerate}[label = (\arabic*), itemsep = 0.3ex]
        \item A (fine) log point $S^\dagger = (S, M_S)\in\LogPt$; see Definition~\ref{def:log-point}.
        
        \item A family $\xi^\dagger = (\xi,\xi^\flat)$ of log curves over $S^\dagger$, with underlying curve $\xi = (C/S,{\bf x})$; see Definition~\ref{def:family-of-log-curves}. Define the log scheme $C^\dagger = (C,M_C)$ using $\xi^\dagger$ as in Definition~\ref{def:family-of-log-curves-total-space}.
        
        \item A log morphism $f^\dagger = (f,f^\flat):C^\dagger\to X^\dagger$, with underlying morphism $f:C\to X$. 
    \end{enumerate}
    
    We refer to $(C,{\bf x},f)$ as the \emph{underlying prestable map} of $(\xi^\dagger/S^\dagger,f^\dagger)$. 
    
    We call $(\xi^\dagger/S^\dagger,f^\dagger)$ an \emph{algebraic stable log map} if $(C,{\bf x},f)$ is a stable map.
\end{Definition}

\begin{Remark}
    Definition~\ref{def:alg-log-map} agrees with \cite[Definition~2.1.2]{Chen-logDF1} for (fine) log points and, using Remark~\ref{rem:log-curve-defs}, it is equivalent to \cite[Definition~1.6]{GS-log} for fs log points.
\end{Remark}

The following result provides a concrete description of algebraic log maps with a fixed underlying prestable map. We freely use Notation~\ref{nota:prestable-map} in its statement and proof.

\begin{Proposition}\label{prop:alg-log-map}
    Fix $S = \Spec\bC$ and a prestable map $(C,{\bf x},f)$ to $X$. For a fine log structure $M_S$ on $S$, write $S^\dagger = (S,M_S)$ and consider the following two types of objects.
    \begin{enumerate}[itemsep = 0.3ex, leftmargin = 9ex]
        \myitem[(ALM$\phantom{}_1$)]\label{abstract-log-map} Algebraic log map $(\xi^\dagger/S^\dagger,f^\dagger)$ to $X^\dagger$ with underlying prestable map $(C,{\bf x},f)$.

        \myitem[(ALM$\phantom{}_2$)]\label{concrete-log-map} Tuple $(\{\rho_e,\Theta_e\},\bmu,\{\delta_{v,j},\Theta_{v,j}\})$ consisting of        
        \begin{enumerate}[label = (\roman*), itemsep = 0.15ex]            
            \item elements $0\ne\rho_e\in\ol M_S$ and linear isomorphisms $\Theta_e:\Lambda_{C,y_e}\xrightarrow{\sim}\clO_{S^\dagger}(\rho_e)$, where $e$ ranges over the edges of the dual graph $\plC$ of $C$,

            \item a locally realizable contact datum $\bmu = \{\mu_{i,j}\in\bN,\mu_{\vec{e},j}\in\bZ\}$ for $(C,{\bf x},f)$ with respect to the divisor $D$ in the sense of Definition~\ref{def-part:loc-realizable}, and
            
            \item elements $\delta_{v,j}\in\ol M_S$ and linear isomorphisms $\Theta_{v,j}:H^0(\tilde C_v,L_{v,j})\xrightarrow{\sim}\clO_{S^\dagger}(\delta_{v,j})$, where $v$ ranges over the vertices of $\plC$, $j$ ranges over $\{1,\ldots,r\}$ and $L_{v,j}$ is the line bundle on $\tilde C_v$ associated to $\bmu$ as in Definition~\ref{def-part:pullback-twisted-by-contact-divisor},
        \end{enumerate}
        which satisfy the following properties.
        \begin{enumerate}[ref = (ALM$\phantom{}_2$)(\alph*), label = (\alph*), itemsep = 0.15ex]
            \item\label{concrete-alg-log-map-vcd-iso} For each vertex $v$ of $\plC$ and $1\le j\le r$, $\Theta_{v,j}$ maps $\tilde f_v^*\taut_{D_j}\otimes\taut_{-D_{v,j}}\in H^0(\tilde C_v,L_{v,j})$ to $\sigma_{S^\dagger,\delta_{v,j}}\in\clO_{S^\dagger}(\delta_{v,j})$. In particular, $j\not\in I_v$ if and only if $\delta_{v,j} = 0$.

            \item\label{concrete-alg-log-map-matching-condition} For each oriented edge $\vec{e}:v\to v'$ of $\plC$ and $1\le j\le r$, we have the identity $\delta_{v',j} = \delta_{v,j} + \mu_{\vec{e},j}\rho_e$ and the following diagram commutes.
            \begin{equation}\label{eqn:twisted-pullback-at-node}
            \begin{tikzcd}
                H^0(\tilde C_{v'},L_{v',j}) \arrow[rr,"\varphi_{\vec{e},j}","\sim"'] \arrow[dd,"\Theta_{v',j}","\mathbin{\rotatebox[origin=c]{90}{$\sim$}}"'] & & H^0(\tilde C_v,L_{v,j})\otimes\Lambda_{C,y_e}^{\otimes\mu_{\vec{e},j}} \arrow[dd,"\Theta_{v,j}\,\otimes\,\Theta_e^{\otimes\mu_{\vec{e},j}}","\mathbin{\rotatebox[origin=c]{90}{$\sim$}}"'] \\ & & \\
                \clO_{S^\dagger}(\delta_{v',j}) \arrow[rr,equal] & & \clO_{S^\dagger}(\delta_{v,j})\otimes\clO_{S^\dagger}(\rho_e)^{\otimes\mu_{\vec{e},j}}
            \end{tikzcd}
            \end{equation}
            Here, $\varphi_{\vec{e},j}$ is the map from \eqref{eqn:fine-basic-slb-edge-iso} and the bottom horizontal identification is induced by the isomorphisms appearing in the global slb on $S^\dagger$.
        \end{enumerate}
    \end{enumerate}
    As the log structure $M_S$ varies, each of \ref{abstract-log-map} and \ref{concrete-log-map} determines a category fibred in groupoids over $\LogPt$, using the obvious notion of pullbacks along log morphisms. In particular, pullbacks in \ref{concrete-log-map} preserve the contact datum $\bmu$.

    These two categories fibred in groupoids over $\LogPt$ are equivalent.
\end{Proposition}

\begin{proof}
    We define functors \ref{abstract-log-map} $\to$ \ref{concrete-log-map} $\to$ \ref{abstract-log-map} but omit the verification that their composites \ref{abstract-log-map} $\to$ \ref{abstract-log-map} and \ref{concrete-log-map} $\to$ \ref{concrete-log-map} are isomorphic to identity functors. We will draw heavily on Proposition~\ref{prop:log-mor-to-snc-pair} and Theorem~\ref{thm:log-curve-over-log-pt} and so, the reader is advised to review their statements first.

    \ref{abstract-log-map} $\to$ \ref{concrete-log-map}: 
    Let $(\xi^\dagger/S^\dagger,f^\dagger)$ be a given algebraic log map to $X^\dagger$ with underlying prestable map $(C,{\bf x},f)$. Let $C^\dagger = (C,M_C)$ denote the log scheme obtained from $\xi^\dagger$ via Definition~\ref{def:family-of-log-curves-total-space}. Define the tuple $(\{\rho_e,\Theta_e\},\bmu,\{\delta_{v,j},\Theta_{v,j}\})$ as follows.

    \begin{enumerate}[label = (\roman*), itemsep = 0.3ex]
        \item Apply the equivalence \ref{abstract-log-curve} $\to$ \ref{concrete-log-curve} from Theorem~\ref{thm:log-curve-over-log-pt} to get, for each edge $e$ of $\plC$, an element $0\ne \rho_e\in\ol M_S$ and an isomorphism $\Theta_e:\Lambda_{C,y_e}\xrightarrow{\sim}\clO_{S^\dagger}(\rho_e)$.

        \item Apply the bijection \ref{abstract-log-mor-to-snc-pair} $\to$ \ref{concrete-log-mor-to-snc-pair} from Proposition~\ref{prop:log-mor-to-snc-pair} to get, for each $1\le j\le r$, an element $\delta_j\in\Gamma(C,\ol M_C)$ and an isomorphism of virtual Cartier divisors
        \begin{align}\label{eqn:alg-log-map-vcd-iso}
            \Psi_j:f^*(\clO_X(D_j),\taut_{D_j})\xrightarrow{\sim}(\clO_{C^\dagger}(\delta_j),\sigma_{C^\dagger,\delta_j}).
        \end{align}
        For each $1\le j\le r$, apply the isomorphism \eqref{eqn:ghost-sections-log-curve} to $\delta_j\in\Gamma(C,\ol M_C)$ to get elements 
        $\delta_{v,j}\in\ol M_S$ for each vertex $v$ of $\plC$, $\mu_{i,j}\in\bN$ for each $1\le i\le n$, and $\mu_{\vec{e},j}\in\bZ$ for each oriented edge $\vec{e}:v\to v'$ of $\plC$ satisfying $\delta_{v',j} = \delta_{v,j} + \mu_{\vec{e},j}\rho_e$.
        \vspace{0.3ex}
        \myitem[] Note that, since $\rho_e\ne 0$, we must have $\mu_{\vec{e},j} + \mu_{\cev{e},j} = 0$ for each $1\le j\le r$. It follows that $\bmu:=\{\mu_{i,j},\mu_{\vec{e},j}\}$ is a contact datum for $(C,{\bf x},f)$ with respect to $D$.

        \item Fix a vertex $v$ of $\plC$ and $1\le j\le r$. Recall from Definition~\ref{def-part:pullback-twisted-by-contact-divisor} that we have a Cartier divisor $D_{v,j}$ and a line bundle $L_{v,j}$ on $\tilde C_v$ associated to the contact datum $\bmu$. From Theorem~\ref{thm-subpart:log-curve-slb-on-irr-comps}, we obtain a natural isomorphism
        \begin{align}\label{eqn:alg-log-map-vcd-on-irr-comp}
            \Phi_{v,j}:(\nu|_{\tilde C_v})^*(\clO_{C^\dagger}(\delta_j),\sigma_{C^\dagger,\delta_j})\xrightarrow{\sim} (\clO_{S^\dagger}(\delta_{v,j})\otimes_\bC\clO_{\tilde C_v}(D_{v,j}),\sigma_{S^\dagger,\delta_{v,j}}\otimes_\bC\taut_{D_{v,j}})\!\!\!\!\!\!\!\!\!\!\!
        \end{align}
        of virtual Cartier divisors on $\tilde C_v$, where $\nu:\tilde C\to C$ is the normalization morphism. Using this, define an isomorphism $\Psi_{v,j}$ of line bundles on $\tilde C_v$ by
        \begin{align*}
            \Psi_{v,j}:L_{v,j} &\xrightarrow{\sim} \clO_{S^\dagger}(\delta_{v,j})\otimes_\bC\clO_{\tilde C_v}, \\
            \Psi_{v,j}&:=(\Phi_{v,j}\circ(\nu|_{\tilde C_v})^*\Psi_j)\otimes\text{id}_{\clO_{\tilde C_v}(-D_{v,j})}.
        \end{align*}
        The existence of the isomorphism $\Psi_{v,j}$ shows that $\bmu$ satisfies Condition~\ref{def-subpart:loc-realizable-noncanonically-trivial}. By definition, the isomorphism $\Psi_{v,j}$ maps the rational section $\tilde f_v^*\taut_{D_j}\otimes\taut_{-D_{v,j}}$ to the constant section $\sigma_{S^\dagger,\delta_{v,j}}$. Thus, we have $\tilde f_v^*\taut_{D_j}\otimes\taut_{-D_{v,j}}\in H^0(\tilde C_v,L_{v,j})$ and, moreover, $\sigma_{S^\dagger,\delta_{v,j}}\ne 0$ if and only if $j\not\in I_v$. This shows that $\bmu$ satisfies Condition~\ref{def-subpart:loc-realizable-canonically-trivial} as well. Thus, the contact datum $\bmu$ is locally realizable.
        \vspace{0.3ex}
        \myitem[] Let $\Theta_{v,j}:H^0(\tilde C_v,L_{v,j})\xrightarrow{\sim}\clO_{S^\dagger}(\delta_{v,j})$ be the map induced by $\Psi_{v,j}$ on global sections.
    \end{enumerate}
    We now verify that $(\{\rho_e,\Theta_e\},\bmu,\{\delta_{v,j},\Theta_{v,j}\})$ is an object of \ref{concrete-log-map} lying over $S^\dagger$.
    \begin{enumerate}[label = (\alph*), itemsep = 0.3ex]
        \item Fix a vertex $v$ of $\plC$ and $1\le j\le r$. We have already seen that $\Theta_{v,j}$ maps $\tilde f_v^*\taut_{D_j}\otimes\taut_{-D_{v,j}}$ to $\sigma_{S^\dagger,\delta_{v,j}}$ and that $j\not\in I_v$ if and only if $\sigma_{S^\dagger,\delta_{v,j}}\ne 0$. Moreover, by the equivalence $\LogPt\simeq\SlbPt$ in Proposition~\ref{prop:cat-of-log-pts}, we have $\sigma_{S^\dagger,\delta_{v,j}} \ne 0$ if and only if $\delta_{v,j} = 0$.

        \item Fix an oriented edge $\vec{e}:v\to v'$ of $\plC$ and $1\le j\le r$. We have already seen that the identity $\delta_{v',j} = \delta_{v,j} + \mu_{\vec{e},j}\rho_e$ holds. By Theorem~\ref{thm-subpart:log-curve-slb-at-nodes}, the elements
        \begin{align*}
            0\ne\Psi_{v,j}|_{y_{\vec{e}}}\in&\Hom(\clO_X(D_j)|_{f(y_e)},\clO_{S^\dagger}(\delta_{v,j})\otimes T_{\tilde C,y_{\vec{e}}}^{\otimes\mu_{\vec{e},j}}), \quad\text{and}\\
            0\ne\Psi_{v',j}|_{y_{\cev{e}}}\in&\Hom(\clO_X(D_j)|_{f(y_e)},\clO_{S^\dagger}(\delta_{v',j})\otimes T_{\tilde C,y_{\cev{e}}}^{\otimes\mu_{\cev{e},j}})
        \end{align*}
        agree under the isomorphism
        \begin{align}\label{eqn:alg-log-map-matching-iso}
            \clO_{S^\dagger}(\delta_{v,j})\otimes T_{\tilde C,y_{\vec{e}}}^{\otimes\mu_{\vec{e},j}} \xrightarrow{\sim} \clO_{S^\dagger}(\delta_{v',j})\otimes T_{\tilde C,y_{\cev{e}}}^{\otimes\mu_{\cev{e},j}}
        \end{align}
        induced by $\Theta_e^{\otimes\mu_{\vec{e},j}}$ and the isomorphisms appearing in the global slb on $S^\dagger$. Here, we are using Remark~\ref{rem:curve-divisor-fibre} to compute the fibre $L_{v,j}$ (resp. $L_{v',j}$) at $y_{\vec{e}}$ (resp. $y_{\cev{e}}$).
        \vspace{0.3ex}
        \myitem[] The agreement of $\Psi_{v,j}|_{y_{\vec{e}}}$ and $\Psi_{v',j}|_{y_{\cev{e}}}$ under \eqref{eqn:alg-log-map-matching-iso} is equivalent to the commutativity of the following diagram, where the top horizontal identification is the composition of the three `equalities' appearing in \eqref{eqn:fine-basic-slb-edge-iso-diagram} and the bottom horizontal identification is induced by the isomorphisms appearing in the global slb on $S^\dagger$.
        \begin{equation}\label{eqn:twisted-pullback-at-node-2}
        \begin{tikzcd}
            L_{v',j}|_{y_{\cev{e}}} \arrow[rr,equals] \arrow[dd,"\Psi_{v',j}|_{y_{\cev{e}}}","\mathbin{\rotatebox[origin=c]{90}{$\sim$}}"'] & & L_{v,j}|_{y_{\vec{e}}} \otimes \Lambda_{C,y_e}^{\otimes\mu_{\vec{e},j}} \arrow[dd,"\Psi_{v,j}|_{y_{\vec{e}}}\,\otimes\,\Theta_e^{\otimes\mu_{\vec{e},j}}","\mathbin{\rotatebox[origin=c]{90}{$\sim$}}"'] \\ & & \\
            \clO_{S^\dagger}(\delta_{v',j}) \arrow[rr,equal] & & \clO_{S^\dagger}(\delta_{v,j})\otimes\clO_{S^\dagger}(\rho_e)^{\otimes\mu_{\vec{e},j}}
        \end{tikzcd}
        \end{equation}
        Finally, examining the definitions of $\Theta_{v,j}$, $\Theta_{v',j}$ and $\varphi_{\vec{e},j}$, we see that the commutativity of \eqref{eqn:twisted-pullback-at-node-2} is equivalent to the commutativity of \eqref{eqn:twisted-pullback-at-node}.
    \end{enumerate}
    Therefore, we obtain a well-defined assignment $(\xi^\dagger/S^\dagger,f^\dagger) \mapsto (\{\rho_e,\Theta_e\},\bmu,\{\delta_{v,j},\Theta_{v,j}\})$. This assignment is functorial by Proposition~\ref{prop:log-mor-to-snc-pair} and Theorem~\ref{thm:log-curve-over-log-pt}.

    \ref{concrete-log-map} $\to$ \ref{abstract-log-map}: Let $(\{\rho_e,\Theta_e\},\bmu,\{\delta_{v,j},\Theta_{v,j}\})$ be an object of \ref{concrete-log-map} lying over the log point $S^\dagger = (S,M_S)$. Define the algebraic log map $(\xi^\dagger/S^\dagger,f^\dagger)$ to $X^\dagger$ as follows. (The discussion that follows is very similar to the one appearing in the construction of the functor \ref{abstract-log-map} $\to$ \ref{concrete-log-map} and so, we will be somewhat brief.)

    \begin{enumerate}[label = (\arabic*), itemsep = 0.3ex]
        \item Apply the equivalence \ref{concrete-log-curve} $\to$ \ref{abstract-log-curve} from Theorem~\ref{thm:log-curve-over-log-pt} to the collection $\{\rho_e,\Theta_e\}$ to get a family $\xi^\dagger$ of log curves over $S^\dagger$ whose underlying curve is $\xi = (C/S,{\bf x})$. Let $C^\dagger = (C,M_C)$ be the log scheme obtained from $\xi^\dagger$ via Definition~\ref{def:family-of-log-curves-total-space}.

        \item Recall that the contact datum $\bmu$ is locally realizable. For each vertex $v$ of $\plC$ and $1\le j\le r$, the triviality of  $L_{v,j}$ provides a unique line bundle isomorphism
        \begin{align*}
            \Psi_{v,j}:L_{v,j}\xrightarrow{\sim}\clO_{S^\dagger}(\delta_{v,j})\otimes_\bC\clO_{\tilde C_v}
        \end{align*}
        which induces the map $\Theta_{v,j}:H^0(\tilde C_v,L_{v,j})\xrightarrow{\sim}\clO_{S^\dagger}(\delta_{v,j})$ on global sections. From Lemma~\ref{lem-subpart:rational-section-is-section}, we get $\tilde f_v^*\taut_{D_j}\otimes\taut_{-D_{v,j}}\in H^0(\tilde C_v,L_{v,j})$. By Condition~\ref{concrete-alg-log-map-vcd-iso}, the global section $\tilde f_v^*\taut_{D_j}\otimes\taut_{-D_{v,j}}$ maps to the constant section $\sigma_{S^\dagger,\delta_{v,j}}$ under $\Psi_{v,j}$. Thus, we obtain an isomorphism of virtual Cartier divisors on $\tilde C_v$ defined by
        \begin{align*}
            \hat\Psi_{v,j}:\tilde f_v^*(\clO_X(D_j),\taut_{D_j})&\xrightarrow{\sim}(\clO_{S^\dagger}(\delta_{v,j})\otimes_\bC\clO_{\tilde C_v}(D_{v,j}),\sigma_{S^\dagger,\delta_{v,j}}\otimes_\bC\taut_{D_{v,j}}),\!\!\!\!\!\!\!\!\!\!\!
            \\
            \hat\Psi_{v,j} &:= \Psi_{v,j}\otimes\text{id}_{\clO_{\tilde C_v}(D_{v,j})}.
        \end{align*}
        
        \item For each oriented edge $\vec{e}:v\to v'$ and $1\le j\le r$, Condition~\ref{concrete-alg-log-map-matching-condition} implies that the elements $\Psi_{v,j}|_{y_{\vec{e}}}$ and $\Psi_{v',j}|_{y_{\cev{e}}}$ agree under the isomorphism \eqref{eqn:alg-log-map-matching-iso}.
        
        \item For each $1\le j\le r$, define $\delta_j\in\Gamma(C,\ol M_C)$ to be the unique section which maps to $(\{\delta_{v,j}\},\{\mu_{i,j}\})$ under the isomorphism \eqref{eqn:ghost-sections-log-curve}. By Theorem~\ref{thm-part:global-slb-log-curve}, there is a unique isomorphism of virtual Cartier divisors on $C$
        \begin{align*}
            \Psi_j:f^*(\clO_X(D_j),\taut_{D_j})\xrightarrow{\sim}(\clO_{C^\dagger}(\delta_j),\sigma_{C^\dagger,\delta_j})
        \end{align*}
        such that we have $\hat\Psi_{v,j} = \Phi_{v,j}\circ(\nu|_{\tilde C_v})^*\Psi_j$ for all vertices $v$ of $\plC$. Here, $\Phi_{v,j}$ is the natural isomorphism in \eqref{eqn:alg-log-map-vcd-on-irr-comp} provided by Theorem~\ref{thm-subpart:log-curve-slb-on-irr-comps}.

        \item Apply the equivalence \ref{concrete-log-mor-to-snc-pair} $\to$ \ref{abstract-log-mor-to-snc-pair} from Proposition~\ref{prop:log-mor-to-snc-pair} to the collection $\{\delta_j,\Psi_j\}$ to obtain a log morphism $f^\dagger:C^\dagger\to X^\dagger$ with underlying morphism $f:C\to X$.
    \end{enumerate}
    The assignment $(\{\rho_e,\Theta_e\},\bmu,\{\delta_{v,j},\Theta_{v,j}\}) \mapsto (\xi^\dagger/S^\dagger,f^\dagger)$ is clearly well-defined. Once again, the functoriality of this assignment follows from Proposition~\ref{prop:log-mor-to-snc-pair} and Theorem~\ref{thm:log-curve-over-log-pt}.
\end{proof}

\begin{Remark}
    Recall the setup from Discussion~\ref{disc:log-map-motivation} consisting of a $1$-parameter family $(\esC\to\esS,{\bf x},F)$ of prestable maps to $X$, whose general fibre is smooth and has prescribed contact orders with the irreducible components of $D$ at its marked points.
    
    Define the log point $S^\dagger = (S,M_S)$ by $M_S := \iota_s^*M_{(\esS,s)}$, where $\iota_s:S = \Spec\bC\to\esS$ is the morphism corresponding to the distinguished point $s\in\esS$ and $M_{(\esS,s)}$ is the divisorial log structure corresponding to $\{s\}\subset\esS$. In other words, we have $\ol M_S = \bN$ and $\clO_{S^\dagger}(1) = T_{\esS,s}$. By Proposition~\ref{prop:alg-log-map}, the output of Discussion~\ref{disc:log-map-motivation} is equivalent to an algebraic log map $(\xi^\dagger/S^\dagger,f^\dagger)$ to $X^\dagger$ whose underlying prestable map is the fibre of $(\esC\to\esS,{\bf x},F)$ over $s$. 
\end{Remark}

\begin{Definition} The following notions arise naturally from Proposition~\ref{prop:alg-log-map}.
    \begin{enumerate}[ref = \theDefinition(\arabic*), label = (\arabic*), itemsep = 0.3ex]
        \item\label{def-part:alg-log-map-trop-data} Let $(\xi^\dagger/S^\dagger,f^\dagger)$ be an algebraic log map to $X^\dagger$. Let $(C,{\bf x},f)$ denote its underlying prestable map and consider the corresponding categories \ref{abstract-log-map} and \ref{concrete-log-map}. 
        Let $(\{\rho_e,\Theta_e\},\bmu,\{\delta_{v,j},\Theta_{v,j}\})$ be the tuple corresponding to the object $(\xi^\dagger/S^\dagger,f^\dagger)$ under the equivalence \ref{abstract-log-map} $\simeq$ \ref{concrete-log-map} from Proposition~\ref{prop:alg-log-map}.
        
        \begin{enumerate}[ref = \theenumi(\roman*), label = (\roman*), itemsep = 0.15ex]
            \item We refer to $\bmu$ as the \emph{associated contact datum} of $(\xi^\dagger/S^\dagger,f^\dagger)$.

            \item\label{def-part:alg-log-map-tropicalization} Recall the fine monoid $\ol M(C,{\bf x},f,\bmu)$ together with its reduced generating set from Construction~\ref{cons:fine-basic-monoid}. The \emph{associated monoid map} of $(\xi^\dagger/S^\dagger,f^\dagger)$, denoted by
            \begin{align}\label{eqn:alg-log-map-tropicalization}
                \Trop(\xi^\dagger/S^\dagger,f^\dagger):\ol M(C,{\bf x},f,\bmu)\to\ol M_S,
            \end{align}
            is defined to be the monoid map given on generators by $m_{v,j}\mapsto\delta_{v,j}$ and $m_e\mapsto\rho_e$. 
            Conditions~\ref{concrete-alg-log-map-vcd-iso} and \ref{concrete-alg-log-map-matching-condition} imply that $\Trop(\xi^\dagger/S^\dagger,f^\dagger)$ is well-defined and that it maps each element of the reduced generating set to $\ol M_S\setminus\{0\}$.
        \end{enumerate}

        \item\label{def:cat-of-log-lifts} Let $(C,{\bf x},f)$ be a prestable map to $X$ and consider the corresponding category \ref{abstract-log-map}. Let $\bmu$ be a locally realizable contact datum for $(C,{\bf x},f)$. Define 
        \begin{align*}
            \fM_{(C,{\bf x},f,\bmu)}(X^\dagger)\subset\text{\ref{abstract-log-map}}
        \end{align*}
        to be the full subcategory consisting of algebraic log maps whose underlying prestable map is $(C,{\bf x},f)$ and whose associated contact datum is $\bmu$. Proposition~\ref{prop:alg-log-map} implies that $\fM_{(C,{\bf x},f,\bmu)}(X^\dagger)$ is a (possibly empty) category fibred in groupoids over $\LogPt$.

        \item\label{def-part:fine-basic} In the setting of part (1), we call $(\xi^\dagger/S^\dagger,f^\dagger)$ a \emph{fine-basic} algebraic log map if it is a final object of the category $\fM_{(C,{\bf x},f,\bmu)}(X^\dagger)$.
    \end{enumerate}
\end{Definition}

\begin{Remark}\label{rem:why-fine-vs-fs-basic}
    In \cite{GS-log} and \cite{Chen-logDF1}, algebraic log maps are primarily studied in the setting of fs log schemes. 
    Therefore, they define the notion of \emph{basic} log maps by a universal property involving only fs log points; \cite{Chen-logDF1} uses the term \emph{minimal} instead of \emph{basic}.

    For the comparison with symplectic log maps, it is necessary to distinguish between the fine and the fs case. Accordingly, we will discuss two notions: \emph{fine-basic} and \emph{fs-basic}. The former has been introduced in Definition~\ref{def-part:fine-basic}; the latter will be introduced in Definition~\ref{def-part:fs-basic} and is equivalent to \emph{basic} in the sense of \cite{GS-log}, as explained in Remark~\ref{rem:fs-basic-agrees-with-GS-basic}.
\end{Remark}

\subsection{Comparison in the fine case}\label{subsec:log-lifts-of-stable-maps}

\begin{Theorem}\label{thm:fine-log-map-comparison}
    Let $(C,{\bf x},f)$ be a prestable map to $X$ and let $\bmu$ be a locally realizable contact datum for it with respect to $D$. Then the following statements are equivalent.
    \begin{enumerate}[label = (\arabic*), itemsep = 0.3ex]
        \item The category $\fM_{(C,{\bf x},f,\bmu)}(X^\dagger)$ has a final object.
        \item The category $\fM_{(C,{\bf x},f,\bmu)}(X^\dagger)$ is non-empty.
        \item The tuple $(C,{\bf x},f,\bmu)$ is a symplectic stable log map.
    \end{enumerate}

    Moreover, if any (and therefore each) of these three statements is true, then an object $(\xi^\dagger/S^\dagger,f^\dagger)$ of $\fM_{(C,{\bf x},f,\bmu)}(X^\dagger)$ is a fine-basic algebraic log map if and only if the associated monoid map $\Trop(\xi^\dagger/S^\dagger,f^\dagger):\ol M(C,{\bf x},f,\bmu)\to\ol M_S$ is an isomorphism.
\end{Theorem}
\begin{proof}
    We will draw heavily on Definition~\ref{def:slb-presentation} and the statement of Proposition~\ref{prop:slb-presentation} and so, the reader is advised to review these first.

    (1) $\Rightarrow$ (2): This is clear.

    (2) $\Rightarrow$ (3): Choose an object $(\xi^\dagger/S^\dagger,f^\dagger)$ of $\fM_{(C,{\bf x},f,\bmu)}(X^\dagger)$ and let $(\{\rho_e,\Theta_e\},\bmu,\{\delta_{v,j},\Theta_{v,j}\})$ denote the tuple which corresponds to it under the equivalence \ref{abstract-log-map} $\simeq$ \ref{concrete-log-map} from Proposition~\ref{prop:alg-log-map}. Let $\clL_{S^\dagger}$ denote the global slb on $S^\dagger$. By Definition~\ref{def-part:slb-presentation-realization-category}, the tuple
    \begin{align*}
        (\ol M_S,\Trop(\xi^\dagger/S^\dagger,f^\dagger),\clL_{S^\dagger},\{\Theta_{v,j},\Theta_e\})
    \end{align*}
    is an object of the category of realizations $\|\fP(C,{\bf x},f,\bmu)\|$ of the slb presentation $\fP(C,{\bf x},f,\bmu)$. The implication (2) $\Rightarrow$ (3) in Proposition~\ref{prop:slb-presentation} now shows that $\fP(C,{\bf x},f,\bmu)$ is consistent. We conclude from Lemma~\ref{lem:basic-fine-slb-consistent-symplectic-log-map} that $(C,{\bf x},f,\bmu)$ is a symplectic log map.

    (3) $\Rightarrow$ (1): By Lemma~\ref{lem:basic-fine-slb-consistent-symplectic-log-map}, the associated slb presentation $\fP(C,{\bf x},f,\bmu)$ on $S = \Spec\bC$ is consistent. The implication (3) $\Rightarrow$ (1) in Proposition~\ref{prop:slb-presentation} produces an initial object $(P,\psi,\clL,\{\Theta_{v,j},\Theta_e\})$ of the category $\|\fP(C,{\bf x},f,\bmu)\|$. This has the following features.    
    \begin{enumerate}[ref = F\arabic*, label = F\arabic*., itemsep = 0.3ex]
        \item\label{slb-sharp-monoid} By the last assertion in Proposition~\ref{prop:slb-presentation}, $\psi:\ol M(C,{\bf x},f,\bmu)\xrightarrow{\sim} P$ is an isomorphism of fine monoids. From Condition~\ref{def-part:tropical-map-existence-condition} and Lemma~\ref{lem:fine-basic-monoid-sharp}, it follows that $P$ is sharp.

        \item\label{slb-all-nontrivial-sections-0} Let us write $\clL = \{L_p,\Phi_{p_1,p_2},\sigma_p\}$. Define $\delta_{v,j}\in P$ (resp. $\rho_e\in P$) to be the images under $\psi$ of $m_{v,j}$ (resp. $m_e$). The virtual Cartier divisor isomorphisms
        \begin{align*}
            \Theta_{v,j}:(\fL_{v,j},\fs_{v,j})\xrightarrow{\sim}(L_{\delta_{v,j}},\sigma_{\delta_{v,j}}) \quad\text{and}\quad \Theta_e:(\fL_e,\fs_e)\xrightarrow{\sim}(L_{\rho_e},\sigma_{\rho_e})
        \end{align*}
        imply that, for each vertex $v$ of $\plC$ and $j\in I_v$, we have $\sigma_{\delta_{v,j}} = 0$ and, for each edge $e$ of $\plC$, we have $\sigma_{\rho_e} = 0$. It follows that we have $\sigma_p = 0$ for all $0\ne p\in P$.
    \end{enumerate}
    
    Using \ref{slb-sharp-monoid}, \ref{slb-all-nontrivial-sections-0} and the equivalence $\SlbPt\simeq\LogPt$ from Proposition~\ref{prop:cat-of-log-pts}, we get a log point $S^\dagger = (S,M_S)$ for which $P$ (resp. $\clL$) is identified with the ghost sheaf $\ol M_S$ (resp. the global slb on $S^\dagger$). Using $\delta_{v,j},\rho_e\in P\simeq\ol M_S$, we get the object $(\{\rho_e,\Theta_e\},\bmu,\{\delta_{v,j},\Theta_{v,j}\})$ of \ref{concrete-log-map} lying over the log point $S^\dagger$. Let $(\xi^\dagger/S^\dagger,f^\dagger)$ denote the object of $\fM_{(C,{\bf x},f,\bmu)}(X^\dagger)$ corresponding to it under the equivalence \ref{abstract-log-map} $\simeq$ \ref{concrete-log-map} from Proposition~\ref{prop:alg-log-map}. Note that we have $\psi = \Trop(\xi^\dagger/S^\dagger,f^\dagger)$ by construction. Using \ref{slb-all-nontrivial-sections-0}, Proposition~\ref{prop:cat-of-log-pts} and the fact that $(P,\psi,\clL,\{\Theta_{v,j},\Theta_e\})$ is an initial object of the category $\|\fP(C,{\bf x},f,\bmu)\|$, one checks that $(\xi^\dagger/S^\dagger,f^\dagger)$ is a final object of the category $\fM_{(C,{\bf x},f,\bmu)}(X^\dagger)$.

    Lastly, an object $(\xi'^\dagger/S'^\dagger,f'^\dagger)$ of $\fM_{(C,{\bf x},f,\bmu)}(X^\dagger)$ is a final object if and only if the (unique) morphism $(\xi'^\dagger/S'^\dagger,f'^\dagger) \to (\xi^\dagger/S^\dagger,f^\dagger)$ is an isomorphism. By the pullback diagram \eqref{eqn:map-of-log-str-cartesian} and the equivalence \ref{abstract-log-map} $\simeq$ \ref{concrete-log-map} from Proposition~\ref{prop:alg-log-map}, this happens if and only if the induced map $\ol M_S\to\ol M\phantom{}'_S$ on ghost sheaves is an isomorphism. Since $\Trop(\xi^\dagger/S^\dagger,f^\dagger) = \psi$ is an isomorphism, we conclude that $(\xi'^\dagger/S'^\dagger,f'^\dagger)$ is a final object if and only if the map $\Trop(\xi'^\dagger/S'^\dagger,f'^\dagger):\ol M(C,{\bf x},f,\bmu)\to\ol M\phantom{}'_S$ is an isomorphism.
\end{proof}

\begin{Definition}
    Given an algebraic log map $(\xi^\dagger/S^\dagger,f^\dagger)$ to $X^\dagger$, denote its underlying prestable map and its associated contact datum by $(C,{\bf x},f)$ and $\bmu$ respectively.
    
    We call $(C,{\bf x},f,\bmu)$ the \emph{underlying symplectic log map} of $(\xi^\dagger/S^\dagger,f^\dagger)$.

    This terminology is justified by Theorem~\ref{thm:fine-log-map-comparison}. 
\end{Definition}

\begin{Corollary}\label{cor:fine-log-map-comparison}
    For any symplectic log map $(C,{\bf x},f,\bmu)$ to $(X,D)$, there is a unique (up to unique isomorphism) fine-basic algebraic log map $(\xi^\dagger/S^\dagger,f^\dagger)$ to $X^\dagger$ whose underlying symplectic log map is given by $(C,{\bf x},f,\bmu)$.
\end{Corollary}

\subsection{Saturated log maps}\label{subsec:fs-stable-log-maps}

\begin{Notation}
    For a finitely generated abelian group $A$, we denote its torsion subgroup by $A_\text{tor}\subset A$. By the structure theorem for finitely generated abelian groups, $A_\text{tor}$ is the direct sum of finitely many finite cyclic groups.
\end{Notation}

\begin{Construction}\label{cons:saturation-data}
    Consider a symplectic log map $(C,{\bf x},f,\bmu)$ to $(X,D)$.
    
    \begin{enumerate}[ref = \theConstruction(\arabic*), label = (\arabic*), itemsep = 0.3ex]
        \item\label{cons-part:sat-data-before-quotient} Define $\widetilde\fS(C,{\bf x},f,\bmu)$ to be the (non-empty) set consisting of all possible choices of $0\ne\sigma_{v,j}\in H^0(\tilde C_v,L_{v,j})$ and $0\ne\theta_e\in\Lambda_{C,y_e}$ satisfying Condition~\ref{def-part:analytic-condition}, i.e.,
        \begin{align}\label{eqn:sat-data-before-quotient}
            \widetilde\fS(C,{\bf x},f,\bmu) := \left\{ 
            \{\sigma_{v,j},\theta_e\}
            \,\middle|\!
            \begin{array}{rcl}
                j\not\in I_v & \!\!\!\Rightarrow & \!\!\! \varphi_{v,j}(\sigma_{v,j}) = 1 \\
                \vec{e}:v\to v' & \!\!\!\Rightarrow & \!\!\! \varphi_{\vec{e},j}(\sigma_{v',j}) = \sigma_{v,j}\otimes\theta_e^{\otimes\mu_{\vec{e},j}}
            \end{array}\!\!\!
            \right\}.
        \end{align}
        
        \item\label{cons-part:sat-data-groups} Using the abbreviation $Q = \ol M(C,{\bf x},f,\bmu)$, define the abelian groups 
        \begin{align*}
            T_{(C,{\bf x},f,\bmu)}&:=\Hom_{\Ab}(Q^\text{gp},\bC^\times),\\ T^0_{(C,{\bf x},f,\bmu)}&:=\Hom_{\Ab}(Q^\text{gp}/Q^\text{gp}_\text{tor},\bC^\times),\quad\text{and}\\
            \pi_0(T_{(C,{\bf x},f,\bmu)})&:=\Hom_{\Ab}(Q^\text{gp}_\text{tor},\bC^\times).
        \end{align*}
        The restriction map $T_{(C,{\bf x},f,\bmu)}\to \pi_0(T_{(C,{\bf x},f,\bmu)})$ is surjective and its kernel is $T^0_{(C,{\bf x},f,\bmu)}$. Evaluating a homomorphism $Q^\text{gp}\to\bC^\times$ at $m_{v,j}$ and $m_e$ yields an identification
        \begin{align}\label{eqn:sat-data-groups}
            T_{(C,{\bf x},f,\bmu)} \xrightarrow{\sim} \left\{
            \{t_{v,j}\in\bC^\times,\varepsilon_e\in\bC^\times\}
            \,\middle|\!
            \begin{array}{rcl}
                j\not\in I_v & \!\!\!\Rightarrow & \!\!\! t_{v,j} = 1 \\
                \vec{e}:v\to v' & \!\!\!\Rightarrow & \!\!\! t_{v',j} = t_{v,j}\varepsilon_e^{\mu_{\vec{e},j}}
            \end{array}\!\!\!
            \right\}.
        \end{align}
        
        \item\label{cons-part:sat-data-group-action} Using the notation from \eqref{eqn:sat-data-before-quotient} and \eqref{eqn:sat-data-groups}, define a transitive free action of the group $T_{(C,{\bf x},f,\bmu)}$ on the set $\widetilde\fS(C,{\bf x},f,\bmu)$ by the formula
        \begin{align*}
            \{\sigma_{v,j},\theta_e\}\cdot\{t_{v,j},\varepsilon_e\} := \{t_{v,j}\sigma_{v,j},\varepsilon_e\theta_e\}.
        \end{align*}

        \item Define the \emph{set of saturation data} for $(C,{\bf x},f,\bmu)$ to be the quotient 
        \begin{align*}
            \fS(C,{\bf x},f,\bmu):=\widetilde\fS(C,{\bf x},f,\bmu)/T^0_{(C,{\bf x},f,\bmu)}.
        \end{align*} 
        By construction, $\fS(C,{\bf x},f,\bmu)$ is a non-empty set and the action in part (3) descends to a transitive free action of $\pi_0(T_{(C,{\bf x},f,\bmu)})$ on $\fS(C,{\bf x},f,\bmu)$.
    \end{enumerate}
    
    An element of $\fS(C,{\bf x},f,\bmu)$ is called a \emph{saturation datum} for $(C,{\bf x},f,\bmu)$.
\end{Construction}

\begin{Lemma}\label{lem:number-of-saturations}
    In the setting of Construction~\ref{cons:saturation-data}, we have
    \begin{align*}
        |\fS(C,{\bf x},f,\bmu)| = |\ol M(C,{\bf x},f,\bmu)^\text{gp}_\text{tor}|.
    \end{align*}
\end{Lemma}
\begin{proof}
    This is immediate from the fact that $\fS(C,{\bf x},f,\bmu)$ is a $\pi_0(T_{(C,{\bf x},f,\bmu)})$-torsor.
\end{proof}

\begin{Remark}
    By examining Construction~\ref{cons:fine-basic-monoid}, it is not hard to see that $\ol M(C,{\bf x},f,\bmu)^\text{gp}_\text{tor}$ is always contained in the subgroup of $\ol M(C,{\bf x},f,\bmu)^\text{gp}$ generated by the elements $\{m_e\}_e$.
\end{Remark}

\begin{Definition}\label{def:sat-log-map}
    We define saturated log maps as follows.
    \begin{enumerate}[ref = \theDefinition(\arabic*), label = (\arabic*), itemsep = 0.3ex]
        \item\label{def-part:sat-symp-log-map} A \emph{saturated symplectic log map} to $(X,D)$ is a tuple $(C,{\bf x},f,\bmu,\eta)$ such that $(C,{\bf x},f,\bmu)$ is a symplectic log map to $(X,D)$ and $\eta$ is a saturation datum for it.
        
        \item A \emph{saturated algebraic log map} to $X^\dagger$ is an algebraic log map $(\xi^\dagger/S^\dagger,f^\dagger)$ to $X^\dagger$ such that $S^\dagger = (S,M_S)$ is a fine saturated (fs) log point.
    \end{enumerate}
\end{Definition}

\begin{Remark}\label{rem:sat-symp-log-map}
    It is not hard to extend Definition~\ref{def-part:sat-symp-log-map} to the more general setting of symplectic manifolds and snc symplectic divisors considered in \cite{Tehrani-log}; see Remark~\ref{rem:symp-log-map}.
\end{Remark}

\begin{Remark}
    In the symplectic case, being saturated is an additional structure on a log map whereas, in the algebraic case, it is an additional property.
\end{Remark}

\begin{Example}\label{exa:simplest-sat-example}
    Let $C$ be the prestable curve obtained by gluing two smooth genus $0$ curves at two points. More precisely, suppose that the dual graph of $C$ consists of two vertices $v,v'$ and two edges $e_1,e_2$ between them and we have $C_v\simeq C_{v'}\simeq\bP^1$. Choose orientations $\vec{e}_1,\vec{e}_2$ going from $v$ to $v'$. We then have $y_{\vec{e}_1},y_{\vec{e}_2}\in C_v$ and $y_{\cev{e}_1},y_{\cev{e}_2}\in C_{v'}$.

    Let $X$ be a smooth variety and let $D=D_1\subset X$ be a smooth divisor. Let $(C,f)$ be a stable map (without marked points) such that $f_v(C_v)\not\subset D$ and $f_{v'}(C_{v'})\subset D$. Let $\bmu$ be a locally realizable contact datum for $(C,f)$ with respect to $D$ and, for simplicity, write $\mu_1,\mu_2$ instead of $\mu_{\vec{e}_1,1},\mu_{\vec{e}_2,1}$ respectively. Concretely, this amounts to the following properties.
    \begin{enumerate}[label = (\roman*), itemsep = 0.3ex]
        \item We have an equality $f_v^{-1}(D) = \mu_1\cdot y_{\vec{e}_1} + \mu_2\cdot y_{\vec{e}_2}$ of Cartier divisors on $C_v$.

        \item The line bundle $f_{v'}^*\clN_{D/X}\otimes\clO_{C_{v'}}(\mu_1\cdot y_{\cev{e}_1} + \mu_2\cdot y_{\cev{e}_2})$ on $C_{v'}$ is trivial.
    \end{enumerate}
    
    The fine monoid associated to $(C,f,\bmu)$ is given by the integral presentation
    \begin{align*}
        \ol M(C,f,\bmu) = \left\langle 
        \!\!
        \begin{array}{ll}
            m_v,m_{v'}, \\
            m_{e_1},m_{e_2}
        \end{array}
        \middle|
        \begin{array}{l}
             m_v = 0\\
             m_{v'} = m_v + \mu_1m_{e_1}\\
             m_{v'} = m_v + \mu_2m_{e_2}
        \end{array}
        \!\!
        \right\rangle^\text{int} = 
        \left\langle
            m_{e_1},m_{e_2}
            \,\middle|\,
            \mu_1m_{e_1} = \mu_2m_{e_2}
        \right\rangle^\text{int},
    \end{align*}
    where we have written $m_v,m_{v'}$ respectively in place of $m_{v,1},m_{v',1}$. One verifies that $(C,f,\bmu)$ is a symplectic stable log map. The associated slb presentation from Construction~\ref{cons:fine-basic-slb-presentation} induces a natural isomorphism $\varphi:\Lambda_{C,y_{e_1}}^{\otimes\mu_1}\xrightarrow{\sim}\Lambda_{C,y_{e_2}}^{\otimes\mu_2}$ using the relation $\mu_1m_{e_1} = \mu_2m_{e_2}$.
    
    The group $\ol M(C,f,\bmu)^\text{gp}_\text{tor}$ is cyclic of order $d := \gcd(\mu_1,\mu_2)\ge 1$ and is generated by the element $(\mu_1/d)m_{e_1} - (\mu_2/d)m_{e_2}$. One verifies that any element $\{\sigma_{v,1},\sigma_{v',1},\theta_{e_1},\theta_{e_2}\}$ of the set $\widetilde\fS(C,f,\bmu)$ from Construction~\ref{cons-part:sat-data-before-quotient} yields an isomorphism
    \begin{align*}
        \psi:\Lambda_{C,y_{e_1}}^{\otimes(\mu_1/d)} \xrightarrow{\sim}\Lambda_{C,y_{e_2}}^{\otimes(\mu_2/d)}, \quad\text{given by} \quad\psi(\theta_{e_1}^{\otimes(\mu_1/d)}) = \theta_{e_2}^{\otimes(\mu_2/d)},
    \end{align*}
    satisfying $\psi^{\otimes d} = \varphi$ and that this assignment descends to a bijection between $\fS(C,f,\bmu)$ and the set of $d^\text{th}$ roots of the isomorphism $\varphi$. In particular, the number of saturation data for the symplectic stable log map $(C,f,\bmu)$ is exactly $d$, as predicted by Lemma~\ref{lem:number-of-saturations}.
\end{Example}

\begin{Construction}\label{cons:underlying-sat-symp-log-map}
    Given a saturated algebraic log map $(\xi^\dagger/S^\dagger,f^\dagger)$ to $X^\dagger$, denote its underlying symplectic log map by $(C,{\bf x},f,\bmu)$.
    \begin{enumerate}[ref = \theConstruction(\arabic*), label = (\arabic*), itemsep = 0.3ex]
        \item Let $(\{\rho_e,\Theta_e\},\bmu,\{\delta_{v,j},\Theta_{v,j}\})$ denote the tuple obtained by applying the equivalence \ref{abstract-log-map} $\simeq$ \ref{concrete-log-map} from Proposition~\ref{prop:alg-log-map} to the algebraic log map $(\xi^\dagger/S^\dagger,f^\dagger)$. 
        
        \item For convenience, we use the abbreviations $Q = \ol M(C,{\bf x},f,\bmu)$ and $\psi = \Trop(\xi^\dagger/S^\dagger,f^\dagger)$. Since $\ol M_S$ is a sharp fs monoid, Remark~\ref{rem:sharp-fs-implies-no-torsion} shows that the map $\psi^\text{gp}:Q^\text{gp}\to\ol M\phantom{}_S^\text{gp}$ vanishes on the torsion subgroup $Q^\text{gp}_\text{tor}\subset Q^\text{gp}$. In particular, the image of the induced map $\Hom_{\Ab}(\ol M\phantom{}_S^\text{gp},\bC^\times)\to T_{(C,{\bf x},f,\bmu)}$ is contained in the subgroup $T^0_{(C,{\bf x},f,\bmu)}$, where the groups $T_{(C,{\bf x},f,\bmu)}$ and $T_{(C,{\bf x},f,\bmu)}^0$ are as in Construction~\ref{cons-part:sat-data-groups}.
        
        \item\label{cons-part:sat-datum-from-fs-alg-log-map} Using $(\xi^\dagger/S^\dagger,f^\dagger)$, we construct a saturation datum $\eta$ for $(C,{\bf x},f,\bmu)$ as follows.
        \begin{enumerate}[ref = \theenumi(\roman*), label = (\roman*), itemsep = 0.15ex]
            \item By Lemma~\ref{lem:local-lift-from-ghost}, the quotient map $M_S\to\ol M\phantom{}_S$ admits a section. Moreover, the set of such sections is a $\Hom_{\Ab}(\ol M\phantom{}_S^\text{gp},\bC^\times)$-torsor; see Remark~\ref{rem:splitting-log-pt}.
            
            \item\label{cons-subpart:sat-datum-from-split-fs-log-map} Fix a section $\beta$ of $M_S\to\ol M_S$. For each $m\in\ol M_S$, we have $\beta(m)\in\clO_{S^\dagger}(-m)^\times$. Define the elements $\sigma_{v,j}\in H^0(\tilde C_v,L_{v,j})$ and $\theta_e\in\Lambda_{C,y_e}$ by the equations
            \begin{align}\label{eqn:lifted-sat-datum-from-splitting}
                \langle\beta(\delta_{v,j}),\Theta_{v,j}(\sigma_{v,j})\rangle = 1 \quad\text{and}\quad \langle\beta(\rho_e),\Theta_e(\theta_e)\rangle = 1,
            \end{align}
            where $\langle\cdot,\cdot\rangle$ is the duality pairing induced by the isomorphisms appearing in the global slb on $S^\dagger$; see Discussion~\ref{disc-part:log-lb-vcd-dual}. Using \ref{concrete-alg-log-map-vcd-iso}, \ref{concrete-alg-log-map-matching-condition} and the fact that $\beta:\ol M_S\to M_S$ is a monoid map, we conclude that the element $\{\sigma_{v,j},\theta_e\}$ belongs to the set $\widetilde\fS(C,{\bf x},f,\bmu)$ from Construction~\ref{cons-part:sat-data-before-quotient}.

            \item Given any other section $\beta'$ of $M_S\to\ol M_S$,  write it as $\beta' = \beta\cdot\lambda$ for a unique homomorphism $\lambda:\ol M\phantom{}_S^\text{gp}\to\bC^\times$. From this, we obtain elements
            \begin{align*}
                t_{v,j} := \lambda(\delta_{v,j})\in\bC^\times \quad\text{and}\quad \varepsilon_e := \lambda(\rho_e)\in\bC^\times.
            \end{align*}
            
            \item From \eqref{eqn:lifted-sat-datum-from-splitting}, we see that the elements $\{\sigma_{v,j},\theta_e\}$ and $\{\sigma_{v,j}',\theta_e'\}$ of $\widetilde\fS(C,{\bf x},f,\bmu)$ determined by the sections $\beta$ and $\beta'$ respectively are related by the equations 
            \begin{align*}
                \sigma'_{v,j} = t_{v,j}^{-1}\sigma_{v,j} \quad\text{and}\quad \theta_e' = \varepsilon_e^{-1}\theta_e.
            \end{align*}
            Note that, under the identification \eqref{eqn:sat-data-groups}, the element $\{t_{v,j},\varepsilon_e\}$ corresponds to the element $\lambda\circ\psi^\text{gp}:Q^\text{gp}\to\bC^\times$ of the group $T_{(C,{\bf x},f,\bmu)}$. From part (2), we conclude that the element $\lambda\circ\psi^\text{gp}$ in fact lies in the subgroup $T^0_{(C,{\bf x},f,\bmu)}$.

            \item Define $\eta$ to be the image of $\{\sigma_{v,j},\theta_e\}$ in the set $\fS(C,{\bf x},f,\bmu)$ of saturation data, obtained using the chosen section $\beta$ of $M_S\to\ol M_S$. It follows from the preceding discussion that $\eta$ is well-defined independent of the choice of $\beta$.
        \end{enumerate}
    \end{enumerate}

    We call $(C,{\bf x},f,\bmu,\eta)$ the underlying \emph{saturated symplectic log map} of $(\xi^\dagger/S^\dagger,f^\dagger)$.
\end{Construction}

\begin{Lemma}\label{lem:sat-data-pullback}
    Let $(C,{\bf x},f,\bmu)$ be a symplectic log map to $(X,D)$. 

    Consider a morphism $(\xi'^\dagger/S'^\dagger,f'^\dagger)\to(\xi^\dagger/S^\dagger,f^\dagger)$ in the category $\fM_{(C,{\bf x},f,\bmu)}(X^\dagger)$ between saturated algebraic log maps. If $\eta$ and $\eta'$ are the elements of $\fS(C,{\bf x},f,\bmu)$ obtained from $(\xi^\dagger/S^\dagger,f^\dagger)$ and $(\xi'^\dagger/S'^\dagger,f'^\dagger)$ respectively via Construction~\ref{cons-part:sat-datum-from-fs-alg-log-map}, then we have $\eta = \eta'$.
\end{Lemma}
\begin{proof}
    By \eqref{eqn:map-of-log-str-cartesian}, we have $M_S\xrightarrow{\sim}\ol M_S\times_{\ol M\phantom{}_S'}M_S'$. Thus, sections of $M_S'\to\ol M\phantom{}_S'$ pull back to sections of $M_S\to\ol M_S$. The result is now immediate from definitions.
\end{proof}

\begin{Definition}
    The following notions arise naturally from Construction~\ref{cons:underlying-sat-symp-log-map}.
    \begin{enumerate}[ref = \theDefinition(\arabic*), label = (\arabic*), itemsep = 0.3ex]
        \item\label{def-part:fs-basic-monoid} Given a symplectic log map $(C,{\bf x},f,\bmu)$ to $(X,D)$, define the sharp fs monoid
        \begin{align*}
            \ol M\phantom{}^\text{fs}(C,{\bf x},f,\bmu) := \ol{\left(\ol M(C,{\bf x},f,\bmu)\right)^\text{sat}}
        \end{align*}
        to be the sharpening of the saturation of the sharp fine monoid $\ol M(C,{\bf x},f,\bmu)$; see Definitions~\ref{def:monoid-props} and \ref{def:integralization-saturation} as well as Lemma~\ref{lem-part:fine-sharp-to-fs-sharp}.

        \item In the setting of part (1), let $(\xi^\dagger/S^\dagger,f^\dagger)$ be a saturated algebraic log map belonging to the category $\fM_{(C,{\bf x},f,\bmu)}(X^\dagger)$. Since $\ol M_S$ is a sharp fs monoid, the associated monoid map $\Trop(\xi^\dagger/S^\dagger,f^\dagger):\ol M(C,{\bf x},f,\bmu)\to\ol M_S$ induces a monoid map
        \begin{align*}
            \Trop^\text{fs}(\xi^\dagger/S^\dagger,f^\dagger):\ol M\phantom{}^\text{fs}(C,{\bf x},f,\bmu)\to\ol M_S.
        \end{align*}
        
        \item Given a saturated symplectic log map $(C,{\bf x},f,\bmu,\eta)$ to $(X,D)$, define
        \begin{align*}
            \fM_{(C,{\bf x},f,\bmu,\eta)}^\text{fs}(X^\dagger)\subset \fM_{(C,{\bf x},f,\bmu)}(X^\dagger)
        \end{align*}
        to be the full subcategory consisting of saturated algebraic log maps $(\xi^\dagger/S^\dagger,f^\dagger)$ which have $(C,{\bf x},f,\bmu,\eta)$ as their underlying saturated symplectic log map.  Lemma~\ref{lem:sat-data-pullback} implies that $\fM_{(C,{\bf x},f,\bmu,\eta)}^\text{fs}(X^\dagger)$ is a category fibred in groupoids over $\LogPt^\text{fs}$, where $\LogPt^\text{fs}\subset\LogPt$ denotes the full subcategory consisting of fs log points.

        \item\label{def-part:fs-basic} Given a saturated algebraic log map $(\xi^\dagger/S^\dagger,f^\dagger)$, let $(C,{\bf x},f,\bmu,\eta)$ be its underlying saturated symplectic log map. We call $(\xi^\dagger/S^\dagger,f^\dagger)$ a \emph{fs-basic} algebraic log map if it is a final object of the category $\fM_{(C,{\bf x},f,\bmu,\eta)}^\text{fs}(X^\dagger)$.
    \end{enumerate}
\end{Definition}

\begin{Remark}\label{rem:fs-basic-agrees-with-GS-basic}
    The universal property in \cite[Proposition~1.24]{GS-log} shows that our notion of fs-basic log maps agrees with \emph{basic} (fs) log maps in the sense of \cite[Definition~1.20]{GS-log}.
\end{Remark}

\begin{Remark}
    Lemma~\ref{lem-part:fine-sharp-to-fs-sharp} shows that, in the setting of Definition~\ref{def-part:fs-basic-monoid}, we have a natural isomorphism $\ol M\phantom{}^{\text{fs}}(C,{\bf x},f,\bmu)\xrightarrow{\sim}\ol M(C,{\bf x},f,\bmu)^\vee\phantom{}^\vee$.
\end{Remark}

\subsection{Comparison in the fine saturated case}

\begin{Theorem}\label{thm:fs-log-map-comparison}
    Consider a symplectic log map $(C,{\bf x},f,\bmu)$ to $(X,D)$. Then, for each saturation datum $\eta\in\fS(C,{\bf x},f,\bmu)$, the category $\fM_{(C,{\bf x},f,\bmu,\eta)}^\text{fs}(X^\dagger)$ has a final object.
    
    Moreover, a saturated algebraic log map $(\xi^\dagger/S^\dagger,f^\dagger)$ belonging to $\fM_{(C,{\bf x},f,\bmu)}(X^\dagger)$ is fs-basic if and only if the map $\Trop^\text{fs}(\xi^\dagger/S^\dagger,f^\dagger):\ol M\phantom{}^\text{fs}(C,{\bf x},f,\bmu)\to\ol M_S$ is an isomorphism.
\end{Theorem}
\begin{proof}
    The argument is divided into three steps.

    \begin{step}
        Using Theorem~\ref{thm:fine-log-map-comparison}, choose a final object $(\xi^\dagger/S^\dagger,f^\dagger)$ of $\fM_{(C,{\bf x},f,\bmu)}(X^\dagger)$. Write $S^\dagger = (S,M)$ with $S = \Spec\bC$. Let $(\{\delta_{v,j},\Theta_{v,j}\},\bmu,\{\rho_e,\Theta_e\})$ be the tuple corresponding to $(\xi^\dagger/S^\dagger,f^\dagger)$ under the equivalence \ref{abstract-log-map} $\simeq$ \ref{concrete-log-map} from Proposition~\ref{prop:alg-log-map}. 
        
        For convenience, we use the abbreviations $\sC = \fM_{(C,{\bf x},f,\bmu)}(X^\dagger)$, $Q = \ol M(C,{\bf x},f,\bmu)$ and $\psi = \Trop(\xi^\dagger/S^\dagger,f^\dagger)$. Theorem~\ref{thm:fine-log-map-comparison} implies that the map $\psi:Q\xrightarrow{\sim}\ol M$ is an isomorphism.        
        Write $\sC^\text{fs}\subset\sC$ for the full subcategory consisting of saturated algebraic log maps. For each saturation datum $\eta\in\fS(C,{\bf x},f,\bmu)$, we use the abbreviation $\sC^\text{fs}_\eta = \fM^\text{fs}_{(C,{\bf x},f,\bmu,\eta)}(X^\dagger)$.
    \end{step}

    \begin{step}
        In this step, we find a simpler but equivalent description of the category $\sC^\text{fs}$.
        
        Define the category $M/\Log(S)$ as follows. Objects of $M/\Log(S)$ are pairs $(M',\varphi')$, where $M'$ is a (fine) log structure on the point $S = \Spec\bC$ and $\varphi':M\to M'$ is a map of log structures. Morphisms $(M',\varphi')\to(M'',\varphi'')$ in $M/\Log(S)$ are maps $F:M'\to M''$ of log structures on $S$ for which we have $F\circ\varphi' = \varphi''$.
        
        Define the functor $\clR:\sC^\text{op}\to M/\Log(S)$ from the opposite category of $\sC$ to the category $M/\Log(S)$ as follows. For any object $(\xi'^\dagger/S'^\dagger,f'^\dagger)$ of $\sC$, its image under $\clR$ is the pair $(M',\varphi')$, where $S'^\dagger = (S,M')$ and $\varphi':M\to M'$ denotes the map induced by the unique morphism $(\xi'^\dagger/S'^\dagger,f'^\dagger)\to(\xi^\dagger/S^\dagger,f^\dagger)$ in $\sC$. For any morphism $(\xi''^\dagger/S''^\dagger,f''^\dagger)\to(\xi'^\dagger/S'^\dagger,f'^\dagger)$ in $\sC$, its image under $\clR$ is the induced map $M'\to M''$, where $S'^\dagger = (S,M')$ and $S''^\dagger = (S,M'')$. 
        The functor $\clR$ is an equivalence of categories since $(\xi^\dagger/S^\dagger,f^\dagger)$ is a final object of $\sC$. 
        
        Define $M/\Log^\text{fs}(S)\subset M/\Log(S)$ to be the full subcategory consisting of pairs $(M',\varphi')$ for which $M'$ is an fs log structure. Then $\clR$ restricts to an equivalence of categories
        \begin{align*}
            \clR^\text{fs}:(\sC^\text{fs})^\text{op}\xrightarrow{\sim}M/\Log^\text{fs}(S).
        \end{align*}
    \end{step}

    \begin{step}
        In this step, we determine the subcategories of $M/\Log^\text{fs}(S)$ which correspond to the subcategories $\sC^\text{fs}_\eta\subset\sC^\text{fs}$ under $\clR^\text{fs}$. In the process, we also construct and characterize initial objects for each of these subcategories of $M/\Log^\text{fs}(S)$.

        Recall that $Q^\text{gp}_\text{tor}\subset Q^\text{gp}$ denotes the torsion subgroup. Define $\widehat\fS(M)$ to be the set of group homomorphisms $\hat\eta:Q^\text{gp}_\text{tor}\to M^\text{gp}$ which make the following diagram commute.
        \begin{equation*}
        \begin{tikzcd}
            Q^\text{gp}_\text{tor} \arrow[d,hook] \arrow[r,"\hat\eta"] & M^\text{gp} \arrow[d] \\
             Q^\text{gp} \arrow[r,"\psi^\text{gp}","\sim"'] & \ol M\phantom{}^\text{gp}
        \end{tikzcd}
        \end{equation*}
        
        By Lemma~\ref{lem:local-lift-from-ghost}, $\psi^\text{gp}$ admits a lift $Q^\text{gp}\to M^\text{gp}$ and restricting this lift to $Q^\text{gp}_\text{tor}$ shows that $\widehat\fS(M)$ is non-empty. Moreover, the group $\Hom_{\Ab}(Q^\text{gp}_\text{tor},\bC^\times) = \pi_0(T_{(C,{\bf x},f,\bmu)})$ acts on $\widehat\fS(M)$ by multiplication; it is clear that this action is transitive and free.

        For any object $(M',\varphi')$ of $M/\Log^\text{fs}(S)$, define an element $\hat\eta(M',\varphi')\in\widehat\fS(M)$ via the following commutative diagram, where ${\bf 1}:Q^\text{gp}_\text{tor}\to M'^{\text{gp}}$ denotes the trivial homomorphism, the bottom right square is a pullback diagram by \eqref{eqn:map-of-log-str-cartesian} and the homomorphism given by the composition $Q^\text{gp}_\text{tor}\subset Q^\text{gp}\to\ol M^\text{gp}\to\ol M\phantom{}'^{\text{gp}}$ is trivial using Remark~\ref{rem:sharp-fs-implies-no-torsion}.
        \begin{equation}\label{eqn:torsion-embedding-induced-by-fs}
        \begin{tikzcd}
            Q^\text{gp}_\text{tor} \arrow[dr,"\hat\eta{(M',\varphi')}"] \arrow[drr,bend left,"{\bf 1}"] \arrow[dd,hook] & & \\
            & M^{\text{gp}} \arrow[r,"{\varphi'^{\text{gp}}}"] \arrow[d] \arrow[dr,phantom,"\lrcorner",pos=0.001] & M\phantom{}'^{\text{gp}} \arrow[d] \\
            Q^\text{gp} \arrow[r,"\psi^\text{gp}","\sim"'] & \ol M^{\text{gp}} \arrow[r,"{\ol\varphi\phantom{}'^{\text{gp}}}"] & \ol M\phantom{}'^{\text{gp}}
        \end{tikzcd}
        \end{equation}
        
        For any morphism $(M',\varphi')\to(M'',\varphi'')$ in the category $M/\Log^\text{fs}(S)$, it follows from the definition that we have the equality $\hat\eta(M',\varphi') = \hat\eta(M'',\varphi'')$. 
        
        For each $\hat\eta\in\widehat\fS(M)$, define $(M/\Log^\text{fs}(S))_{\hat\eta}\subset M/\Log^\text{fs}(S)$ to be the full subcategory consisting of the pairs $(M',\varphi')$ for which we have $\hat\eta(M',\varphi') = \hat\eta$.

        \begin{claim}\label{claim:fs-initial-object}
            For each $\hat\eta\in \widehat\fS(M)$, the category $(M/\Log^\text{fs}(S))_{\hat\eta}$ has an initial object. Moreover, an object $(M',\varphi')$ of $(M/\Log^\text{fs}(S))_{\hat\eta}$ is an initial object if and only if the map $\ol{Q^\text{sat}}\to\ol M\phantom{}'$ induced by $\ol{\varphi}\phantom{}'\circ\psi:Q\to\ol M\phantom{}'$ is an isomorphism.
        \end{claim}
        \begin{proof}[Proof of claim]
            Fix $\hat\eta\in\widehat\fS(M)$ and recall that $\ol M\phantom{}^\text{gp}_\text{tor}\subset\ol M\phantom{}^\text{gp}$ denotes the torsion subgroup. 
            We then have the following commutative diagram with exact rows, where we use $0$ (resp. $1$) respectively to denote the trivial group written additively (resp. multiplicatively).
            \begin{equation*}
            \begin{tikzcd}
                1 \arrow[r] & \bC^\times \arrow[d,equal] \arrow[r] & M^{\text{gp}} \arrow[r] \arrow[d] & \ol M\phantom{}^{\text{gp}} \arrow[r] \arrow[d] & 0 \\
                1 \arrow[r] & \bC^\times \arrow[r] & M^{\text{gp}}/\hat\eta(Q^{\text{gp}}_{\text{tor}}) \arrow[r] & \ol M\phantom{}^{\text{gp}}/\ol M\phantom{}^{\text{gp}}_{\text{tor}} \arrow[r] & 0
            \end{tikzcd}
            \end{equation*}
            
            By Lemma~\ref{lem-part:saturation-units-torsion}, the group of units in $(\ol M)^\text{sat}$ is $\ol M\phantom{}^\text{gp}_\text{tor}$. Using this, define
            \begin{align*}
                M_{\hat\eta}\subset M^\text{gp}/\hat\eta(Q^\text{gp}_\text{tor})
            \end{align*}
            to be the inverse image of the submonoid $\ol{(\ol M)^\text{sat}} = (\ol M)^\text{sat}/\ol M\phantom{}^\text{gp}_\text{tor}\subset\ol M\phantom{}^\text{gp}/\ol M\phantom{}^\text{gp}_\text{tor}$. Define the monoid map $M_{\hat\eta}\to\clO_S=\bC$ to be the identity on $\bC^\times\subset M_{\hat\eta}$ and zero on its complement; by construction, this is an fs log structure on $S$. Note that we have $M_{\hat\eta}^\text{gp} = M^\text{gp}/\hat\eta(Q^\text{gp}_\text{tor})$ and $\ol M_{\hat\eta} = (\ol M)^\text{sat}/\ol M\phantom{}^\text{gp}_\text{tor}$. The natural map $M\to M^\text{gp}/\hat\eta(Q^\text{gp}_\text{tor})$ factors through $M_{\hat\eta}$ to yield a map of log structures $\varphi_{\hat\eta}:M\to M_{\hat\eta}$. Using \eqref{eqn:ghost-groupification-cartesian}, \eqref{eqn:map-of-log-str-cartesian} and \eqref{eqn:torsion-embedding-induced-by-fs}, one verifies that the pair $(M_{\hat\eta},\varphi_{\hat\eta})$ is an initial object of $(M/\Log^\text{fs}(S))_{\hat\eta}$.
            
            Any other object $(M',\varphi')$ of $(M/\Log^\text{fs}(S))_{\hat\eta}$ is an initial object if and only if the unique morphism $(M_{\hat\eta},\varphi_{\hat\eta})\to (M',\varphi')$ in $(M/\Log^\text{fs}(S))_{\hat\eta}$ is an isomorphism. By \eqref{eqn:map-of-log-str-cartesian}, this happens if and only if the induced map $\ol M_{\hat\eta}\to\ol M\phantom{}'$ is an isomorphism. Composing with the isomorphism $\ol{Q^\text{sat}}\xrightarrow{\sim}\ol M_{\hat\eta}$ induced by $\ol\varphi_{\hat\eta}\circ\psi$, we conclude that $(M',\varphi')$ is an initial object if and only if the map $\ol{Q^\text{sat}}\to\ol M\phantom{}'$ induced by $\ol{\varphi}\phantom{}'\circ\psi:Q\to\ol M\phantom{}'$ is an isomorphism.
        \end{proof}
        
        \begin{claim}\label{claim:fs-cat-op-equiv}
            For each saturation datum $\eta\in\fS(C,{\bf x},f,\bmu)$, there exists a unique element $\hat\eta\in\widehat\fS(M)$ such that $\clR^\text{fs}$ restricts to an equivalence of categories
            \begin{align*}
                \clR^\text{fs}_\eta:(\sC^\text{fs}_\eta)^\text{op}\xrightarrow{\sim}(M/\Log^\text{fs}(S))_{\hat\eta}.
            \end{align*}
            Moreover, the resulting 
            assignment $\eta\mapsto\hat\eta$
            defines a bijection $\fS(C,{\bf x},f,\bmu)\xrightarrow{\sim}\widehat\fS(M)$.
        \end{claim}
        \begin{proof}[Proof of claim]
            We describe the bijection $\eta\mapsto\hat\eta$ as follows.
            \begin{enumerate}[label = (\roman*), itemsep = 0.3ex]
                \item 
                Given any element $\tilde\eta$ of $\widetilde\fS(C,{\bf x},f,\bmu)$, define a monoid map $\psi_{\tilde\eta}:Q\to M$ lifting $\psi:Q\xrightarrow{\sim}\ol M$ as follows. Using the notation from \eqref{eqn:sat-data-before-quotient}, write $\tilde\eta = \{\sigma_{v,j},\theta_e\}$ and define $\tilde\eta_{v,j}\in\clO_{S^\dagger}(-\delta_{v,j})^\times\subset M$ and $\tilde\eta_e\in\clO_{S^\dagger}(-\rho_e)^\times\subset M$ by the equations
                \begin{align}\label{eqn:splitting-from-lifted-sat-datum}
                    \langle\tilde\eta_{v,j},\Theta_{v,j}(\sigma_{v,j})\rangle = 1 \quad\text{and}\quad \langle\tilde\eta_e,\Theta_e(\theta_e)\rangle = 1,
                \end{align}
                where $\langle\cdot,\cdot\rangle$ is the duality pairing induced by the isomorphisms appearing in the global slb on $S^\dagger$; see Discussion~\ref{disc-part:log-lb-vcd-dual}. By \ref{concrete-alg-log-map-vcd-iso} and \ref{concrete-alg-log-map-matching-condition}, the assignment 
                \begin{align*}
                    m_{v,j}\mapsto\tilde\eta_{v,j},\quad m_e\mapsto\tilde\eta_e
                \end{align*}
                defines a homomorphism $\psi_{\tilde\eta}:Q\to M$ lifting $\psi$. One checks that, given any group element $\lambda\in\Hom_{\Ab}(Q^\text{gp},\bC^\times) = T_{(C,{\bf x},f,\bmu)}$, we have the equality
                \begin{align}\label{eqn:anti-equivariance}
                    \psi_{\tilde\eta\cdot\lambda} = \psi_{\tilde\eta}\cdot \lambda^{-1}.
                \end{align}

                \item 
                Given $\eta\in\fS(C,{\bf x},f,\bmu)$, lift it to an element $\tilde\eta\in\widetilde\fS(C,{\bf x},f,\bmu)$. Define $\hat\eta\in\widehat\fS(M)$ to be the restriction of $\psi^\text{gp}_{\tilde\eta}:Q^\text{gp}\to M^\text{gp}$ to the subgroup $Q^\text{gp}_\text{tor}$. 
                
                \item The fact that $\fS(C,{\bf x},f,\bmu)$ and $\widehat\fS(M)$ are $\pi_0(T_{(C,{\bf x},f,\bmu)})$-torsors, together with the equivariance property \eqref{eqn:anti-equivariance}, shows that the assignment $\eta\mapsto\hat\eta$ gives rise to a well-defined bijection $\fS(C,{\bf x},f,\bmu)\xrightarrow{\sim}\widehat\fS(M)$.
            \end{enumerate}

            Fix $\eta\in\fS(C,{\bf x},f,\bmu)$ and let $\hat\eta\in\widehat\fS(M)$ be its image. Given any object $(\xi'^\dagger/S'^\dagger,f'^\dagger)$ of $\sC^\text{fs}_\eta$, we argue that its image $(M',\varphi')$ under $\clR^\text{fs}$ lies in $(M/\Log^\text{fs}(S))_{\hat\eta}$ as follows.
            \begin{enumerate}[label = (\roman*), itemsep = 0.3ex]
                \item Following Construction~\ref{cons:underlying-sat-symp-log-map}, choose a section $\beta':\ol M\phantom{}'\to M'$ of the projection $M'\to\ol M\phantom{}'$. The section $\beta'$ determines an element $\tilde\eta\in\widetilde\fS(C,{\bf x},f,\bmu)$ lifting $\eta$, as in Construction~\ref{cons-subpart:sat-datum-from-split-fs-log-map}. We thus get a monoid map $\psi_{\tilde\eta}:Q\to M$ lifting $\psi$.
                
                \item On the other hand, by \eqref{eqn:map-of-log-str-cartesian}, the section $\beta'$ pulls back to a section $\beta:\ol M\to M$ of the projection $M = \ol M\times_{\ol M\phantom{}'}M'\to\ol M$. By comparing \eqref{eqn:lifted-sat-datum-from-splitting} and \eqref{eqn:splitting-from-lifted-sat-datum}, we get the equality $\psi_{\tilde\eta} = \beta\circ\psi:Q\to M$. The argument up to this point is summarized by the following commutative diagram.
                \begin{equation}\label{eqn:fs-splitting-pullback}
                \begin{tikzcd}
                    & & M \arrow[r,"{\varphi'}"] \arrow[d] \arrow[dr,phantom,"\lrcorner", pos = 0.001] & M' \arrow[d] \\
                    Q \arrow[rr,"\psi"] \arrow[urr,"\psi_{\tilde\eta}"] & & \ol M \arrow[r,"{\ol\varphi\phantom{}'}"] \arrow[u, bend left,"\beta", pos = 0.4] & \ol M\phantom{}' \arrow[u, bend left,"{\beta'}", pos = 0.4]
                \end{tikzcd}
                \end{equation}
                
                \item The map $\ol\varphi\phantom{}'^{\text{gp}}\circ\psi^\text{gp}$ vanishes on $Q^\text{gp}_\text{tor}$ by Remark~\ref{rem:sharp-fs-implies-no-torsion}. Using \eqref{eqn:fs-splitting-pullback}, this means that $\varphi'^{\text{gp}}\circ\psi_{\tilde\eta}^\text{gp} = \beta'^{\text{gp}}\circ\ol\varphi\phantom{}'^{\text{gp}}\circ\psi^\text{gp}$ restricts to the trivial homomorphism on $Q^\text{gp}_\text{tor}$, i.e., we have $\hat\eta(M',\varphi') = \hat\eta$ as desired.
            \end{enumerate}

            Thus, $\clR^\text{fs}$ restricts to a functor $\clR^\text{fs}_\eta:(\sC^\text{fs}_\eta)^\text{op}\to(M/\Log^\text{fs}(S))_{\hat\eta}$. Since $\clR^\text{fs}$ is fully faithful, it follows that $\clR^\text{fs}_\eta$ is also fully faithful. Since $\eta\mapsto\hat\eta$ is a bijection and $\clR^\text{fs}$ is essentially surjective, it follows that $\clR^\text{fs}_\eta$ is also essentially surjective. Thus, $\clR^\text{fs}_\eta$ is an equivalence.
        \end{proof}
    \end{step}

    Taken together, Claims~\ref{claim:fs-initial-object} and \ref{claim:fs-cat-op-equiv} complete the proof.
\end{proof}

\begin{Corollary}\label{cor:fs-log-map-comparison}
    For any saturated symplectic log map $(C,{\bf x},f,\bmu,\eta)$ to $(X,D)$, there is a unique (up to unique isomorphism) fs-basic algebraic log map $(\xi^\dagger/S^\dagger,f^\dagger)$ to $X^\dagger$ whose underlying saturated symplectic log map is given by $(C,{\bf x},f,\bmu,\eta)$.
\end{Corollary}

\begin{Corollary}[{\cite[Conjecture 1.1]{Tehrani-log}}]\label{cor:number-of-saturations}
    Each fs-basic algebraic log map to $X^\dagger$ has a well-defined underlying symplectic log map to $(X,D)$. Moreover, for each symplectic log map $(C,{\bf x},f,\bmu)$ to $(X,D)$, there are precisely $|\ol M(C,{\bf x},f,\bmu)^\text{gp}_\text{tor}|$ distinct fs-basic algebraic log maps to $X^\dagger$ (up to isomorphism) with underlying symplectic log map $(C,{\bf x},f,\bmu)$.
\end{Corollary}
\begin{proof}
    This is immediate from Lemma~\ref{lem:number-of-saturations} and Corollary~\ref{cor:fs-log-map-comparison}; see also Remark~\ref{rem:number-of-saturations}.
\end{proof}

\begin{Remark}\label{rem:number-of-saturations}
    In the setting of Corollary~\ref{cor:number-of-saturations}, \cite[Equation~(5.15)]{Tehrani-log} gives a seemingly different (conjectural) expression for the number of distinct fs-basic algebraic log maps (up to isomorphism) lying over the given symplectic log map $(C,{\bf x},f,\bmu)$. We will recall this expression and explain why it is in fact equal to $|\ol M(C,{\bf x},f,\bmu)^\text{gp}_\text{tor}|$.
    \begin{enumerate}[label = (\arabic*), itemsep = 0.3ex]
        \item Denote the set of vertices (resp. edges) of $\plC$ by $\bV$ (resp. $\bE$). Choose an \emph{orientation} on $\plC$, i.e., a collection $\esO$ of oriented edges of $\plC$ such that, for each $e\in\bE$, exactly one of the two oriented edges lying over $e$ belongs to the collection $\esO$.
        \vspace{0.3ex}
        
        \myitem[] For each $e\in\bE$ and $1\le j\le r$, define the integer $\mu_{e,j}^{\esO}:=\mu_{\vec{e},j}$, where ${\vec{e}}$ is the unique lift of $e$ in $\esO$. 
        For each $v\in\bV$, define the set $\esO^-_v$ (resp. $\esO^+_v$) to consist of all $e\in\bE$ whose unique lift to $\esO$ is of the form $\vec{e}:v\to v'$ (resp. $\vec{e}:v'\to v$) for some $v'\in\bV$.
        \vspace{0.3ex}

        \myitem[] For each $v\in\bV$, we have previously defined the subset $I_v\subset\{1\,\ldots,r\}$ to consist of all $j$ for which $f(C_v)\subset D_j$. Likewise, for each $e\in\bE$, define the subset $I_e\subset\{1,\ldots,r\}$ to consist of all $j$ for which $f(y_e)\in D_j$. We conclude from \cite[Equation~(2.21)]{Tehrani-log} that, for each oriented edge $\vec{e}:v\to v'$ of $\plC$, we have 
        \begin{align}\label{eqn:trivial-mu-values}
            I_e = I_v\cup I_{v'}\quad\text{and}\quad \mu_{\vec{e},j} = 0\text{ for all }j\not\in I_e.
        \end{align}
        With these conventions, \cite[Equation~(2.26)]{Tehrani-log} defines a homomorphism
        \begin{align*}
            \varrho:\textstyle\bZ^{\bE}\oplus\bigoplus_{v\in\bV}\bZ^{I_v}\to\bigoplus_{e\in\bE}\bZ^{I_e}
        \end{align*}
        by the formulas
        \begin{align}\label{eqn:varrho-map}
            \varrho({\bf 1}_e) := \sum_{j\in I_e}\mu_{e,j}^\esO{\bf 1}_{e,j} \quad\text{and}\quad \varrho({\bf 1}_{v,j}) := \sum_{e\in\esO^-_v}{\bf 1}_{e,j} - \sum_{e\in\esO^+_v}{\bf 1}_{e,j},
        \end{align}
        where ${\bf 1}_e$, ${\bf 1}_{v,j}$ and ${\bf 1}_{e,j}$ denote the standard basis vectors of $\bZ^\bE$, $\bZ^{I_v}$ and $\bZ^{I_e}$ respectively. In terms of $\varrho$, the expression in \cite[Equation~(5.15)]{Tehrani-log} equals $|(\ker\varrho)^\perp/\img(\varrho^\vee)|$.

        \item Using \eqref{eqn:trivial-mu-values}, \eqref{eqn:varrho-map} and Remark~\ref{rem:fine-basic-monoid}, we get a natural identification
        \begin{align}\label{eqn:alt-presentation-fine-basic-monoid}
            \textstyle(\bZ^{\bE}\oplus\bigoplus_{v\in\bV}\bZ^{I_v})/\img(\varrho^\vee) = \ol M(C,{\bf x},f,\bmu)^\text{gp}
        \end{align}
        of abelian groups. By elementary linear algebra, the inclusion $\img(\varrho^\vee)\subset(\ker\varrho)^\perp$ becomes an equality after applying the functor $-\otimes_\bZ\bQ$. Therefore, we obtain an inclusion $(\ker\varrho)^\perp/\img(\varrho^\vee)\subset \ol M(C,{\bf x},f,\bmu)^\text{gp}_\text{tor}$ from the identification \eqref{eqn:alt-presentation-fine-basic-monoid}. 
        \vspace{0.3ex}
        
        \myitem[] This inclusion is in fact an equality, since the quotient $(\bZ^{\bE}\oplus\bigoplus_{v\in\bV}\bZ^{I_v})/(\ker\varrho)^\perp$ is torsion-free. Indeed, if ${\bf x}\in \bZ^{\bE}\oplus\bigoplus_{v\in\bV}\bZ^{I_v}$ is an element such that some nonzero integer multiple of ${\bf x}$ lies in $(\ker\varrho)^\perp$, then ${\bf x}$ itself lies in $(\ker\varrho)^\perp$.
    \end{enumerate}
    
    In conclusion, we have an identification $(\ker\varrho)^\perp/\img(\varrho^\vee) = \ol M(C,{\bf x},f,\bmu)^\text{gp}_\text{tor}$. Thus, the expression in \cite[Equation~(5.15)]{Tehrani-log} agrees with $|\ol M(C,{\bf x},f,\bmu)^\text{gp}_\text{tor}|$.
\end{Remark}

\subsection{Saturation in log Gromov convergence}\label{subsec:sat-log-gromov-conv}

In this subsection, we recall the notion of \emph{log Gromov convergence} for sequences of symplectic log maps (\cite[Definition~3.7]{Tehrani-log}) and then explain how to extend it to the setting of saturated symplectic log maps.

For simplicity, we work in the setting of a smooth projective snc pair $(X,D)$, rather than the more general setting of symplectic manifolds and snc symplectic divisors considered in \cite{Tehrani-log}. Unlike in the rest of this article, we use the ``usual" (i.e., complex analytic) topology on varieties and schemes in this subsection. To illustrate the main idea, we will consider only the following situation, which is analogous to the setup in Discussion~\ref{disc:log-map-motivation}.

\begin{Situation}\label{situ:log-convergence}
    Consider a sequence $\{(C^{(k)},{\bf x}^{(k)},f^{(k)})\}_{k\ge 1}$ of stable maps to $X$ of class $(g,n,\beta)$ together with a collection of numbers $\lambda_{i,j}\in\bN$, where $i$ and $j$ range over $\{1,\ldots,n\}$ and $\{1,\ldots,r\}$ respectively. For each $k\ge 1$, assume that the curve $C^{(k)}$ is smooth and that, for each $1\le j\le r$, we have the following equality of divisors on $C^{(k)}$:
    \begin{align*}
        (f^{(k)})^{-1}(D_j) = \sum_{1\le i\le n}\lambda_{i,j}\cdot x_i^{(k)}.
    \end{align*}
    Writing $\blambda := \{\lambda_{i,j}\}$, we get a sequence of symplectic stable log maps $\{(C^{(k)},{\bf x}^{(k)},f^{(k)},\blambda)\}_{k\ge 1}$. In addition, also consider a symplectic stable log map $(C,{\bf x},f,\bmu)$ to $(X,D)$ such that the stable map $(C,{\bf x},f)$ is of class $(g,n,\beta)$ and we have $\mu_{i,j} = \lambda_{i,j}$ for each $i$ and $j$.
\end{Situation}

\begin{Definition}
    In Situation~\ref{situ:log-convergence}, denote the complex structure on $\tilde C$ by $\fj$. 
    \begin{enumerate}[ref = \theDefinition(\arabic*), label = (\arabic*), itemsep = 0.3ex]
        \item\label{def-part:gromov-conv} We say that $\{(C^{(k)},{\bf x}^{(k)},f^{(k)})\}_{k\ge 1}$ \emph{Gromov converges} to $(C,{\bf x},f)$ if there exist
        \begin{enumerate}[label = (\roman*), itemsep = 0.15ex]
            \item a sequence $\{\fj^{(k)}\}_{k\ge 1}$ of complex structures on (the underlying smooth manifold of) $\tilde C$ converging uniformly with all derivatives to the complex structure $\fj$,

            \item a sequence $\{z_{\vec{e}}^{(k)}\}_{k\ge 1}$ of $\fj^{(k)}$-holomorphic local coordinates centred at $y_{\vec{e}}\in\tilde C$ converging uniformly with all derivatives to a $\fj$-holomorphic local coordinate $z_{\vec{e}}$ centred at $y_{\vec{e}}\in\tilde C$, for each oriented edge ${\vec{e}}$ of $\plC$, and

            \item a sequence $\{\varepsilon^{(k)}_e\}_{k\ge 1}$ of elements of $\bC^\times$ converging to $0$, for each edge $e$ of $\plC$,
        \end{enumerate}
        \vspace{0.3ex}
        with the following properties.
        \vspace{0.3ex}
        \begin{enumerate}[label = (\alph*), itemsep = 0.15ex]
            \item For sufficiently large $k$, there exists a biholomorphic map $\Phi^{(k)}$ from $(C^{(k)},{\bf x}^{(k)})$ to the pointed Riemann surface obtained by smoothing $(C,{\bf x})$, equipped with the complex structure $\fj^{(k)}$, using the equations $z_{\vec{e}}^{(k)}z_{\cev{e}}^{(k)} = \varepsilon^{(k)}_e$ for each $e$, and

            \item The sequence of biholomorphic maps $\{\Phi^{(k)}\}_{k\gg 1}$ can be chosen such that the sequence of maps $\{f^{(k)}\circ(\Phi^{(k)})^{-1}\}_{k\gg 1}$ converges uniformly with all derivatives to the map $f$ on all compact subsets of $C_\text{gen}$.
        \end{enumerate}

        \item\label{def-part:log-gromov-conv} We say that $\{(C^{(k)},{\bf x}^{(k)},f^{(k)},\blambda)\}_{k\ge 1}$ \emph{log Gromov converges} to $(C,{\bf x},f,\bmu)$ if we have the following two properties.
        \begin{enumerate}[label = (\roman*), itemsep = 0.15ex]
            \item The sequence $\{(C^{(k)},{\bf x}^{(k)},f^{(k)})\}_{k\ge 1}$ Gromov converges to $(C,{\bf x},f)$.

            \item For each vertex $v$ of $\plC$ and $1\le j\le r$, there exist a sequence $\{t^{(k)}_{v,j}\}_{k\ge 1}$ of elements of $\bC^\times$ and a rational section $\zeta_{v,j}$ of $\tilde f_v^*(\clO_X(D_j))$ such that 
            \begin{enumerate}[label = (\alph*), itemsep = 0.15ex]
                \item the divisor of zeros and poles of the rational section $\zeta_{v,j}$ is given by the divisor $D_{v,j}$ from Definition~\ref{def-part:pullback-twisted-by-contact-divisor},
                
                \item we have $t^{(k)}_{v,j} = 1$ for all $k\ge 1$ if $j\not\in I_v$, and

                \item the sequence of sections $\{(t^{(k)}_{v,j})^{-1}\cdot (f^{(k)}\circ(\Phi^{(k)})^{-1})^*\taut_{D_j}\}_{k\gg 1}$ of (pullbacks of) $\clO_X(D_j)$ converges uniformly with all derivatives to $\zeta_{v,j}$ on all compact subsets of $C_v\cap C_\text{gen}$, for some choices of $\{{\fj}^{(k)}\}_{k\ge 1}$, $\{z_{\vec{e}}^{(k)}\}_{k\ge 1}$, $\{\varepsilon_e^{(k)}\}_{k\ge 1}$ and $\{\Phi^{(k)}\}_{k\gg 1}$ as in the definition of Gromov convergence.
            \end{enumerate}
        \end{enumerate}
        \vspace{0.3ex}
        In fact, by \cite[Proposition~3.20]{Tehrani-log}, the choices $\{z_{\vec{e}}^{(k)}\}_{k\ge 1}$, $\{\varepsilon_e^{(k)}\}_{k\ge 1}$ and $\{t^{(k)}_{v,j}\}_{k\ge 1}$ can be assumed to have the following two additional properties.
        \vspace{0.3ex}
        \begin{enumerate}[resume, ref = \theenumi(\roman*), label = (\roman*), itemsep = 0.15ex]
            \item\label{def-subpart:rescaling-lies-in-group} For each $k\ge 1$, the collection $\{t^{(k)}_{v,j},\varepsilon_e^{(k)}\}$ defines an element of the group $T_{(C,{\bf x},f,\bmu)}$ from Construction~\ref{cons-part:sat-data-groups}, under the identification \eqref{eqn:sat-data-groups}.

            \item\label{def-subpart:rescaled-limit-smoothing-parameter-compatible} The formulas \eqref{eqn:rescaled-limits-for-analytic-matching-condition} and \eqref{eqn:smoothing-parameter-for-analytic-matching-condition} together define an element $\{\sigma_{v,j},\theta_e\}$ of the set $\widetilde\fS(C,{\bf x},f,\bmu)$ from Definition~\ref{cons-part:sat-data-before-quotient}.
        \end{enumerate}
    \end{enumerate}
\end{Definition}

\begin{Remark}
    Definition~\ref{def-part:gromov-conv} is standard in symplectic topology. Definition~\ref{def-part:log-gromov-conv} is a reformulation of \cite[Definition~3.7 and Proposition~3.20]{Tehrani-log}; see also Remark~\ref{rem:symp-log-map}.
\end{Remark}

We propose the following definition of \emph{saturated log Gromov convergence}, as a refinement of the notion of log Gromov convergence described in Definition~\ref{def-part:log-gromov-conv}.

\begin{Definition}\label{def:sat-log-gromov-conv}
    In Situation~\ref{situ:log-convergence}, observe that $(C^{(k)},{\bf x}^{(k)},f^{(k)},\blambda)$ has a unique saturation datum $\eta_k$ for each $k\ge 1$. Fix a saturation datum $\eta\in\fS(C,{\bf x},f,\bmu)$ for $(C,{\bf x},f,\bmu)$.

    We say that $\{(C^{(k)},{\bf x}^{(k)},f^{(k)},\blambda,\eta_k)\}_{k\ge 1}$ \emph{saturated log Gromov converges} to $(C,{\bf x},f,\bmu,\eta)$ if we have the following properties.
    \begin{enumerate}[label = (\arabic*), itemsep = 0.3ex]
        \item The sequence $\{(C^{(k)},{\bf x}^{(k)},f^{(k)},\blambda)\}_{k\ge 1}$ log Gromov converges to $(C,{\bf x},f,\bmu)$.

        \item There exist choices of $\{z_{\vec{e}}^{(k)}\}_{k\ge 1}$, $\{\varepsilon_e^{(k)}\}_{k\ge 1}$ and $\{t^{(k)}_{v,j}\}_{k\ge 1}$ with the four properties in Definition~\ref{def-part:log-gromov-conv} and the following two additional properties.
        \vspace{0.3ex}
        \begin{enumerate}[label = (\roman*$'$), itemsep = 0.15ex]
            \setcounter{enumii}{2}
            \item For each $k\ge 1$, the group element $\{t^{(k)}_{v,j},\varepsilon_e^{(k)}\}$ from Definition~\ref{def-subpart:rescaling-lies-in-group} lies in the subgroup $T^0_{(C,{\bf x},f,\bmu)}\subset T_{(C,{\bf x},f,\bmu)}$ from Construction~\ref{cons-part:sat-data-groups}.

            \item The element $\{\sigma_{v,j},\theta_e\}$ of $\widetilde\fS(C,{\bf x},f,\bmu)$ from Definition~\ref{def-subpart:rescaled-limit-smoothing-parameter-compatible} is a lift of the given saturation datum $\eta\in\fS(C,{\bf x},f,\bmu)$.
        \end{enumerate}
    \end{enumerate}
\end{Definition}

\begin{Remark}\label{rem:forgetful-map-from-sat-to-unsat-is-separated}
    In the setting of Definition~\ref{def:sat-log-gromov-conv}, if $\eta,\eta'\in\fS(C,{\bf x},f,\bmu)$ are saturation data such that the sequence $\{(C^{(k)},{\bf x}^{(k)},f^{(k)},\blambda,\eta_k)\}_{k\ge 1}$ saturated log Gromov converges to both $(C,{\bf x},f,\bmu,\eta)$ and $(C,{\bf x},f,\bmu,\eta')$, then it is not hard to show that we have $\eta = \eta'$.
\end{Remark}

We conclude with the following observation on saturated log Gromov convergence.

\begin{Proposition}\label{prop:forgetful-map-from-sat-to-unsat-is-proper}
    In Situation~\ref{situ:log-convergence}, assume that the sequence $\{(C^{(k)},{\bf x}^{(k)},f^{(k)},\blambda)\}_{k\ge 1}$ log Gromov converges to $(C,{\bf x},f,\bmu)$. Then the sequence $\{(C^{(k)},{\bf x}^{(k)},f^{(k)},\blambda,\eta_k)\}_{k\ge 1}$ has a saturated log Gromov convergent subsequence, where $\eta_k$ is as in Definition~\ref{def:sat-log-gromov-conv}.
\end{Proposition}
\begin{proof}
    Fix choices of data $\{{\fj}^{(k)}\}_{k\ge 1}$, $\{z_{\vec{e}}^{(k)}\}_{k\ge 1}$, $\{\varepsilon_e^{(k)}\}_{k\ge 1}$, $\{t^{(k)}_{v,j}\}_{k\ge 1}$ and $\{\Phi^{(k)}\}_{k\gg 1}$ satisfying the four properties in Definition~\ref{def-part:log-gromov-conv}. Define $z_{\vec{e}}$ and $\zeta_{v,j}$ be the resulting limits of $z_{\vec{e}}^{(k)}$ and $(t^{(k)}_{v,j})^{-1}\cdot (f^{(k)}\circ(\Phi^{(k)})^{-1})^*\taut_{D_j}$ on $C_v\cap C_\text{gen}$ respectively as $k\to\infty$. Also define $\sigma_{v,j}$ and $\theta_e$ via the formulas \eqref{eqn:rescaled-limits-for-analytic-matching-condition} and \eqref{eqn:smoothing-parameter-for-analytic-matching-condition} respectively.
    
    After passing to a subsequence, we may assume that the image of the group elements $\gamma_k:=\{t^{(k)}_{v,j},\varepsilon_e^{(k)}\}$ in the finite abelian group $\pi_0(T_{(C,{\bf x},f,\bmu)})$ is a constant independent of $k\ge 1$; see Construction~\ref{cons-part:sat-data-groups}. In other words, the group element $\gamma_k\cdot\gamma_1^{-1}$ belongs to $T^0_{(C,{\bf x},f,\bmu)}$ for all $k\ge 1$. Define the saturation datum $\eta\in\fS(C,{\bf x},f,\bmu)$ to be the image of the element $\{\sigma_{v,j},\theta_e\}\cdot\gamma_1^{-1}$ of the set $\widehat\fS(C,{\bf x},f,\bmu)$; see Construction~\ref{cons-part:sat-data-group-action}.
    
    For each edge $e$ of $\plC$, choose elements $\varepsilon_{\vec{e}}^{(1)},\varepsilon_{\cev{e}}^{(1)}\in\bC^\times$ whose product is $\varepsilon_e^{(1)}$. Replace the original choices $\{z_{\vec{e}}^{(k)}\}_{k\ge 1}$, $\{\varepsilon_e^{(k)}\}_{k\ge 1}$ and $\{t^{(k)}_{v,j}\}_{k\ge 1}$ by the new choices
    \begin{align*}
        \left\{\frac{z_{\vec{e}}^{(k)}}{\varepsilon_{\vec{e}}^{(1)}}\right\}_{k\ge 1}, \left\{\frac{\varepsilon_e^{(k)}}{\varepsilon_{e}^{(1)}}\right\}_{k\ge 1}\quad\text{and}\quad \left\{\frac{t_{v,j}^{(k)}}{t_{v,j}^{(1)}}\right\}_{k\ge 1}\quad\text{respectively.}
    \end{align*}
    With these new choices, one directly verifies the conditions in Definition~\ref{def:sat-log-gromov-conv} to conclude that $\{(C^{(k)},{\bf x}^{(k)},f^{(k)},\blambda,\eta_k)\}_{k\ge 1}$ saturated log Gromov converges to $(C,{\bf x},f,\bmu,\eta)$.
\end{proof}

\begin{Remark}
    Definition~\ref{def:sat-log-gromov-conv}, Remark~\ref{rem:forgetful-map-from-sat-to-unsat-is-separated} and Proposition~\ref{prop:forgetful-map-from-sat-to-unsat-is-proper} extend to the more general setting of symplectic manifolds and snc symplectic divisors considered in \cite{Tehrani-log} in a straightforward fashion.
\end{Remark}
    \bibliography{references}

@article{AMS21,
  title={Complex cobordism, {H}amiltonian loops and global {K}uranishi charts},
  author={Abouzaid, M. and McLean, M. and Smith, I.},
  journal={arXiv:2110.14320},
  year={2021},
}

@article {AMS24,
    AUTHOR = {Abouzaid, M. and McLean, M. and Smith, I.},
     TITLE = {Gromov-{W}itten invariants in complex and {M}orava-local
              {$K$}-theories},
   JOURNAL = {Geom. Funct. Anal.},
  FJOURNAL = {Geometric and Functional Analysis},
    VOLUME = {34},
      YEAR = {2024},
    NUMBER = {6},
     PAGES = {1647--1733},
      ISSN = {1016-443X,1420-8970},
   MRCLASS = {53D45 (14N35 55N15)},
  MRNUMBER = {4823210},
       DOI = {10.1007/s00039-024-00697-4},
       URL = {https://doi.org/10.1007/s00039-024-00697-4},
}

@article {AC-logDF2,
    AUTHOR = {Abramovich, D. and Chen, Q.},
     TITLE = {Stable logarithmic maps to {D}eligne-{F}altings pairs {II}},
   JOURNAL = {Asian J. Math.},
  FJOURNAL = {Asian Journal of Mathematics},
    VOLUME = {18},
      YEAR = {2014},
    NUMBER = {3},
     PAGES = {465--488},
      ISSN = {1093-6106,1945-0036},
   MRCLASS = {14D23 (14A20 14H10 14N35)},
  MRNUMBER = {3257836},
MRREVIEWER = {Jonathan\ Wise},
       DOI = {10.4310/AJM.2014.v18.n3.a5},
       URL = {https://doi.org/10.4310/AJM.2014.v18.n3.a5},
}

@incollection {log-geometry-and-moduli,
    AUTHOR = {Abramovich, D. and Chen, Q. and Gillam, D. and Huang, Y. and Olsson, M. and Satriano, M. and Sun, S.},
     TITLE = {Logarithmic geometry and moduli},
 BOOKTITLE = {Handbook of moduli. {V}ol. {I}},
    SERIES = {Adv. Lect. Math. (ALM)},
    VOLUME = {24},
     PAGES = {1--61},
 PUBLISHER = {Int. Press, Somerville, MA},
      YEAR = {2013},
      ISBN = {978-1-57146-257-2},
   MRCLASS = {14D20 (14A20)},
  MRNUMBER = {3184161},
MRREVIEWER = {Howard\ M.\ Thompson},
}

@article {ACGS-decomp,
    AUTHOR = {Abramovich, D. and Chen, Q. and Gross, M. and Siebert,
              B.},
     TITLE = {Decomposition of degenerate {G}romov-{W}itten invariants},
   JOURNAL = {Compos. Math.},
  FJOURNAL = {Compositio Mathematica},
    VOLUME = {156},
      YEAR = {2020},
    NUMBER = {10},
     PAGES = {2020--2075},
      ISSN = {0010-437X,1570-5846},
   MRCLASS = {14N35 (14D23 14T20)},
  MRNUMBER = {4177284},
MRREVIEWER = {William\ Liu},
       DOI = {10.1112/s0010437x20007393},
       URL = {https://doi.org/10.1112/s0010437x20007393},
}

@article {borne-vistoli-DF,
    AUTHOR = {Borne, N. and Vistoli, A.},
     TITLE = {Parabolic sheaves on logarithmic schemes},
   JOURNAL = {Adv. Math.},
  FJOURNAL = {Advances in Mathematics},
    VOLUME = {231},
      YEAR = {2012},
    NUMBER = {3-4},
     PAGES = {1327--1363},
      ISSN = {0001-8708,1090-2082},
   MRCLASS = {14F05 (14A20 18D10)},
  MRNUMBER = {2964607},
MRREVIEWER = {Jon\ Eivind\ Vatne},
       DOI = {10.1016/j.aim.2012.06.015},
       URL = {https://doi.org/10.1016/j.aim.2012.06.015},
}

@article {moduli-trop-curves,
    AUTHOR = {Cavalieri, R. and Chan, M. and Ulirsch, M. and
              Wise, J.},
     TITLE = {A moduli stack of tropical curves},
   JOURNAL = {Forum Math. Sigma},
  FJOURNAL = {Forum of Mathematics. Sigma},
    VOLUME = {8},
      YEAR = {2020},
     PAGES = {Paper No. e23, 93},
      ISSN = {2050-5094},
   MRCLASS = {14T20 (14D23)},
  MRNUMBER = {4091085},
MRREVIEWER = {Hannah\ Markwig},
       DOI = {10.1017/fms.2020.16},
       URL = {https://doi.org/10.1017/fms.2020.16},
}

@article{Chen-Church-Zhao,
  title={Curves on complete intersections and measures of irrationality},
  author={Chen, N. and Church, B. and Zhao, J.},
  journal={arXiv:2406.12101},
  year={2024},
}

@article {Chen-logDF1,
    AUTHOR = {Chen, Q.},
     TITLE = {Stable logarithmic maps to {D}eligne-{F}altings pairs {I}},
   JOURNAL = {Ann. of Math. (2)},
  FJOURNAL = {Annals of Mathematics. Second Series},
    VOLUME = {180},
      YEAR = {2014},
    NUMBER = {2},
     PAGES = {455--521},
      ISSN = {0003-486X,1939-8980},
   MRCLASS = {14N35 (14D23)},
  MRNUMBER = {3224717},
MRREVIEWER = {Sergiy\ Koshkin},
       DOI = {10.4007/annals.2014.180.2.2},
       URL = {https://doi.org/10.4007/annals.2014.180.2.2},
}

@article {Tehrani-log,
    AUTHOR = {{Farajzadeh-Tehrani}, M.},
     TITLE = {Pseudoholomorphic curves relative to a normal crossings
              symplectic divisor: compactification},
   JOURNAL = {Geom. Topol.},
  FJOURNAL = {Geometry \& Topology},
    VOLUME = {26},
      YEAR = {2022},
    NUMBER = {3},
     PAGES = {989--1075},
      ISSN = {1465-3060,1364-0380},
   MRCLASS = {14N35 (53D45)},
  MRNUMBER = {4466645},
MRREVIEWER = {William\ Liu},
       DOI = {10.2140/gt.2022.26.989},
       URL = {https://doi.org/10.2140/gt.2022.26.989},
}

@article {Tehrani-semistable,
    AUTHOR = {{Farajzadeh-Tehrani}, M.},
     TITLE = {Limits of stable maps in a semi-stable degeneration},
   JOURNAL = {Geom. Dedicata},
  FJOURNAL = {Geometriae Dedicata},
    VOLUME = {216},
      YEAR = {2022},
    NUMBER = {6},
     PAGES = {Paper No. 66, 42},
      ISSN = {0046-5755,1572-9168},
   MRCLASS = {53D45 (14D06 14N35)},
  MRNUMBER = {4482071},
MRREVIEWER = {Johan\ Asplund},
       DOI = {10.1007/s10711-022-00731-5},
       URL = {https://doi.org/10.1007/s10711-022-00731-5},
}

@article{FRTU-chip-firing,
    AUTHOR = {Foster, T. and Ranganathan, D. and Talpo, M. and
              Ulirsch, M.},
     TITLE = {Logarithmic {P}icard groups, chip firing, and the
              combinatorial rank},
   JOURNAL = {Math. Z.},
  FJOURNAL = {Mathematische Zeitschrift},
    VOLUME = {291},
      YEAR = {2019},
    NUMBER = {1-2},
     PAGES = {313--327},
      ISSN = {0025-5874,1432-1823},
   MRCLASS = {14T05 (14H10)},
  MRNUMBER = {3936072},
MRREVIEWER = {Andreas\ Gross},
       DOI = {10.1007/s00209-018-2085-2},
       URL = {https://doi.org/10.1007/s00209-018-2085-2},
}

@article{gillam-minimal,
    AUTHOR = {Gillam, W. D.},
     TITLE = {Logarithmic stacks and minimality},
   JOURNAL = {Internat. J. Math.},
  FJOURNAL = {International Journal of Mathematics},
    VOLUME = {23},
      YEAR = {2012},
    NUMBER = {7},
     PAGES = {1250069, 38},
      ISSN = {0129-167X,1793-6519},
   MRCLASS = {14A20 (18D30)},
  MRNUMBER = {2945649},
MRREVIEWER = {Jon\ Eivind\ Vatne},
       DOI = {10.1142/S0129167X12500693},
       URL = {https://doi.org/10.1142/S0129167X12500693},
}

@article {GS-log,
    AUTHOR = {Gross, M. and Siebert, B.},
     TITLE = {Logarithmic {G}romov-{W}itten invariants},
   JOURNAL = {J. Amer. Math. Soc.},
  FJOURNAL = {Journal of the American Mathematical Society},
    VOLUME = {26},
      YEAR = {2013},
    NUMBER = {2},
     PAGES = {451--510},
      ISSN = {0894-0347,1088-6834},
   MRCLASS = {14N35 (14D23)},
  MRNUMBER = {3011419},
MRREVIEWER = {Hsian-Hua\ Tseng},
       DOI = {10.1090/S0894-0347-2012-00757-7},
       URL = {https://doi.org/10.1090/S0894-0347-2012-00757-7},
}

@book{Har77,
	AUTHOR = {Hartshorne, R.},
	TITLE = {Algebraic geometry},
	SERIES = {Graduate Texts in Mathematics, No. 52},
	PUBLISHER = {Springer-Verlag, New York-Heidelberg},
	YEAR = {1977},
	PAGES = {xvi+496},
	ISBN = {0-387-90244-9},
	MRCLASS = {14-01},
	MRNUMBER = {0463157},
	MRREVIEWER = {Robert Speiser},
}

@article {HS24,
    AUTHOR = {Hirschi, A. and Swaminathan, M.},
     TITLE = {Global {K}uranishi charts and a product formula in symplectic
              {G}romov-{W}itten theory},
   JOURNAL = {Selecta Math. (N.S.)},
  FJOURNAL = {Selecta Mathematica. New Series},
    VOLUME = {30},
      YEAR = {2024},
    NUMBER = {5},
     PAGES = {Paper No. 87, 74},
      ISSN = {1022-1824,1420-9020},
   MRCLASS = {53D45 (14N35)},
  MRNUMBER = {4807086},
       DOI = {10.1007/s00029-024-00982-y},
       URL = {https://doi.org/10.1007/s00029-024-00982-y},
}

@article{Ionel-Parker-relGW,
    AUTHOR = {Ionel, E.-N. and Parker, T. H.},
     TITLE = {Relative {G}romov-{W}itten invariants},
   JOURNAL = {Ann. of Math. (2)},
  FJOURNAL = {Annals of Mathematics. Second Series},
    VOLUME = {157},
      YEAR = {2003},
    NUMBER = {1},
     PAGES = {45--96},
      ISSN = {0003-486X,1939-8980},
   MRCLASS = {53D45 (57R17 57R57)},
  MRNUMBER = {1954264},
MRREVIEWER = {Ignasi\ Mundet-Riera},
       DOI = {10.4007/annals.2003.157.45},
       URL = {https://doi.org/10.4007/annals.2003.157.45},
}

@article {FKato-logcurve,
    AUTHOR = {Kato, F.},
     TITLE = {Log smooth deformation and moduli of log smooth curves},
   JOURNAL = {Internat. J. Math.},
  FJOURNAL = {International Journal of Mathematics},
    VOLUME = {11},
      YEAR = {2000},
    NUMBER = {2},
     PAGES = {215--232},
      ISSN = {0129-167X,1793-6519},
   MRCLASS = {14D22 (14H10)},
  MRNUMBER = {1754621},
MRREVIEWER = {Samuel\ Dalalyan},
       DOI = {10.1142/S0129167X0000012X},
       URL = {https://doi.org/10.1142/S0129167X0000012X},
}

@incollection {KKato-log,
    AUTHOR = {Kato, K.},
     TITLE = {Logarithmic structures of {F}ontaine-{I}llusie},
 BOOKTITLE = {Algebraic analysis, geometry, and number theory ({B}altimore,
              {MD}, 1988)},
     PAGES = {191--224},
 PUBLISHER = {Johns Hopkins Univ. Press, Baltimore, MD},
      YEAR = {1989},
      ISBN = {0-8018-3841-X},
   MRCLASS = {14F30 (14G20)},
  MRNUMBER = {1463703},
MRREVIEWER = {Adolfo\ Quir\'os},
}

@article {KKato-toric-singularities,
    AUTHOR = {Kato, K.},
     TITLE = {Toric singularities},
   JOURNAL = {Amer. J. Math.},
  FJOURNAL = {American Journal of Mathematics},
    VOLUME = {116},
      YEAR = {1994},
    NUMBER = {5},
     PAGES = {1073--1099},
      ISSN = {0002-9327,1080-6377},
   MRCLASS = {14M25 (11G99 14B05)},
  MRNUMBER = {1296725},
MRREVIEWER = {G.\ K.\ Sankaran},
       DOI = {10.2307/2374941},
       URL = {https://doi.org/10.2307/2374941},
}

@article{kleiman-rational-functions,
    AUTHOR = {Kleiman, S. L.},
     TITLE = {Misconceptions about {$K_X$}},
   JOURNAL = {Enseign. Math. (2)},
  FJOURNAL = {L'Enseignement Math\'ematique. Revue Internationale. 2e
              S\'erie},
    VOLUME = {25},
      YEAR = {1979},
    NUMBER = {3-4},
     PAGES = {203--206 (1980)},
      ISSN = {0013-8584},
   MRCLASS = {14A25},
  MRNUMBER = {570309},
MRREVIEWER = {Takao\ Fujita},
}

@article {Li-Ruan-relGW,
    AUTHOR = {Li, A.-M. and Ruan, Y.},
     TITLE = {Symplectic surgery and {G}romov-{W}itten invariants of
              {C}alabi-{Y}au 3-folds},
   JOURNAL = {Invent. Math.},
  FJOURNAL = {Inventiones Mathematicae},
    VOLUME = {145},
      YEAR = {2001},
    NUMBER = {1},
     PAGES = {151--218},
      ISSN = {0020-9910,1432-1297},
   MRCLASS = {53D45 (14J32 14N35 32Q25 32Q65)},
  MRNUMBER = {1839289},
MRREVIEWER = {Ignasi\ Mundet-Riera},
       DOI = {10.1007/s002220100146},
       URL = {https://doi.org/10.1007/s002220100146},
}

@article{Jun-Li-1,
    AUTHOR = {Li, J.},
     TITLE = {Stable morphisms to singular schemes and relative stable
              morphisms},
   JOURNAL = {J. Differential Geom.},
  FJOURNAL = {Journal of Differential Geometry},
    VOLUME = {57},
      YEAR = {2001},
    NUMBER = {3},
     PAGES = {509--578},
      ISSN = {0022-040X,1945-743X},
   MRCLASS = {14N35},
  MRNUMBER = {1882667},
MRREVIEWER = {Andreas\ Gathmann},
       URL = {http://projecteuclid.org/euclid.jdg/1090348132},
}

@article {Jun-Li-2,
    AUTHOR = {Li, J.},
     TITLE = {A degeneration formula of {GW}-invariants},
   JOURNAL = {J. Differential Geom.},
  FJOURNAL = {Journal of Differential Geometry},
    VOLUME = {60},
      YEAR = {2002},
    NUMBER = {2},
     PAGES = {199--293},
      ISSN = {0022-040X,1945-743X},
   MRCLASS = {14N35},
  MRNUMBER = {1938113},
       URL = {http://projecteuclid.org/euclid.jdg/1090351102},
}

@book{ogus-log-book,
    AUTHOR = {Ogus, A.},
     TITLE = {Lectures on logarithmic algebraic geometry},
    SERIES = {Cambridge Studies in Advanced Mathematics},
    VOLUME = {178},
 PUBLISHER = {Cambridge University Press, Cambridge},
      YEAR = {2018},
     PAGES = {xviii+539},
      ISBN = {978-1-107-18773-3},
   MRCLASS = {14D06 (14A20 14M25)},
  MRNUMBER = {3838359},
MRREVIEWER = {Howard\ M.\ Thompson},
       DOI = {10.1017/9781316941614},
       URL = {https://doi.org/10.1017/9781316941614},
}

@article {olsson-stack-of-log-str,
    AUTHOR = {Olsson, M.},
     TITLE = {Logarithmic geometry and algebraic stacks},
   JOURNAL = {Ann. Sci. \'{E}cole Norm. Sup. (4)},
  FJOURNAL = {Annales Scientifiques de l'\'{E}cole Normale Sup\'{e}rieure.
              Quatri\`eme S\'{e}rie},
    VOLUME = {36},
      YEAR = {2003},
    NUMBER = {5},
     PAGES = {747--791},
      ISSN = {0012-9593},
   MRCLASS = {14D20 (14A20)},
  MRNUMBER = {2032986},
MRREVIEWER = {Ivan\ S.\ Kausz},
       DOI = {10.1016/j.ansens.2002.11.001},
       URL = {https://doi.org/10.1016/j.ansens.2002.11.001},
}

@article {olsson-log-twisted-curves,
    AUTHOR = {Olsson, M.},
     TITLE = {({L}og) twisted curves},
   JOURNAL = {Compos. Math.},
  FJOURNAL = {Compositio Mathematica},
    VOLUME = {143},
      YEAR = {2007},
    NUMBER = {2},
     PAGES = {476--494},
      ISSN = {0010-437X,1570-5846},
   MRCLASS = {14D22 (14A20 14N35)},
  MRNUMBER = {2309994},
MRREVIEWER = {Charles\ D.\ Cadman},
       DOI = {10.1112/S0010437X06002442},
       URL = {https://doi.org/10.1112/S0010437X06002442},
}

@book {olsson-stacks-book,
    AUTHOR = {Olsson, M.},
     TITLE = {Algebraic spaces and stacks},
    SERIES = {American Mathematical Society Colloquium Publications},
    VOLUME = {62},
 PUBLISHER = {American Mathematical Society, Providence, RI},
      YEAR = {2016},
     PAGES = {xi+298},
      ISBN = {978-1-4704-2798-6},
   MRCLASS = {14D23 (14D22)},
  MRNUMBER = {3495343},
MRREVIEWER = {Stefan\ Schr\"oer},
       DOI = {10.1090/coll/062},
       URL = {https://doi.org/10.1090/coll/062},
}

@article {Parker-exploded-vs-log,
    AUTHOR = {Parker, B.},
     TITLE = {Log geometry and exploded manifolds},
   JOURNAL = {Abh. Math. Semin. Univ. Hambg.},
  FJOURNAL = {Abhandlungen aus dem Mathematischen Seminar der Universit\"at
              Hamburg},
    VOLUME = {82},
      YEAR = {2012},
    NUMBER = {1},
     PAGES = {43--81},
      ISSN = {0025-5858,1865-8784},
   MRCLASS = {53D45 (14N35)},
  MRNUMBER = {2922725},
MRREVIEWER = {Mark\ Gross},
       DOI = {10.1007/s12188-012-0065-8},
       URL = {https://doi.org/10.1007/s12188-012-0065-8},
}

@article {RSPW,
    AUTHOR = {Ranganathan, D. and {Santos-Parker}, K. and Wise, J.},
     TITLE = {Moduli of stable maps in genus one and logarithmic geometry,
              {I}},
   JOURNAL = {Geom. Topol.},
  FJOURNAL = {Geometry \& Topology},
    VOLUME = {23},
      YEAR = {2019},
    NUMBER = {7},
     PAGES = {3315--3366},
      ISSN = {1465-3060,1364-0380},
   MRCLASS = {14N35 (14D23)},
  MRNUMBER = {4046967},
MRREVIEWER = {Amin\ Gholampour},
       DOI = {10.2140/gt.2019.23.3315},
       URL = {https://doi.org/10.2140/gt.2019.23.3315},
}

@misc{Ranganathan-personal-comm,
  author       = {Ranganathan, D.},
  title        = {Personal communication},
  howpublished = {},
  year         = {2024},
}

@misc{stacks-project,
  author       = {The {Stacks project authors}},
  title        = {The Stacks project},
  howpublished = {\url{https://stacks.math.columbia.edu}},
  year         = {2025},
}

@article{vistoli-lci,
  title={The deformation theory of local complete intersections},
  author={Vistoli, A.},
  journal={arXiv:alg-geom\slash 9703008v2},
  year={1999}
}

\end{document}